\documentclass[11pt]{article}
\usepackage{verbatim,amsmath,amssymb}
\usepackage{epsfig,float,color}
\usepackage{geometry}
\usepackage{setspace}
\usepackage{subfig}
\usepackage{natbib}
\usepackage{placeins}
\usepackage{amsthm,mathrsfs,amsfonts,dsfont} 
\usepackage{hyperref}
\usepackage{enumitem}
\usepackage{xcolor}
\usepackage[normalem]{ulem}
\usepackage[toc,page]{appendix}
\usepackage{tcolorbox}

\usepackage[abs]{overpic}
\usepackage{color, colortbl}

\definecolor{darkblue}{rgb}{0,0,1}
\definecolor{dgreen}{rgb}{0,0.5,0}
\definecolor{cgn}{rgb}{1,0,1}
\hypersetup{pdftex=true, colorlinks=true, breaklinks=true, linkcolor=darkblue, menucolor=darkblue, pagecolor=darkblue, citecolor=darkblue, urlcolor=darkblue}

\newcommand\crule[3][white]{\textcolor{#1}{\rule{#2}{#3}}}

\newcommand {\eqb}[1]{\begin{equation}\begin{array}{#1}}
\newcommand {\eqe}{\end{array}\end{equation}}

\newcommand {\esb}[1]{\begin{equation*}\begin{array}{#1}}
\newcommand {\ese}{\end{array}\end{equation*}}
\newcommand {\ds}{\displaystyle}

\newcommand {\pa}[2]{\frac{\partial{#1}}{\partial{#2}}}

\newcommand {\back}{\! \! \!}
\newcommand {\is}{\back &=& \back}
\newcommand {\dis}{\back &:=& \back}

\newcommand {\plus}{\back &+& \back}

\newcommand {\norm}[1]{\|#1\|}

\newcommand {\dif}{\mathrm{d}}

\newcommand {\II}{{I\kern-.3em I}}
\newcommand {\III}{{I\kern-.3em I\kern-.3em I}}

\newcommand {\mrb}{\mathrm{b}}
\newcommand {\mrc}{\mathrm{c}}
\newcommand {\mrd}{\mathrm{d}}

\newcommand {\mrf}{\mathrm{f}}

\newcommand {\mri}{\mathrm{i}}

\newcommand {\mrm}{\mathrm{m}}
\newcommand {\mrn}{\mathrm{n}}

\newcommand {\mrp}{\mathrm{p}}

\newcommand {\mrr}{\mathrm{r}}
\newcommand {\mrs}{\mathrm{s}}

\newcommand {\mrx}{\mathrm{x}}

\newcommand {\mf}{\mathbf{f}}

\newcommand {\mk}{\mathbf{k}}

\newcommand {\mm}{\mathbf{m}}

\newcommand {\mx}{\mathbf{x}}

\newcommand {\ba}{\boldsymbol{a}}

\newcommand {\be}{\boldsymbol{e}}
\newcommand {\bff}{\boldsymbol{f}}

\newcommand {\bj}{\boldsymbol{j}}

\newcommand {\bn}{\boldsymbol{n}}

\newcommand {\bq}{\boldsymbol{q}}

\newcommand {\bv}{\boldsymbol{v}}
\newcommand {\bw}{\boldsymbol{w}}
\newcommand {\bx}{\boldsymbol{x}}

\newcommand {\bnu}{\mbox{\boldmath$\nu$}}
\newcommand {\bmu}{\mbox{\boldmath$\mu$}}

\newcommand {\bphi}{\mbox{\boldmath$\phi$}}

\newcommand {\mC}{\mathbf{C}}

\newcommand {\mK}{\mathbf{K}}

\newcommand {\mM}{\mathbf{M}}
\newcommand {\mN}{\mathbf{N}}

\newcommand {\mX}{\mathbf{X}}

\newcommand {\bA}{\boldsymbol{A}}
\newcommand {\bB}{\boldsymbol{B}}
\newcommand {\bC}{\boldsymbol{C}}

\newcommand {\bM}{\boldsymbol{M}}
\newcommand {\bN}{\boldsymbol{N}}

\newcommand {\bT}{\boldsymbol{T}}

\newcommand {\bX}{\boldsymbol{X}}

\newcommand {\sig}{\sigma}

\newcommand {\bsig}{\mbox{\boldmath$\sigma$}}

\newcommand {\bone}{\mathbf{1}}

\newcommand {\bbR}{\mathbb{R}}

\newcommand {\IR}{{\rm\kern.24em
   \vrule width.02em height1.53ex depth-.05ex
   \kern-.3em R}}
\newcommand {\ic}{{\rm\kern.20em
   \vrule width.02em height1.0ex depth-.05ex
   \kern-.22em c}}
\newcommand {\ia}{{\rm\kern.20em
   \vrule width.02em height1.05ex depth-.0ex
   \kern-.25em a}}
\newcommand {\IC}{{\rm\kern.24em
   \vrule width.02em height1.4ex depth-.05ex
   \kern-.26em C}}
\newcommand {\ID}{{\rm\kern.34em
   \vrule width.02em height1.5ex depth-.05ex
   \kern-.36em D}}
\newcommand {\IS}{{\rm\kern.24em
   \vrule width.02em height1.6ex depth.05ex
   \kern-.26em S}}
\newcommand {\IT}{{\rm\kern.50em
   \vrule width.02em height1.55ex depth-.05ex
   \kern-.52em T}}

\newcommand {\IE}{{\rm\kern.24em
   \vrule width.02em height1.55ex depth-.05ex
   \kern-.33em E}}
\newcommand {\IEa}{{\rm\kern.24em
   \vrule width.02em height1.55ex depth-.05ex
   \kern-.33em E}^{1}_{ijkl}}
\newcommand {\IEb}{{\rm\kern.24em
   \vrule width.02em height1.55ex depth-.05ex
   \kern-.33em E}^{2}_{ijkl}}

\newcommand {\sS}{\mathcal{S}}

\newcommand {\sU}{\mathcal{U}}
\newcommand {\sV}{\mathcal{V}}

\newcommand{\splSpace}[1]{\mathbb{S}_{#1}}
\newcommand{\splBasis}{N^A}
\newcommand{\splBasisWOA}{N}
\newcommand{\splBasisMod}{N^D}
\newcommand{\mbf}[1]{\mathbf{#1}}

\newcommand {\Ass}[2]{\kern 0.9ex \vrule width0.45em height0.2ex depth0ex \kern -2.1ex \bigwedge_{#1}^{#2}}
\newcommand {\ASS}[2]{\kern 1.45ex \vrule width0.5em height0.2ex depth0ex \kern -2.65ex \bigwedge_{#1}^{#2}}

\newcommand{\mrT}{\mathrm{T}}

\newcommand {\mrB}{\mathrm{B}}

\newcommand{\nablas}{\nabla_{\!\!\mrs}}

\begin{document}

\begin{center}
\Large{\bf{An isogeometric finite element formulation for phase transitions on deforming surfaces}}\\

\end{center}

\begin{center}
\large{Christopher Zimmermann$^\ast$, Deepesh Toshniwal$^\ddagger$, Chad M. Landis$^\ddagger$, \\[1mm] Thomas J.R. Hughes$^\ddagger$, Kranthi K. Mandadapu$^{\dagger\S}$\footnote{corresponding author, email: kranthi@berkeley.edu}, Roger A. Sauer$^\ast$\footnote{corresponding author, email: sauer@aices.rwth-aachen.de}}
\vspace{4mm}

\small{\textit{$^\ast$Aachen Institute for Advanced Study in Computational Engineering Science (AICES), \\ RWTH Aachen University, Templergraben 55, 52056 Aachen, Germany}}

\small{\textit{$^\ddagger$Institute for Computational Engineering and Sciences,\\ The University of Texas at Austin, 1 University Station, C0200,\\ 201 E. 24th Street, Austin, TX 78712, USA}}

\small{\textit{$^\dagger$Department of Chemical and Biomolecular Engineering,\\ University of California at Berkeley,
110A Gilman Hall, Berkeley, CA 94720-1460, USA}}

\small{\textit{$^\S$Chemical Sciences Division, Lawrence Berkeley National Laboratory, CA 94720, USA}}



\vspace{4mm}

Published\footnote{This pdf is the personal version of an article whose final publication is available at \href{http://dx.doi.org/10.1016/j.cma.2019.03.022}{www.sciencedirect.com}} 
in \textit{Computer Methods in Applied Mechanics and Engineering}, \\
\href{http://dx.doi.org/10.1016/j.cma.2019.03.022}{DOI: 10.1016/j.cma.2019.03.022} \\
Submitted on 01.~June 2018, Revised on 15.~February 2019, Accepted on 8.~March 2019

\end{center}



\rule{\linewidth}{.15mm}
{\bf Abstract}\\
This paper presents a general theory and isogeometric finite element implementation for studying mass conserving phase transitions on deforming surfaces.
The mathematical problem is governed by two coupled fourth-order nonlinear partial differential equations (PDEs) that live on an evolving two-dimensional manifold.
For the phase transitions, the PDE is the Cahn-Hilliard equation for curved surfaces, which can be derived from surface mass balance in the framework of irreversible thermodynamics. 
For the surface deformation, the PDE is the (vector-valued) Kirchhoff-Love thin shell equation.
Both PDEs can be efficiently discretized using $C^1$-continuous interpolations without derivative degrees-of-freedom (dofs).
Structured NURBS and unstructured spline spaces with pointwise $C^1$-continuity are utilized for these interpolations.
The resulting finite element formulation is discretized in time by the generalized-$\alpha$ scheme with adaptive time-stepping, and it is fully linearized within a monolithic Newton-Raphson approach.
A curvilinear surface parameterization is used throughout the formulation to admit general surface shapes and deformations.
The behavior of the coupled system is illustrated by several numerical examples exhibiting phase transitions on deforming spheres, tori and double-tori.


\vspace{-4mm}
\rule{\linewidth}{.15mm}


\tableofcontents

\section{Introduction}

A wide range of biological, chemical, electro- and thermo-mechanical applications are governed by phase transitions, which include de-mixing of a well-mixed phase into two separate phases. 
For example, in electro-chemical devices such as batteries \citep{tang10,ebner13}, phase transitions can affect the resulting mechanical and kinetic behavior. 
In biology, it is known that lipid membranes can separate into two distinct phases when quenched from high temperatures to low temperatures depending on the mole fraction of the constituents that make up the membrane \citep{veatch03}. 
Under temperature quenches, these two-dimensional lipid membranes can undergo severe shape changes as a result of the coupling between in-plane phase transitions and out-of-plane bending \citep{baumgart03}. 
This interplay between in-plane phase transitions and out-of-plane bending has not been explored in its entirety, except for simple situations where the membrane deformations are either axi-symmetric or small. 
Recently, \citet{sahu17} presented a general theory to describe the coupling between in-plane phase transitions and out-of-plane bending for arbitrarily curved surfaces, employing the framework of irreversible thermodynamics. 
Specifically, this new theory introduces Korteweg stresses induced by in-plane phase transitions in the context of deformable surfaces and shows how they couple to out of plane deformations. 
This theory can be regarded as an extension of the Cahn-Hillard theory \citep{cahn58-1, cahn61} to arbitrarily curved surfaces. 
To study the coupling between in-plane phase transitions and surface deformations governed by the theory of \citet{sahu17} requires the development of suitable numerical methods.

Modeling phase transitions requires defining an order parameter 
that distinguishes the 
phases. 
The evolution of the phases is described by the Cahn-Hilliard theory that results in a partial differential equation (PDE) that is of fourth order in the order parameter. 
Deforming surfaces are commonly described by the Kirchhoff-Love thin shell equation, which is a vector-valued PDE that is of fourth order in the out-of-plane deformation. 
The standard weak forms of these fourth-order PDEs involve products of second-order derivatives.
Such weak forms require either using globally $C^1$-continuous discretizations \citep{gomez08-1, bartezzaghi15, kastner16-1}, mixed formulations \citep{elliott89,barrett99} or discontinuous Galerkin methods \citep{wells06,xia07}.
The latter two avoid the necessity of global $C^1$-continuity. 
They lead, however, to an increase of the computational cost, since additional dofs or additional operators are required.
Further, mixed methods have to satisfy additional stability requirements.
$C^1$-continuous formulations, on the other hand, avoid this overhead and thus provide a more direct numerical approach.

A very powerful methodology that allows for $C^1$-continuous discretizations within the finite element (FE) method is isogeometric analysis (IGA) \citep{hughes05}.
This stems from the fact that the high-order discretizations of IGA also provide much better spectral behavior  \citep{hughes05,cottrell2006isogeometric,cottrell2007studies}, efficiency \citep{akkerman2008role,morganti2015patient} and robustness \citep{lipton2010robustness} when compared to their $C^0$-continuous FE counterparts.
Within IGA, global B-spline- and NURBS-patches are the most widely used basis functions \citep{cottrell}. 
In recent years, these have been extended to local refinement techniques using T-splines \citep{scott12}, hierarchical B-splines \citep{hollig2003finite,schill12}, truncated hierarchical B-splines \citep{giannelli2012thb}, locally refinable (LR) B-splines \citep{dokken13,johannessen14} and LR NURBS \citep{zimmermann17}.

IGA on any sufficiently complex geometry of arbitrary topology requires parametric representations containing isolated parameterization singularities. With regard to quadrilateral meshes, the two types of singularities employed are corner singularities, called \emph{extraordinary points} \citep{scott2013isogeometric,toshniwal2017smooth}, and collapsed-edge singularities, called \emph{polar points} \citep{myles2011c2polar,toshniwal2017multi}. 
While the latter can be used for surfaces of genus zero, the former can be used to handle surfaces of arbitrary genii.
The construction of smooth splines on meshes containing such singularities must follow special rules.
In this work, we employ the bi-cubic splines construction presented in \citet{toshniwal2017smooth}.

Recent works have demonstrated the benefit of using IGA in the context of phase transitions. 
Examples are the study of spinodal decompositions of binary mixtures \citep{gomez08-1, bartezzaghi15, kastner16-1}, 
spinodal decompositions under shear flow \citep{liu13-1}, 
topology optimization \citep{dede12}, 
phase segregation in Li-ion electrodes \citep{stein2014,dileo14-1,zhao15-1,zhao2016,Xu2016} 
and fracture mechanics \citep{borden201277, borden14-1, borden16-1}.
IGA and other techniques have been used to study phase transitions on fixed surfaces \citep{mercker12,bartezzaghi15}. 
General phase transitions on deforming surfaces, however, have not yet been studied with IGA:
The approaches that exist use spring-based network models \citep{mcwhirter04}, mixed FE methods \citep{elliot10}, 2D and axi-symmetric formulations \citep{embar13}, or use a second phase-field in order to describe the surface in a diffuse manner \citep{wang08,lowengrub09}. 

Other approaches, which have been used for surface PDEs are spectral finite element methods, e.g. \citet{taylor97}, trace finite element methods, e.g. \citet{reusken14}, level-set methods, e.g. \citet{sethian99,bertalimo01},
and evolving surface finite element methods, e.g. \citet{dziuk07,dziuk12}. The latter is applied to the Cahn-Hilliard equation by \citet{eilks08} and analyzed by \citet{elliott15}.
There are also related works on PDEs on rigidly rotating \citep{taylor97} and moving surfaces \citep{elliott09}.

Since a general IGA formulation for deforming surfaces is still lacking, it is studied in the present work.
The proposed formulation is based on the theory of \citet{sahu17}, which is combined with the isogeometric shell model of \citet{solidshell2}.
A monolithic and fully implicit time integration scheme is used to solve the coupled system based on the generalized-$\alpha$ method of \citet{chung93}. 
The proposed formulation features the following novelties:
\begin{itemize}
	\item it couples phase transitions with general surface deformations,
	\item it accounts for geometrical and material nonlinearities,
	\item it is implemented within a monolithic and fully implicit finite element formulation,
	\item it uses an automatic, adaptive time-stepping scheme,
	\item it uses isogeometric surface discretizations based on unstructured spline spaces, and
	\item it is used to determine and study the surface Korteweg stresses.
\end{itemize}

The remainder of this paper is organized as follows.
Sec.~\ref{Sec:def_s} summarizes the description of deforming surfaces. 
The balance laws for mass and momentum are presented in Sec.~\ref{s:bal}, while Sec.~\ref{s:consti} presents the corresponding constitutive equations. 
Those lead to the weak form of Sec.~\ref{s:wf}. 
The spatial and temporal discretization of the coupled problem is then presented in Sec.~\ref{sec:disc}. 
Sec.~\ref{s:num_ex} then shows several numerical examples that illustrate the coupled model behavior. 
The paper concludes with  Sec.~\ref{s:concl}.

\section{Deforming surfaces} 
\label{Sec:def_s}
This section gives a brief summary of the general description of curved surfaces and their deformation according to Kirchhoff-Love kinematics. 
A more detailed description can be found for example in \citet{cism}.

\subsection{Surface description}

In general, a curved surface can be denoted by a set of surface points $\bx\in\sS$. Their motion can be described by the mapping
\eqb{l}
\bx = \bx(\xi^\alpha,t)\,,
\label{e:bx}\eqe
where $\xi^\alpha$, $\alpha=1,2$ denote the coordinates (or parameters) associated with a material point on the surface. 
Such coordinates are also termed convected coordinates.\footnote{See \cite{sahu17} for a description of the surface using different coordinate parametrizations.}
The tangent vectors at $\bx$ then follow from
\eqb{l}
\ba_\alpha := \ds\pa{\bx}{\xi^\alpha}\,.
\label{e:ba}\eqe
They define the surface metric 
\eqb{l}
a_{\alpha\beta} := \ba_\alpha\cdot\ba_\beta\,,
\eqe
the surface normal
\eqb{l}
\bn := \ds\frac{\ba_1\times\ba_2}{\norm{\ba_1\times\ba_2}}\,,
\label{e:bn}\eqe
and the contravariant tangent vectors 
\eqb{l}
\ba^\alpha = a^{\alpha\beta}\,\ba_\beta\,,
\eqe
through $[a^{\alpha\beta}] = [a_{\alpha\beta}]^{-1}$. 
Here, all Greek indices run from 1 to 2 and obey the Einstein summation convention. 
The second parametric derivative $\ba_{\alpha,\beta} := \partial\ba_\alpha/\partial\xi^\beta$ defines the curvature components
\eqb{l}
b_{\alpha\beta} := \ba_{\alpha,\beta}\cdot\bn
\label{e:bab}\eqe
and the mean surface curvature
\eqb{l}
H := a^{\alpha\beta}\,b_{\alpha\beta}/2\,.
\label{e:H}\eqe

Given the parametrization in \eqref{e:bx}, the surface gradient, surface divergence and surface Laplacian can be defined, respectively, as
\eqb{lllll}
\mathrm{grad}_\mrs\phi \dis \nablas\phi \dis \phi_{;\alpha}\,\ba^\alpha\,, \\[1mm]
\mathrm{div}_\mrs\bv \dis \nablas\cdot\bv \dis \bv_{;\alpha}\cdot\ba^\alpha = v^\alpha_{;\alpha} - 2Hv\,, \\[1mm]
\Delta_\mrs\phi \dis \nablas\cdot\nablas\phi \is \phi_{;\alpha\beta}\,a^{\alpha\beta}\,,
\label{e:diff}\eqe
where $\phi$ and $\bv$ denote general scalars and vectors and $v^\alpha:=\bv\cdot\ba^\alpha$ and $v:=\bv\cdot\bn$ are the vector components corresponding to the $\{\ba_1,\ba_2,\bn\}$ basis.
The symbol `;' denotes the covariant derivative. 
It is equal to the parametric derivative for general scalars and vectors, i.e.
$\phi_{;\alpha} = \phi_{,\alpha}:=\partial\phi/\partial\xi^\alpha$ and
$\bv_{;\alpha} = \bv_{,\alpha}:=\partial\bv/\partial\xi^\alpha$.
However, $v^\alpha_{;\beta} \ne v^\alpha_{,\beta}$ and $\phi_{;\alpha\beta} \ne \phi_{,\alpha\beta}$. 
Instead
\eqb{lll}
v^\alpha_{;\beta} \is v^\alpha_{,\beta} + \Gamma^\alpha_{\beta\gamma}\,v^\gamma\,, \\[2mm]
\phi_{;\alpha\beta} \is \phi_{,\alpha\beta} - \Gamma^\gamma_{\alpha\beta}\,\phi_{,\gamma}\,,
\label{e:phiab}\eqe
where $\Gamma^\gamma_{\alpha\beta} := \ba_{\alpha,\beta}\cdot\ba^\gamma$ are the Christoffel symbols of the second kind on surface $\sS$.

\subsection{Surface kinematics}
\label{s:kin}
Given the motion of the surface over time $t$ in \eqref{e:bx}, we can define the surface at $t=0$ as a reference configuration and denote it $\sS_0$.
The set of surface points $\bX\in\sS_0$ follow from $\bX := \bx(\xi^\alpha,0)$.
Analogous to Eqs.~\eqref{e:ba}--\eqref{e:H}, the surface quantities
$\bA_\alpha := \partial\bX/\partial\xi^\alpha$,
$A_{\alpha\beta} :=  \bA_\alpha\cdot\bA_\beta$,
$\bN := \bA_1\times\bA_2/\norm{\bA_1\times\bA_2}$,
$\bA^\alpha := A^{\alpha\beta}\,\bA_\beta$,
$[A^{\alpha\beta}] := [A_{\alpha\beta}]^{-1}$,
$B_{\alpha\beta} := \bA_{\alpha,\beta}\cdot\bN$,
and
$H_0 := A^{\alpha\beta}\,B_{\alpha\beta}/2$
are introduced.
The surface kinematics are then characterized by the relation between corresponding objects on $\sS_0$ and $\sS$.
An example is the left surface Cauchy-Green tensor 
\eqb{l}
\bB = A^{\alpha\beta}\,\ba_\alpha\otimes\ba_\beta\,,
\eqe
that has the two invariants
\eqb{l}
I_1 := A^{\alpha\beta}\,a_{\alpha\beta}
\eqe
and
\eqb{l}
J := \sqrt{\det[A^{\alpha\beta}]\det[a_{\alpha\beta}]}\,,
\label{e:J}\eqe
which characterizes the change in surface area between $\sS_0$ and $\sS$. 
\\
The material velocity at $\bx$ is given by
\eqb{l}\label{e:vel}
\bv := \dot\bx~,
\eqe
where the material time derivative is defined by
\eqb{l}
\dot{(...)} := \ds\pa{...}{t}\Big|_{\xi^\alpha=\,\mathrm{fixed}}\,.
\eqe
The velocity vector in \eqref{e:vel} can be used to define the material time derivatives of various surface quantities such as
\eqb{l}
\dot\ba_\alpha = \bv_{,\alpha} = \ds\pa{\bv}{\xi^\alpha}
\eqe
and
\eqb{l}
\dot a_{\alpha\beta} = \ba_\alpha\cdot\dot\ba_\beta + \dot\ba_\alpha\cdot\ba_\beta\,.
\eqe

\subsection{Surface variations}

In order to formulate the weak form of the governing PDEs for thin shells, the variations of various surface measures are needed. 
For example, considering a kinematically admissible variation of the deformation, denoted $\delta\bx$, we can write
\eqb{lll}
\delta a_{\alpha\beta} \is \ba_\alpha\cdot\delta\ba_\beta + \delta\ba_\alpha\cdot\ba_\beta\,, \\[2mm]
\delta b_{\alpha\beta} \is \big(\delta\ba_{\alpha,\beta}-\Gamma^\gamma_{\alpha\beta}\,\delta\ba_\gamma\big)\cdot\bn\,, \\[2mm]
\delta\bn \is -(\ba^\alpha\otimes\bn)\,\delta\ba_\alpha\,,
\eqe
where $\delta\ba_\alpha = \partial(\delta\bx)/\partial\xi^\alpha$ and $\delta\ba_{\alpha,\beta} = \partial(\delta\ba_\alpha)/\partial\xi^\beta$.
The variation of further measures related to deforming surfaces can be found in \citet{shelltheo2}.

\section{Balance laws}\label{s:bal}

This section gives a brief summary of the equations that govern the physical behavior of thin shells. 
They follow from the balance laws of mass and momentum and describe the evolution of the surface concentration and shape, respectively. 
A detailed derivation of the surface balance laws for multicomponent systems in the framework of irreversible thermodynamics can be found in \cite{sahu17}.

\subsection{Balance of mass}

Consider that surface $\sS$ consists of two species with the mass densities per unit area $\rho_1$ and $\rho_2$. 
The total mass of each species is assumed to be conserved. This implies that the total density $\rho=\rho_1+\rho_2$ satisfies $\rho = \hat\rho/J$, where $\hat\rho$ denotes the initial density and $J$ is the area change defined in~\eqref{e:J}.
The dimensionless concentration $\phi:=\rho_{1}/\rho$ is sufficient to model the local density fractions of both species. 
The concentration field $\phi$ is sometimes also denoted as \textit{order parameter field} or \textit{phase field}.
The rate of change of $\phi$ follows as 
\eqb{l}
\rho\,\dot\phi = -j^\alpha_{;\alpha}~\quad\forall\,\bx\in\sS
\label{e:sf}\eqe
\citep{sahu17}, where
\eqb{l}
j^\alpha = a^{\alpha\beta}\,j_\beta 
\eqe
and $j_\alpha=\bj\cdot\ba_\alpha$ are the contra- and covariant components of the diffusive surface flux vector $\bj$, respectively. 
They follow from the constitutive equations discussed in Sec.~\ref{s:consti}. As long as there is no mass inflow from the boundary, such as is considered here, the mass of each species is conserved by Eq.~\eqref{e:sf}.

\subsection{Balance of momentum}\label{s:mombal}

From the balance of linear momentum for an arbitrarily deforming surface $\sS$ follows the equation of motion
\eqb{l}
\rho\,\dot{\bv} = \bT^{\alpha}_{;\alpha} + \bff~\quad\forall\,\bx\in\sS\,,
\label{e:sfm}\eqe
where $\bff$ is a body force and 
\eqb{l}
\bT^{\alpha} = N^{\alpha\beta}\,\ba_\beta  + S^\alpha\,\bn\,
\label{e:Ta}\eqe 
($\alpha=1,2$) are the stress vectors that have the in-plane membrane components $N^{\alpha\beta}$ and the out-of-plane shear components $S^\alpha$ \citep{naghdi1971theory,steigmann99b,shelltheo2}. 
The stress vectors are related to the stress tensor  
\begin{equation}
\bsig = N^{\alpha \beta} \ba_\alpha \otimes \ba_\beta + S^\alpha \ba_\alpha \otimes \bn
\label{e:bsig}\end{equation}
through Cauchy's formula $\bT^\alpha = \bsig^\mrT \ba^\alpha$.
From this, the traction $\bT$, acting on any cut through the surface with outward normal $\bnu=\nu_\alpha\ba^\alpha$, 
follows as $\bT = \bsig^\mrT \bnu=\bT^\alpha\nu_\alpha$.
\\
Similarly, the moment on the cut can be written as $\bM  = \bmu^\mrT \bnu$, where
\begin{equation}
\bmu = -M^{\alpha \beta} \ba_\alpha \otimes \ba_\beta\,, 
\end{equation}
is the moment tensor that has the in-plane components $M^{\alpha\beta}$ \citep{shelltheo2,sahu17}.
The balance of angular momentum dictates that 
\eqb{lll} 
S^\alpha \is - M^{ \beta \alpha}_{;\beta}\,, \\[1mm]
\sigma^{\alpha \beta} \is \sigma^{\beta\alpha}\,,
\label{e:amom}\eqe
where $\sigma^{\alpha \beta} := N^{\alpha \beta} - b^\beta_{\gamma} M^{\gamma\alpha}$.
The stress components $\sig^{\alpha\beta}$ and $M^{\alpha\beta}$ follow from constitution, which is discussed in the following section.
\\
Combining Eqs. \eqref{e:sfm}, \eqref{e:bsig} and (\ref{e:amom}.1), the equation of motion can be written in the component form
\eqb{lll}
\rho\,a^\alpha \is f^\alpha + N^{\lambda \alpha}_{; \lambda} - S^\lambda b^\alpha_\lambda\,,\\[1mm]
\rho\,a_\mrn \is p + N^{\alpha \beta} b_{\alpha \beta} + S^\alpha_{; \alpha}\,,
\label{e:sfm2}\eqe
where $a^\alpha := \dot\bv\cdot\ba^\alpha$, $a_\mrn := \dot\bv\cdot\bn$, $f^\alpha := \bff\cdot\ba^\alpha$ and $p := \bff\cdot\bn$.

\section{Constitutive equations}\label{s:consti}

This section presents the constitutive equations for the diffusive flux vector and the stress and moment tensors using a combined elasticity and Cahn-Hilliard model.
The formulation follows the framework of irreversible thermodynamics of curved surfaces \citep{sahu17}.

\subsection{Helmholtz free energy}

\begin{figure}[H]
	\centering
	\includegraphics[width=0.5\linewidth, trim = 0 0 0 0,clip]{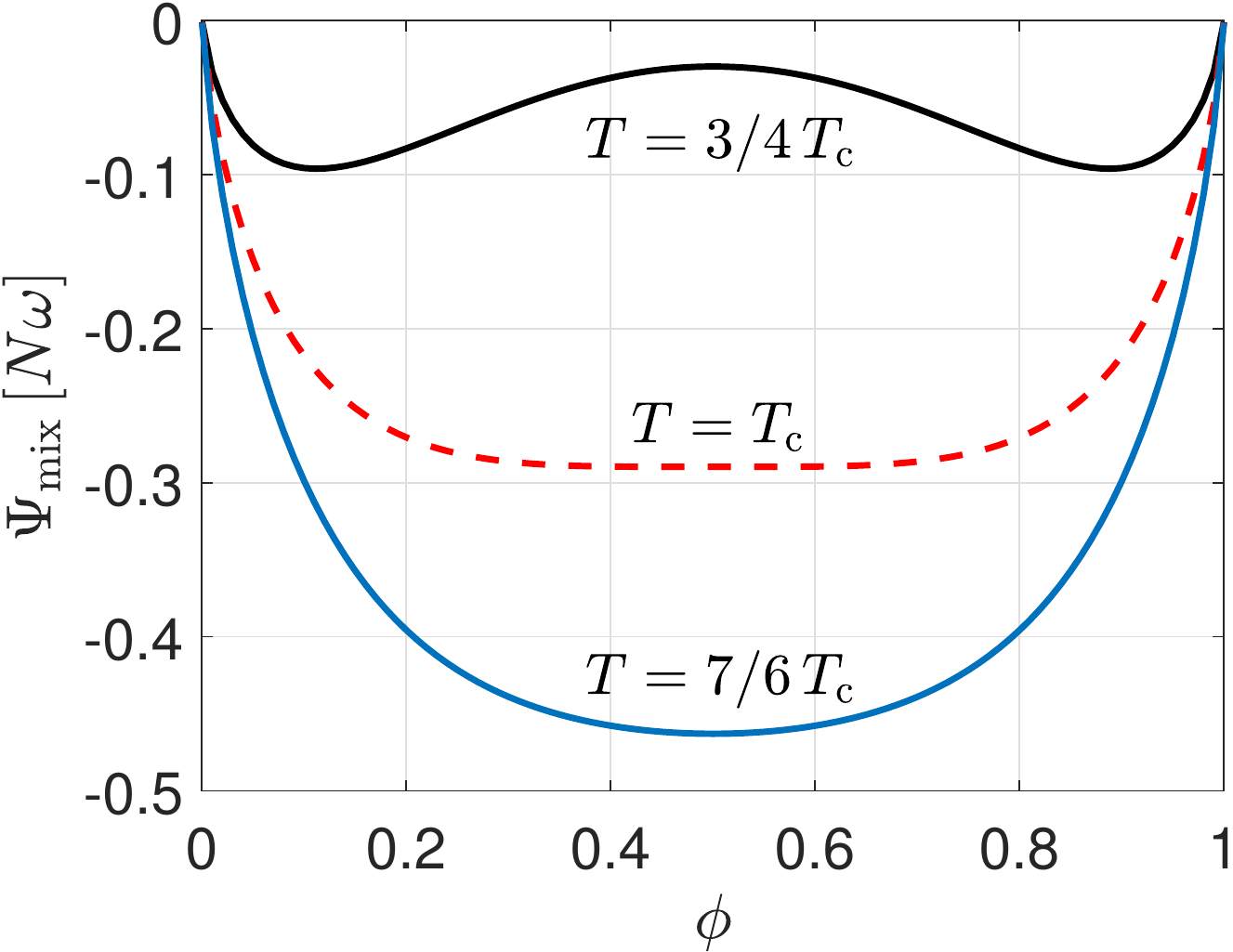}
	\caption{The free energy of mixing $\Psi_\mathrm{mix}$ vs.~$\phi$ for various temperatures $T$.}
	\label{f:Psimix}
\end{figure}

The constitutive description for the system can be obtained given an appropriate description of the Helmholtz free energy. 
In this paper, we consider phase transforming systems that 
are visco-elastic in-plane and elastic out-of-plane. 
In this context, the Helmholtz free energy per unit reference area, $\Psi$, is dependent on the metric tensor $a_{\alpha \beta}$, the curvature tensor $b_{\alpha \beta}$, the concentration field $\phi$, its surface gradient $\nablas\phi$ and the temperature $T$, i.e., 
\eqb{l}
\Psi = \Psi(a_{\alpha \beta}, b_{\alpha \beta}, \phi,\nablas\phi, T)\,.
\eqe
In what follows, we assume that the temperature is uniform across the surface $\sS$. 
The phase transformation is assumed to be governed by the Cahn-Hilliard energy combined with an elastic potential in an additive manner, i.e.,
\eqb{l}
\Psi = \Psi_\mathrm{el} + \Psi_\mathrm{CH}\,.
\label{e:Psi}\eqe
$\Psi_\mathrm{el}$ is taken as an additive composition of dilatational, deviatoric and bending energies in the form
\eqb{l}
\Psi_\mathrm{el} = \Psi_\mathrm{dil}(a_{\alpha\beta},\phi) + \Psi_\mathrm{dev}(a_{\alpha\beta},\phi) + \Psi_\mathrm{bend}(b_{\alpha\beta},\phi)\,.
\label{e:Wel}\eqe
Assuming the in-plane response to be isotropic, a suitable choice for the first two terms is the Neo-Hookean surface material model
\eqb{l}\label{e:Wdil}
\Psi_\mathrm{dil} = \ds\frac{K(\phi)}{4}\big(J^2-1-2\,\ln J\big)\,,
\eqe
and
\eqb{l}\label{e:Wdev}
\Psi_\mathrm{dev} = \ds\frac{G(\phi)}{2}\big(I_1/J-2\big)
\eqe
\citep{shelltheo2}.
Here, $K$ and $G$ denote the 2D bulk and shear moduli. 
The bending term is taken from the Koiter model
\eqb{l}
\Psi_\mathrm{bend} = \ds\frac{c(\phi)}{2}\big(b_{\alpha\beta}-B_{\alpha\beta}\big)\big(b^{\alpha\beta}_0-B^{\alpha\beta}\big)
\eqe
\citep{ciarlet}, where $b^{\alpha\beta}_0:=A^{\alpha\gamma}b_{\gamma\delta}A^{\beta\delta}$, and $c$ denotes the bending modulus. 
The moduli $K$, $G$ and $c$ are functions of $\phi$ according to the mixtures rules of Sec.~\ref{s:mix}. 
\\
In analogy to 3D problems \citep{cahn58-1}, the Cahn-Hillard energy for surfaces takes the form 
\eqb{l}
\Psi_\mathrm{CH} 
= \Psi_\mathrm{mix}(\phi,T)  + \Psi_\mri(J,\nablas\phi)\,.
\label{e:WCH}\eqe
Here $\Psi_\mathrm{mix}=W_\mathrm{mix}(\phi)- T\,\eta_\mathrm{mix}(\phi)$ is the free energy of mixing that contains the internal mixing energy
\eqb{l} 
W_\mathrm{mix} = N\omega\,\phi\,(1-\phi)
\eqe
and the mixing entropy 
\eqb{l} 
\eta_\mathrm{mix} = -Nk_\mrB\,\big(\phi\ln\phi + (1-\phi)\ln(1-\phi)\big)\,,
\eqe
and 
\eqb{l}
\Psi_\mri = J\,N\,\omega\ds\frac{\lambda}{2}\,\nablas\phi\cdot\nablas\phi
\label{e:Psii}\eqe
is the energy of maintaining an interface between the two species when the system is phase separated \citep{embar13}.
$N$ denotes the number of molecules per reference area, $k_\mrB$ is Boltzmann's constant, $\sqrt{\lambda}$ represents the length scale of the phase interface, and $\omega = 2\,k_\mrB\,T_\mrc$ is a bulk energy related to the critical temperature, $T_\mrc$, below which phase separation occurs.\footnote{In the subsequent examples, the temperature $T = 2T_\mrc/3$ is chosen, such that the minimization of $\Psi_\mathrm{mix}$ drives the phase separation.} 
$N$, $k_\mrB$, $\lambda$ and $\omega$ are treated as constants here.
The area stretch $J$ (see Eq.~\eqref{e:J}) is included in $\Psi_\mathrm{i}$, since $\Psi_\mathrm{CH}$ is an energy w.r.t.~the reference configuration, while $\nablas\phi$ refers to the current configuration (it can be viewed as having units of 1/(current length)).
For this reason the last term in $\Psi_\mathrm{CH}$ explicitly depends on $a_{\alpha\beta}$ apart from depending on $\nablas\phi$.
The first term in $\Psi_\mathrm{CH}$, on the other hand, is only a function of $\phi$ and $T$.
Fig.~\ref{f:Psimix} shows the variation of $\Psi_\mathrm{mix}$ with $\phi$ and $T$.
For $T>T_\mrc$, $\Psi_\mathrm{mix}$ has a single minimum -- indicating that a mixed state is preferred -- while
for $T<T_\mrc$, $\Psi_\mathrm{mix}$ has two minima -- indicating that a phase separated state is preferred. 

\subsection{Diffusive flux}\label{s:flux}

Given the Helmholtz free energy $\Psi$, the components of the diffusive flux can be written as
\eqb{l}
j_\alpha = -\ds\frac{M}{J}\,\mu_{\mrc,\alpha}\,,
\label{e:jaMJ}\eqe
where $M=D\,\phi\,(1-\phi)$, with $D=$ const., is the degenerate mobility\footnote{The mobility $M$ should not be confused with the bending moment components $M^{\alpha\beta}$.} \citep{wells06} and
\eqb{l}
\mu_\mrc = \mu_\mrb + \mu_\mri\,,
\label{e:muc}\eqe
is the chemical potential that has the bulk and interface contributions
\eqb{lll}
\mu_\mrb \dis \ds\pa{\Psi}{\phi}\,, \\[3mm]
\mu_\mri \dis -J\,\ds\bigg(\frac{1}{J}\pa{\Psi}{\phi_{;\alpha}}\bigg)_{\!\!;\alpha}\,,
\label{e:mubi}\eqe  
respectively (see Appendix~\ref{s:thermo}).
Division by $J$ is included in \eqref{e:jaMJ} since $j_\alpha$ relates to the current area, while $\mu_\mrc$ is defined per reference area.
From \eqref{e:Psi}-\eqref{e:Psii} we find
\eqb{lrl}
\mu_\mrb \is \mu_\phi + \mu_\mathrm{el}\,, \\
\mu_{\phi} \is Nk_\mathrm{B}T\ln{\ds\frac{\phi}{1-\phi}} + N\omega\,(1-2\phi)\,, \\[1.5mm]
\mu_\mathrm{el}\is \Psi'_\mathrm{el}\,, \\[1.5mm]
\mu_\mri \is -J\,N\,\omega\,\lambda\,\Delta_\mrs\phi\,,
\label{e:mubi2}\eqe
where $(...)':=\partial.../\partial\phi$.
The elastic contribution to the chemical potential follows from \eqref{e:Wel} as
\eqb{l}
\label{e:muel}
\mu_\mathrm{el} = \ds \frac{K'}{4}\big(J^2-1-2\,\ln J\big) + \frac{G'}{2}\big(I_1/J-2\big) + \frac{c'}{2}\big(b_{\alpha\beta}-B_{\alpha\beta}\big)\big(b^{\alpha\beta}_0-B^{\alpha\beta}\big)~.
\eqe
The diffusive flux can be decomposed as
\eqb{l}
j^\alpha = j^\alpha_{\phi} +  j^\alpha_\mathrm{el} + j^\alpha_\mri\,,
\label{e:ja}\eqe
with
\eqb{l}
j^\alpha_\bullet = - \ds\frac{M}{J}a^{\alpha\beta}\,\mu_{\bullet,\beta}\,,
\label{e:jad}\eqe
for the three different contributions.

\subsection{Stress and moments}\label{s:sigM}

The components of the stress and moment tensors 
follow from the Helmholtz free energy per reference area as
\eqb{lll}
\sig^{\alpha\beta} \is \ds\frac{2}{J}\pa{\Psi}{a_{\alpha\beta}} - \eta\,\dot a^{\alpha\beta}\,, \\[4mm]
M^{\alpha\beta} \is \ds\frac{1}{J}\pa{\Psi}{b_{\alpha\beta}}
\label{e:tauM0}\eqe
\citep{liquidshell,sahu17}.
The second term in (\ref{e:tauM0}.1) accounts for viscous in-plane stress considering finite linear surface shear viscosity \citep{rangamani13,rangamani14,sahu17}.
Here $\eta$ is the dynamic surface viscosity and $\dot a^{\alpha\beta} = -a^{\alpha\gamma}\,\dot a_{\gamma\delta}\,a^{\delta\beta}$ 
corresponds to the components of the surface velocity gradient multiplied by $-$2 \citep{cism}. 
Given the different contributions to the total Helmholtz free energy in \eqref{e:Psi}, the stress components follow as
\eqb{l}
\sig^{\alpha\beta} = \sig^{\alpha\beta}_\mathrm{el} +\sig^{\alpha\beta}_\mathrm{visc} + \sig^{\alpha\beta}_\mathrm{CH}\,, 
\label{e:sigtot}\eqe
where the elastic stress contribution is 
\eqb{l}
\sig^{\alpha\beta}_\mathrm{el} := \ds\frac{2}{J}\pa{\Psi_\mathrm{el}}{a_{\alpha\beta}} 
= \ds\frac{K(\phi)}{2J}\big(J^2-1\big)\,a^{\alpha\beta} + \frac{G(\phi)}{2J^2}\Big(2A^{\alpha\beta}-I_1\,a^{\alpha\beta}\Big)\,, 
\label{e:sigel}\eqe
the viscous stress contribution is 
\eqb{l}
\sig^{\alpha\beta}_\mathrm{visc} := -\ds\eta(\phi)\,\dot a^{\alpha\beta}\,,
\label{e:visc}
\eqe
and the \textit{Korteweg stresses} \citep{sahu17} 
due to the Cahn-Hilliard energy is given by 
\eqb{l}
\sig^{\alpha\beta}_\mathrm{CH} := \ds\frac{2}{J}\pa{\Psi_\mathrm{CH}}{a_{\alpha\beta}} = N\,\omega\frac{\lambda}{2}\big(a^{\alpha\beta}\,a^{\gamma\delta}-2a^{\alpha\gamma}a^{\beta\delta}\big)\,\phi_{;\gamma}\,\phi_{;\delta}\,.
\label{e:sigCH}\eqe
These Korteweg stresses lead to a coupling between in-plane phase transformations and out-of-plane bending according to (\ref{e:sfm2}.2). 
We illustrate the Korteweg stresses in the numerical examples of Sec.~\ref{s:num_ex}. 
\\
The components of the moment tensor only stem from $\Psi_{\mathrm{bend}}$ in \eqref{e:Wel}. 
They follow as 
\eqb{l}
M^{\alpha\beta}  = \ds\frac{c(\phi)}{J}\big(A^{\alpha\gamma}b_{\gamma\delta}A^{\beta\delta} - B^{\alpha\beta}\big)\,.
\label{e:Mab}\eqe
We emphasize that $\sig^{\alpha\beta}$ are the stresses following from constitution, but they are not the total stresses appearing in the equations of motion. 
These are
\eqb{l}
N^{\alpha \beta} = \sig^{\alpha \beta} + b^\beta_{\gamma} M^{\gamma\alpha}\,,
\label{e:Ntot}\eqe
as noted in Sec.~\ref{s:mombal}.

\subsection{Mixture rules}\label{s:mix}

In this section, we propose a model for the dependency of the material parameters on the field variable $\phi \in [0,1]$. Recall that $\phi$ characterizes the current composition of the mixture and the two separate phases are characterized by values close to 0 and close to 1. Due to the characteristics of $\Psi_\mathrm{mix}$ shown in Fig.~\ref{f:Psimix}, $\phi$ does not attain the exact values of 0 and 1. 
\begin{figure}[h]
	\centering
	\includegraphics[width=0.5\linewidth, trim = 0 0 0 0,clip]{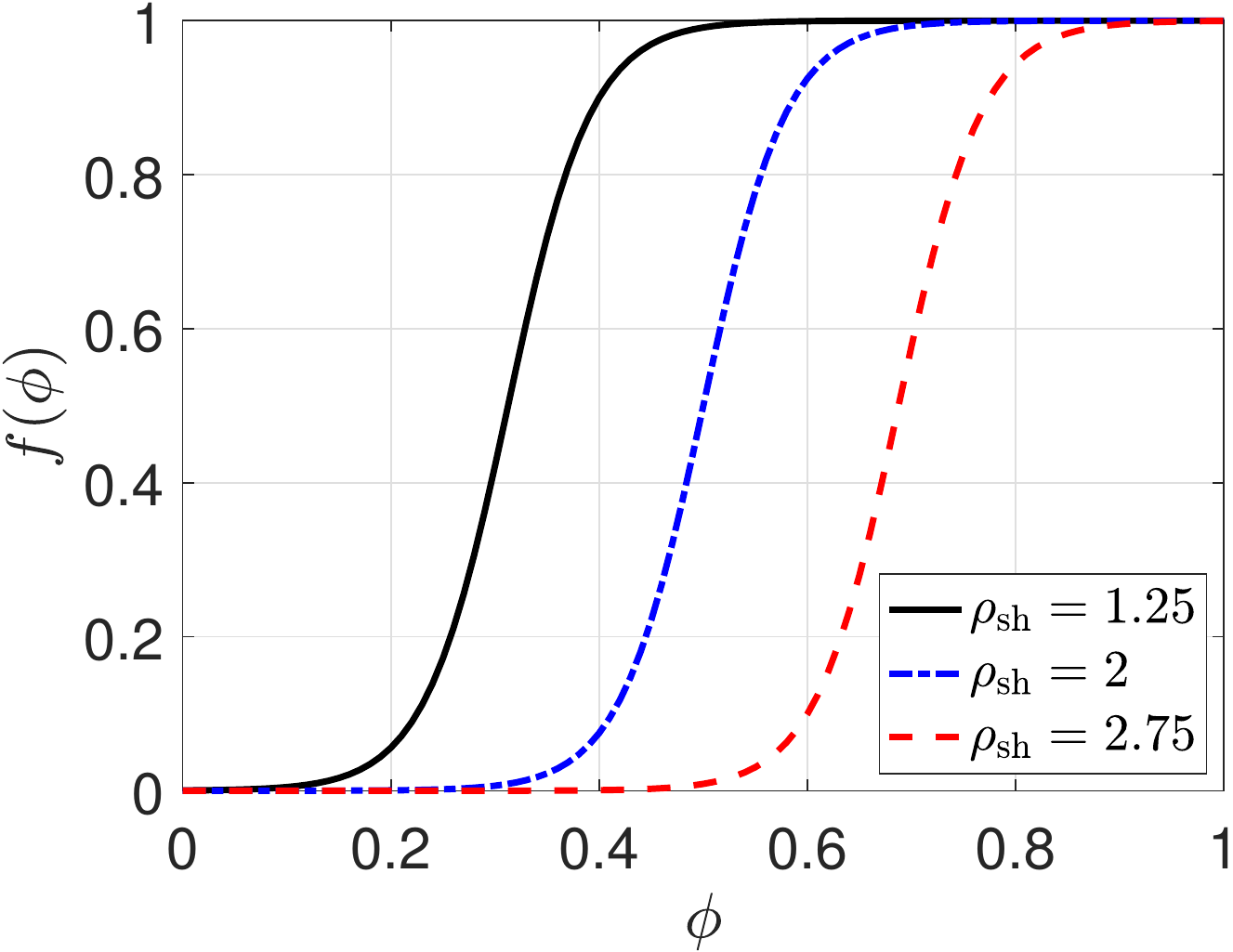}
	\caption{Modeling of mixtures: Characteristics of Eq.~\eqref{e:mFun}.}
	\label{fig:fphi}
\end{figure}
At any point $\bx\in\sS$ there will thus be a mixture of two phases. In this work, we model the behavior of the mixture by proposing the following mixture rule
\eqb{lll}
K(\phi) \is K_1\,f(\phi) + K_0\,(1-f(\phi))\,, \\[1mm]
G(\phi) \is G_1\,f(\phi) + G_0\,(1-f(\phi))\,, \\[1mm]
c(\phi) \is c_1\,f(\phi) + c_0\,(1-f(\phi))\,, \\[1mm]
\eta(\phi) \is \eta_1\,f(\phi) + \eta_0\,(1-f(\phi))\,,
\label{e:mix}
\eqe
with the interpolation function
\eqb{lll}
f(\phi) = \ds \frac{1}{2}\Big(1+\text{tanh}\left(-\rho_\mathrm{sh}\,\pi +4\,\pi\,\phi \right) \Big)\,.
\label{e:mFun}
\eqe
Here, $K_i$, $G_i$, $c_i$ and $\eta_i$ are the material parameters corresponding to $f(\phi) = i$, $i=0,1$.
The constant $\rho_\mathrm{sh}\in\bbR$ prescribes whether a smaller or a larger portion of the phase interface is characterized by material properties corresponding to $\phi=1$.
The function $f(\phi)$ is shown in Fig.~\ref{fig:fphi} for different $\rho_\mathrm{sh}$. 
In the subsequent examples, $\rho_\mathrm{sh}=1.25$ is chosen in order to increase the influence of phase $\phi=1$, which is the softer phase in the examples. 

\section{Weak form}\label{s:wf}

This section presents the weak form of PDEs~\eqref{e:sf} and \eqref{e:sfm}.
Combining \eqref{e:sf} with \eqref{e:ja}, \eqref{e:jad} and \eqref{e:mubi2} and combining \eqref{e:sfm2} with \eqref{e:amom}, \eqref{e:Mab} and \eqref{e:bab},
shows that both are fourth-order PDEs \citep{sahu17}.
Hence, the surface divergence theorem is applied twice in order to obtain second-order weak forms.

\subsection{Weak form for the Kirchhoff-Love thin shell equation}

For Kirchhoff-Love shells the weak form is given by
\eqb{l}
G_\mathrm{in} + G_\mathrm{int} - G_\mathrm{ext} = 0 \quad\forall\,\delta\bx\in\sU\,,
\label{e:wfu}\eqe
with
\eqb{lll}
G_\mathrm{in} 
\dis \ds\int_{\sS} \delta\bx\cdot\rho\,\dot\bv\,\dif a\,, \\[4mm]
G_\mathrm{int} 
\dis \ds\int_{\sS} \frac{1}{2}\,\delta a_{\alpha\beta} \, \sig^{\alpha\beta} \, \dif a  + \int_{\sS} \delta b_{\alpha\beta} \, M^{\alpha\beta} \, \dif a\,, \\[4mm]
G_\mathrm{ext} 
\dis \ds\int_{\sS}\delta\bx\cdot\bff\,\dif a + \int_{\partial_t\sS}\delta\bx\cdot\bT\,\dif s + \int_{\partial_m\sS}\delta\bn\cdot\bM\,\dif s
\label{e:Giie}\eqe
\citep{shelltheo2}.
Here,  $\mathcal{U} = \left\lbrace \delta\bx\in \mathcal{H}^2\big(\sS(\bx,t)^3\big)|\,\delta\bx=0\,\text{on}\,\partial_x\sS\, ,\,\delta\bn=0\,\text{on}\,\partial_n\sS\right\rbrace$ is the space of suitable surface variations, where $\mathcal{H}^2$ is the Sobolev space with square integrable first and second derivatives and $\partial_x\sS$ and $\partial_n\sS$ are the Dirichlet boundaries for displacements and rotations. Further, $\bff = f^\alpha\,\ba_\alpha + p\,\bn$, $\bT = \bsig^\mrT\bnu$ and $\bM = \bmu^\mrT\bnu$ denote prescribed surface forces, edge tractions and edge moments.
The latter act on the boundaries $\partial_t\sS$ and $\partial_m\sS$, respectively. 
For closed surfaces (as used in the examples of Sec.~\ref{s:num_ex}), $\partial_t\sS=\partial_m\sS=\emptyset$.
If desired, $\dif a = J\,\dif A$ 
can be used to map integrals to the reference surface $\sS_0$.

\subsection{Weak form for the Cahn-Hilliard surface equation}

Multiplying field equation \eqref{e:sf} with the test function $\delta\phi$, and applying the surface divergence theorem
\eqb{l}
\ds\int_\sS j^\alpha_{;\alpha}\,\dif a = \ds\int_{\partial\sS}j^\alpha\,\nu_\alpha\,\dif s\,,
\eqe 
gives
\eqb{l}
\ds\int_\sS\delta\phi\,\rho\,\dot\phi\,\dif a - \int_\sS\delta\phi_{;\alpha}\,j^\alpha\,\dif a + \int_{\partial_j\sS}\delta\phi\,\bar\bj\cdot\bnu\,\dif s = 0\quad\forall\,\delta\phi\in\sV\,,
\eqe
where $\bar\bj$ is the prescribed flux on boundary $\partial_j\sS$ with outward unit normal $\bnu=\nu_\alpha\ba^\alpha$.
Here we have assumed that $\delta\phi = 0$ on $\partial\sS\backslash\partial_j\sS$. 
Further, $\mathcal{V} = \left\lbrace \delta\phi\in \mathcal{H}^2\big(\sS(\phi,t)\big)\right|\,\delta\phi=0\,\text{on}\, \partial\sS\backslash\partial_j\sS\rbrace$ is the space of suitable test functions. 
\\
According to \eqref{e:ja} and \eqref{e:jad} the flux $j^\alpha$ has three contributions.
The last of those, the interfacial flux $j^\alpha_\mri$, contains three derivatives, and so we again apply the surface divergence theorem to this term to reduce it to second order.
In order to avoid handling complex expressions for terms arising from $\nablas \mu_\mathrm{el}$ that later need to be linearized,\footnote{Since $a_{\alpha\beta;\gamma}=0$, we can write $\nablas \mu_\mathrm{el}=\mu'_\mathrm{el}\,\nablas\phi+\partial \mu_\mathrm{el}/\partial b_{\alpha\beta}\,b_{\alpha\beta;\gamma}\ba^\gamma$.} we will apply the surface divergence theorem also to $j^\alpha_\mathrm{el}$.
Doing so, we obtain,
\eqb{l}
\ds\int_\sS\delta\phi_{;\alpha}\,(j^\alpha_\mri+j^\alpha_\mathrm{el})\,\dif a = 
\int_{\sS_0}\big(\nablas\delta\phi\cdot\nablas M + \Delta_\mrs\delta\phi\,M\big)\big(\mu_\mri+\mu_\mathrm{el}\big)\,\dif A
+\int_{\partial_r\sS}\nablas\delta\phi\cdot\bnu\,\bar r\,\dif s\,,
\eqe
where $\bar r$ is the prescribed boundary value for the quantity $r := M(\lambda\,\omega N \Delta_\mrs\phi - \mu_\mathrm{el}/J)$. 
Here we have assumed that $\nablas\delta\phi = \mathbf{0}$ on $\partial\sS\backslash\partial_r\sS$, and transformed integrals using $\dif a = J\,\dif A$. 
Writing $\nablas M = M'\,\nablas\phi$ and $\delta\phi_{;\alpha}\,j^\alpha_{\phi}=-\nablas\delta\phi\cdot\nablas\phi\,M\mu'_{\phi}/J$ then leads to the weak form
\eqb{l}
\bar G_\mathrm{in} + \bar G_\mathrm{int} - \bar G_\mathrm{ext} = 0 \quad\forall\,\delta\phi\in\sV\,,
\label{e:CHwf}
\eqe 
with
\eqb{lll}
\bar G_\mathrm{in} \dis \ds\int_\sS\delta\phi\,\rho\,\dot\phi\,\dif a\,,  \\[4mm]
\bar G_\mathrm{int} \dis \ds \int_{\sS_0}\nablas\delta\phi\cdot\Big( M\mu'_\phi - M'\big(\mu_\mri+\mu_\mathrm{el}\big) \Big)\,\nablas\phi\,\dif A 
- \int_{\sS_0}\Delta_\mrs\delta\phi\,M\big(\mu_\mri+\mu_\mathrm{el}\big)\,\dif A\,, \\[4mm]
\bar G_\mathrm{ext} \dis -\ds\int_{\partial_j\sS}\delta\phi\,\bar\bj\cdot\bnu\,\dif s + \int_{\partial_r\sS}\nablas\delta\phi\cdot\bnu\,\bar r\,\dif s\,.
\label{e:wfCH}\eqe
For closed surfaces (as in the examples of Sec.~\ref{s:num_ex}), $\partial_j\sS=\partial_r\sS=\emptyset$ and hence $\bar G_\mathrm{ext}=0$.

\textbf{Remark:} As an alternative to using the surface divergence theorem on $\nablas \mu_\mathrm{el}$, this term can be expanded as in footnote 7. 
This results in
\eqb{l}
\bar G_\mathrm{int} := \ds \int_{\sS_0}\nablas\delta\phi\cdot\bigg[\Big( M\big(\mu'_\phi+\mu'_\mathrm{el}\big) - M'\mu_\mri\Big)\,\nablas\phi
+ M\pa{\mu_\mathrm{el}}{b_{\alpha\beta}}b_{\alpha\beta;\gamma}\,\ba^\gamma\bigg]\,\dif A
- \int_{\sS_0}\Delta_\mrs\delta\phi\,M\mu_\mri\,\dif A\,.
\label{e:wfCHa}\eqe
This avoids having $\mu_\mathrm{el}$ inside $\bar r$, but it requires dealing with $\partial \mu_\mathrm{el}/\partial b_{\alpha\beta}\,b_{\alpha\beta;\gamma}\,\ba^\gamma$ in the discretization and linearization of the coupled system.
Therefore we will use \eqref{e:wfCH} instead of \eqref{e:wfCHa}.

\subsection{Dimensionless form}
\label{sec:dimles}
The preceding equations can be normalized by defining dimensionless quantities for position, time and the Helmholtz free energy as
\eqb{l}
\bx^\star := \ds \frac{\bx}{L_0}\,, \quad 
t^\star := \ds\frac{t}{T_0}\,,\quad
\Psi^\star := \ds\frac{\Psi}{\Psi_0}\,,
\label{e:dimles}\eqe
where $L_0$, $T_0$ and $\Psi_0$ are chosen scales for length, time and energy density, respectively.
From this, the normalization of surface stress, surface moment, chemical potential, mobility, density and mass flux follow as\footnote{Considering that $\xi^\alpha$ has units of length, and so $\bA_\alpha$, $\ba_{\alpha}$, $A_{\alpha\beta}$ and $a_{\alpha\beta}$ are dimensionless.}
\eqb{l}
\sig^{\alpha\beta}_\star := \ds\frac{\sig^{\alpha\beta}}{\Psi_0}\,,\quad
M^{\alpha\beta}_\star := \ds\frac{M^{\alpha\beta}}{\Psi_0L_0}\,, \quad
\mu_\mrc^\star := \ds\frac{\mu_\mrc}{\Psi_0}\,, \quad
M^\star := \ds\frac{M}{T_0}\,, \quad
\rho^\star := \ds\frac{\rho}{\rho_0}\,, \quad 
j^\alpha_\star := \ds\frac{j^\alpha}{j_0}\,,
\eqe
where $\rho_0 :=T_0^2\Psi_0/L_0^2$ and $j_0 := \rho_0/T_0$.
Further, the normalizations of the temporal and spatial derivative operators yield
\eqb{l}
\ds\pa{...}{t^\star} := T_0\pa{...}{t}\,, \quad
\nablas^\star := L_0\nablas\,,\quad
\Delta_\mrs^\star := L_0^2\Delta_\mrs\,.
\eqe
With these definitions, the weak forms in Eq.~\eqref{e:wfu} and \eqref{e:CHwf} can be fully normalized as
\eqb{l}
G^\star := \ds\frac{G}{\Psi_0L_0^2}\,,\quad
\bar G^\star := \ds\frac{\bar G}{\Psi_0T_0}\,.
\eqe
Likewise,
\eqb{l}
\bar\Psi^\star := \ds\frac{\bar\Psi}{\Psi_0L_0^2}\,,\quad
\bar\Psi := \ds\int_{\sS_0} \Psi\, \dif A
\label{e:enI}\eqe
is the normalization of the total energy in the system.
In the following, we will only work with the dimensionless form of all equations and will omit the superscript $\star$ for notational simplicity.
In the examples we use $N \omega = \Psi_0$ and $N k_\mrB T = \Psi_0/3$. 
\section{Discretization of the coupled system}\label{sec:disc}

This section presents the discretization of the governing equations in the framework of isogeometric finite elements.
Due to the smoothness of spline basis functions, we can directly discretize the fields within the two second-order weak forms that describe the coupled problem.
That is, we do not need to employ rotational degrees-of-freedom (dofs) or resort to mixed formulations.
The spatial discretization used here is based on the unstructured spline construction presented in Sec.~\ref{sec:spaceDisc_splines}, which is at least $C^1$-continuous at all points $\bx(t)$ for all $t$.
This is then used in Secs.~\ref{sec:spaceDisc_mech} and \ref{sec:spaceDisc_phase} to discretize the two governing weak forms.
In Sec.~\ref{sec:timeDisc} we discuss their temporal discretization using an adaptive time-stepping scheme.

\subsection{Unstructured spline spaces}\label{sec:spaceDisc_splines}

The numerical examples presented in this work utilize structured NURBS meshes and unstructured quadrilateral meshes to describe the surface geometry. 
The construction of unstructured spline spaces for the latter is based on the approach of \citet{toshniwal2017smooth}. 
This approach is advantageous since it allows the description of surfaces that are point-wise $C^1$-continuous even during deformation.
The approach is briefly summarized here.

The tasks of geometric modeling and computational analysis place differing requirements on the spaces of spline functions to be used. 
Acknowledging these differences, \cite{toshniwal2017smooth} built separate spline spaces for these tasks, $\splSpace{D}$ and $\splSpace{A}$, respectively.
The following sections give a conceptual overview of the construction and properties of the spline basis functions spanning $\splSpace{D}$ and $\splSpace{A}$.

\subsubsection{Construction of spline spaces}

The construction of spline spaces is explained using the concept of \emph{extraction operators} \citep{borden2011isogeometric}. 
Those allow to write the IGA formulation in classical FE notation.
We explain this concept with piecewise-polynomial (abbreviated as p-w-p) splines in mind (B-/T-/LR-/HBS-splines, for instance). 
The p-w-p splines restricted to any parametric element of the mesh $\Omega_e$ are tensor-product polynomials\footnote{For rational polynomial splines, simply consider homogeneous coordinates.}. 
Then, the extraction operator is the map from the local tensor-product polynomial basis, typically chosen as the Bernstein polynomial basis, to the element-local polynomial representation of the spline basis.

\begin{figure}[h]
	\centering
	\includegraphics[width=0.28\textwidth]{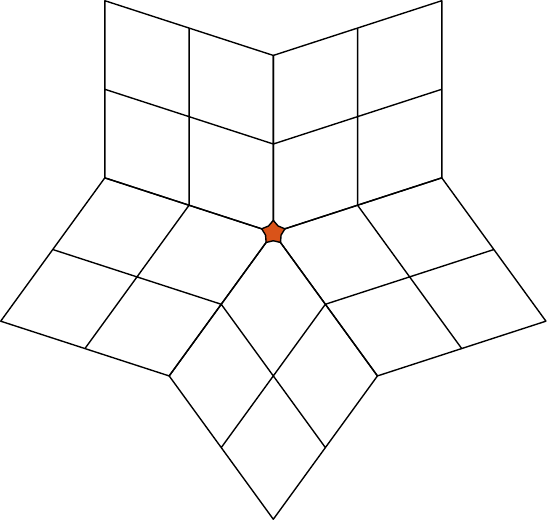}\hspace{0.5cm}
	\includegraphics[width=0.28\textwidth]{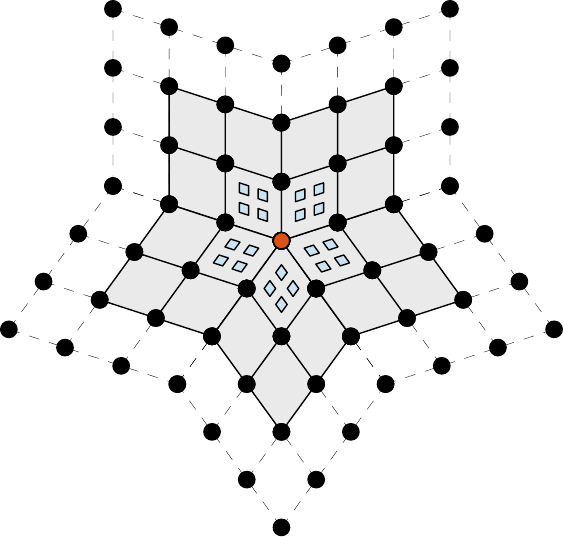}\hspace{0.5cm}
	\includegraphics[width=0.28\textwidth]{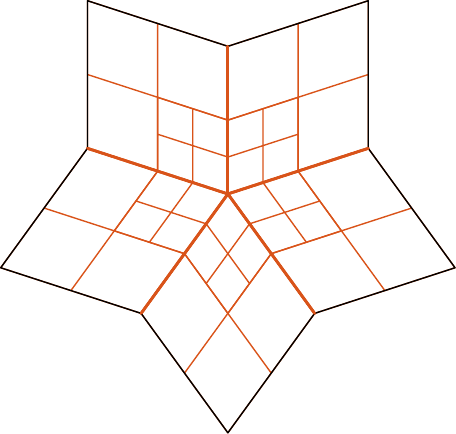}
	\caption{The above figures show a neighborhood of an extraordinary point. The left side displays the 2-ring elements that surround the extraordinary point, while the middle shows the local dof structure around the extraordinary point for $\splSpace{A}$. Instead of all mesh vertices being identified as dofs, some mesh vertices are labelled \emph{inactive} (red disk), and in their place dofs on the adjacent faces are introduced (blue squares); then, the face-based dofs and the mesh vertices not labelled inactive collectively form the full set of \emph{active} dofs (black disks and blue squares) -- the ones used for performing computations. Each dof is associated to a $C^k$ smooth spline function, $k \geq 1$. The rightmost figure elucidates the smoothness of an arbitrary spline in the space spanned by these functions -- smoothness across the red (resp. black) edges is $C^1$ (resp. $C^2$), while it is $C^\infty$ in the white space as the splines are simply polynomials in the element interiors.}
	\label{fig:splineDOFConfig}
\end{figure}

Keeping the above basics in mind, the constructions of splines reduces to defining suitable extraction operators for each element. We do this for the basis functions spanning $\splSpace{D}$ and $\splSpace{A}$, $\splBasisMod_i$ and $\splBasis_i$, respectively, in the following manner:
\begin{enumerate}[label=(\alph*)]
	\item \underline{Initial, macro extractions}:
	First, the p-w-p forms of $\splBasisMod_i$ and $\splBasis_i$ on the underlying quadrilateral mesh are initialized. 
	This amounts to initialization of the extraction operators for these splines on each element; see Appendix \ref{App:ExtractionInitialization} for details.
	\item \underline{Smoothed, micro extraction}:
	After initialization, $\splBasisMod_i$ and $\splBasis_i$ are only $C^0$-smooth on the elements containing extraordinary points (EPs), i.e., vertices where the number of edges that meet is not equal to $4$, like the central vertex in Fig.~\ref{fig:splineDOFConfig} left. 
	Then, the splines are smoothed by (a) splitting their p-w-p forms on the elements containing EPs \citep{nguyen2016refinable} using the de Casteljau algorithm \citep{piegl2012nurbs}, and then (b) by a smoothing of the p-w-p forms using a \emph{smoothing matrix} and the theory of \emph{D-patches} \citep{reif1997refineable}.
\end{enumerate}

One of the salient features of the above construction is that each step is carried out while ensuring satisfaction of \emph{isogeometric compatibility}, $\splSpace{D} \subset \splSpace{A}$. 
This is a sufficient condition for allowing exact representation of geometries built using $\splSpace{D}$ as members of $\splSpace{A}$. In other words, at each step of the construction, we ensure that the following holds,
\begin{equation}
\begin{bmatrix}
\splBasisMod_1\\
\vdots\\
\splBasisMod_{n_D}
\end{bmatrix}
=
\mbf{C}^{A\rightarrow D}
\begin{bmatrix}
\splBasis_1\\
\vdots\\
\splBasis_{n_A}
\end{bmatrix},
\end{equation}
where $\mbf{C}^{A\rightarrow D}$ is an explicitly computable matrix, and $n_A$ and $n_D$ are the numbers of control points (or nodes) used for the analysis and design, respectively. Then, isogeometric compatibility follows trivially,
\begin{equation}
	\splSpace{D} := \text{span}(\splBasisMod_1,\dots,\splBasisMod_{n_D}) \subset \text{span}(\splBasis_1,\dots,\splBasis_{n_A}) =: \splSpace{A}\,.
\end{equation}

Initial geometries at time $t = 0$, $\sS(0)$, are built using $\splSpace{D}$ and, because of isogeometric compatibility, we can express them exactly as members of $\splSpace{A}$. 
In the subsequent analysis only $\splSpace{A}$ is needed. 
Therefore, 
we restrict the remaining discussion to
the usage of $\splSpace{A}$ and omit index $A$ to simplify notation, i.e.
\begin{equation}
	\splBasisWOA_i := \splBasis_i\,,\qquad\splSpace{} := \splSpace{A}\,.	
\end{equation}
The dof structure corresponding to $\splSpace{}$ in the vicinity of extraordinary points is shown in the middle of Fig.~\ref{fig:splineDOFConfig}. 
The smoothness of an arbitrary spline $s \in \splSpace{}$ is illustrated on the right of Fig.~\ref{fig:splineDOFConfig}. 
As shown, the extraordinary point's neighborhood contains edges across which the smoothness is only $C^1$ (depicted in red in the figure). Also note that this zone of $C^1$-continuity is limited to the 2-ring elements of each extraordinary point (at the coarsest level of refinement), and outside of this zone the splines are maximally smooth, i.e., $C^2$-continuous. 

\subsubsection{Properties of $\splSpace{}$}

The spline space $\splSpace{}$ is built exclusively from bi-cubic polynomial pieces, and is identical to the space of bi-cubic analysis-suitable T-splines (or, AST-splines) \citep{scott2013isogeometric,li2015some} in the regular (locally structured) regions of the mesh. 
In particular, the basis functions spanning $\splSpace{}$ form a convex partition of unity and are locally supported. 
Additionally, the space $\splSpace{}$ was observed to possess good approximation properties as evidenced by the suite of numerical tests presented in \citet{toshniwal2017smooth}, and makes numerical investigation of high-order problems on arbitrary surfaces possible.

\subsubsection{Spatial discretization of primary fields}
In this section, finite dimensional approximations to all primary fields of interest (surface geometry and phase field order parameter) will be expressed as members of $\splSpace{}$. Let $n_e$ spline basis functions, with global indices $i_1,\dots,i_{n_e}$, be supported on parametric element $\Omega_e$. Then, we can express the local element representations of the surface $\sS$, $\sS_0$ and phase field $\phi$ as,
\eqb{l}
	\bx^h = \mN\,\mx_e\,,\qquad\bX^h = \mN\,\mX_e\,,
\label{e:xXh}
\eqe 
and,
\eqb{l}
	\phi^h = \bar\mN\,\bphi_e\,,
\label{e:phih}
\eqe
respectively, where,
\eqb{l}
\mN := [N_{i_1}\bone,\,N_{i_2}\bone,\,...,\,N_{i_{n_e}}\bone]\,,
\eqe
\eqb{l}
	\bar\mN := [N_{i_1},\,N_{i_2},\,...,\,N_{i_{n_e}}]\,.
\eqe
Here, $\bone$ denotes the $(3\times3)$ identity matrix, and $\mX_e$, $\mx_e$ and $\bphi_e$ denote element-level vectors containing the positions and dofs at nodes $i_1,\dots,i_{n_e}$. These local vectors can be extracted from the global vectors $\mX$, $\mx$ and $\bphi$ that contain all nodal positions and dofs. The respective variations are defined analogously, given by
\eqb{l}
	\delta\bx^h = \mN\,\delta\mx_e\,,\qquad\delta\bX^h = \mN\,\delta\mX_e\,,
\label{e:deltaxXh}
\eqe
and,
\eqb{l}
\delta\phi^h = \bar\mN\,\delta\bphi_e\,.
\eqe 
Using the above equations, the weak forms for the surface deformation and the phase field are discretized as described in Secs.~\ref{sec:spaceDisc_mech} and \ref{sec:spaceDisc_phase}, respectively.

\subsection{Spatial discretization of the mechanical weak form}\label{sec:spaceDisc_mech}

Using Eqs.~\eqref{e:xXh} and \eqref{e:phih}, the tangent vectors on the surface are discretized as
\eqb{lll}
\ba_\alpha^h \is \mN_{\!,\alpha}\,\mx_e\,,~\mathrm{and}\quad \bA_\alpha^h = \mN_{\!,\alpha}\,\mX_e\,,
\label{e:bah}\eqe
where $\mN_{\!,\alpha}:=\partial\mN/\partial\xi^\alpha$. 
The discretized tangent vectors of \eqref{e:bah} lead to the discretized normal vectors $\bn^h$ and $\bN^h$ following Eq.~\eqref{e:bn}.\footnote{To avoid confusion we write discrete arrays, such as the shape function array $\mN$, in roman font, whereas continuous tensors, such as the normal vector $\bN$, are written in italic font.} 
The metric tensor and curvature components can then be expressed as
\eqb{lll}
a_{\alpha\beta}^h \is \mx_e^\mrT\mN^\mrT_{\!,\alpha}\,\mN_{\!,\beta}\,\mx_e\,,~\mathrm{and}\quad b_{\alpha\beta}^h = \bn^h\cdot\mN_{\!,\alpha\beta}\,\mx_e\,,
\eqe 
and similarly
\eqb{lll}
A_{\alpha\beta}^h \is \mX_e^\mrT\mN^\mrT_{\!,\alpha}\,\mN_{\!,\beta}\,\mX_e\,, ~\mathrm{and}\quad
B_{\alpha\beta}^h = \bN^h\cdot\mN_{\!,\alpha\beta}\,\mX_e\,.
\eqe 
The contravariant metrics $[a^{\alpha\beta}_h]=[a_{\alpha\beta}^h]^{-1}$ and $[A^{\alpha\beta}_h]=[A_{\alpha\beta}^h]^{-1}$ then follow.
Using Eq.\;\eqref{e:deltaxXh}, the variations of the surface metric and curvature can be obtained as
\eqb{lll}
\delta a_{\alpha\beta}^h \is \delta\mx_e^\mrT\big(\mN^\mrT_{\!,\alpha}\,\mN_{\!,\beta}+\mN^\mrT_{\!,\beta}\,\mN_{\!,\alpha}\big)\,\mx_e\,, ~\mathrm{and}\quad
\delta b_{\alpha\beta}^h =\delta\mx_e^\mrT\,\mN^\mrT_{\!;\alpha\beta}\,\bn^h\,,
\eqe 
with
\eqb{l}
\mN_{\!;\alpha\beta} := \mN_{\!,\alpha\beta} - \Gamma^\gamma_{\alpha\beta}\,\mN_{\!,\gamma}\,.
\eqe
Here,
\eqb{l}
\Gamma^\gamma_{\alpha\beta} = \mx_e^\mrT\,\mN^\mrT_{,\alpha\beta}\,a^{\gamma\delta}_h\,\mN_{,\delta}\,\mx_e
\eqe
denotes the discretized Christoffel symbols.
Using the above expressions, the discretized mechanical weak form becomes
\eqb{l}
\delta\mx^\mrT\,\big[\mf_\mathrm{in} + \mf_\mathrm{int} - \mf_\mathrm{ext}\big] = 0 \quad \forall~\delta\mx\in\sU^h\,,
\label{e:fx}
\eqe
where the global force vectors $\mf_\mathrm{in}$, $\mf_\mathrm{int}$ and $\mf_\mathrm{ext}$ are assembled from their respective elemental contributions
\eqb{lll}
\mf^e_\mathrm{in} \dis \mm_e\,\ddot\mx_e~,\quad\mm_e := \ds\int_{\Omega^e}\rho\,\mN^\mrT\mN\,\dif a\,,  \\[4mm]
\mf^e_\mathrm{int} \dis \ds\int_{\Omega^e}\sig^{\alpha\beta}\,\mN_{\!,\alpha}^\mrT\,\ba^h_\beta\,\dif a
+ \int_{\Omega^e}M^{\alpha\beta}\,\mN^\mrT_{\!;\alpha\beta}\,\bn^h\,\dif a\,, \\[4mm] 
\mf^e_\mathrm{ext} \dis \ds\int_{\Omega^e}\mN^\mrT\,p\,\bn^h\,\dif a + \ds\int_{\Omega^e}\mN^\mrT\,f^\alpha\,\ba^h_\alpha\,\dif a\,.
\label{e:finext}\eqe
Further, $\delta\mx$ denotes the global vector of all nodal variations, and $\mathcal{U}^h := \mathcal{U}\,\cap\, \splSpace{}$ is its corresponding discrete space. The expression of $\mf^e_\mathrm{ext}$, corresponds to the case that there are no boundary loads $\bT$ and $\bM$ acting on $\partial\sS$.
This is the case in all the subsequent examples. 
The extension to boundary loads can be found in \citet{solidshell2}.
$\mf_\mathrm{in}$, $\mf_\mathrm{int}$ and $\mf_\mathrm{ext}$ depend on $\mx(t)$, while $\mf_\mathrm{int}$ also depends on $\bphi(t)$ through the material properties in $\sig^{\alpha\beta}$ and $M^{\alpha\beta}$, and the Korteweg stresses $\sig_\mathrm{CH}^{\alpha\beta}$.
The resulting equations at the free nodes (after application of Dirichlet boundary conditions) can thus be written as
\eqb{l}
\mf(\mx,\bphi) = \mM\,\ddot\mx + \mf_\mathrm{int}(\mx,\bphi) - \mf_\mathrm{ext}(\mx) = \mathbf{0}~,
\label{e:ODEx}
\eqe
where $\mM$ denotes the global mass matrix assembled from $\mm_e$, and $\mx$ and $\bphi$ denote the global vectors of the unknown nodal positions and unknown nodal phase parameters. 

\subsection{Spatial discretization of the phase field equations}\label{sec:spaceDisc_phase}

Using Eq.~\eqref{e:diff} and Eq.~\eqref{e:phih}, we can write
\eqb{lll}
\phi^h_{;\alpha} \is \bar\mN_{\!,\alpha}\,\bphi_e\,, \\[1mm]
\nablas\phi^h \is \ba^\alpha_h\,\bar\mN_{\!,\alpha}\,\bphi_e\,, \\[1mm]
\nablas\delta\phi^h \is \ba^\alpha_h\,\bar\mN_{\!,\alpha}\,\delta\bphi_e\,, \\[1mm]
\Delta_\mrs\phi^h \is \Delta_\mrs\bar\mN\,\bphi_e\,, \\[1mm]
\Delta_\mrs\delta\phi^h \is \Delta_\mrs\bar\mN\,\delta\bphi_e\,,
\label{e:approx_diff}\eqe
where $\bar\mN_{\!,\alpha}:=\partial\bar\mN/\partial\xi^\alpha$
and
\eqb{l}
\Delta_\mrs\bar\mN := a^{\alpha\beta}_h\,\bar\mN_{\!;\alpha\beta}\,, 
\eqe
with
\eqb{l}
\bar\mN_{\!;\alpha\beta} = \bar\mN_{\!,\alpha\beta} - \Gamma^\gamma_{\alpha\beta}\,\bar\mN_{\!,\gamma}
\label{e:bN;ab}\eqe
according to Eq.~\eqref{e:phiab}. 
Here, $\ba^\alpha_h = a^{\alpha\beta}_h\,\ba^h_\beta$ and $\Gamma^\gamma_{\alpha\beta}$ follows from the surface discretization discussed in Sec.~\ref{sec:spaceDisc_mech}.
The discretized weak form of the Cahn-Hilliard Eq.~\eqref{e:CHwf} then becomes
\eqb{l}
\delta\bphi^\mrT\,\big[\bar\mf_\mathrm{in} + \bar\mf_\mathrm{int} - \bar\mf_\mathrm{ext}\big] = 0 \quad\forall~\delta\bphi\in\sV^h\,,
\eqe 
where the global vectors $\bar\mf_\mathrm{in}$, $\bar\mf_\mathrm{int}$ and $\bar\mf_\mathrm{ext}$ are assembled from their respective elemental contributions
\eqb{lll}
\bar\mf^e_\mathrm{in} \dis \bar\mm_e\,\dot\bphi_e\,,\quad\bar\mm_e := \ds\int_{\Omega^e}\rho\,\bar\mN^\mrT\bar\mN\,\dif a\,, \\[4mm]
\bar\mf^e_\mathrm{int} \dis \bar\mk_e\,\bphi_e-\bar\mf^e_\mathrm{el}\,,\quad\!
\bar\mk_e :=\! \ds\int_{\Omega^e_0}\!\Big[\bar\mN_{\!,\alpha}^\mrT\,a^{\alpha\beta}_h\Big(\!M \mu'_\phi - M'\big(\mu_\mri+\mu_\mathrm{el}\big)\!\Big)\,\bar\mN_{\!,\beta}+\Delta_\mrs\bar\mN^\mrT J\lambda M\,\Delta_\mrs\bar\mN\Big]\,\dif A\,, \\[3mm]
\bar\mf^e_\mathrm{el} \dis \ds\int_{\Omega^e_0}\Delta_\mrs\bar\mN^\mrT M\mu_\mathrm{el}\,\dif A\,, \\[4mm]
\bar\mf^e_\mathrm{ext} \dis -\ds\int_{\Gamma_j^e}\bar\mN^\mrT\,\bar\bj\cdot\bnu\,\dif s + \int_{\Gamma_r^e}\bar\mN_{\!,\alpha}^\mrT\,\nu^\alpha\,\bar r\,\dif s\,.
\label{e:fphi_sys}
\eqe
Further, $\delta\bphi$ denotes the variation of global vector $\bphi$, and $\mathcal{V}^h := \mathcal{V}\,\cap\, \splSpace{}$ is its corresponding discrete space. Note that the expressions in \eqref{e:fphi_sys} depend on $\bphi_e$ through $M$, $\mu_\phi$, $\mu_\mri$ and $\mu_\mathrm{el}$. 
They further depend on the geometry $\mx_e$ through  $a^{\alpha\beta}_h$, $\mu_\mri$, $\mu_\mathrm{el}$, $\Delta_\mrs\bar\mN$, $J$ 
and the boundary quantities $\bnu$ and $\dif s$.
The resulting dynamical equations at the free nodes (after application of Dirichlet boundary conditions) can thus be written as
\eqb{l}
\bar\mf(\mx,\bphi) = \bar\mM\,\dot\bphi 
+ \bar\mf_\mathrm{int}(\mx,\bphi) - \bar\mf_\mathrm{ext}(\mx) = \mathbf{0}\,,
\label{e:ODEphi}
\eqe
where $\bar\mM$ denotes the global mass matrix assembled from $\bar\mm_e$. 
This, in conjunction with Eq.~\eqref{e:ODEx}, completes the semi-discrete formulation, and we discuss the temporal discretization next.
The spatially discretized equations of the coupled problem are summarized in Table~\ref{t:gov}.

\begin{table}[h]
\begin{tcolorbox}[colback=white,boxrule=0.15mm]
Governing mechanical ODE (2. order)
\esb{l}
\mf(\mx,\bphi) = \mM\,\ddot\mx + \mf_\mathrm{int}(\mx,\bphi) - \mf_\mathrm{ext}(\mx) = \mathbf{0}~.
\ese
Governing phase field ODE (1. order)
\esb{l}
\bar\mf(\mx,\bphi) = \bar\mM\,\dot\bphi 
+ \bar\mf_\mathrm{int}(\mx,\bphi) = \mathbf{0}\,.
\ese
These are assembled from the elemental contributions
\esb{lll}
\mm_e := \ds\int_{\Omega^e}\rho\,\mN^\mrT\mN\,\dif a\,,\quad
\bar\mm_e := \ds\int_{\Omega^e}\rho\,\bar\mN^\mrT\bar\mN\,\dif a\,, \\[4mm]
\mf^e_\mathrm{int} := \ds\int_{\Omega^e}\sig^{\alpha\beta}\,\mN_{\!,\alpha}^\mrT\,\ba^h_\beta\,\dif a
+ \int_{\Omega^e}M^{\alpha\beta}\,\mN^\mrT_{\!;\alpha\beta}\,\bn^h\,\dif a\,,\quad \mf^e_\mathrm{ext} := \ds\int_{\Omega^e}\mN^\mrT\,p\,\bn^h\,\dif a\,,\\[4mm]
\bar\mf^e_\mathrm{int} := \ds\int_{\Omega^e_0}\!\Big[\bar\mN_{\!,\alpha}^\mrT\,a^{\alpha\beta}_h\Big(\!M \mu'_\phi - M'\big(\mu_\mri+\mu_\mathrm{el}\big)\!\Big)\,\bar\mN_{\!,\beta}\,\bphi_e+\Delta_\mrs\bar\mN^\mrT M \big(J\lambda\Delta_\mrs\bar\mN\,\bphi_e -\mu_\mathrm{el} \big)\Big]\,\dif A\,.\\[3mm]
\ese
Here, $\sig^{\alpha\beta}$ is given by Eqs.~\eqref{e:sigtot},~\eqref{e:sigel},~\eqref{e:visc},~\eqref{e:sigCH},~\eqref{e:mix} and~\eqref{e:mFun}, $M^{\alpha\beta}$ is given by Eqs.~\eqref{e:Mab},~\eqref{e:mix} and~\eqref{e:mFun}. The degenerate mobility is given by $M= D\,\phi\,(1-\phi)$ and its derivative by $M' = D\,(1-2\phi)$. $\mu_\mri$ is given by Eq.~\eqref{e:mubi2}.4 and $\mu_\mathrm{el}$ by Eq.~\eqref{e:muel}, while $\mu'_\phi = 1/\big(3(\phi-\phi^2)\big)-2$. The mass matrices $\mM$ and $\bar{\mM}$ are constant since a mass conserving system in the Lagrangian frame is considered. All variables and integrals are understood to be dimensionless $(\star)$-quantities as introduced in Sec.~\ref{sec:dimles}. 
\end{tcolorbox}
\caption{Summary of the governing discretized equations for closed surfaces and $f^\alpha = 0$, as is used in the following examples.}
\label{t:gov}
\end{table}

\subsection{Temporal discretization of the coupled problem}\label{sec:timeDisc}

In this work, monolithic time integration based on the fully implicit generalized-$\alpha$ scheme \citep{chung93} is used. 
The resulting discrete nonlinear system of equations is solved by the Newton-Raphson iteration at each time step. 
Given the quantities $(\mx_{n},\dot\mx_{n},\ddot\mx_{n},\bphi_{n},\dot\bphi_{n})$ at time $t_n$, the new values $(\mx_{n+1},\dot\mx_{n+1},\ddot\mx_{n+1},\bphi_{n+1},\dot\bphi_{n+1})$ at time $t_{n+1}$ can be computed. The generalized-$\alpha$ method proceeds by requiring the system of equations to be satisfied at intermediate values $(\mx_{n+\alpha_\mrf}, \dot\mx_{n+\alpha_\mrf}, \ddot\mx_{n+\alpha_\mrm}, \bphi_{n+\alpha_\mrf}, \dot\bphi_{n+\alpha_\mrm})$, i.e.
\eqb{lll}
\label{e:gen_a1}
\begin{bmatrix}
\mf\left(\mx_{n+\alpha_\mrf},\dot\mx_{n+\alpha_\mrf},\ddot{\mx}_{n+\alpha_\mrm},\bphi_{n+\alpha_\mrf} \right) \\[2mm]
\bar\mf\left(\mx_{n+\alpha_\mrf},\bphi_{n+\alpha_\mrf},\dot{\bphi}_{n+\alpha_\mrm} \right) 
\end{bmatrix}
\is \mathbf{0}\,.
\eqe
The intermediate quantities, and the quantities at time step $n+1$, are evaluated as described in Appendix~\ref{sec:genA}.
The system of nonlinear equations \eqref{e:gen_a1} is solved at each time step using the iterative Newton-Raphson method, see Appendix~\ref{sec:NR} for details.
Therefore, the linearized system of equations can be expressed as
\eqb{lll}
\ds
\begin{bmatrix}
\mK_\mrx & \mK_\phi \\[2mm]
\bar\mK_\mrx & \bar\mK_\phi
\end{bmatrix}
\begin{bmatrix}
\Delta\mx_{n+1} \\[2mm]
\Delta\bphi_{n+1}
\end{bmatrix}
\is \ds -
\begin{bmatrix}
\mf\left(\mx_{n+\alpha_\mrf},\dot\mx_{n+\alpha_\mrf},\ddot{\mx}_{n+\alpha_\mrm},\bphi_{n+\alpha_\mrf} \right)\\[2mm]
\bar\mf\left(\mx_{n+\alpha_\mrf},\bphi_{n+\alpha_\mrf},\dot{\bphi}_{n+\alpha_\mrm} \right)
\end{bmatrix},
\label{e:linsys}
\eqe
where the tangent matrix blocks are computed from
\eqb{lll}
\mK_\mrx \is \ds \frac{\partial \mf}{\partial \mx_{n+1}} = \alpha_\mrf \frac{\partial \mf}{\partial \mx_{n+\alpha_\mrf}} +
\frac{\alpha_\mathrm{f}\gamma}{\beta\Delta t_{n+1}} \frac{\partial \mf}{\partial \dot{\mx}_{n+\alpha_\mathrm{f}}} + \frac{\alpha_\mrm}{\beta\Delta t^2_{n+1}} \frac{\partial \mf}{\partial \ddot\mx_{n+\alpha_\mrm}}\,,\\[4mm]
\mK_\phi \is \ds \frac{\partial \mf}{\partial \bphi_{n+1}} = \alpha_\mrf \frac{\partial \mf}{\partial \bphi_{n+\alpha_\mrf}}\,,\\[4mm]
\bar\mK_\mrx \is \ds \frac{\partial \bar\mf}{\partial \mx_{n+1}} = \alpha_\mrf \frac{\partial \bar\mf}{\partial \mx_{n+\alpha_\mrf}}\,, \\[4mm]
\bar\mK_\phi \is \ds \frac{\partial \bar\mf}{\partial \bphi_{n+1}} = \alpha_\mrf \frac{\partial \bar\mf}{\partial \bphi_{n+\alpha_\mrf}} + \frac{\alpha_\mrm}{\gamma\Delta t_{n+1}} \frac{\partial \bar\mf}{\partial \dot\bphi_{n+\alpha_\mrm}}\,.
\label{e:tan_m}
\eqe
They are assembled from the elemental contributions reported in Appendix~\ref{sec:lin}. 
Eq.~\eqref{e:linsys} is solved iteratively for $\Delta \mx_{n+1}$ and $\Delta\bphi_{n+1}$. The new surface quantities at $t_{n+1} = t_n + \Delta t_{n+1}$ are then updated as $\mx_{n+1} = \mx_n+\Delta \mx_{n+1}$ and $\bphi_{n+1} = \bphi_n+\Delta \bphi_{n+1}$.
With this, all surface quantities at $t_{n+1}$ can be evaluated as described in Secs.~\ref{Sec:def_s} -~\ref{s:wf}.
At all times, the mass matrices remain constant. 
Systems involving time dependent mass matrices were for example studied in \citet{lubich13}.

\subsection{Adaptive time-stepping}\label{Sec:ATS}

Phase transitions evolve at different time scales, which motivates an adaptive adjustment of the time step.
This section presents an adaptive time-stepping scheme for the proposed coupled system. To begin with, we note that \cite{hulbert95} present an automatic time step control algorithm for studying structural dynamics. Here, we adapt and reformulate their idea in the context of phase fields on deforming surfaces. For this purpose, the local time truncation errors of the phase field, $\be^\mrp_{n+1}$, and the surface deformation, $\be^\mrd_{n+1}$, are introduced and examined. Note that these are estimates occurring in the time step from $t_n$ to  $t_{n+1}$.
An estimate for the local time truncation error of the deformation can be expressed as \citep{hulbert95}
\eqb{lll}
\be^\mrd_{n+1} = \ds \Delta t_{n+1}^2 \bigg(c_1^\mrd \Delta \ddot{\mx}_n + c_2^\mrd \,\sum_{j=1}^n\,(-\rho_\infty)^{j-1} \Delta\ddot{\mx}_{n-j}  \bigg),
\label{eq:ede}
\eqe
where $\Delta \ddot{\mx}_n := \ddot{\mx}_{n+1}-\ddot{\mx}_n$ and $\ddot\mx_{n+1}:=\ddot\mx_{n+1}(t_{n+1})$.\footnote{Here, the Newton-Raphson iteration index is omitted for notational simplicity.} 
Here, $\ddot{\mx}_{n+1}$ is computed using Newmark's formulae (Appendix \ref{App:A}, Eq.~\eqref{eq:NM}) given the solution ${\mx}_{n+1}$ from the current Newton-Raphson iteration.  
Expressions for the constants $c_{1,2}^\mrd$ and $\rho_\infty$ can be found in Appendix \ref{App:B}, Eq.~\eqref{eq:cd}. A detailed derivation of Eq.~\eqref{eq:ede} and further information are provided in Appendix~\ref{App:B}.
The local time truncation error of the phase field can be expressed in a similar way as defined for the deformation in Eq.~\eqref{eq:ede}. In contrast to the deformation, the differential equation for the phase field is only first order in time. An estimate for the local time truncation error is
\eqb{lll}
\be^\mrp_{n+1} = \ds \Delta t_{n+1} \bigg(c_1^\mrp \Delta \dot{\bphi}_n + c_2^\mrp \,\sum_{j=1}^n\,(-\rho_\infty )^{j-1} \Delta\dot{\bphi}_{n-j}  \bigg),
\label{eq:ep4}
\eqe
where $\Delta \dot{\bphi}_n := \dot{\bphi}_{n+1}-\dot{\bphi}_n$. 
Here, $\dot{\bphi}_{n+1}$ is computed using Newmark's formulae (Appendix \ref{App:A}, Eq.~\eqref{eq:NM2}) given the solution ${\bphi}_{n+1}$ from the current Newton-Raphson iteration. 
Expressions for the constants $c_{1,2}^\mrp$ can be found in Appendix \ref{App:B2}, Eq.~\eqref{eq:cp}. A detailed derivation of Eq.~\eqref{eq:ep4} and further information are given in Appendix~\ref{App:B2}.
By using the normalized errors,
\eqb{lll}
\text{err}^\mrp = \ds \frac{\norm{\be^\mrp_{n+1}}}{\norm{\bphi_n}}\,, \quad \text{and} \quad
\text{err}^\mrd = \ds \frac{\norm{\be^\mrd_{n+1}}}{\norm{\mx_n}}\,,
\label{e:err}
\eqe
the time step is then updated according to
\eqb{lll}
\Delta t_{n+1} = \ds \rho_\mathrm{sc} \, \Delta t_n \,\text{min}\left( \sqrt{\left(\frac{\text{tol}^\mrp}{\text{err}^\mrp} \right)}, \sqrt{\left(\frac{\text{tol}^\mrd}{\text{err}^\mrd} \right)}\right).
\label{eq:ts_up}
\eqe
We found that $\text{tol}^\mrp=\text{tol}^\mrd=7.5\cdot 10^{-5}$ and $\rho_\mathrm{sc}=0.8$ are good choices for the tolerances and the safety coefficient, respectively. The time step is rejected and recomputed if either $\text{err}^\mrp > 10^{-4}$ or $\text{err}^\mrd > 10^{-4}$ in all of the following numerical examples.

\section{Numerical examples}\label{s:num_ex}

This section presents several examples in order to demonstrate the numerical behavior of the proposed model.
First, the decoupled model is verified based on existing results from literature. 
Then, coupling is investigated for deforming tori, spheres and double-tori. 
In all examples, the initial condition for the Cahn-Hilliard equation is chosen as
\eqb{l}
\phi(\bx)=\bar{\phi}+\phi_\mrr\,, 
\label{eq:ic}
\eqe
where $\bar{\phi}$ is a constant value representing the volume fraction of the mixtures, and $\phi_\mrr \in [-0.05,0.05]$ is a random perturbation. 
$\bar{\phi}=1/3$ is chosen if not otherwise stated and the density is $\rho = \rho_0$. 
In the case of deformation, the mechanical material parameters (see Sec.~\ref{s:mix}) are chosen as listed in Table~\ref{tab:mat}.
They are expressed in terms of 2D Young's modulus $E$ (force per length) and Poisson's ratio $\nu$, which are chosen as
$E = N\omega$ and $\nu = 0.3$. 

\begin{table}[h]
\begin{center}
\begin{tabular}{| p{1.8cm} | r | r |}
    \hline
    $\,$ & Pure phase state $\phi=0$ (blue color)  & Pure phase state $\phi= 1$ (red color)\\ \hline
    $\ds K_i$ & $K_0 = 1.25\,\ds \frac{E\,\nu}{(1+\nu)(1-2\,\nu)}$ & $K_1 = 0.0375\,\ds \frac{E\,\nu}{(1+\nu)(1-2\,\nu)}$ \\ \hline
    $G_i$ & $G_0 = 6.25\,\ds \frac{E}{2\,(1+\nu)}$ & $G_1 = 0.375\,\ds \frac{E}{2\,(1+\nu)}$ \\ \hline
    $c_i$ & $c_0=0.01\,E\,L_0$ & $c_1 = 0.0001875\,E\,L_0$ \\ \hline
    $\eta_i$ & $\eta_0 = 1.5\,K_0\,T_0$ & $\eta_1 = 1.5\,K_0\,T_0$ \\
    \hline
\end{tabular}
\end{center}
\caption{Material parameters for all the following numerical examples presented in this work.}
\label{tab:mat}
\end{table}

\subsection{Verification}\label{Sec:verif}

We first discuss the verification of the phase field formulation by rerunning examples from the literature. 
The verification of the shell formulation was already demonstrated in \citet{solidshell2} and is not repeated here.

\subsubsection{Phase separation on a 2D square}

\begin{figure}
\centering
\includegraphics[width=0.24\linewidth, trim = 250 70 250 120,clip]{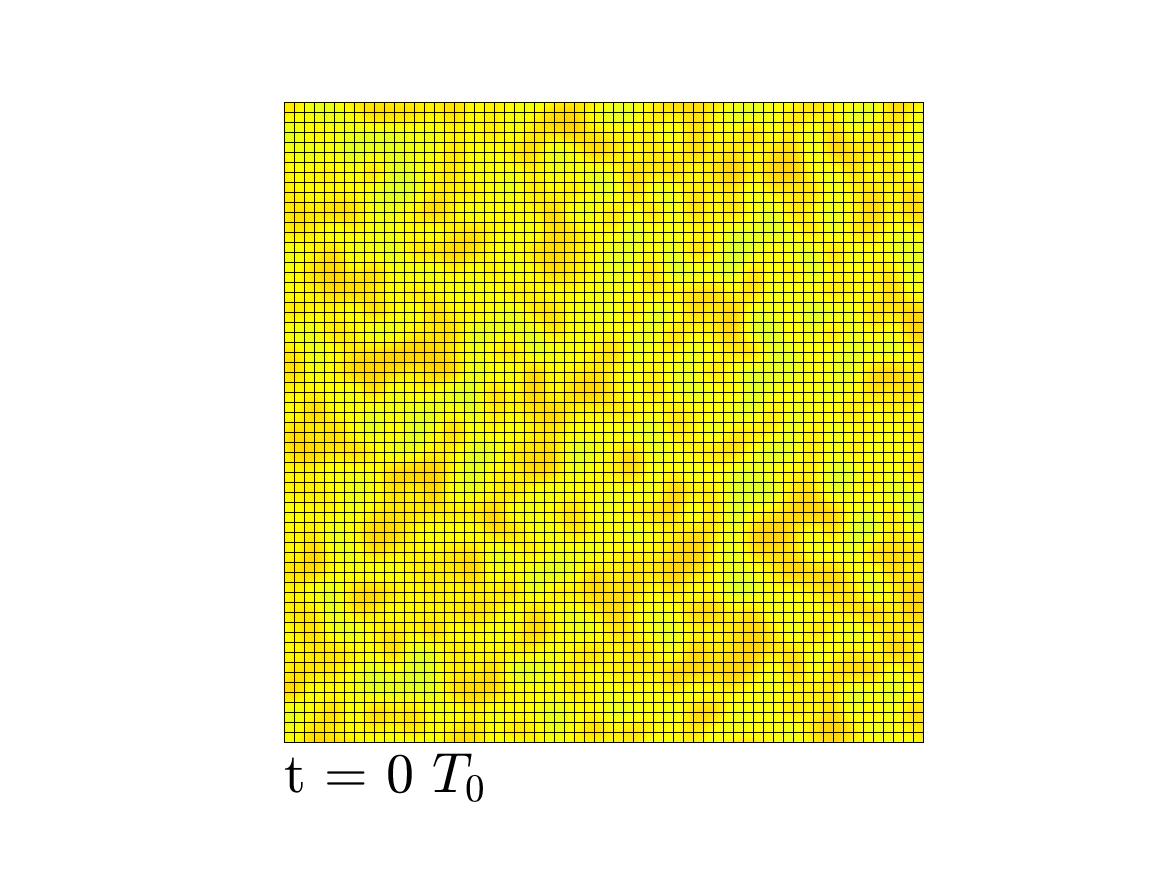}
\includegraphics[width=0.24\linewidth, trim = 250 70 250 120,clip]{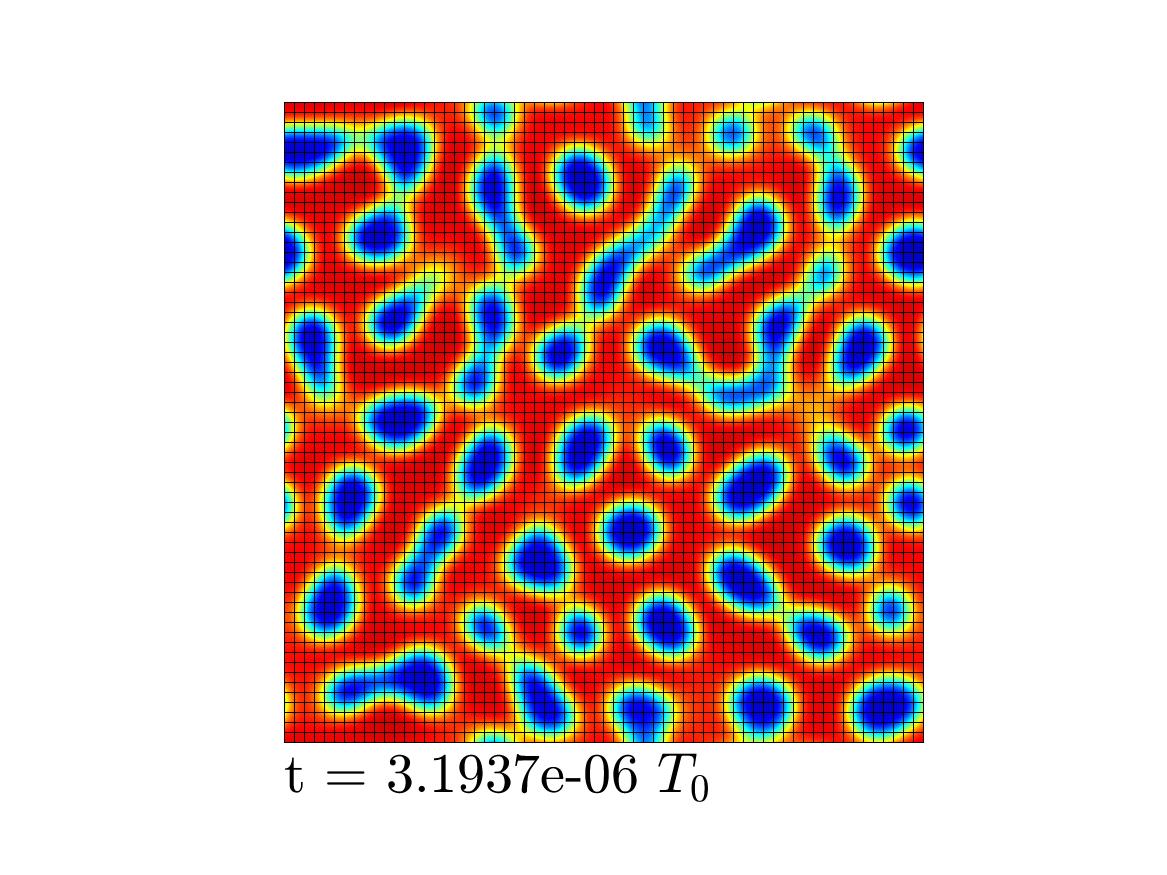}
\includegraphics[width=0.24\linewidth, trim = 250 70 250 120,clip]{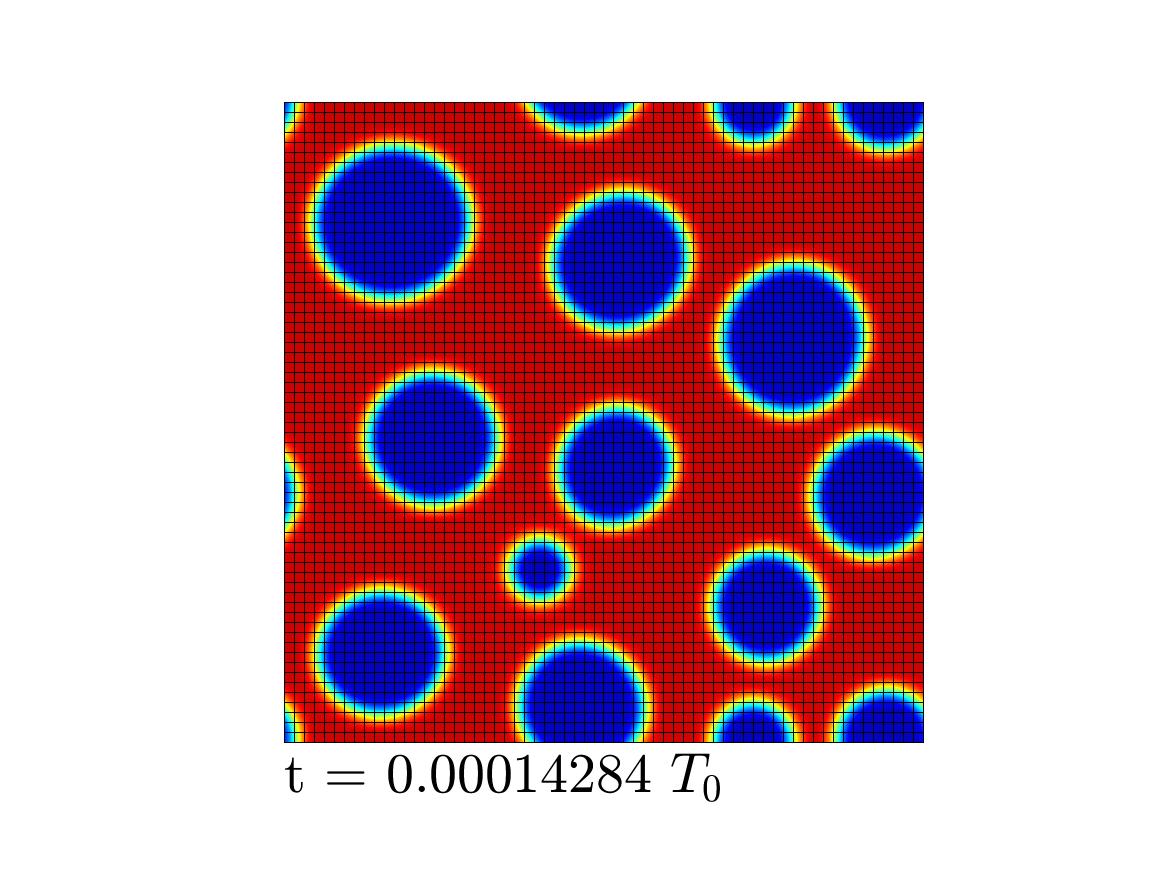}
\includegraphics[width=0.24\linewidth, trim = 250 70 250 120,clip]{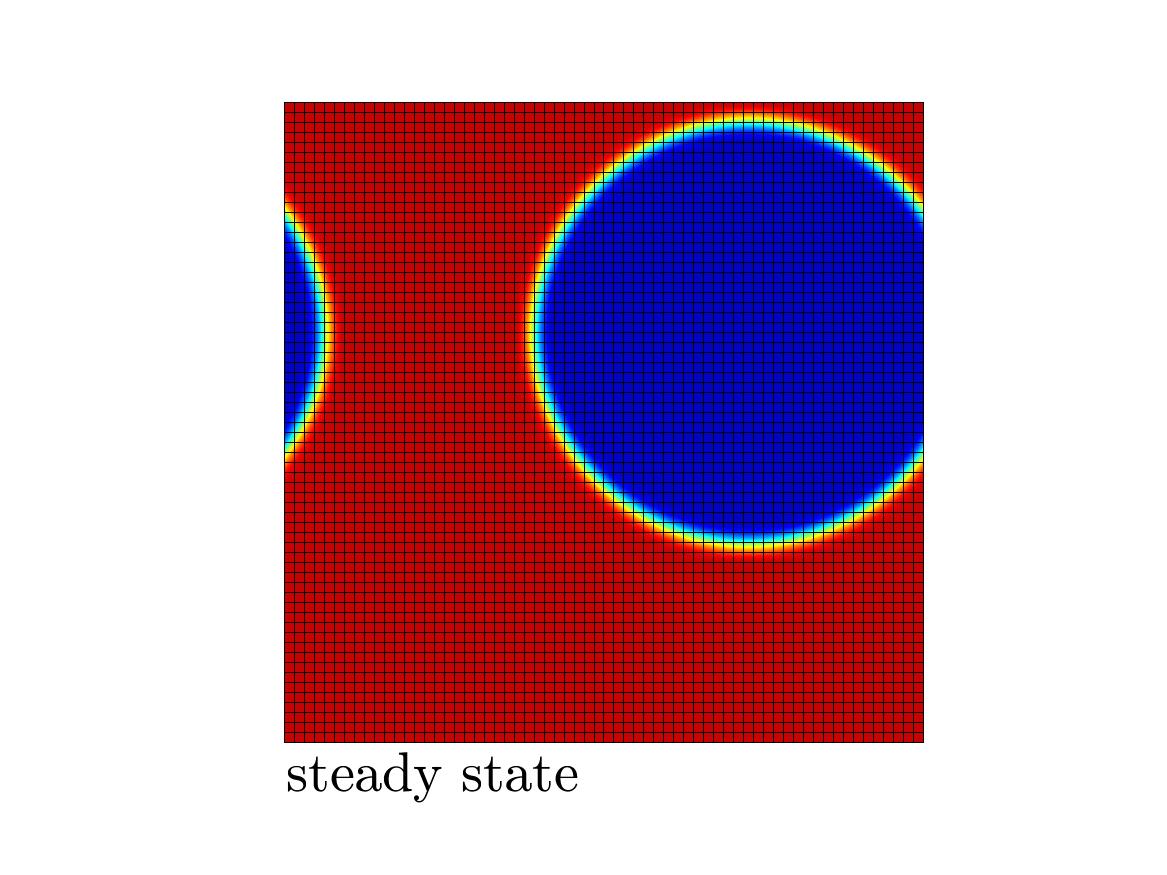}
\caption{Phase separation on a 2D square. Evolution of the phase field $\phi$ for $\lambda=1/9000\,L_0^2$ and volume fraction $\bar{\phi} = 0.63$ on a mesh containing $64 \times 64$ quadratic NURBS elements. The coloring follows Table~\ref{tab:mat}.}
\label{fig:ss_evo}
\end{figure}

\begin{figure}
\centering
\includegraphics[width=0.49\linewidth, trim = 0 0 0 0,clip]{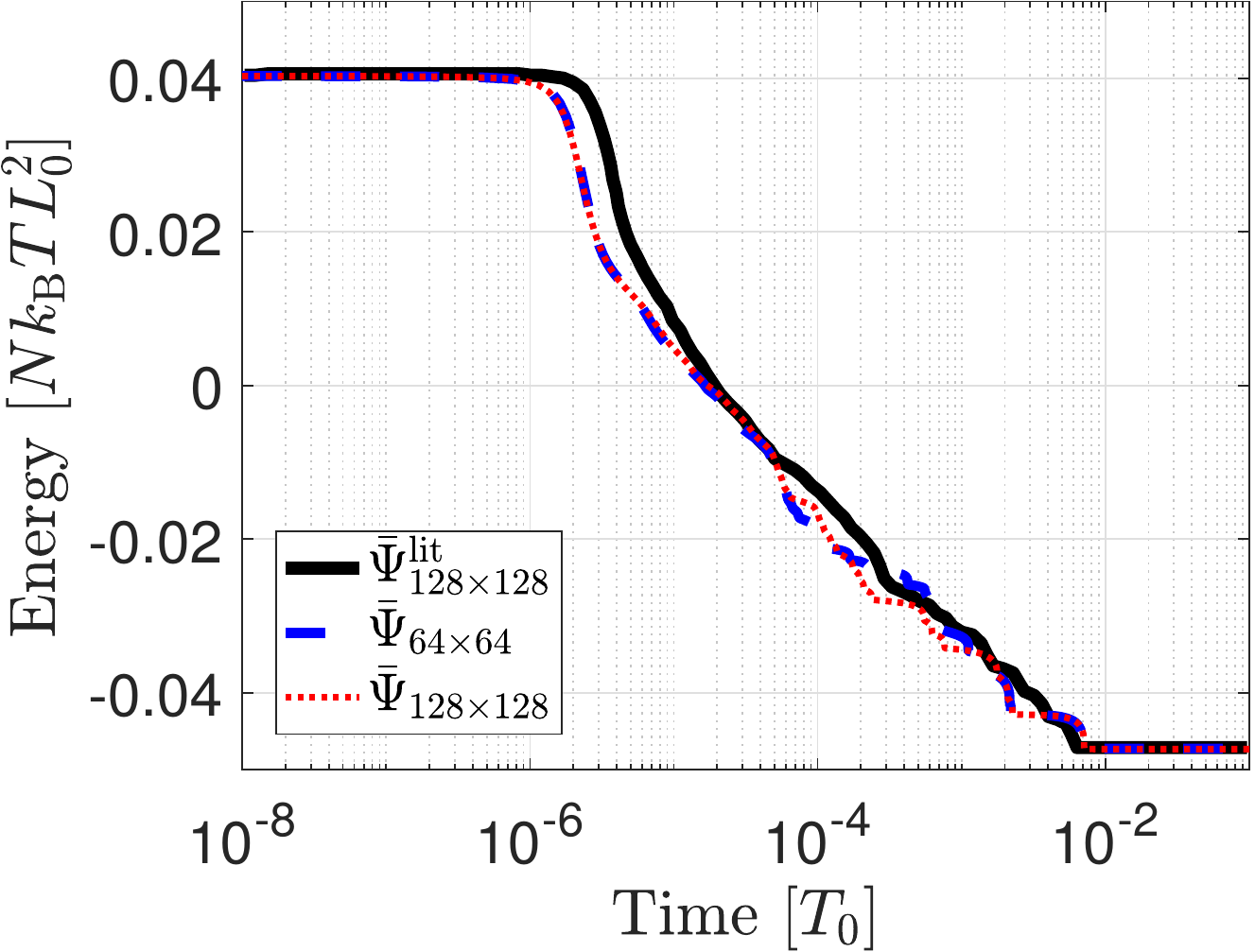}
\includegraphics[width=0.49\linewidth, trim = 0 0 0 0,clip]{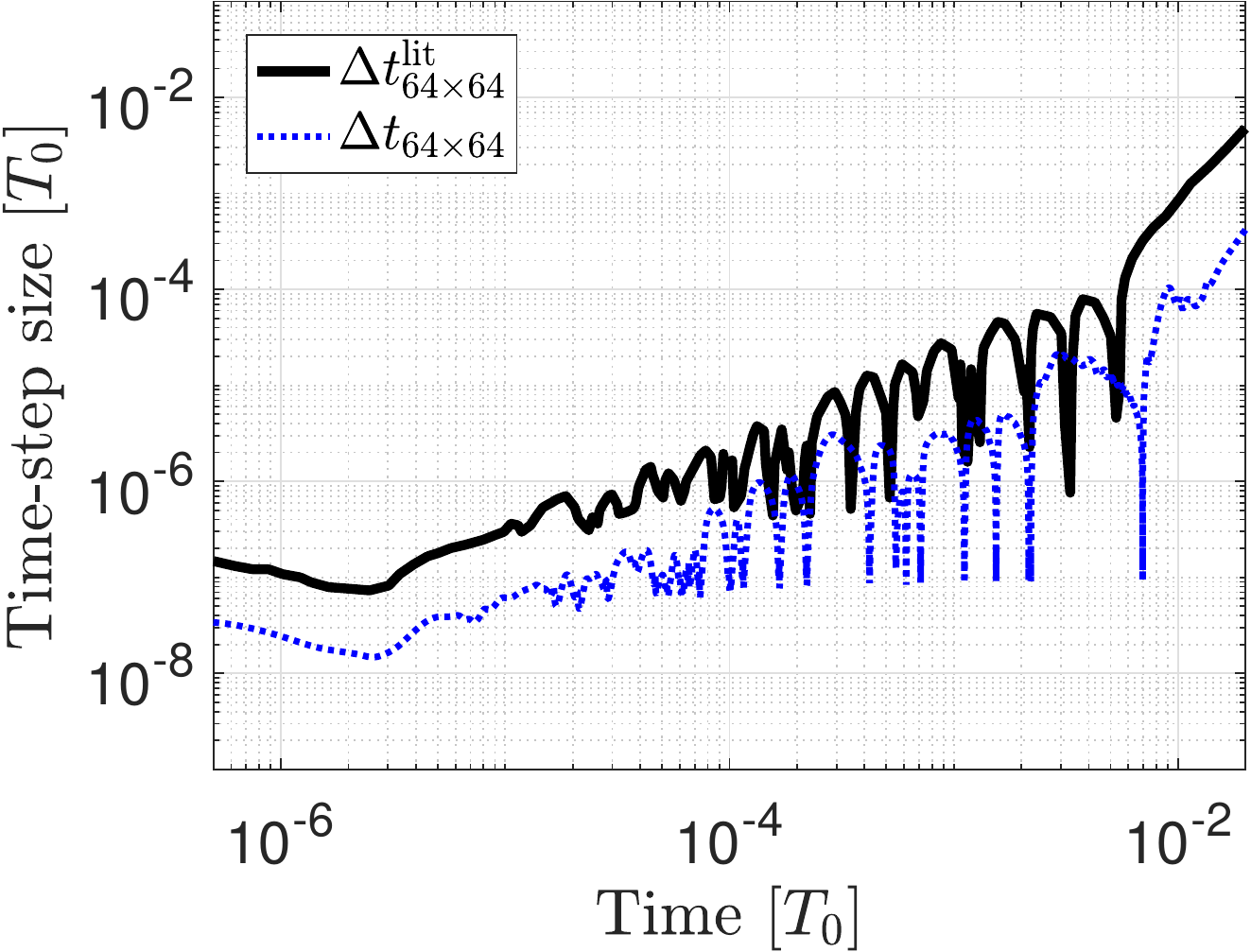}
\caption{Phase separation on a 2D square. Left: Evolution of the Helmholtz free energy defined in Eq.~\eqref{e:enI} for $\lambda=1/9000\,L_0^2$ on meshes containing $64\times 64$ and $128\times 128$ quadratic NURBS elements. 
Right: Evolution of the time step size. 
The results for $\bar\Psi_{128\times 128}^\mathrm{lit}$ and $\Delta t_{64\times 64}^\mathrm{lit}$ are taken from \citet{gomez08-1}.}
\label{fig:ss_dt}
\end{figure}

The first example considers phase separation on a 2D square following the setup of \citet{gomez08-1}.
The square has dimensions $L_0\times L_0$ and periodic boundary conditions. 
The initial volume fraction is $\bar{\phi}=0.63$.
Fig.~\ref{fig:ss_evo} shows the evolution of the phase field as a function of time starting from a random configuration and leading to complete phase separation. 
Fig.~\ref{fig:ss_dt} shows a comparison of the time step size (determined here by Eq.~\eqref{eq:ts_up}) and the free energy $\bar{\Psi}$. 
Both quantities show similar behavior and good agreement with \citet{gomez08-1}. 
Due to the randomness of the initial distribution of $\phi$, our initial condition is not exactly the same as in \citet{gomez08-1}, which results in minor differences for $\bar\Psi$.
The reason for the lower values of our time step size is a smaller tolerance for the adaptive time-stepping. 

\subsubsection{Phase separation on a rigid sphere}

\begin{figure}[H]
\centering
\includegraphics[width=0.24\linewidth, trim = 300 70 300 150,clip]{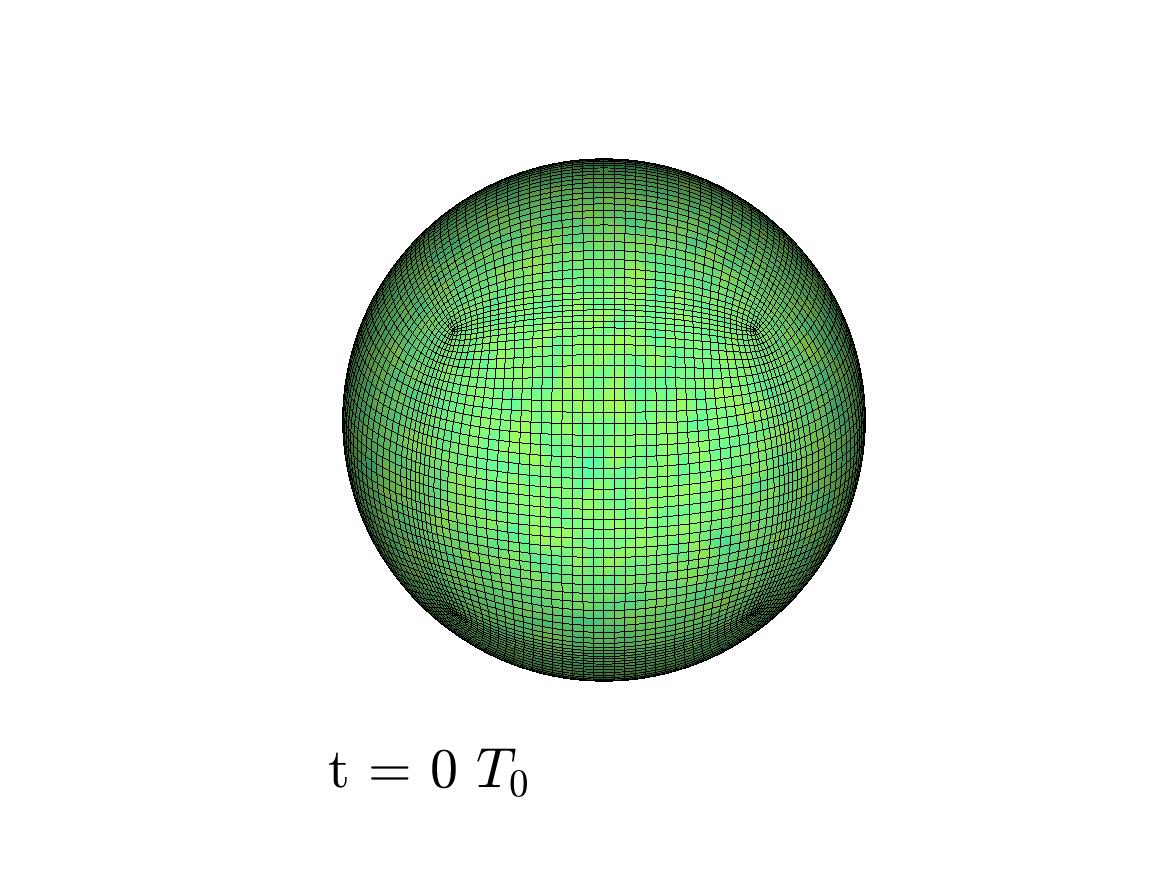}
\includegraphics[width=0.24\linewidth, trim = 300 70 300 150,clip]{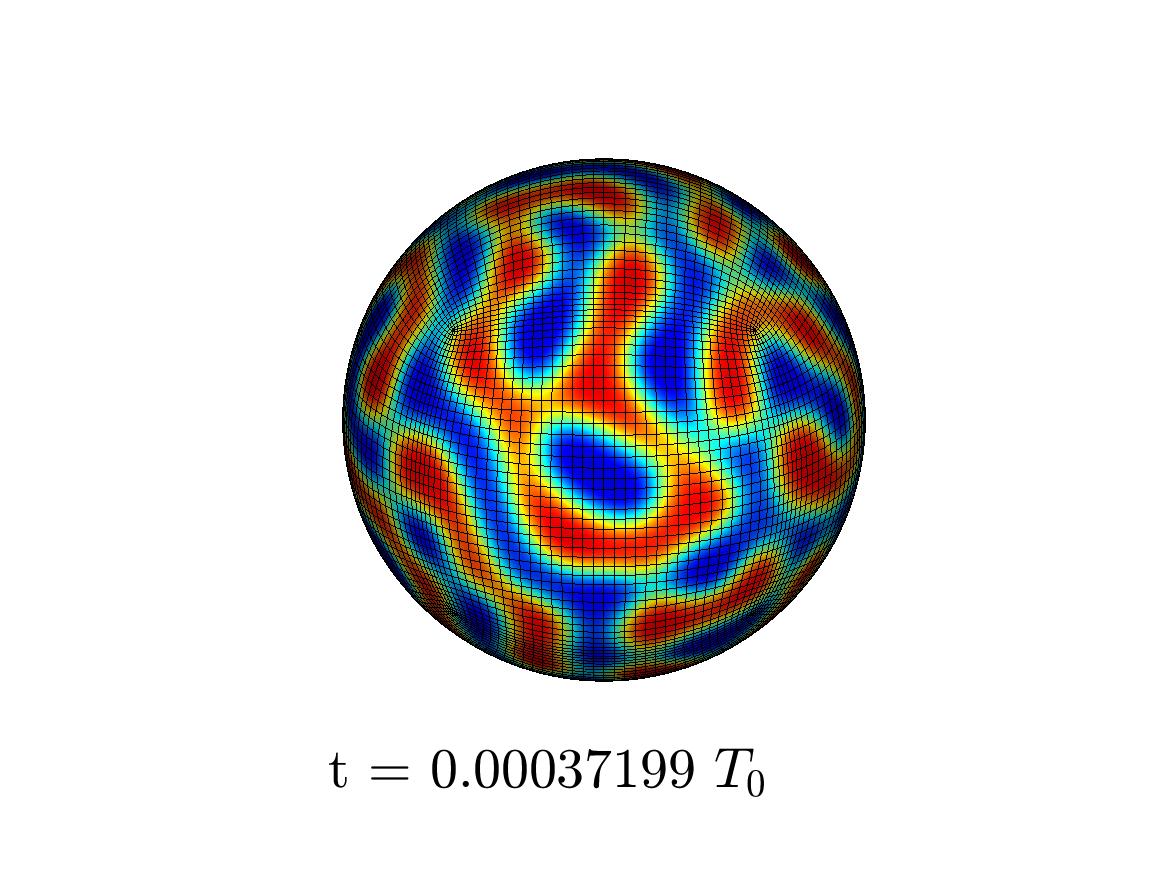}
\includegraphics[width=0.24\linewidth, trim = 300 70 300 150,clip]{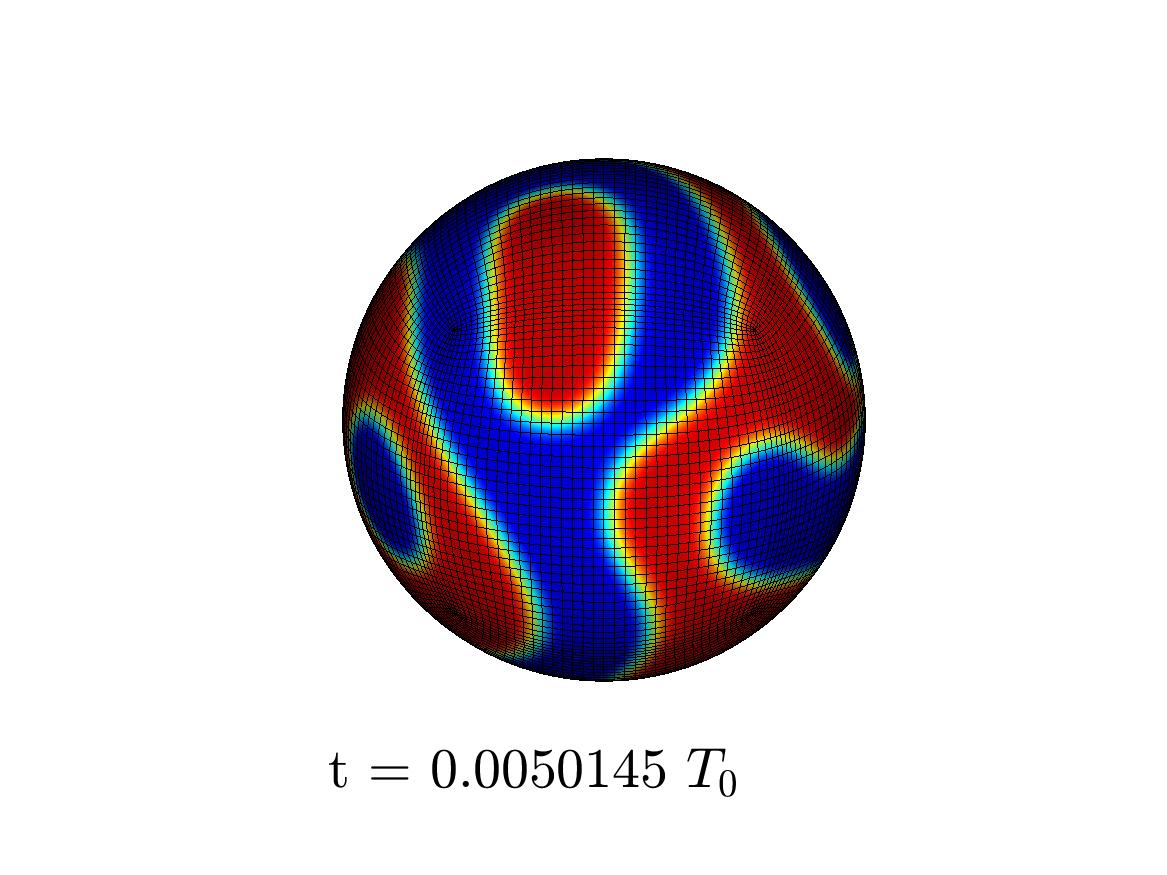}
\includegraphics[width=0.24\linewidth, trim = 300 70 300 150,clip]{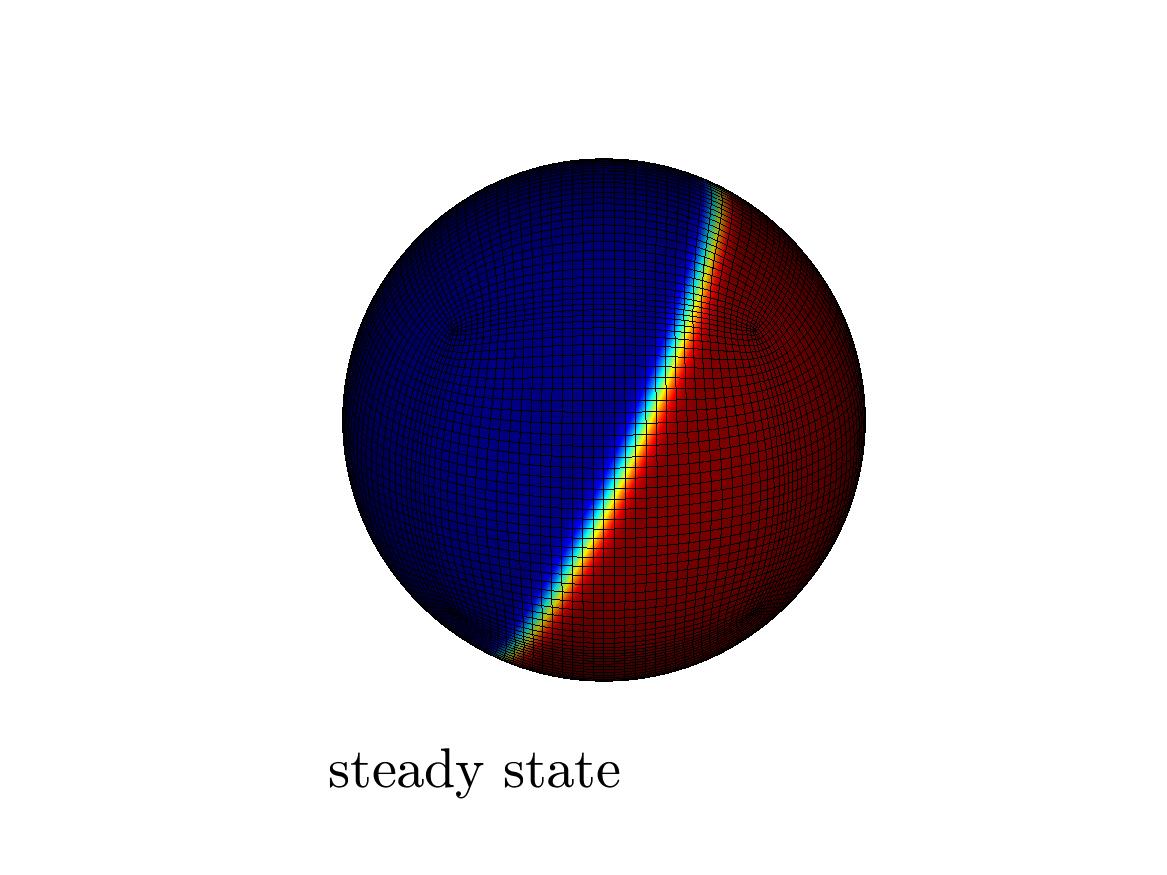}
\caption{Phase separation on a rigid sphere. Evolution of the solution for $\lambda=1.3144\cdot10^{-3}\,L_0^2$ and volume fraction $\bar{\phi} = 0.5$ on a mesh containing $9672$ cubic elements. The coloring follows Table~\ref{tab:mat}.}
\label{fig:sp_evol}
\end{figure}

\begin{figure}[H]
\centering
\includegraphics[width=0.49\linewidth, trim = 0 0 0 0,clip]{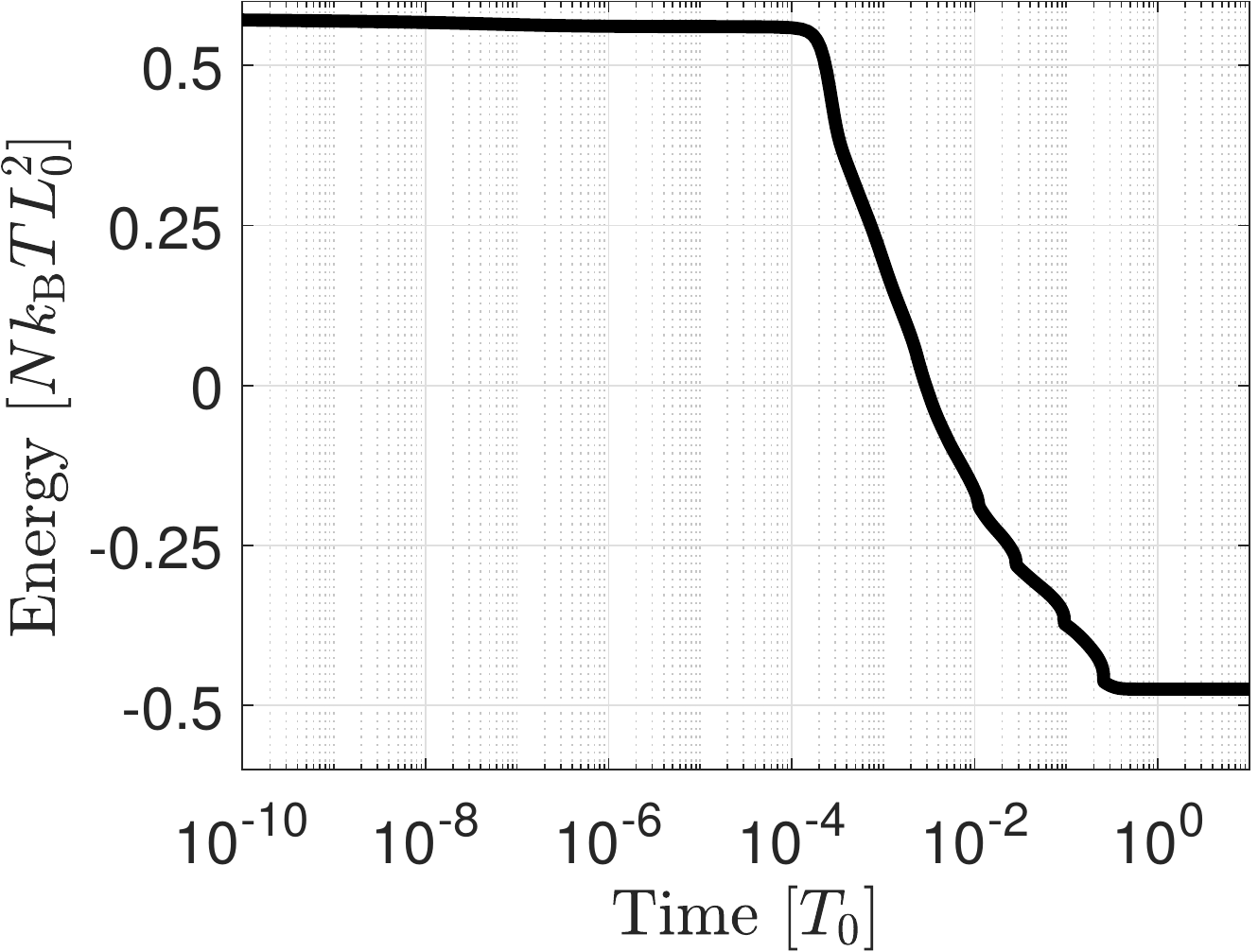}
\includegraphics[width=0.49\linewidth, trim = 0 0 0 0,clip]{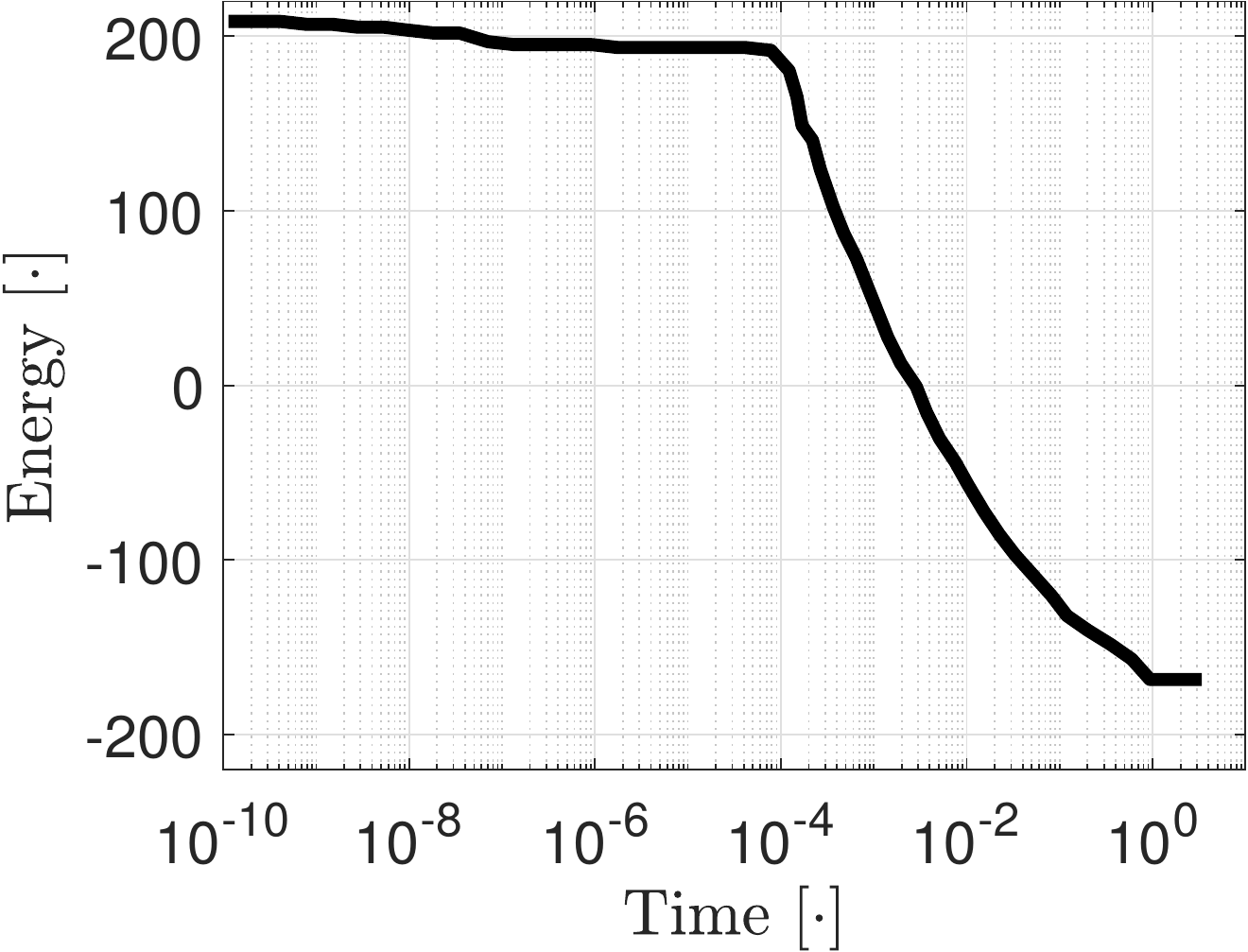}
\caption{Phase separation on a rigid sphere with $\lambda=1.3144\cdot10^{-3}\,L_0^2$ and $\bar{\phi} = 0.5$. 
Left: Evolution of the free energy on a mesh containing $9672$ cubic elements (present result). 
Right: Evolution of the free energy on a mesh containing $8844$ quadratic NURBS elements \citep{bartezzaghi15}.}
\label{fig:sp_lit}
\end{figure}

The second example studies the phase separation on a rigid sphere following the setup of \citet{bartezzaghi15}. 
Fig.~\ref{fig:sp_evol} shows the phase separation over time for this example. 
Fig.~\ref{fig:sp_lit} shows the evolution of the free energy compared to the results from \citet{bartezzaghi15}. 
The comparison shows a similar evolution in time, but the absolute values are different due to a different normalization of the governing equations. 

\subsection{Phase separation on a deforming torus}\label{s:exT}

\begin{figure}[H]
\centering
\includegraphics[width=0.40\linewidth, trim = 290 330 230 300,clip]{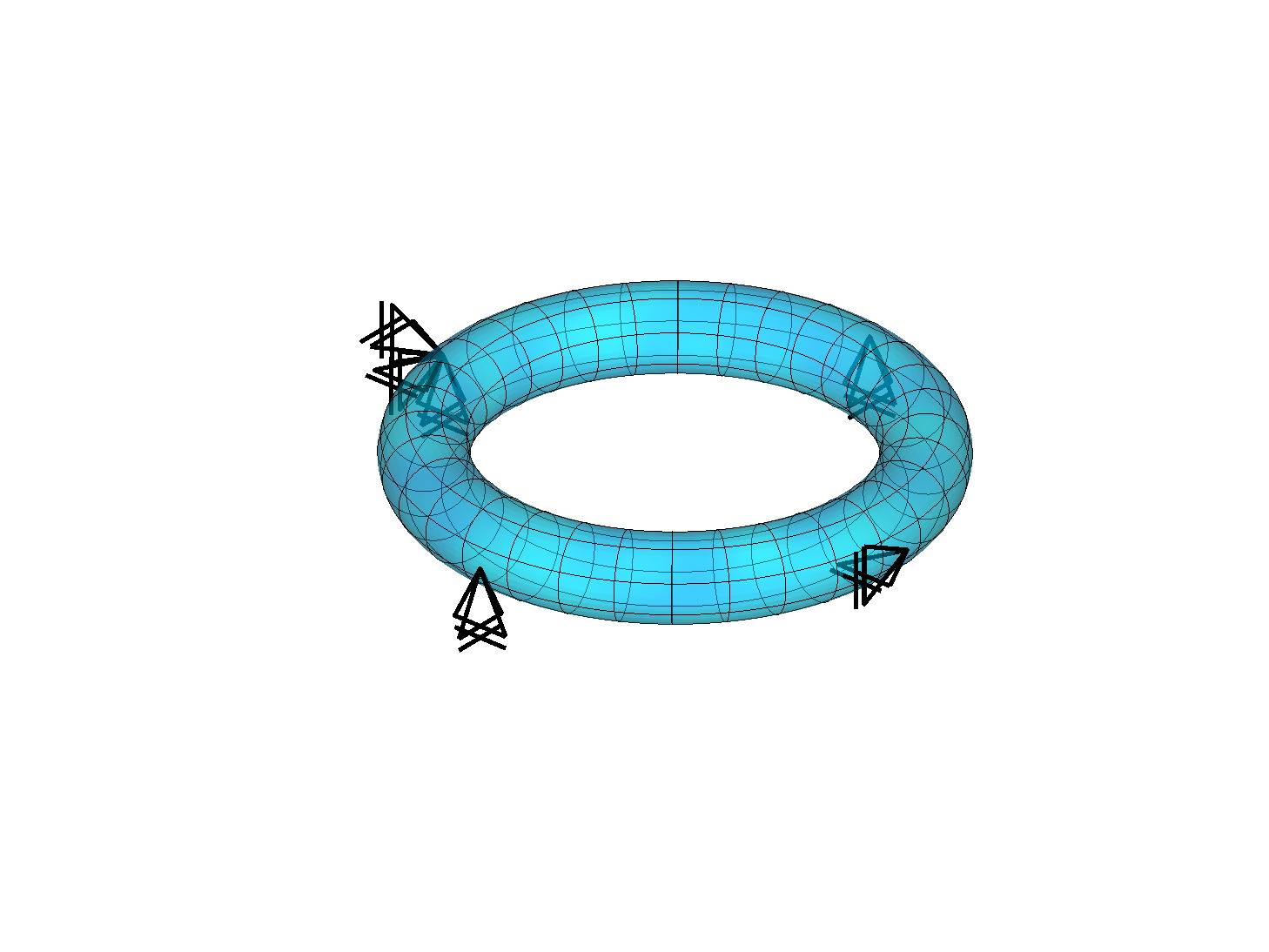}\hspace{4mm}
\includegraphics[width=0.40\linewidth, trim = 290 330 230 300,clip]{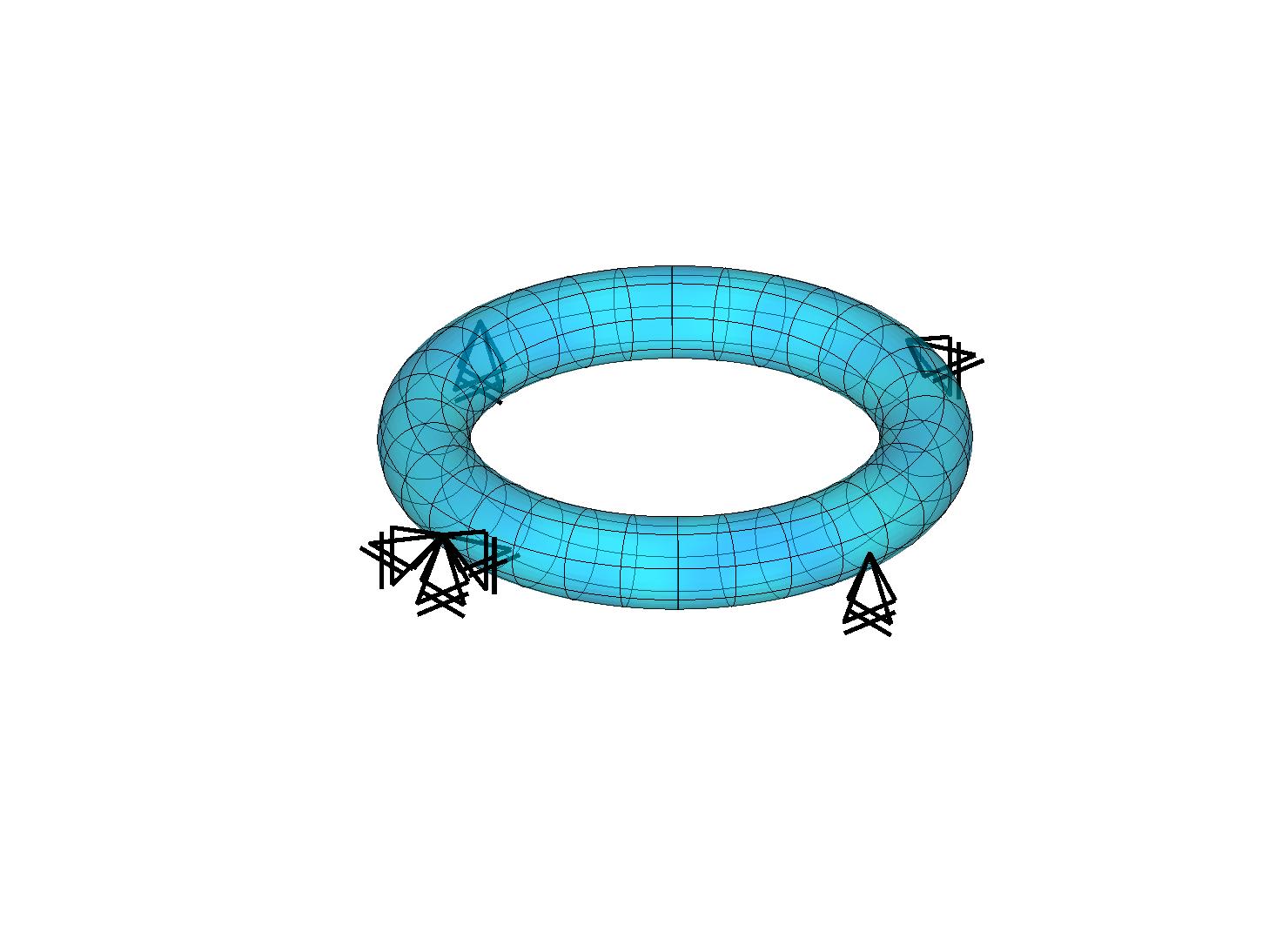}
\caption{Phase separation on a deforming torus: Initial configuration and boundary conditions shown from two different viewpoints. 
The displacement of six dofs is fixed to prevent rigid body motions.
}
\label{fig:t_IC1}
\end{figure}

The following two examples study phase separations on a deformable torus using the proposed material coupling of Sec.~\ref{s:mix} and Table~\ref{tab:mat}. 
A constant internal pressure $p_{\mathrm{int}}=0.1\,EL_0^{-1}$ is prescribed for all $t$ to provide mechanical loading. 
The boundary conditions are illustrated in Fig.~\ref{fig:t_IC1}. 
This is the first non-trivial example, where both the phase field and surface deformations evolve simultaneously. 

\subsubsection{Large phase interface}\label{Sec:ex1}

The first example studies the behavior of different spatial discretizations with identical initial configurations. 
A comparison for meshes containing $8 \times 32$, $16 \times 64$,  $32 \times 128$, $64 \times 256$ and  $128 \times 512$ quadratic NURBS elements is provided. 
The mechanical material parameters are listed in Table~\ref{tab:mat}. 
The mobility constant is selected to be $D=4\,T_0$ and the interfacial thickness parameter $\lambda = 0.075\,L_0^2$ is chosen. 
This is a relatively large value that allows to use coarse meshes: for the coarsest mesh $\sqrt{\lambda}\approx h$ and for the finest mesh $\sqrt{\lambda}\approx 16\,h$, where h is the average element size. 
The constant internal pressure $p_{\mathrm{int}}=0.1\,EL_0^{-1}$ is prescribed for all $t$.
\begin{figure}[H]
\centering
\includegraphics[width=0.32\linewidth, trim = 300 310 200 570,clip]{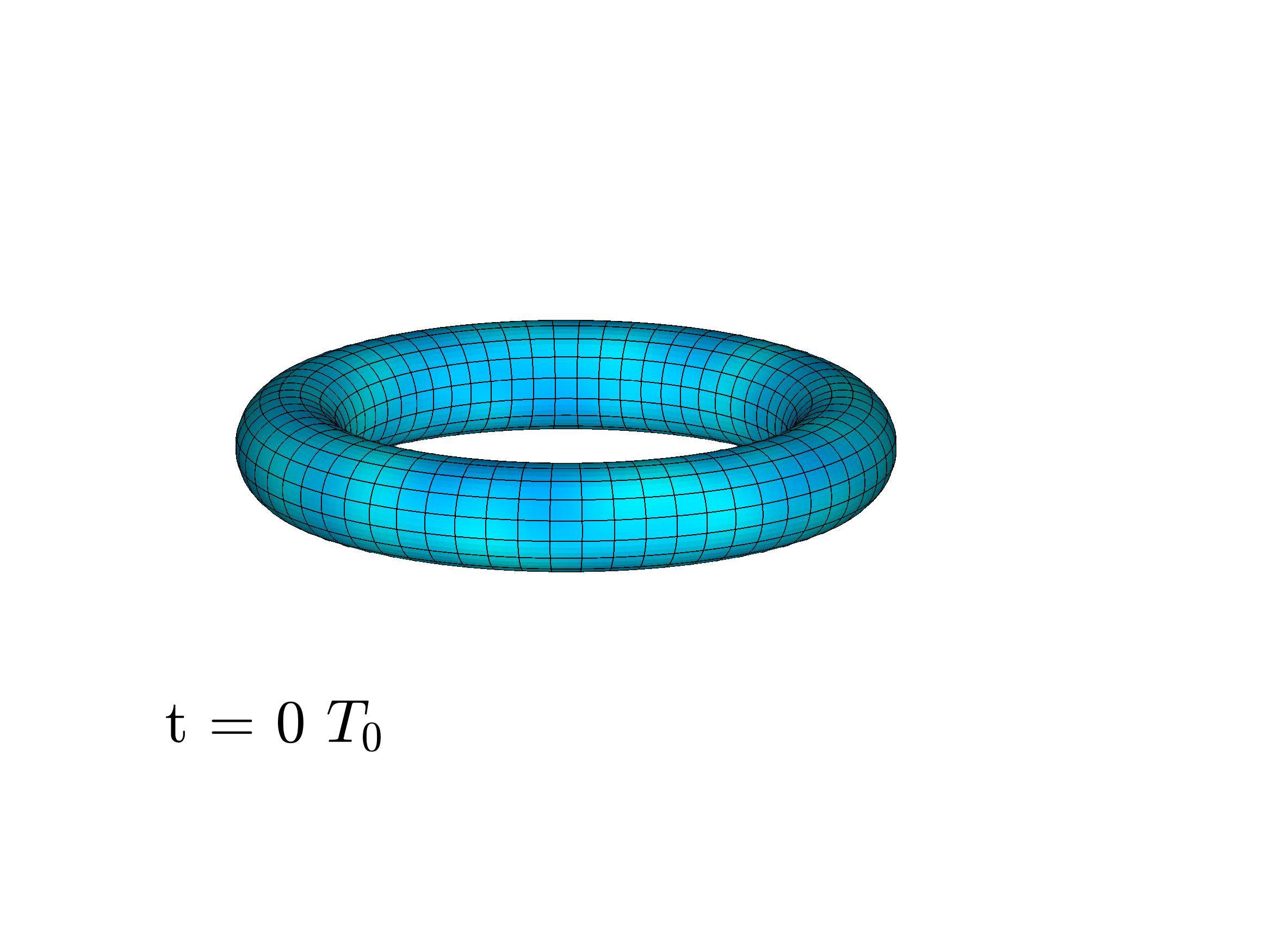}
\includegraphics[width=0.32\linewidth, trim = 300 310 200 570,clip]{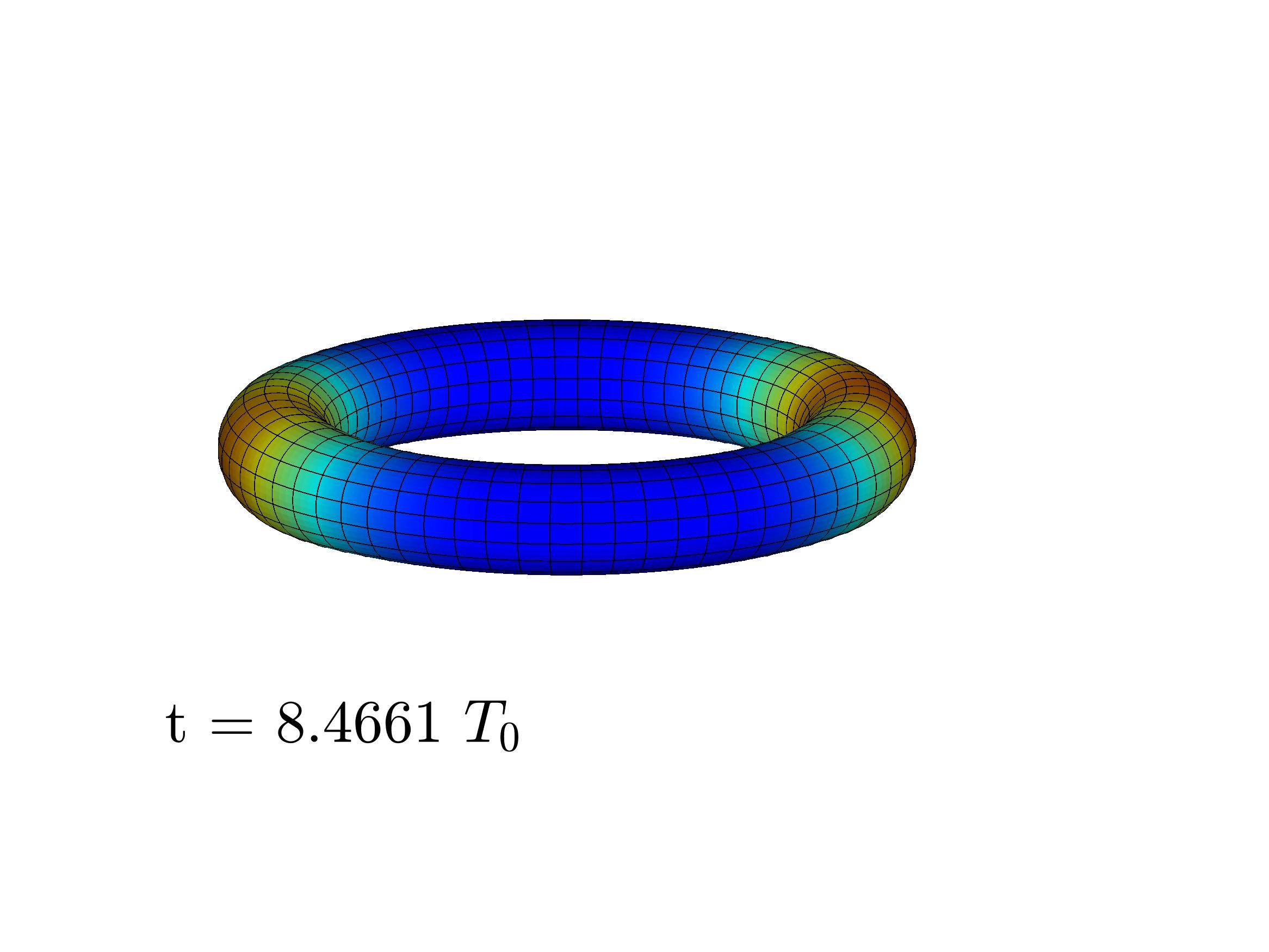}
\includegraphics[width=0.32\linewidth, trim = 300 310 200 570,clip]{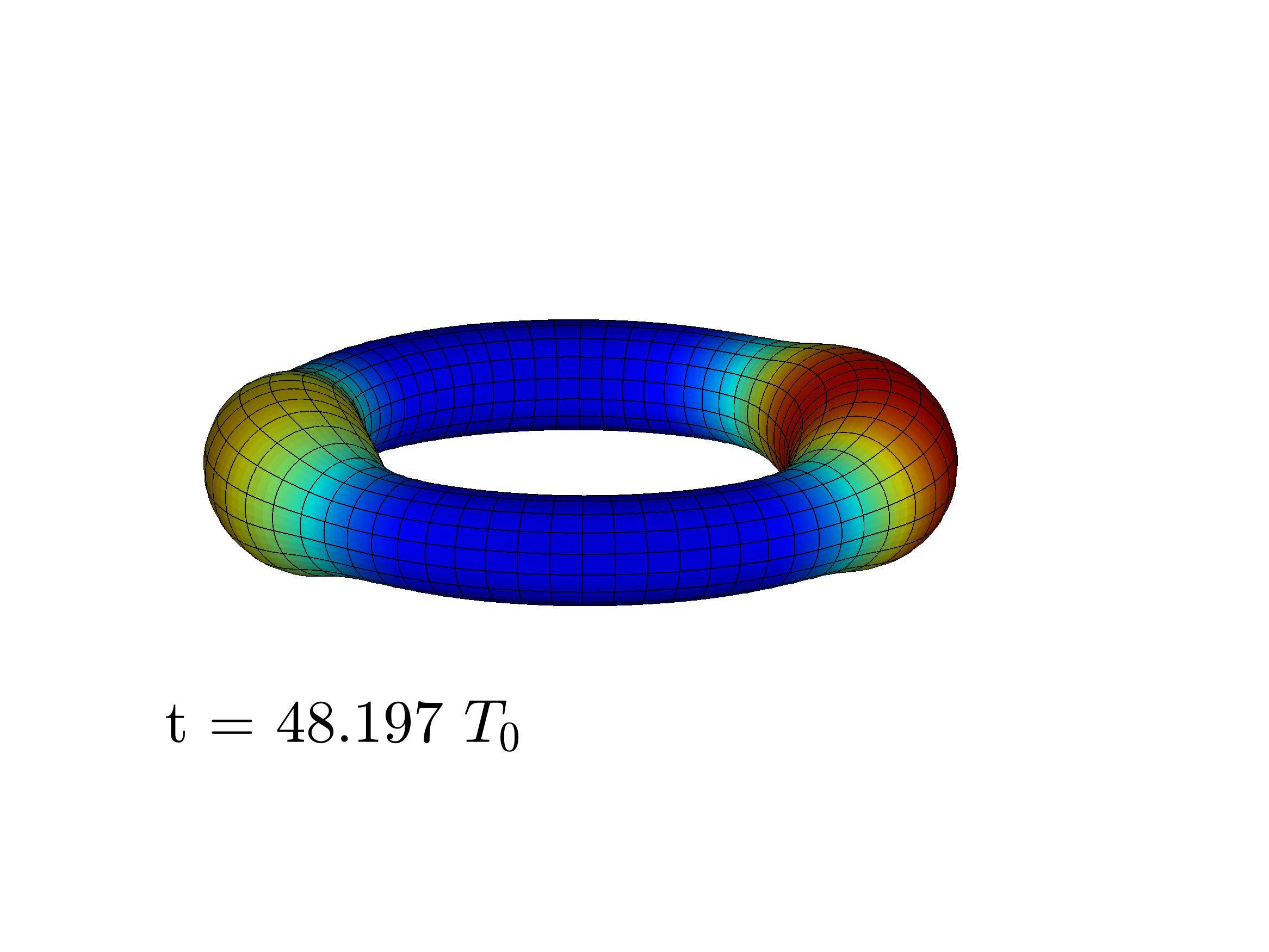}\\
\includegraphics[width=0.32\linewidth, trim = 300 310 200 500,clip]{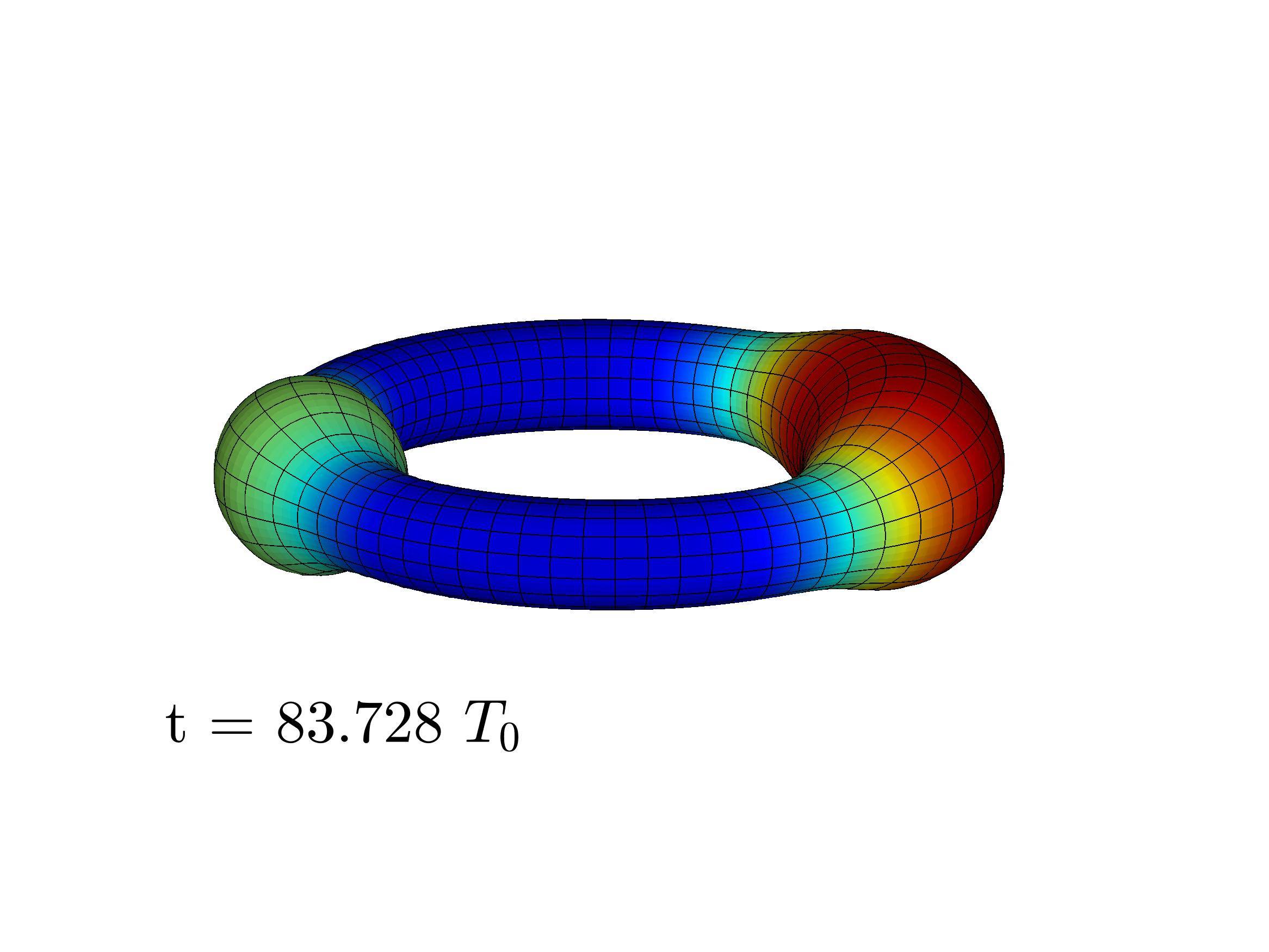}
\includegraphics[width=0.32\linewidth, trim = 300 310 200 500,clip]{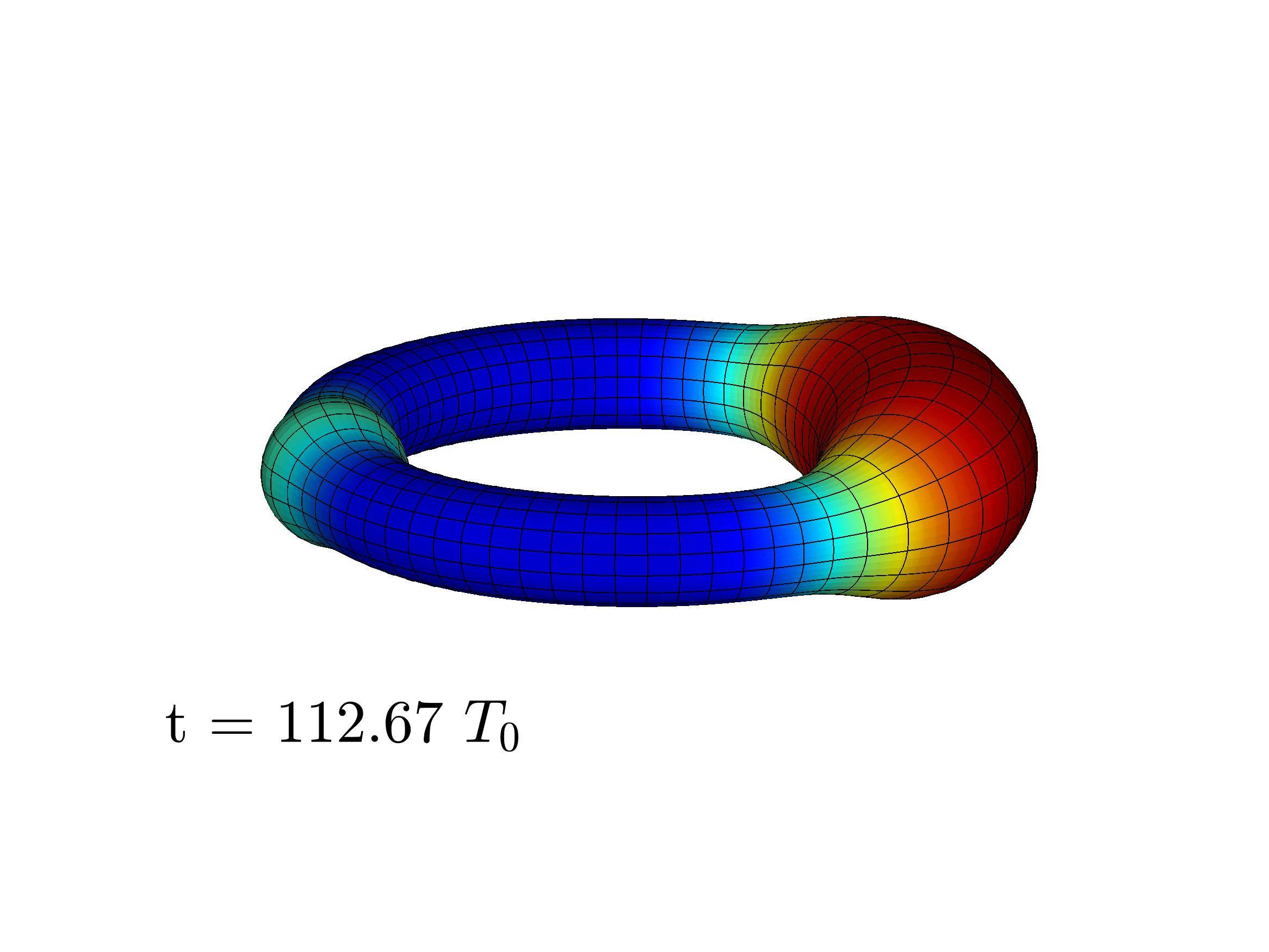}
\includegraphics[width=0.32\linewidth, trim = 300 310 200 500,clip]{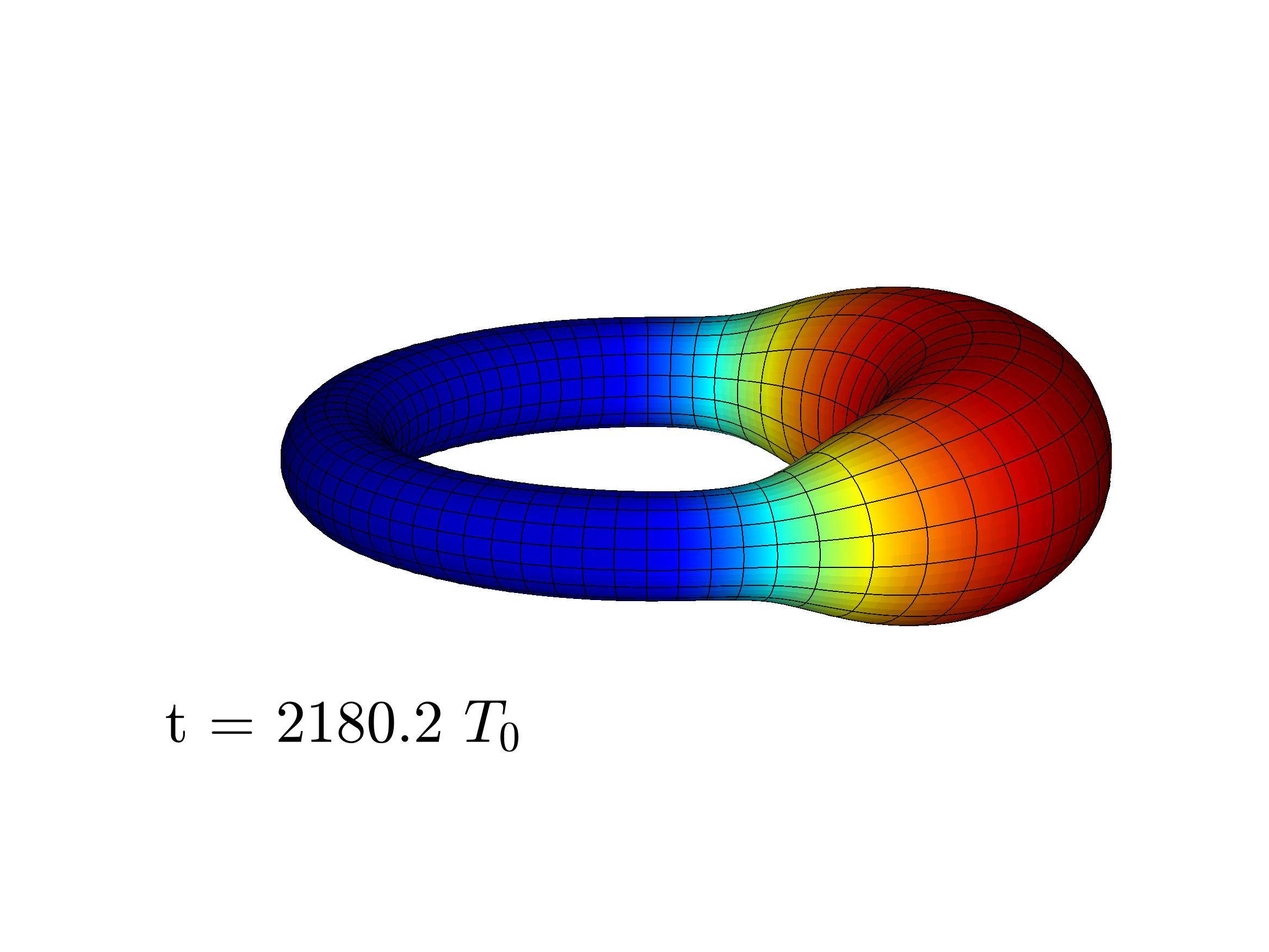}
\caption{Phase separation on a deforming torus: Evolution of the phase field for $\lambda=0.075\,L_0^2$ on a deforming torus containing $16 \times 64$ quadratic NURBS elements. The coloring follows Table~\ref{tab:mat}.}
\label{fig:t_evo1}
\end{figure}
Fig.~\ref{fig:t_evo1} shows the evolution of the phase separation over time. 
The material behavior of phase $\phi=0$ (blue color) is much stiffer than the material behavior of phase $\phi=1$ (red color). 
Therefore, large bulges appear in the red phase that grow in time as the red phase becomes larger\footnote{In all the following figures, the true deformation without any scaling is visualized.}.
Fig.~\ref{fig:t_evo1} also shows that the deformation and phase separation evolve at a similar time scale. 
The mechanical response is strongly affected by viscosity.
Low values of $\eta$ lead to strong oscillations.

\begin{figure}[H]
\centering
\includegraphics[width=0.49\linewidth, trim = 0 0 0 0,clip]{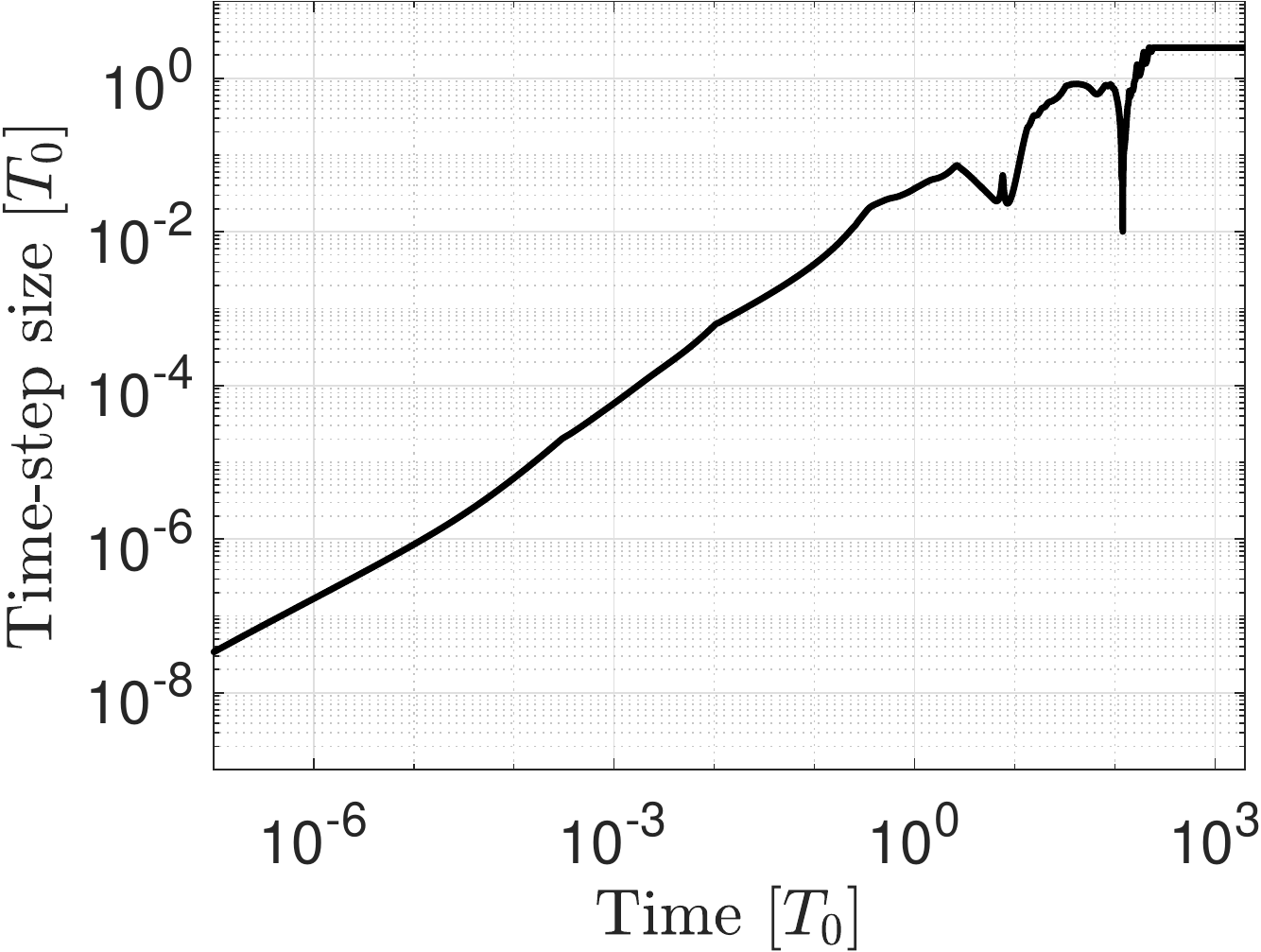}
\includegraphics[width=0.49\linewidth, trim = 0 0 0 0,clip]{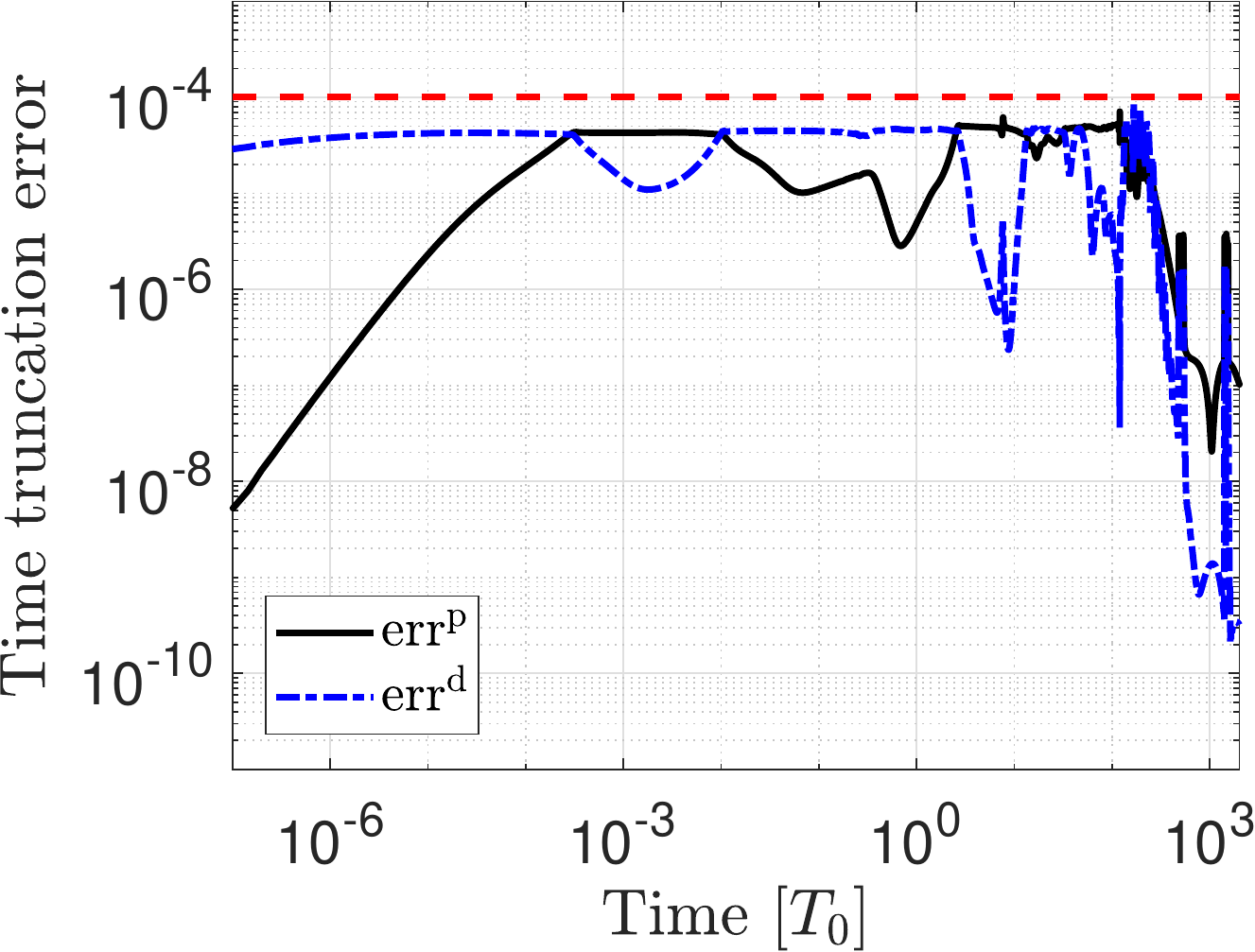}
\caption{Phase separation on a deforming torus: Left: Time step sizes during phase separation. Right: Evolution of the local time truncation errors for the phase field, $\mathrm{err}^\mathrm{p}$, and mechanical field, $\mathrm{err}^\mathrm{d}$. The temporal error bound is shown by a red dashed line.}
\label{fig:t_ts1}
\end{figure}
The left side of Fig.~\ref{fig:t_ts1} shows the time step size resulting from the adaptive time stepping procedure of Sec.~\ref{Sec:ATS}. 
The right side of Fig.~\ref{fig:t_ts1} shows the local time truncation errors $\text{err}^\mrp$ and $\text{err}^\mrd$ defined in Eq.~\eqref{e:err}. 
It can be observed that the time step is restricted in an alternating manner, by either the mechanical error (dot-dashed blue line) or the phase field error (solid black line). 
The temporal error bound for rejecting and recomputing the time step is chosen at $10^{-4}$ (red dashed line). 
The maximum time step size is limited to $\Delta t = 2.5\,T_0$ to ensure sufficient accuracy and stability. 

\begin{figure}[H]
\centering
\includegraphics[width=0.49\linewidth, trim = 0 0 0 0,clip]{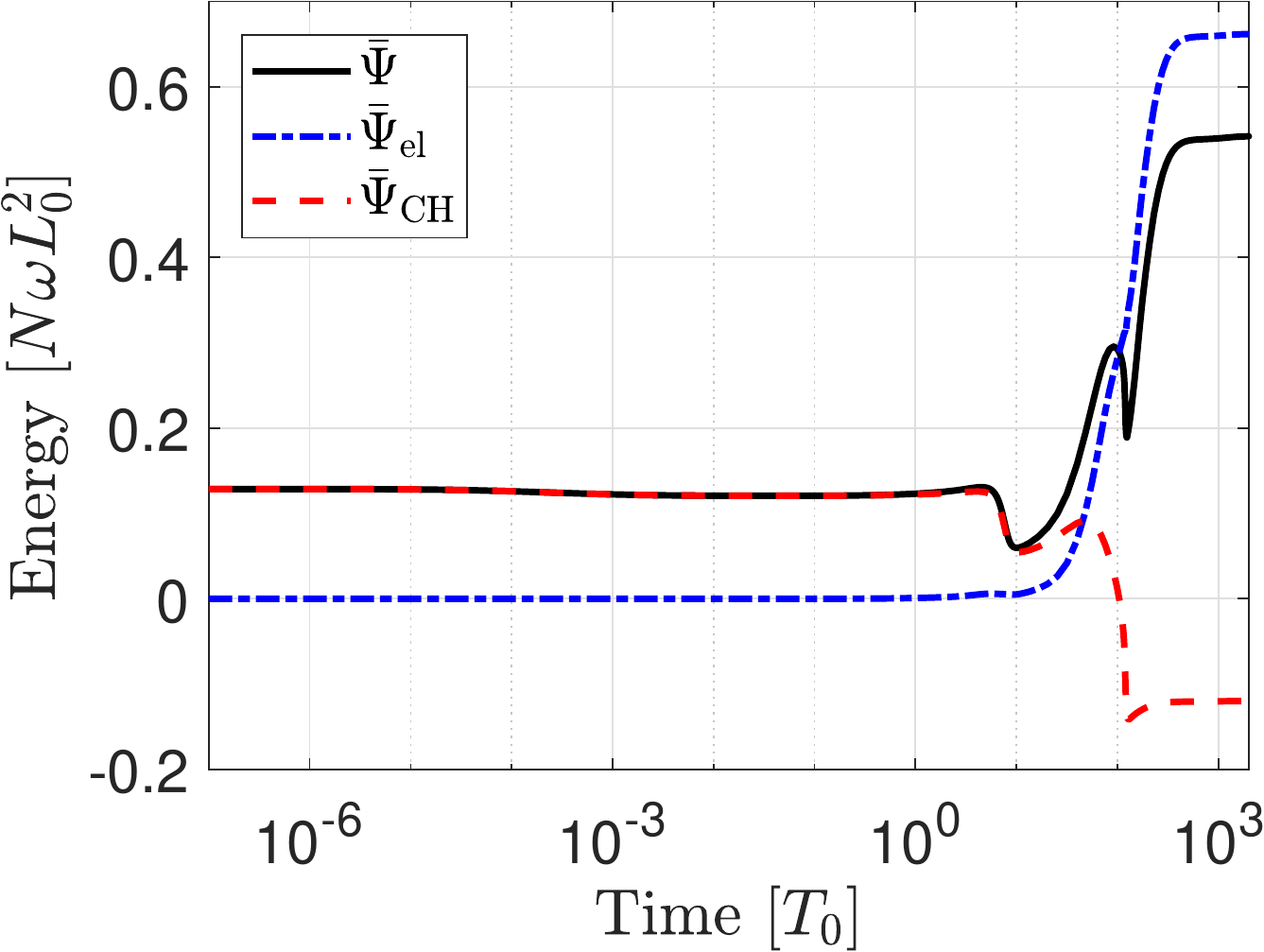}
\includegraphics[width=0.49\linewidth, trim = 0 0 0 0,clip]{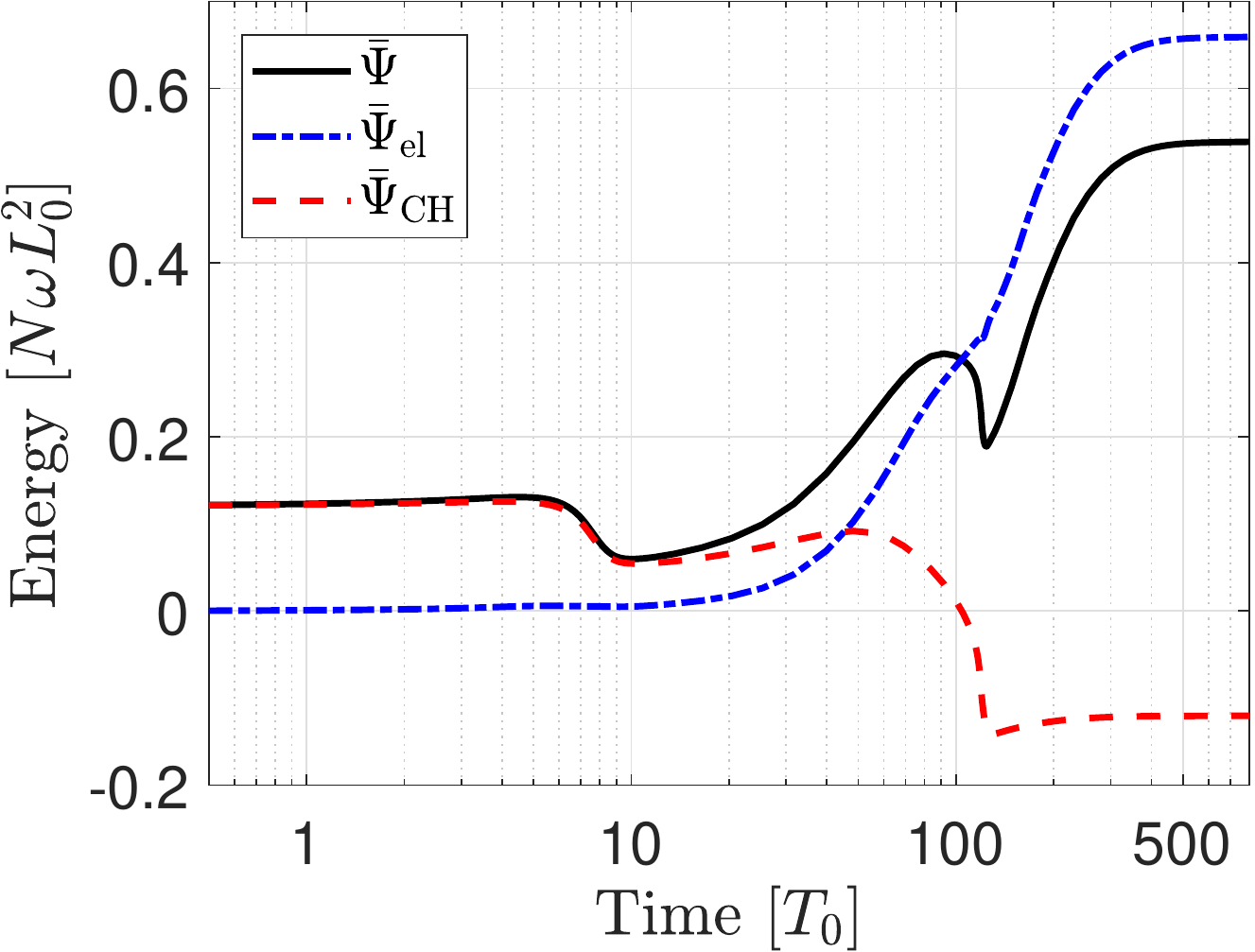}
\caption{Phase separation on a deforming torus: Evolution of the characteristic energies of the system on the left side (quadratic NURBS, mesh: $16\times 64$). Close-up on the right side.}
\label{fig:t_ene1}
\end{figure}
Fig.~\ref{fig:t_ene1} shows the evolution of the characteristic energies of the system. 
Initially, $\bar\Psi_\mathrm{CH}$ is large compared to $\bar\Psi_\mathrm{el}$, but then decreases during phase separation due to lowering of $\Psi_\mathrm{mix}$ shown in Fig.~\ref{f:Psimix}. The kink in Fig.~11 at time $125\,T_0$ reflects the state at which the two phases completely separate. At that time the Cahn-Hilliard energy decreases, while the deformation, and thus $\bar\Psi_\mathrm{el}$, increases with a slight delay due to viscosity.
After the phases are completely separated the system reaches a steady state.

\begin{figure}[H]
\centering
\includegraphics[width=0.49\linewidth, trim = 0 0 0 0,clip]{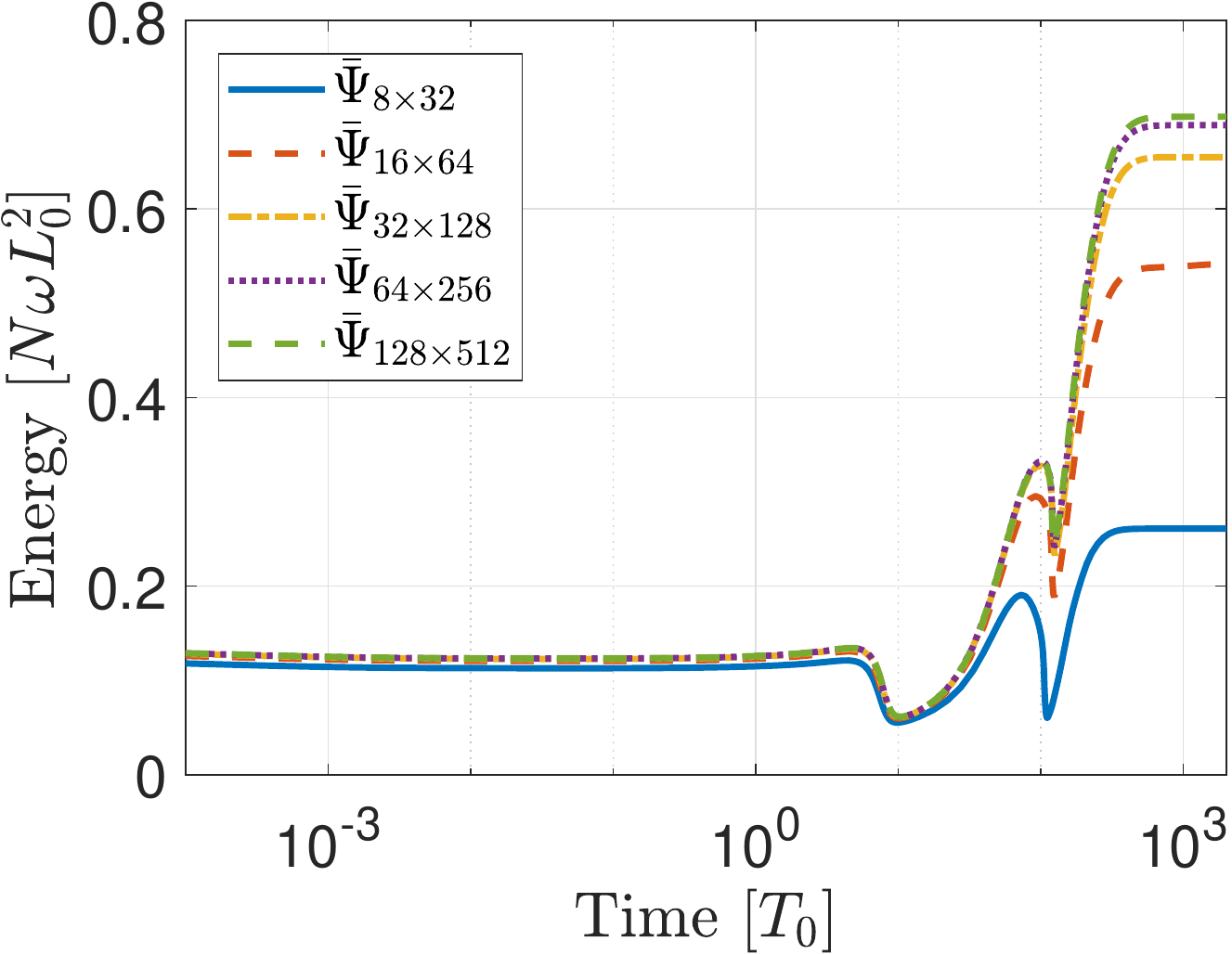}
\includegraphics[width=0.49\linewidth, trim = 0 0 0 0,clip]{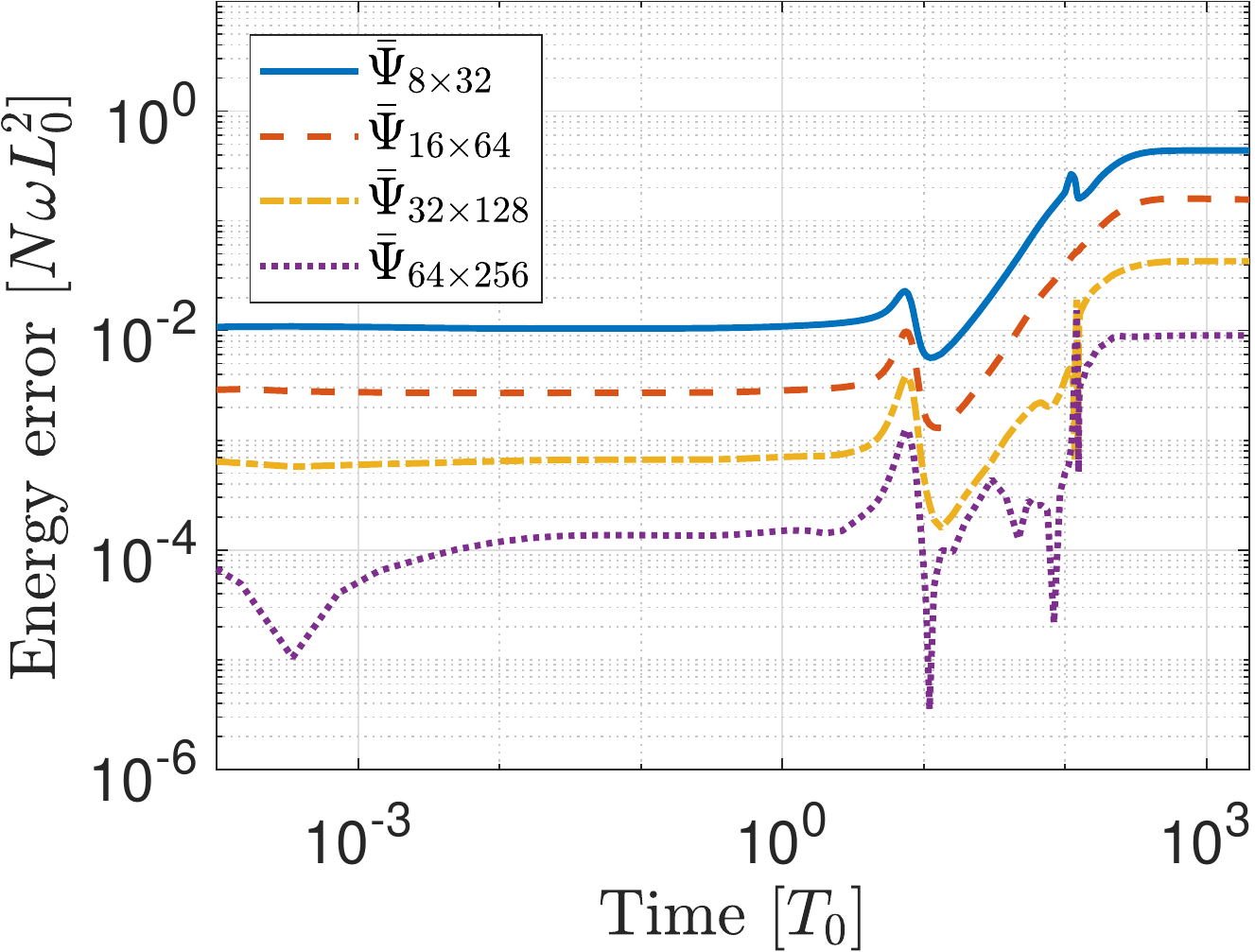}
\caption{Phase separation on a deforming torus: Left: Evolution of the Helmholtz free energy for five different quadratic NURBS discretizations. Right: Error of the Helmholtz free energy for the coarser meshes $8\times 32$, $16\times 64$, $32\times 128$ and $64\times 256$ with respect to the fine mesh $128\times 512$.}
\label{fig:t_ena1}
\end{figure}
The left side of Fig.~\ref{fig:t_ena1} shows the evolution of the Helmholtz free energy for the five different NURBS meshes. 
A good agreement of the evolution of the Helmholtz free energy can be observed for all meshes, except for the coarse meshes $8\times 32$ and $16\times 64$. 
The right side of Fig.~\ref{fig:t_ena1} shows the error of the Helmholtz free energy of the coarser meshes with respect to the finest mesh. 
The error decreases with increasing mesh refinement.
After the steady state is reached, the energy error stays constant for all meshes. 

\subsubsection{Small phase interface}
\label{Sec:ex2}

\begin{figure}[H]
\centering
\includegraphics[width=0.32\linewidth, trim = 300 310 200 400,clip]{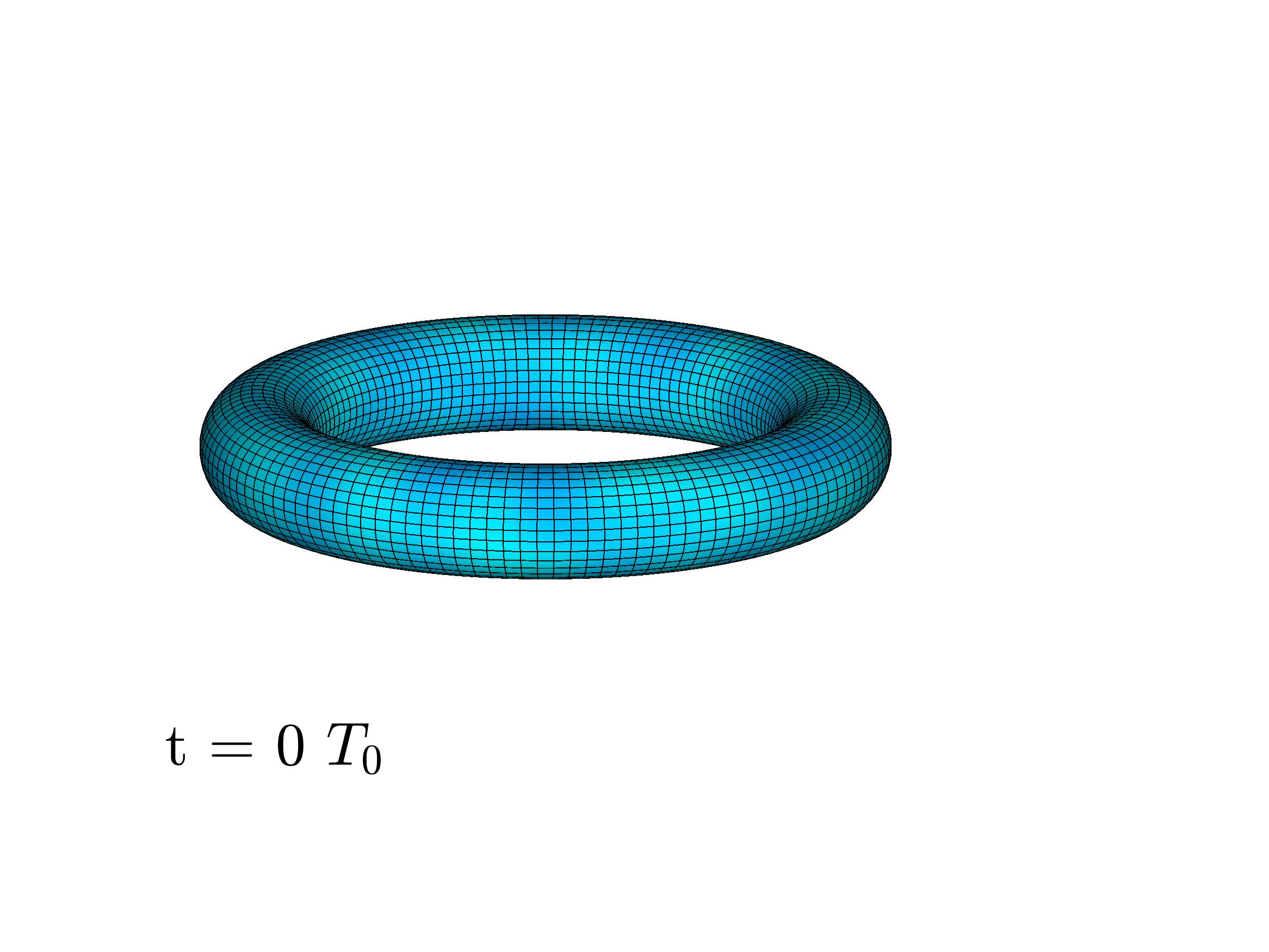}
\includegraphics[width=0.32\linewidth, trim = 300 310 200 400,clip]{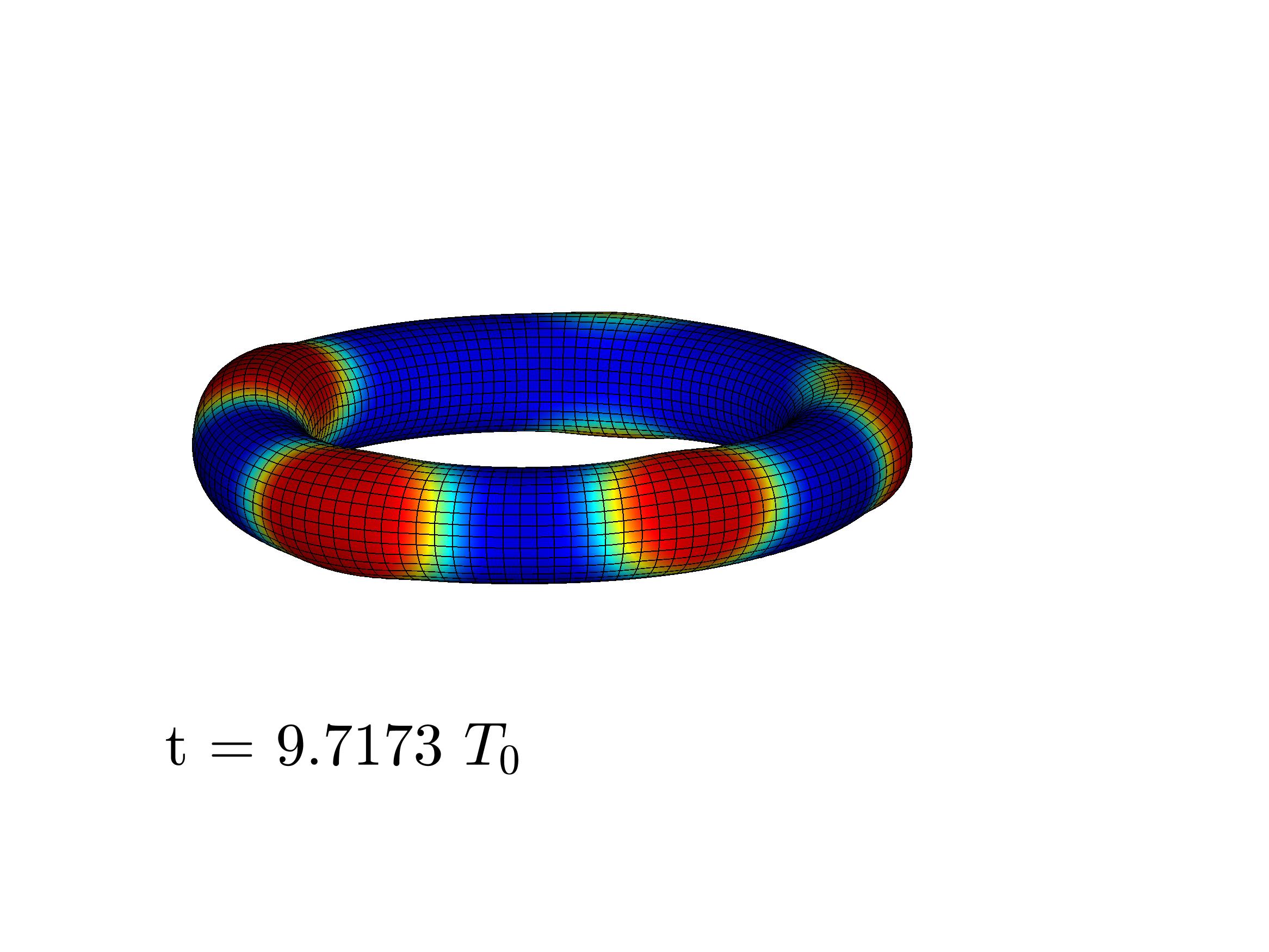}
\includegraphics[width=0.32\linewidth, trim = 300 310 200 400,clip]{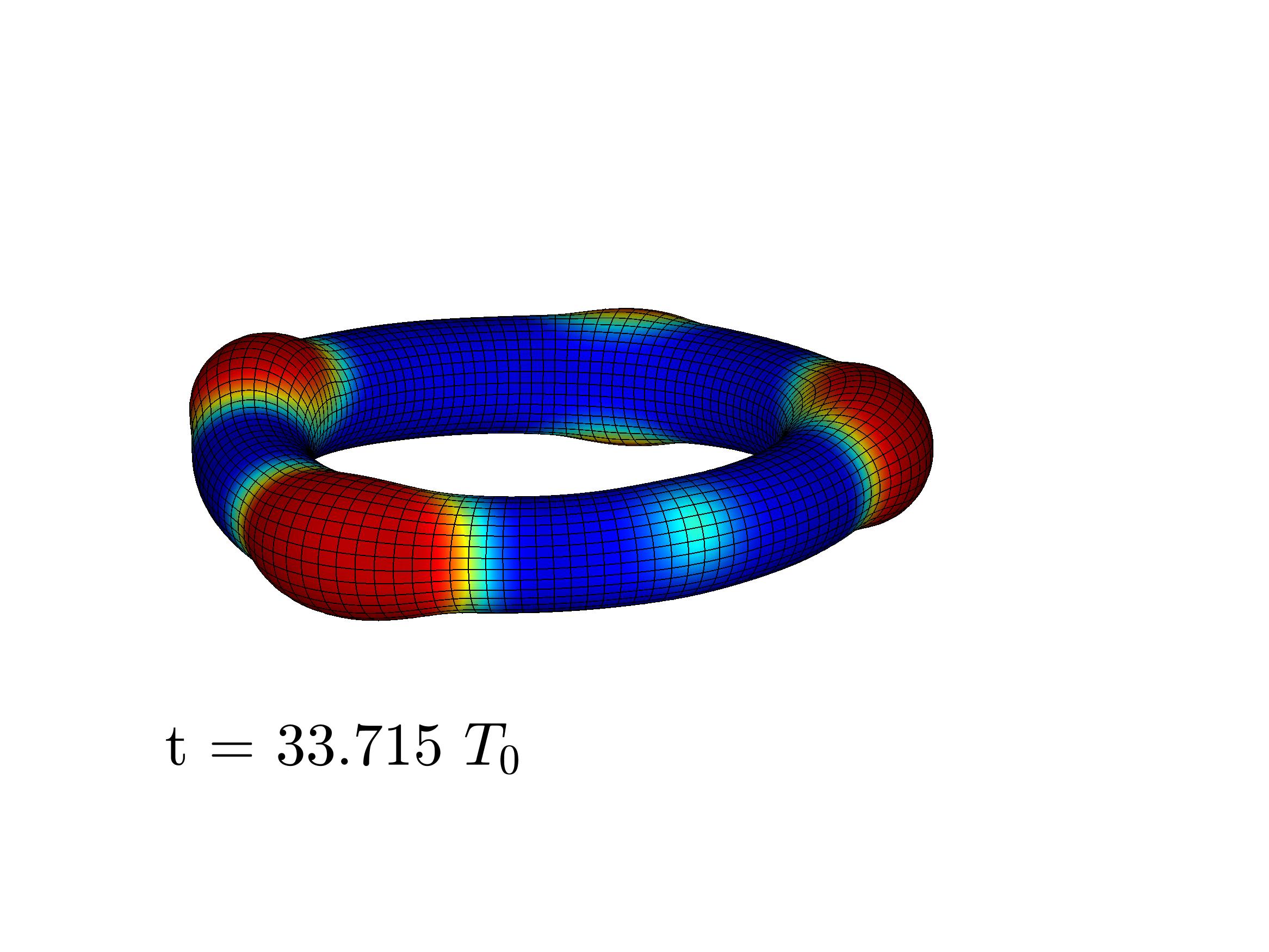}\\
\includegraphics[width=0.32\linewidth, trim = 300 310 200 400,clip]{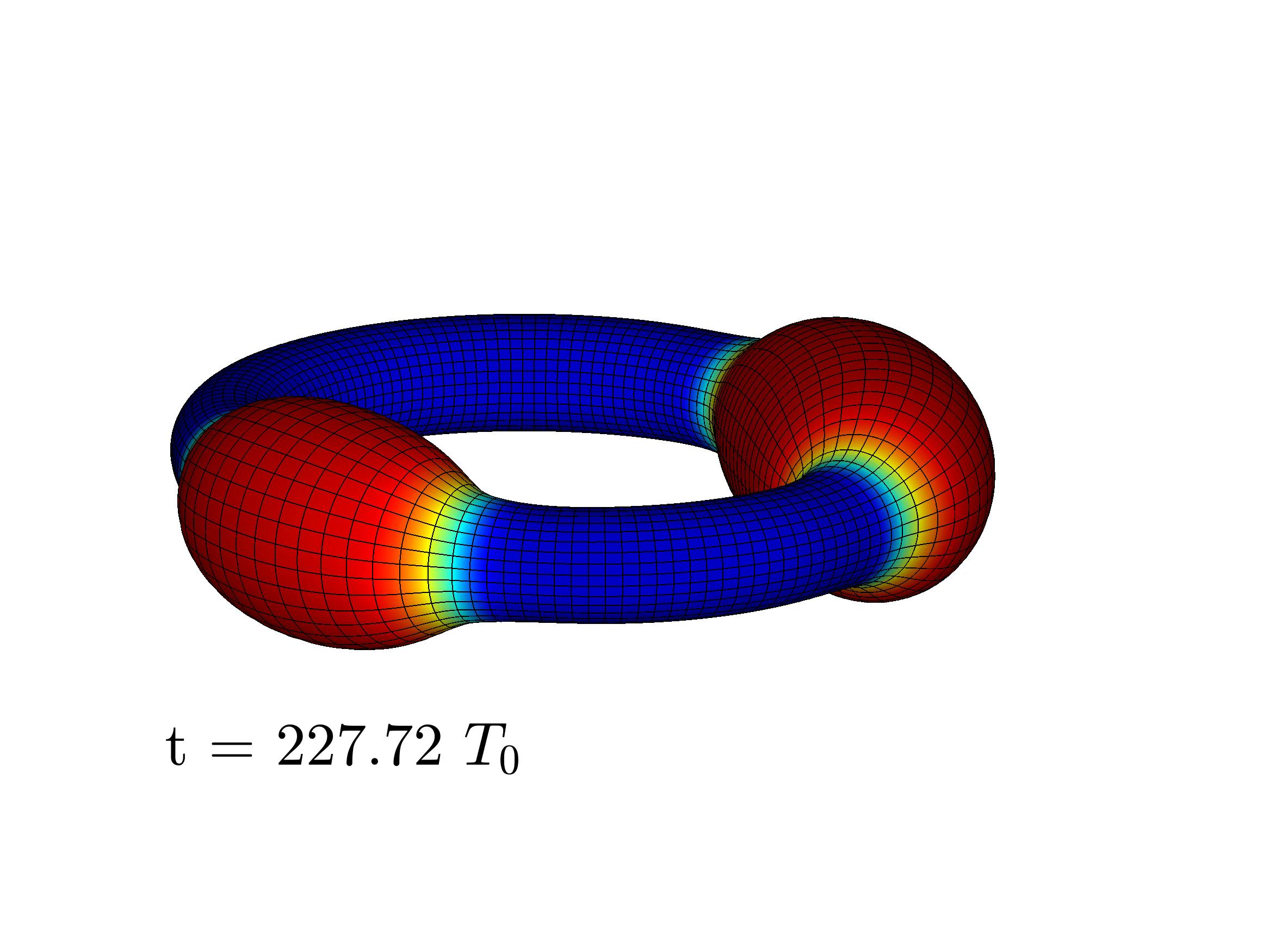}
\includegraphics[width=0.32\linewidth, trim = 300 310 200 400,clip]{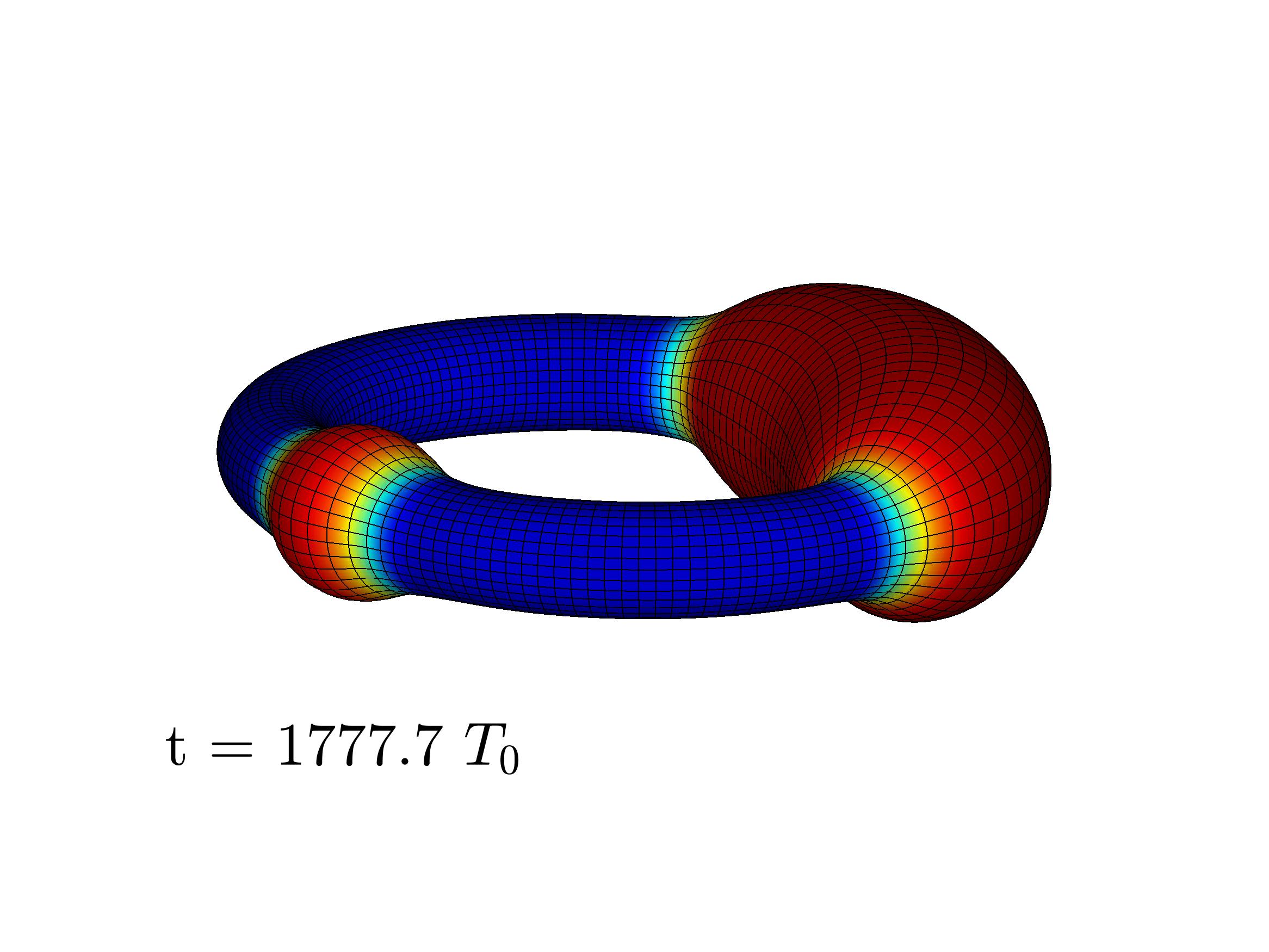}
\includegraphics[width=0.32\linewidth, trim = 300 310 200 400,clip]{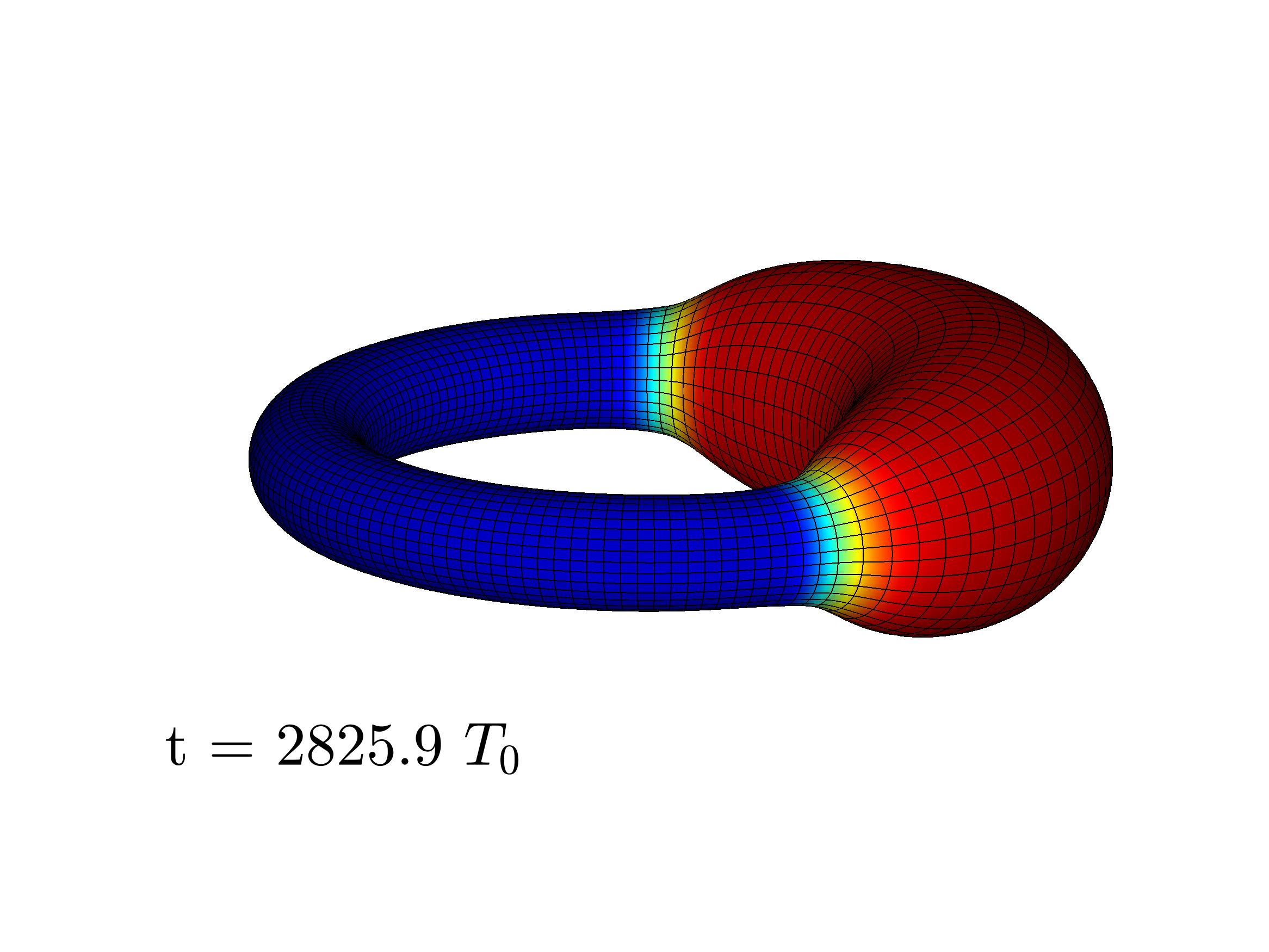}
\caption{Phase separation on a deforming torus: Evolution of the solution with $\lambda=0.0075\,L_0^2$ on a mesh containing $32 \times 128$ quadratic NURBS elements, with $\sqrt{\lambda}\approx 1.75\,h$. The colors follow Table~\ref{tab:mat}.}
\label{fig:t_evo2}
\end{figure}

For the second example, $\lambda = 0.0075\,L_0^2$ is selected and the constant internal pressure $p_{\mathrm{int}}=0.1\,EL_0^{-1}$ is prescribed for all $t$. 
The parameters listed in Table~\ref{tab:mat} and $D=4\,T_0$ are used. 
Fig.~\ref{fig:t_evo2} shows a series of snapshots of the evolution of the phase field on the deforming torus. 
Multiple bulges appear, evolve and merge during the phase separation process. 

\begin{figure}[H]
\centering
\includegraphics[width=0.49\linewidth, trim = 0 0 0 0,clip]{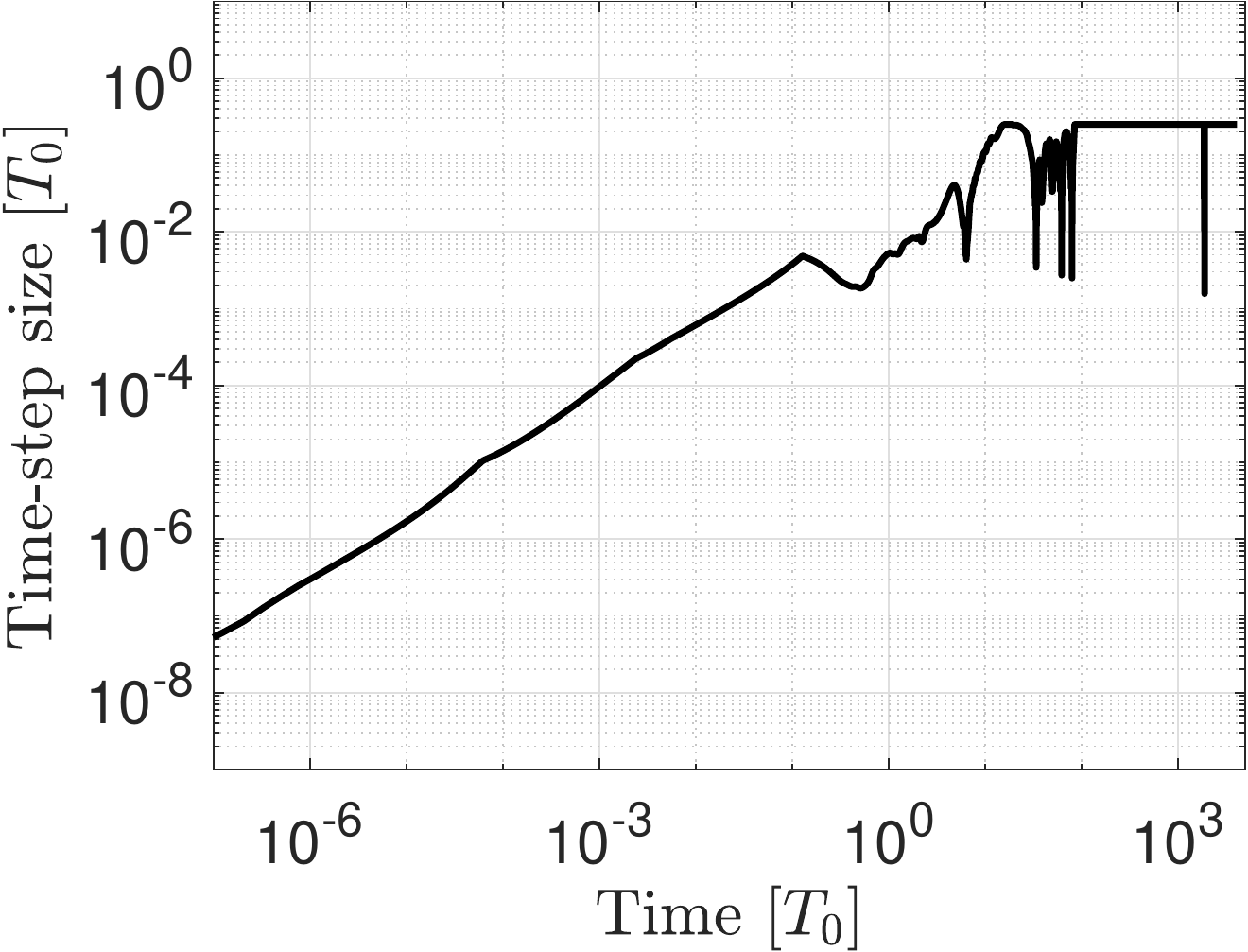}
\includegraphics[width=0.49\linewidth, trim = 0 0 0 0,clip]{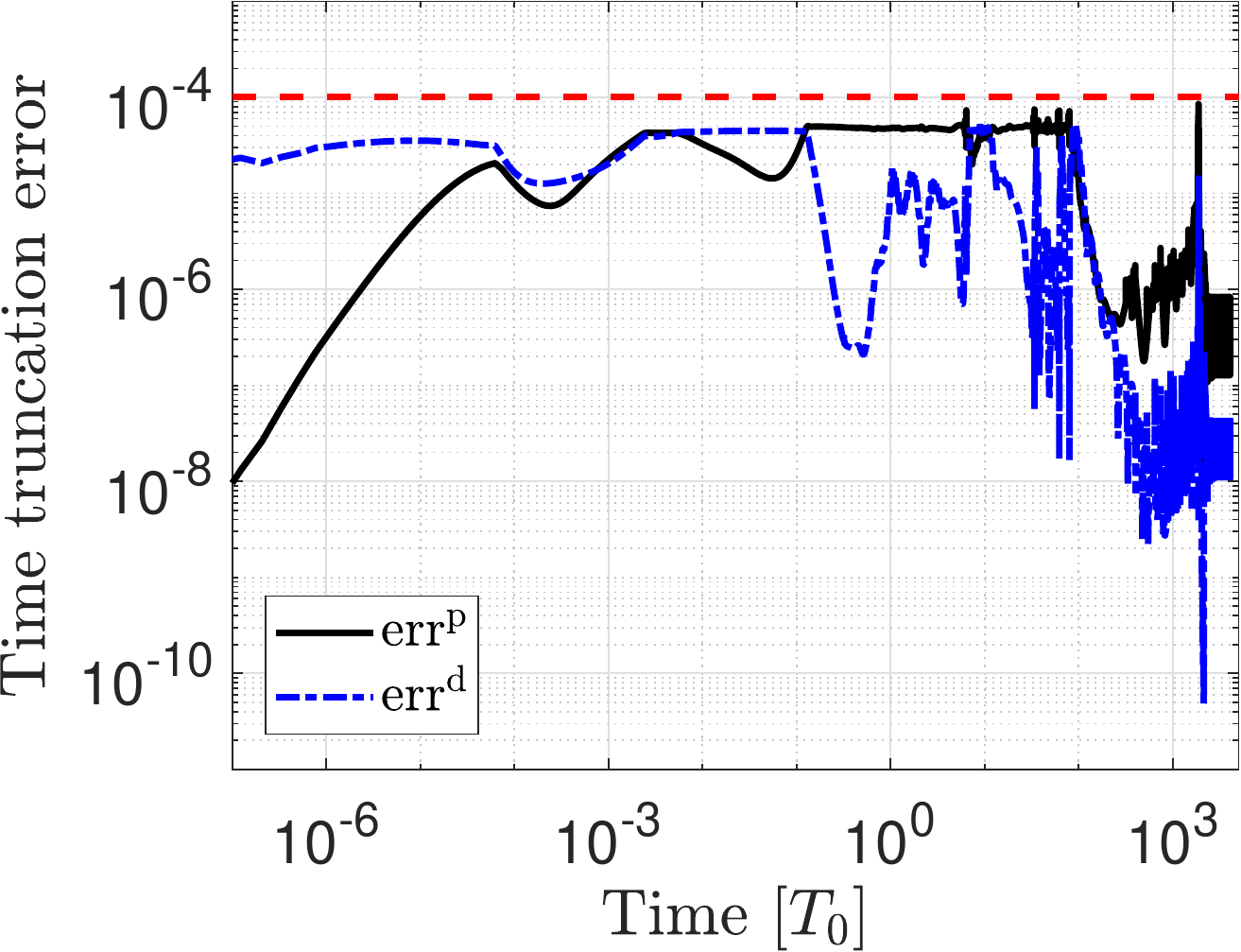}
\caption{Phase separation on a deforming torus: Left: Adaptive time step sizes. Right: Evolution of the local time truncation errors of the phase field, $\mathrm{err}^\mathrm{p}$, and mechanical field, $\mathrm{err}^\mathrm{d}$. The temporal error bound is shown by a dashed red line.}
\label{fig:t_ts2}
\end{figure}
\begin{figure}[H]
\centering
\includegraphics[width=0.49\linewidth, trim = 0 0 0 0,clip]{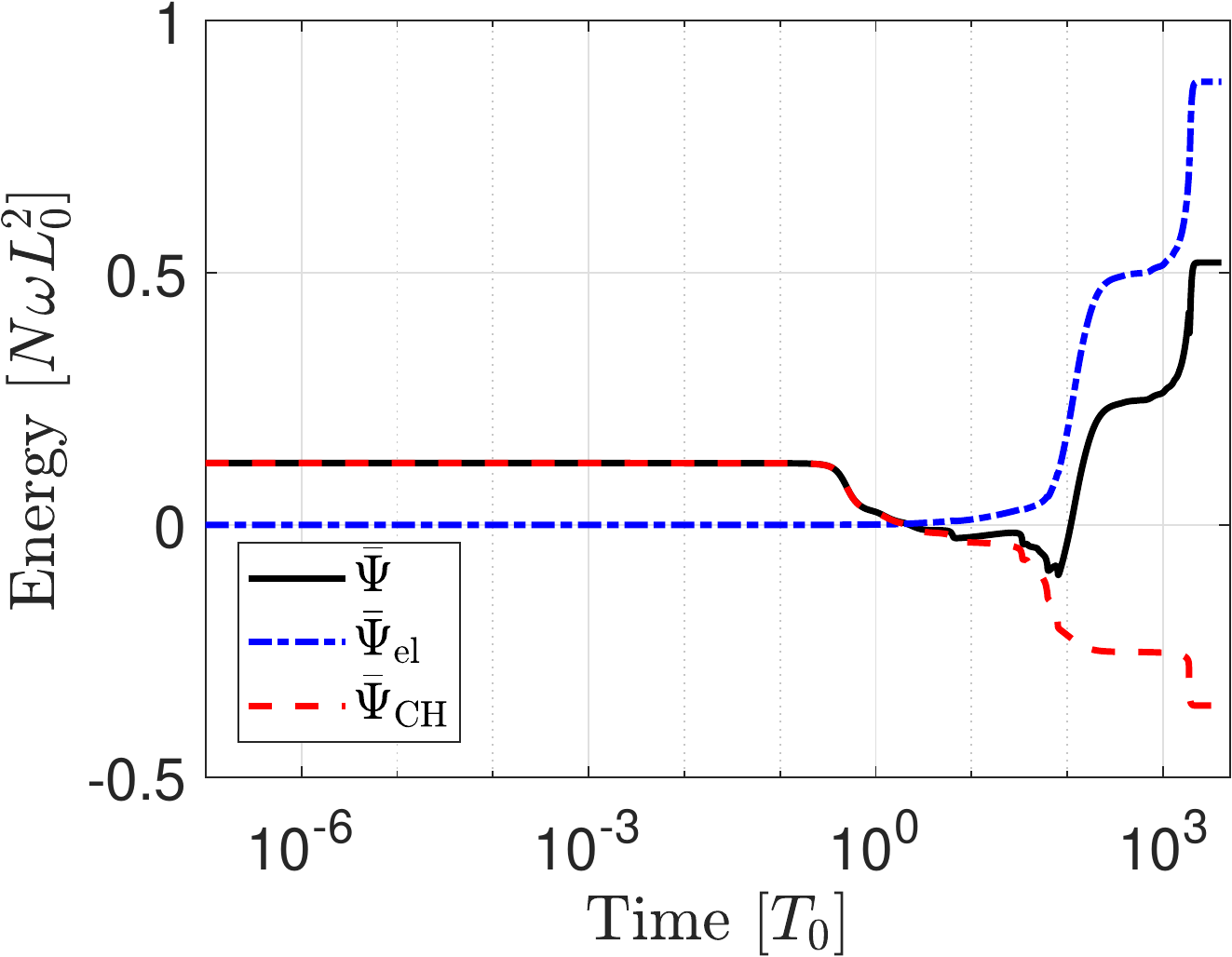}
\includegraphics[width=0.49\linewidth, trim = 0 0 0 0,clip]{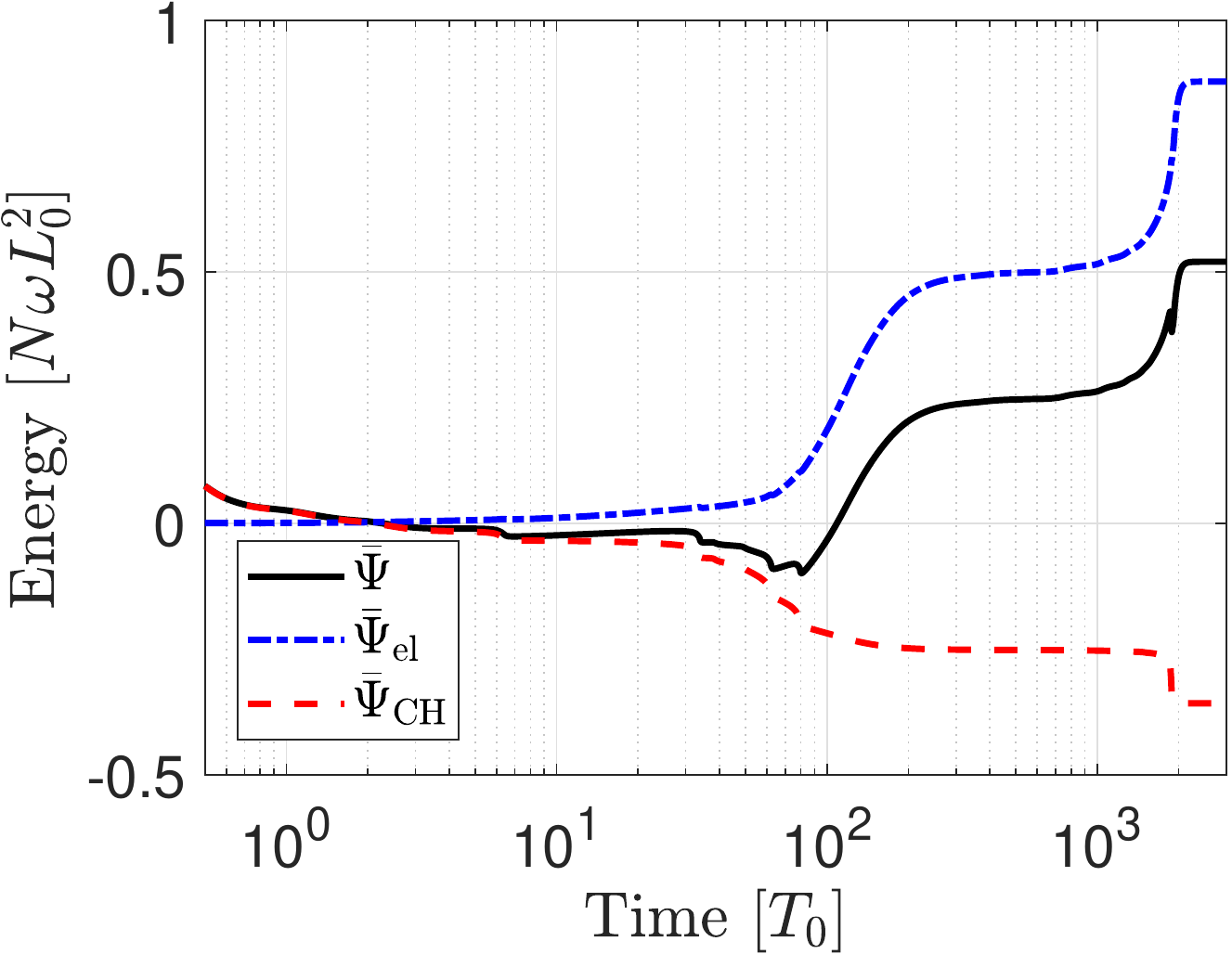}
\caption{Phase separation on a deforming torus: Evolution of the characteristic energies of the system. Close-up on the right side.}
\label{fig:t_ene2}
\end{figure}

Fig.~\ref{fig:t_ts2} shows the evolution of the time step size and the local truncation error. The time truncation error shows oscillations, which result in abrupt changes of the time step size. 
This reflects rapid changes and interactions of the phase field and the mechanical field. 
By choice, the time step size and the local time truncation error are limited to $t=0.25\,T_0$ and $10^{-4}$, respectively. Fig.~\ref{fig:t_ene2} shows the evolution of the characteristic energies of the system. 
The behavior is similar to the previous example (see Fig.~\ref{fig:t_ene1}).

\subsection{Phase separation on a deforming sphere}

\begin{figure}[H]
\centering
\includegraphics[width=0.24\linewidth, trim = 390 150 300 220,clip]{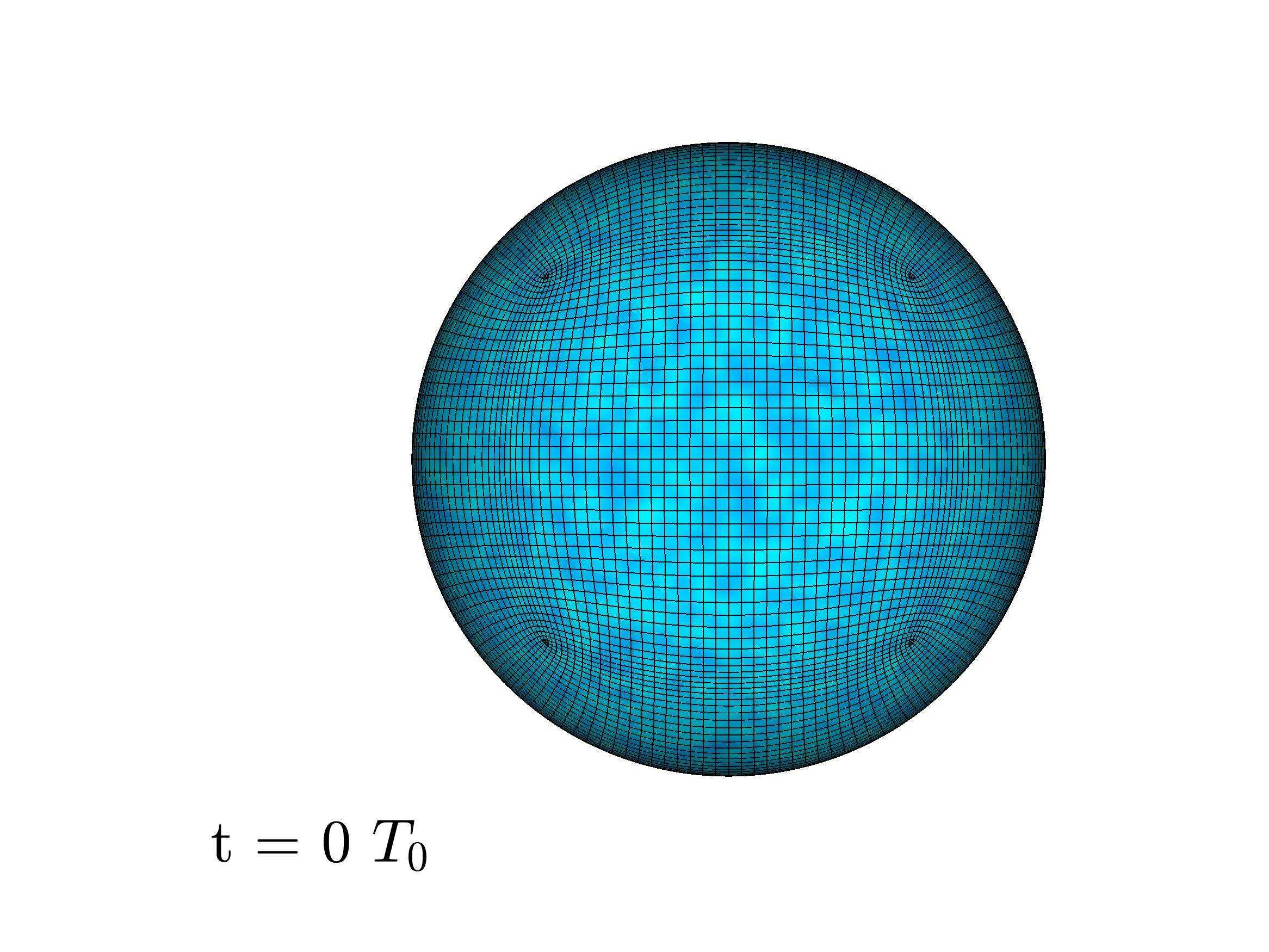}
\includegraphics[width=0.24\linewidth, trim = 390 150 300 220,clip]{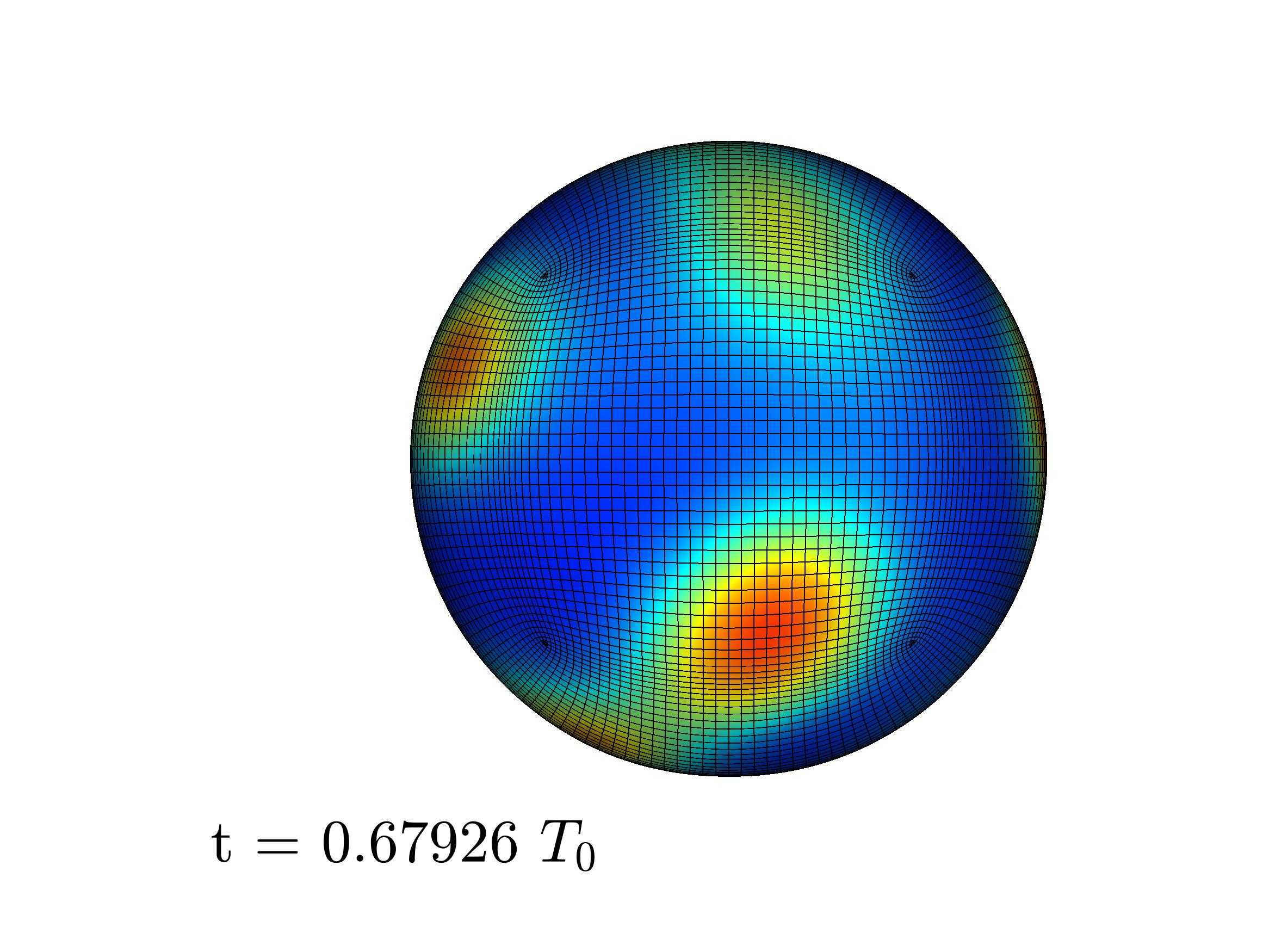}
\includegraphics[width=0.24\linewidth, trim = 390 150 300 220,clip]{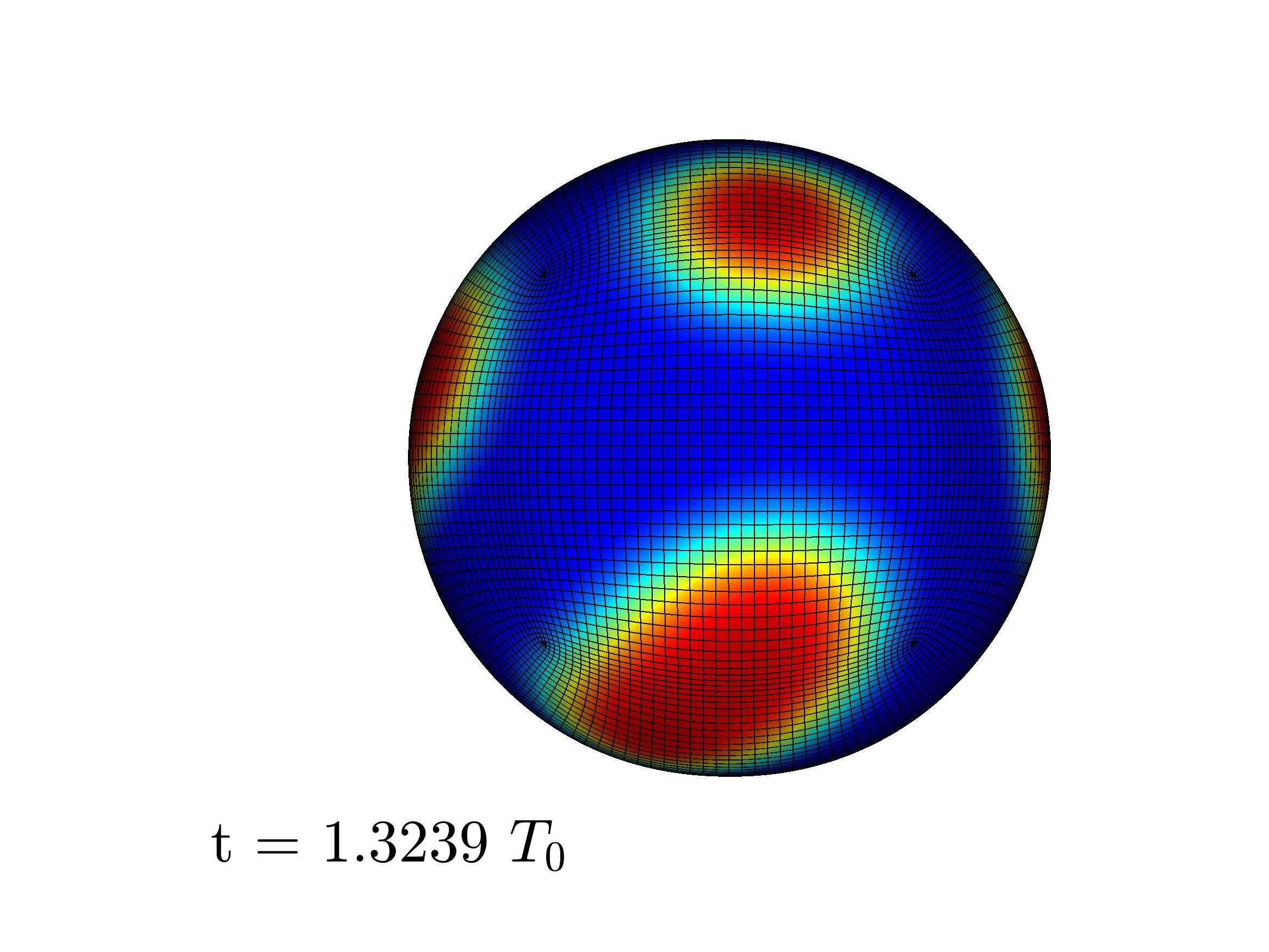}
\includegraphics[width=0.24\linewidth, trim = 390 150 300 220,clip]{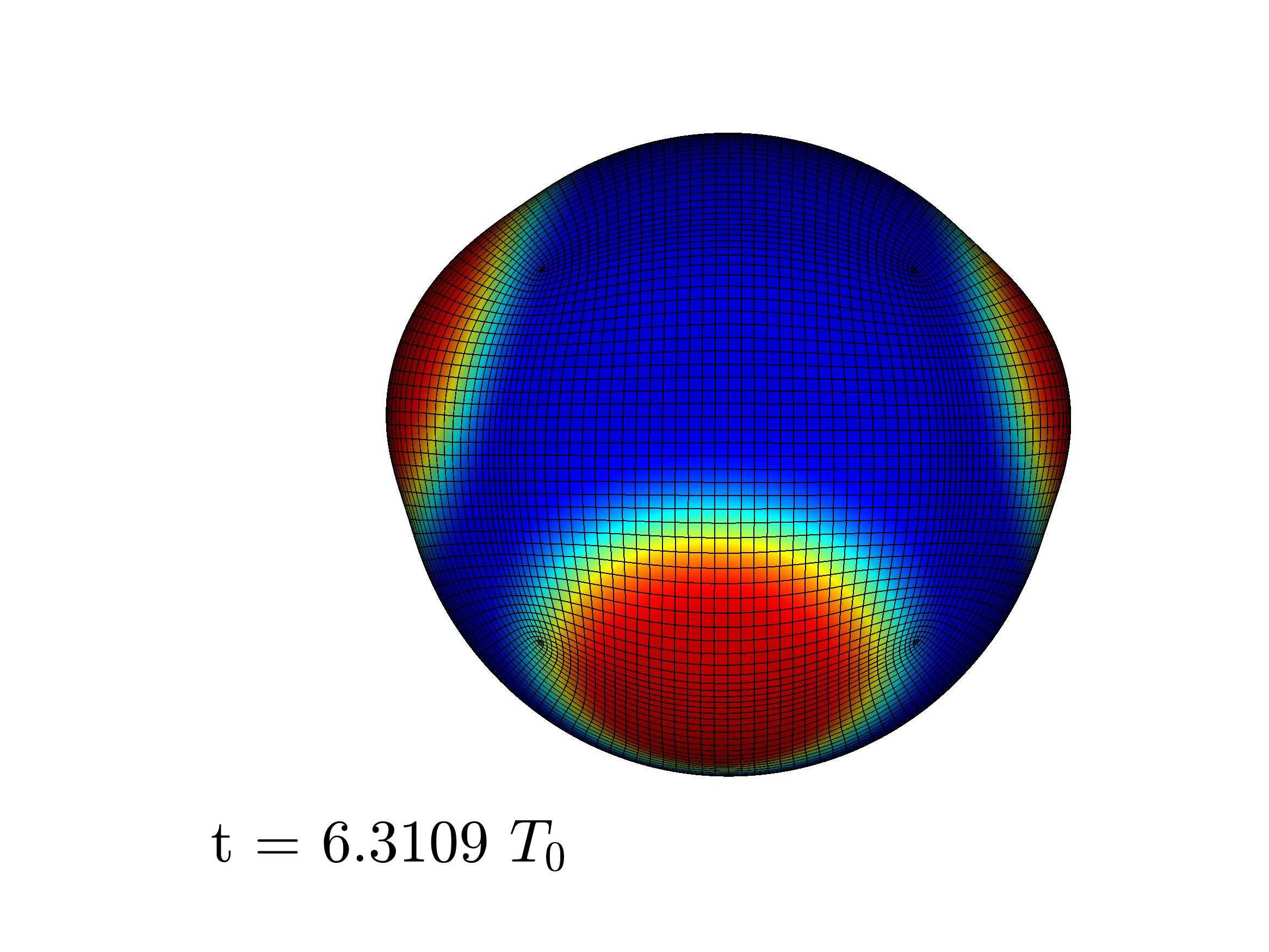}
\\
\includegraphics[width=0.24\linewidth, trim = 390 150 300 220,clip]{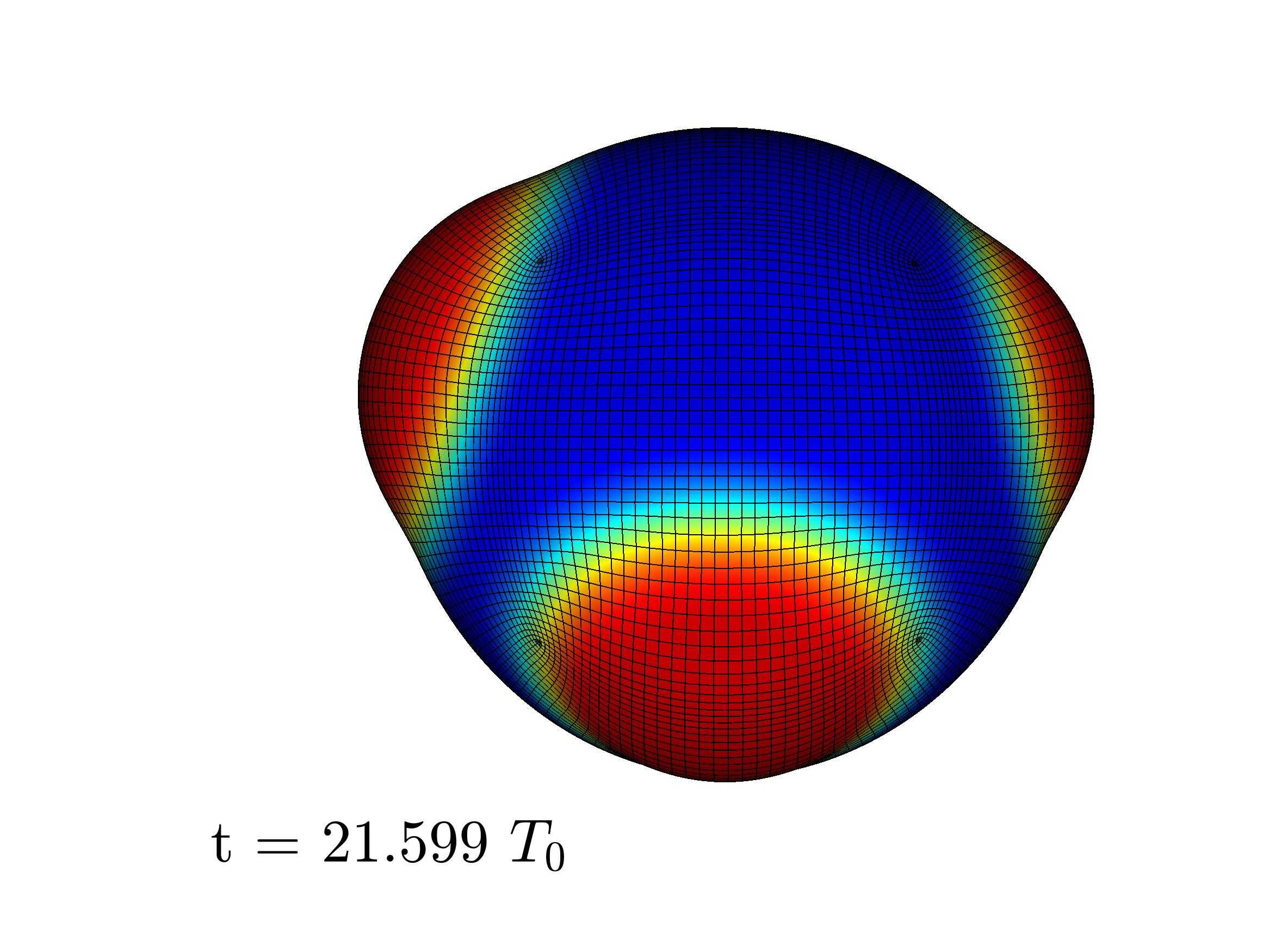}
\includegraphics[width=0.24\linewidth, trim = 390 150 300 220,clip]{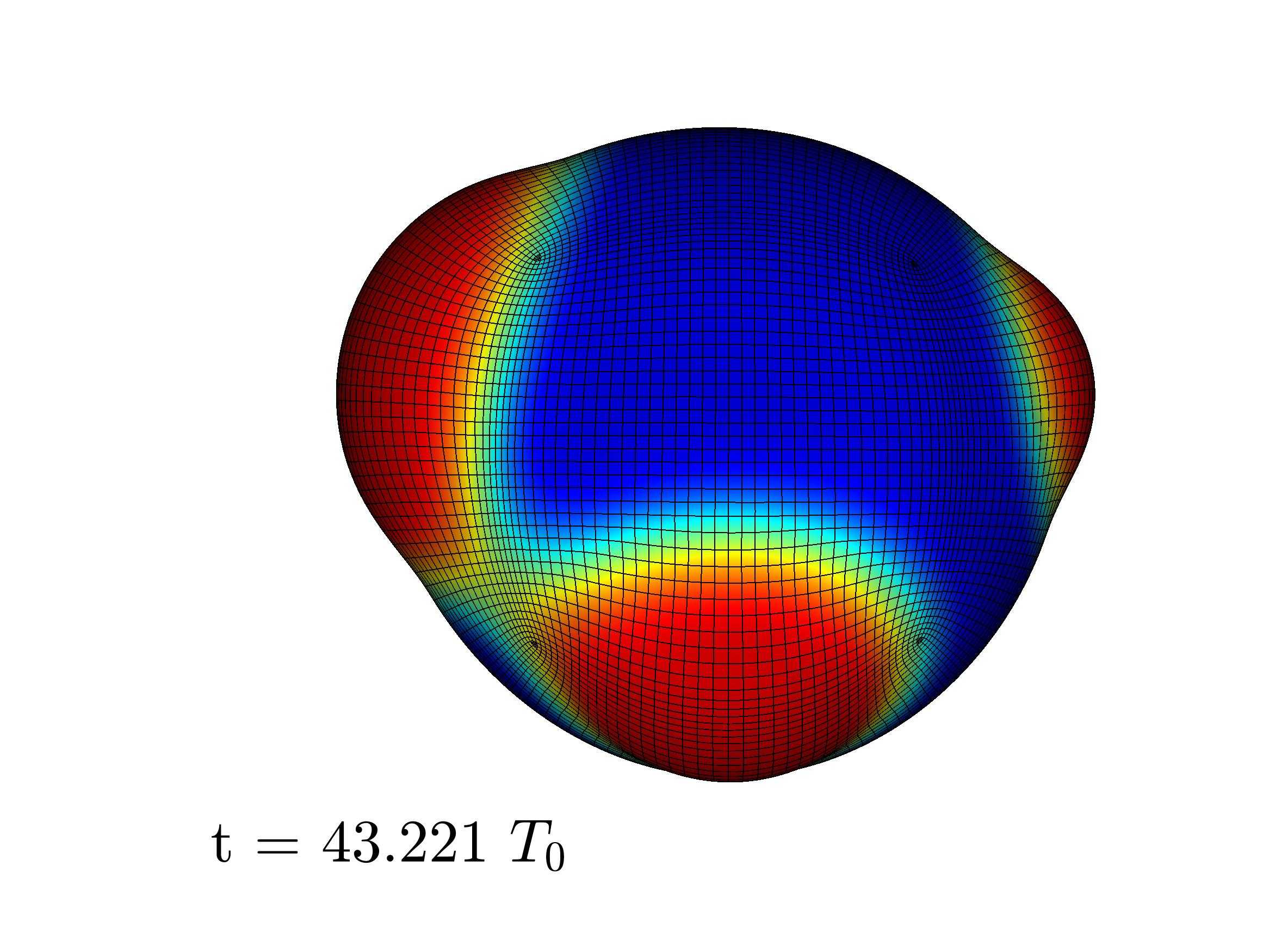}
\includegraphics[width=0.24\linewidth, trim = 390 150 300 220,clip]{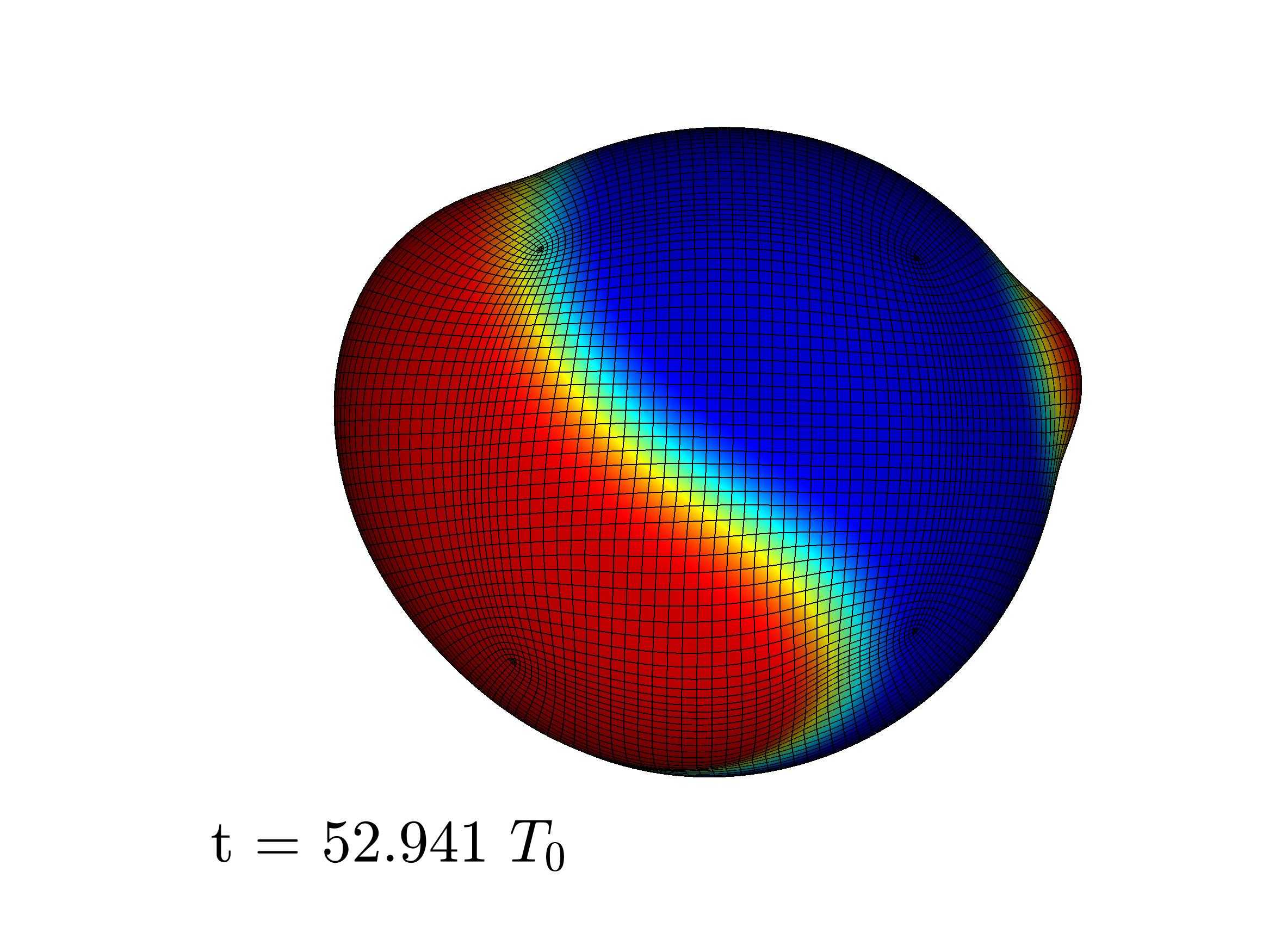}
\includegraphics[width=0.24\linewidth, trim = 390 150 300 220,clip]{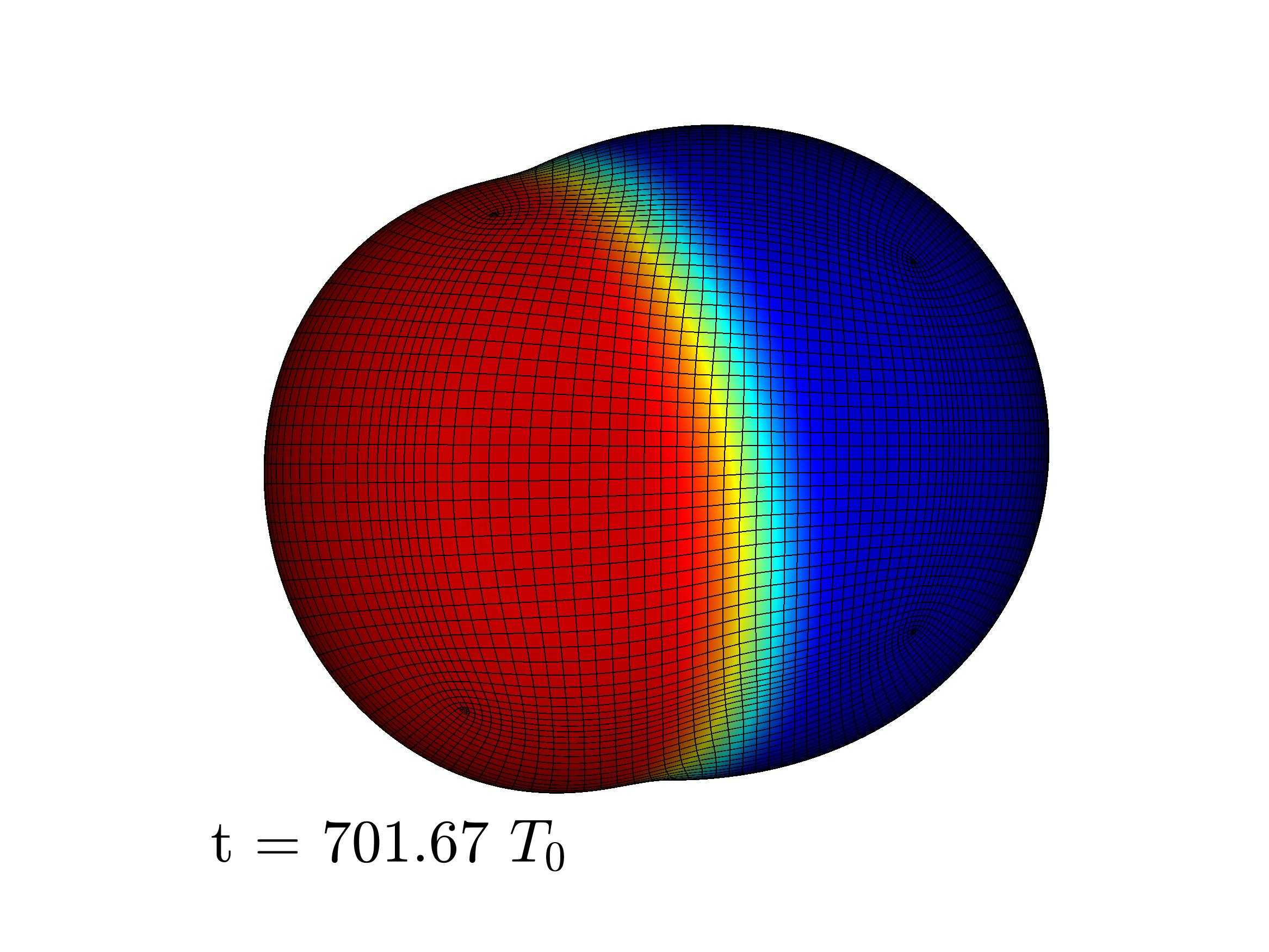}
\caption{Phase separation on a deforming sphere: Evolution of the solution with $\lambda=0.0075\,L_0^2$ on an unstructured mesh containing $9672$ cubic elements, with $\sqrt{\lambda}\approx 2.15\,h$. The colors follow Table~\ref{tab:mat}. See also supplemental movie file at \href{https://doi.org/10.5446/40802}{\texttt{https://doi.org/10.5446/40802}}.}
\label{fig:s_evo2}
\end{figure}

The third example studies phase separation on a deforming sphere that is discretized by the unstructured splines from Sec.~\ref{sec:spaceDisc_splines}. 
The parameters $D=4\,T_0$ and $\lambda=0.0075\,L_0^2$ are used together with those in Table~\ref{tab:mat}.
The constant internal pressure $p_{\mathrm{int}}=0.0425\,EL_0^{-1}$ is prescribed for all $t$. 
The unstructured mesh consists of $9672$ cubic elements and has $8$ extraordinary points.
This mesh provides $C^2$-continuity except for the extraordinary points that are only $C^1$-continuous. 
Rigid body deformations are prevented by analogous boundary conditions to those shown in Fig.~\ref{fig:t_IC1}. 
Fig.~\ref{fig:s_evo2} shows a series of snapshots of the phase separation on the deforming sphere. 
Multiple red phase nuclei appear, bulge, evolve and merge during phase separation.  
As the nuclei grow the deformations become larger.

\begin{figure}[H]
\centering
\includegraphics[width=0.49\linewidth, trim = 0 0 0 0,clip]{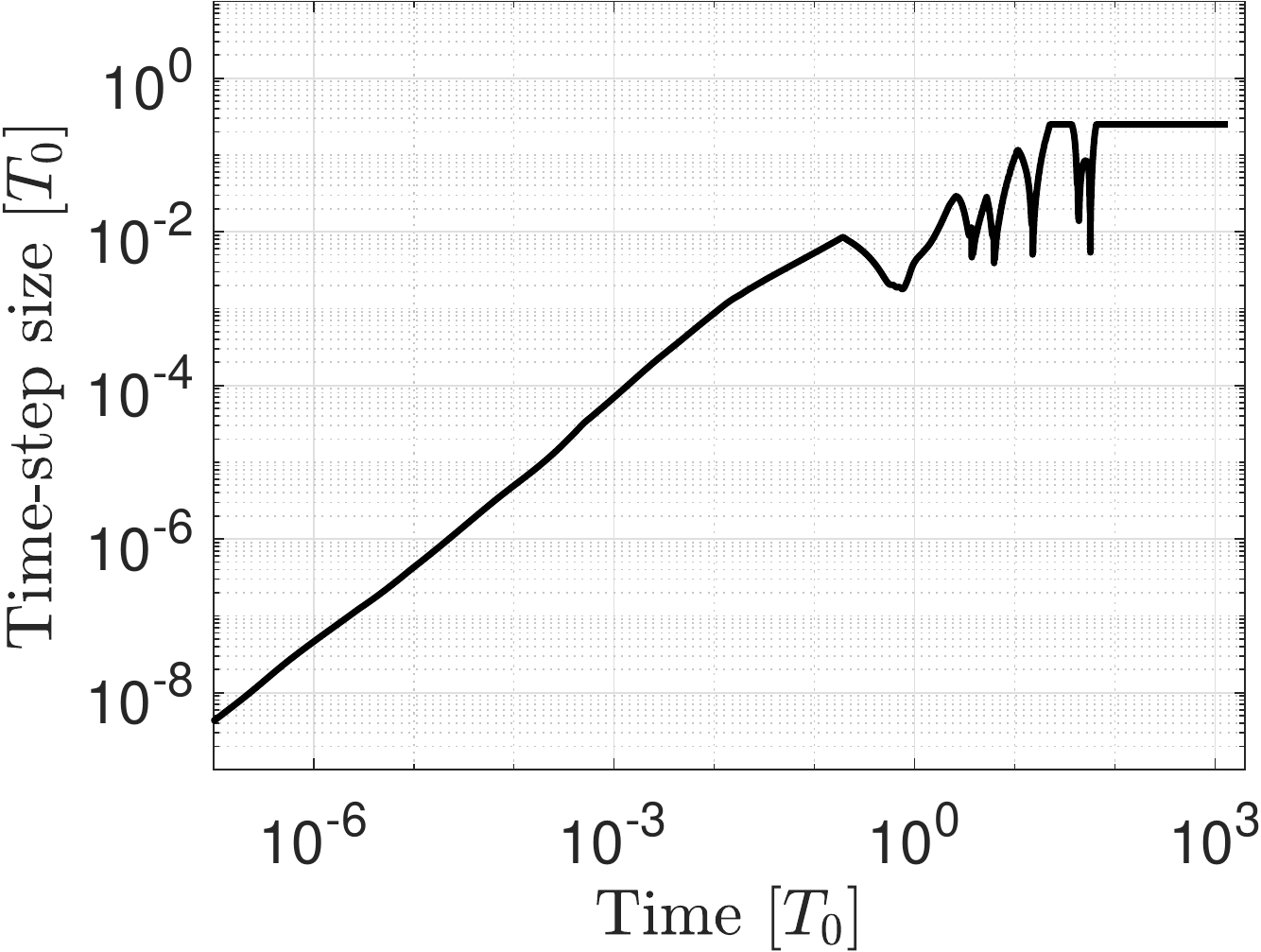}
\includegraphics[width=0.49\linewidth, trim = 0 0 0 0,clip]{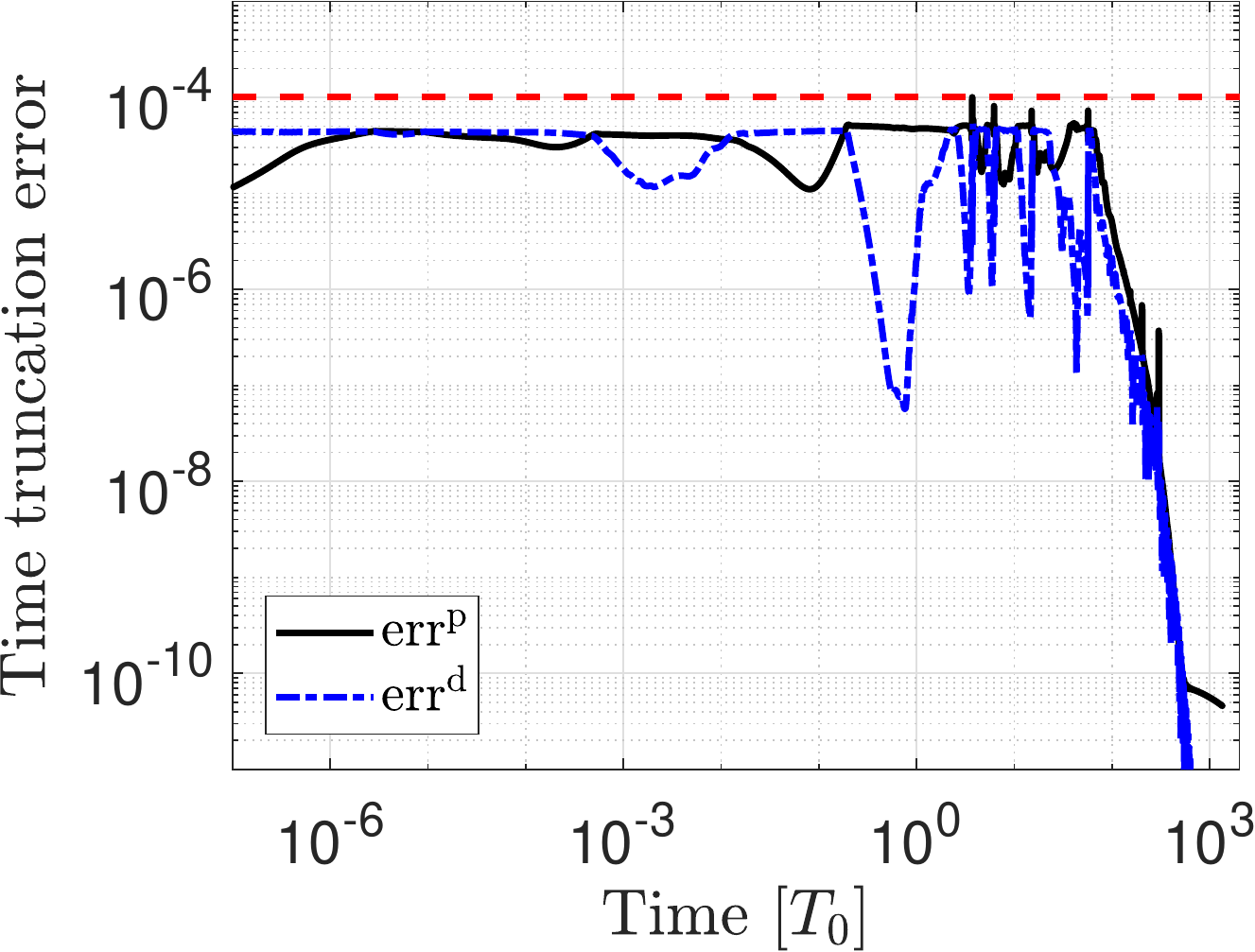}
\caption{Phase separation on a deforming sphere: Left: Adaptive time step size. Right: Evolution of the local time truncation errors of the phase field, $\mathrm{err}^\mathrm{p}$, and mechanical field, $\mathrm{err}^\mathrm{d}$. The chosen temporal error bound is shown by a dashed red line.}
\label{fig:s_ts2}
\end{figure}

\begin{figure}[H]
\centering
\includegraphics[width=0.49\linewidth, trim = 0 0 0 0,clip]{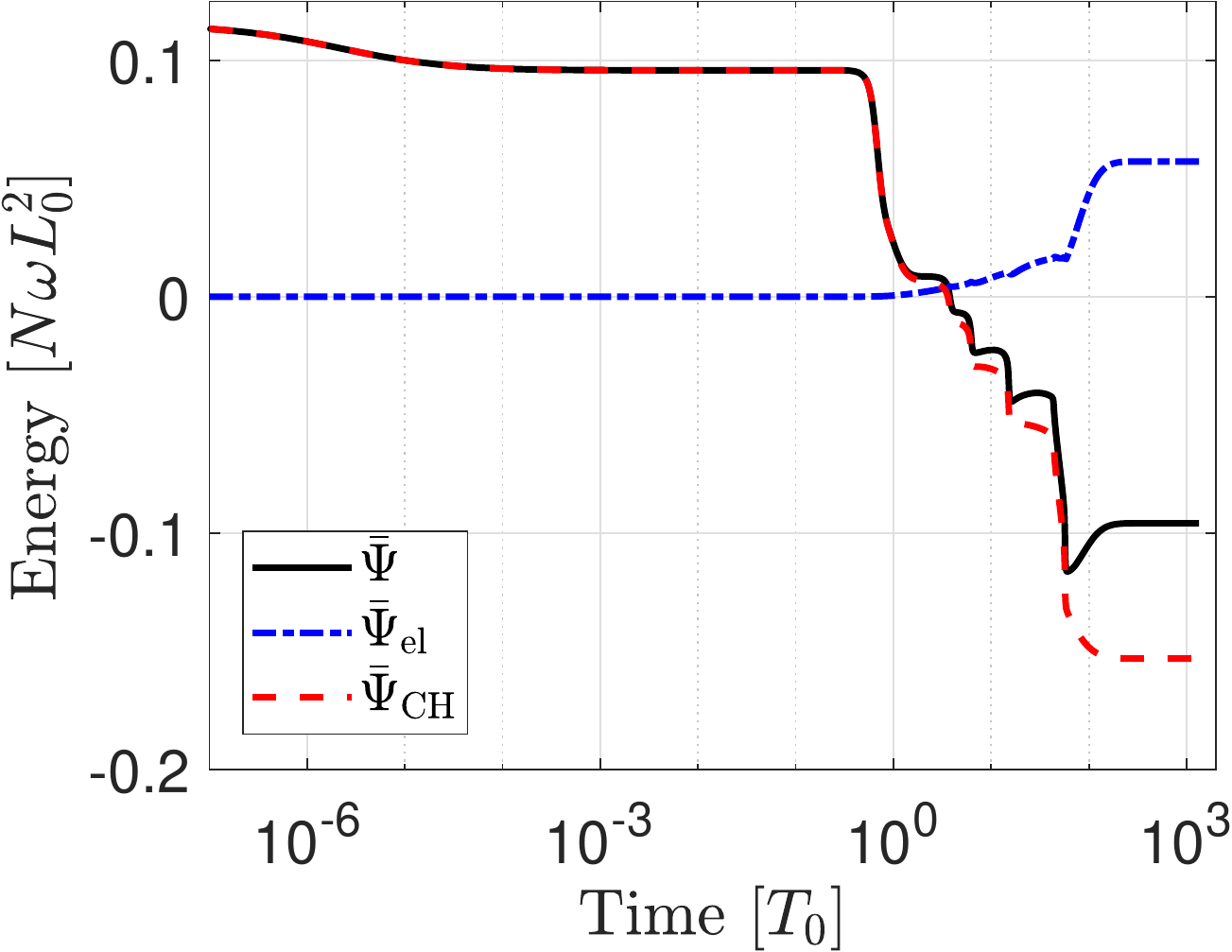}
\includegraphics[width=0.49\linewidth, trim = 0 0 0 0,clip]{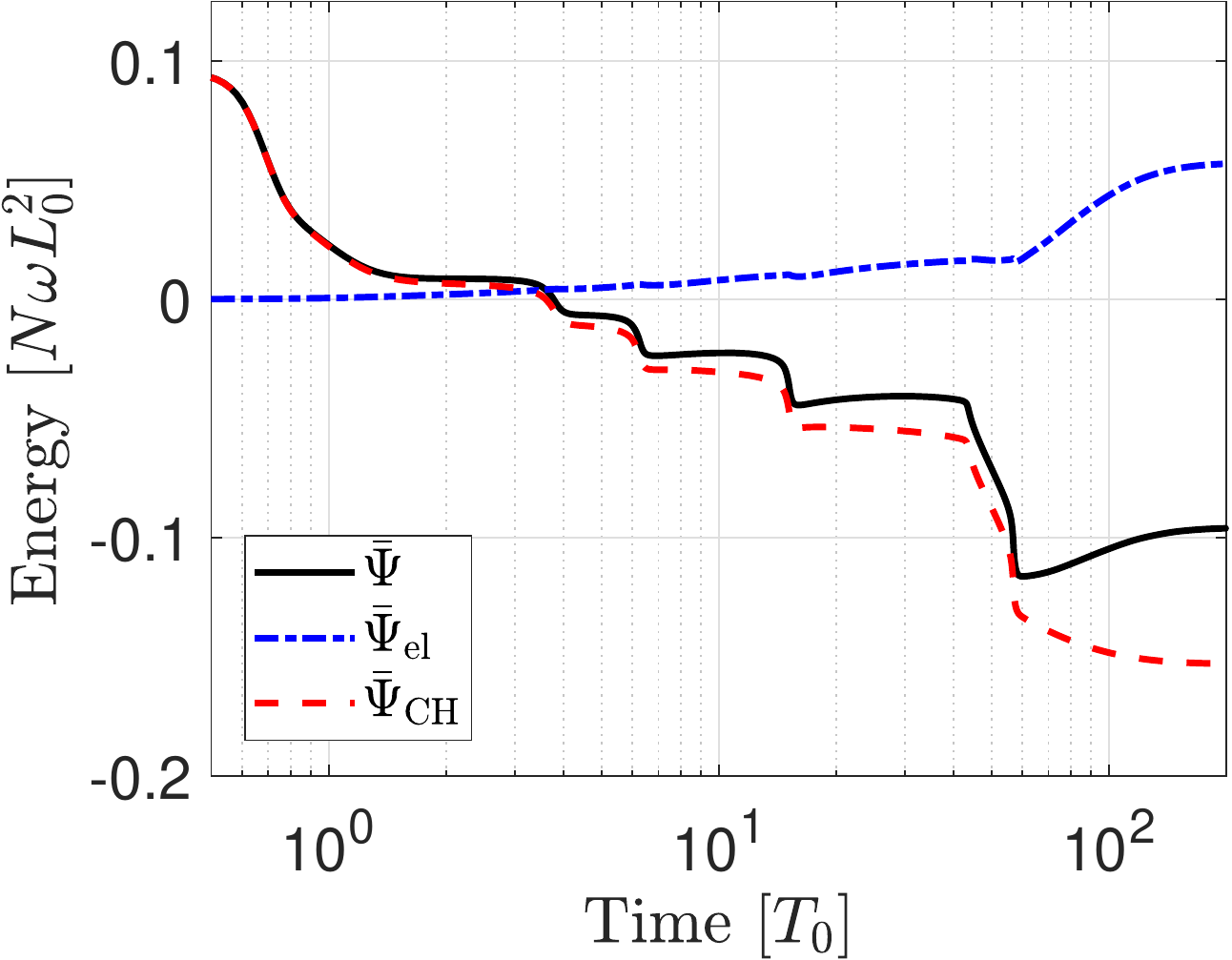}
\caption{Phase separation on a deforming sphere: Evolution of the characteristic energies of the system. Close-up on the right side.}
\label{fig:s_ene2}
\end{figure}

\begin{figure}[H]
\centering
\begin{overpic}[width=0.21\linewidth, trim = 300 145 760 538,clip]{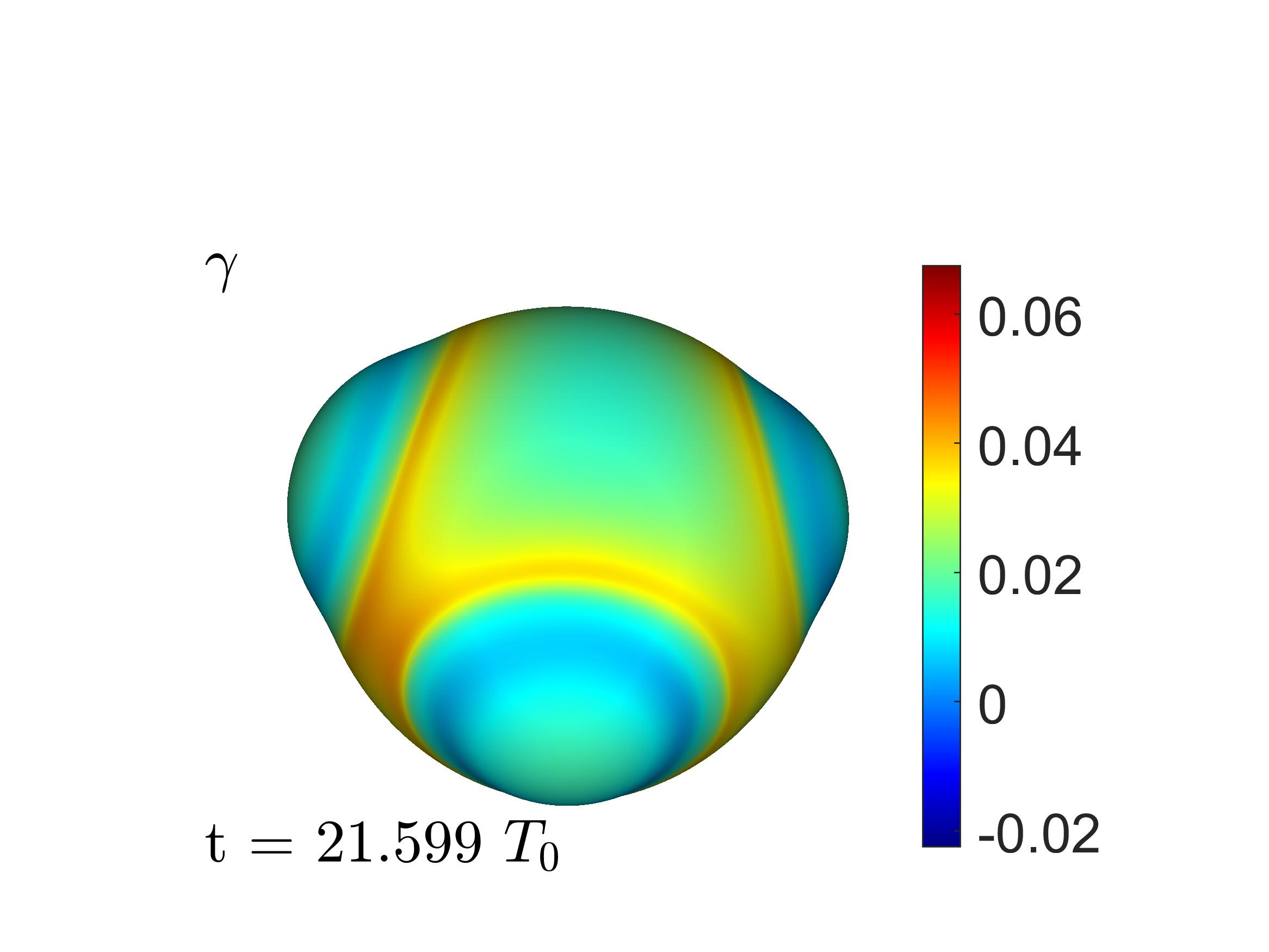}
\put(0,0){\crule[white]{2.1cm}{0.275cm}}
\put(5,80){$\gamma$}
\put(5,0){$t=21.599\,T_0$}
\end{overpic}
\begin{overpic}[width=0.21\linewidth, trim = 300 145 760 538,clip]{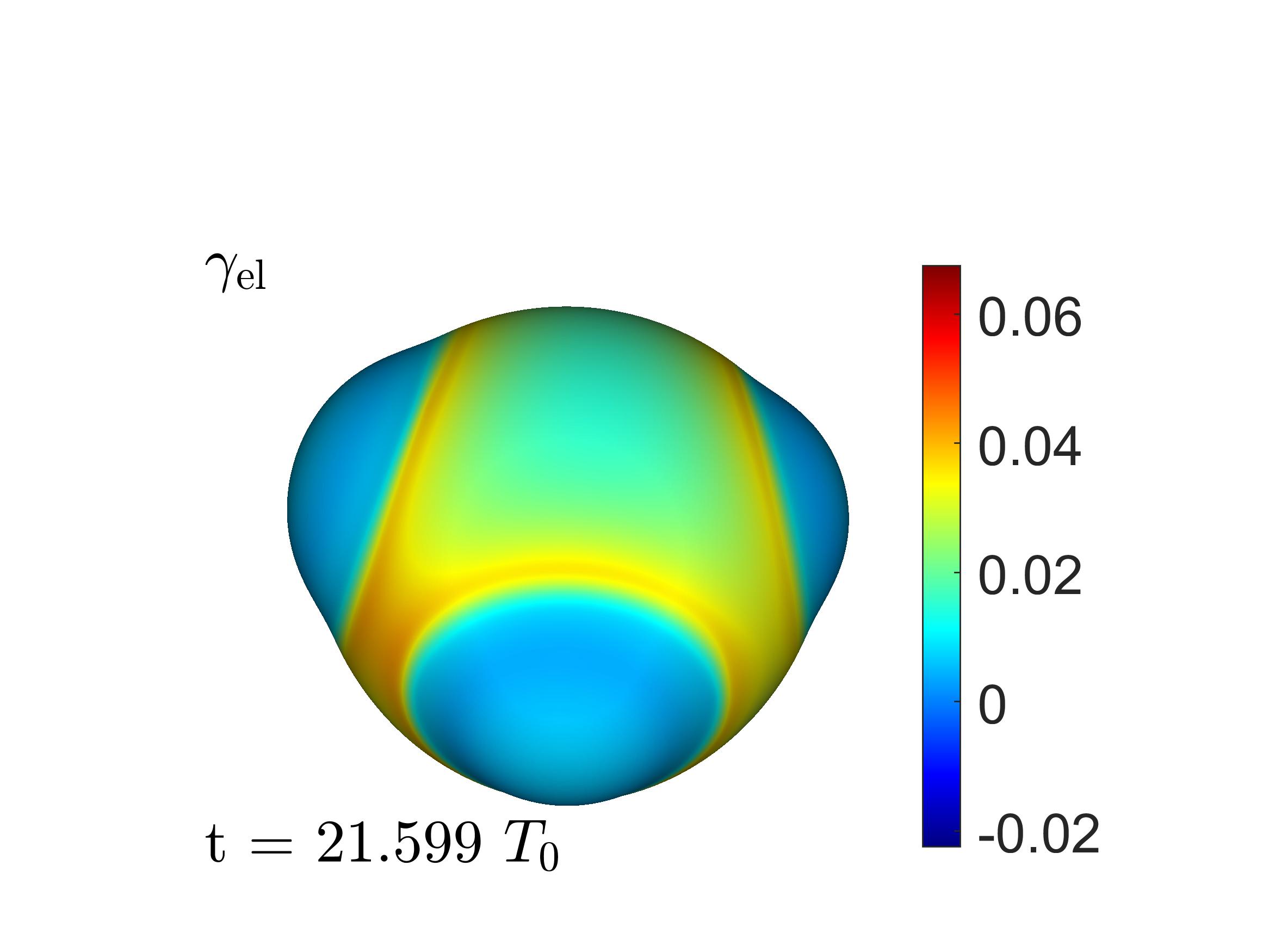}
\put(0,0){\crule[white]{2.1cm}{0.275cm}}
\put(5,80){$\gamma_\mathrm{el}$}
\put(5,0){$t=21.599\,T_0$}
\end{overpic}
\begin{overpic}[width=0.21\linewidth, trim = 300 145 760 538,clip]{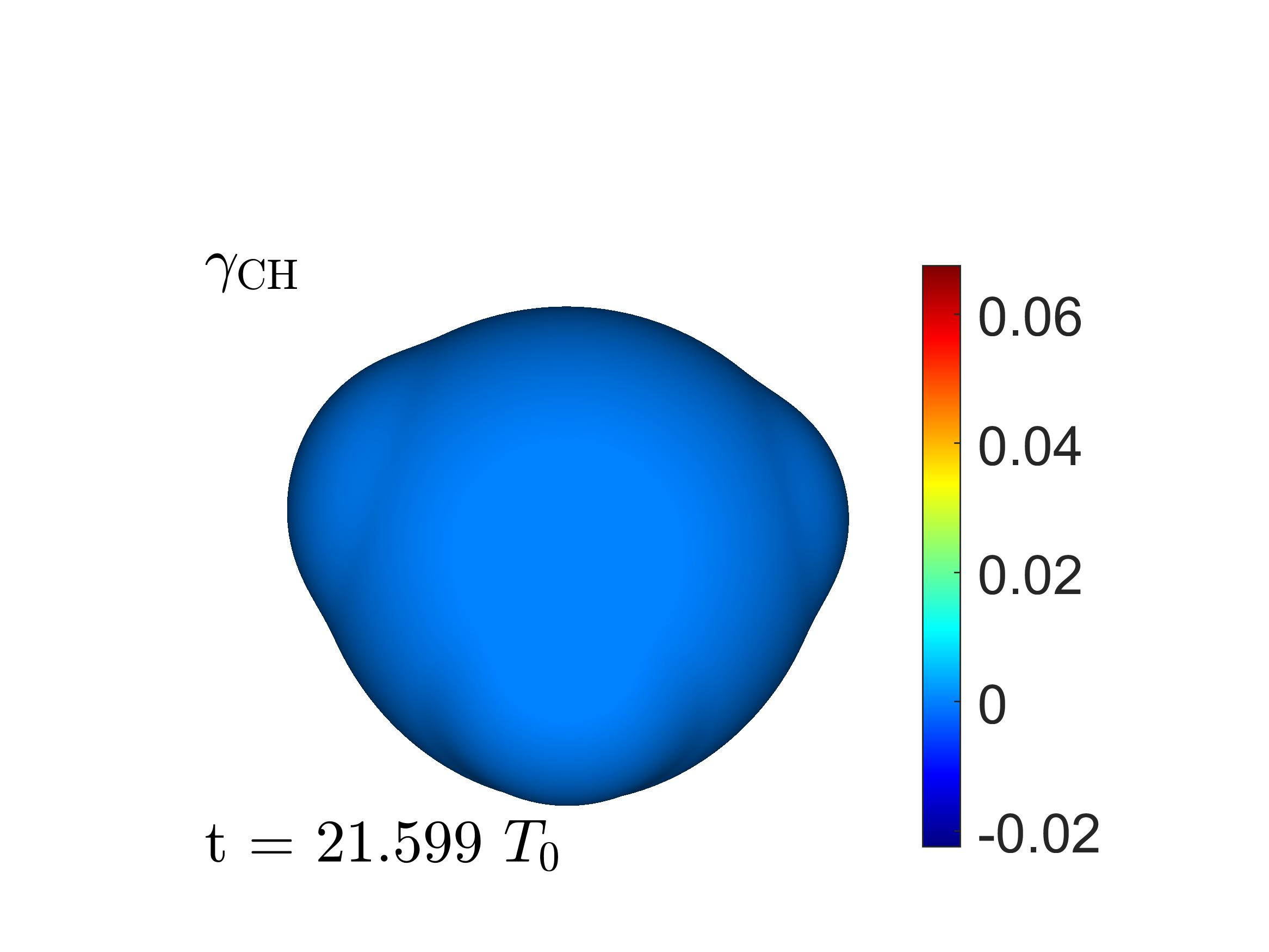}
\put(0,0){\crule[white]{2.1cm}{0.275cm}}
\put(5,80){$\gamma_\mathrm{CH}$}
\put(5,0){$t=21.599\,T_0$}
\end{overpic}
\begin{overpic}[width=0.21\linewidth, trim = 300 145 760 538,clip]{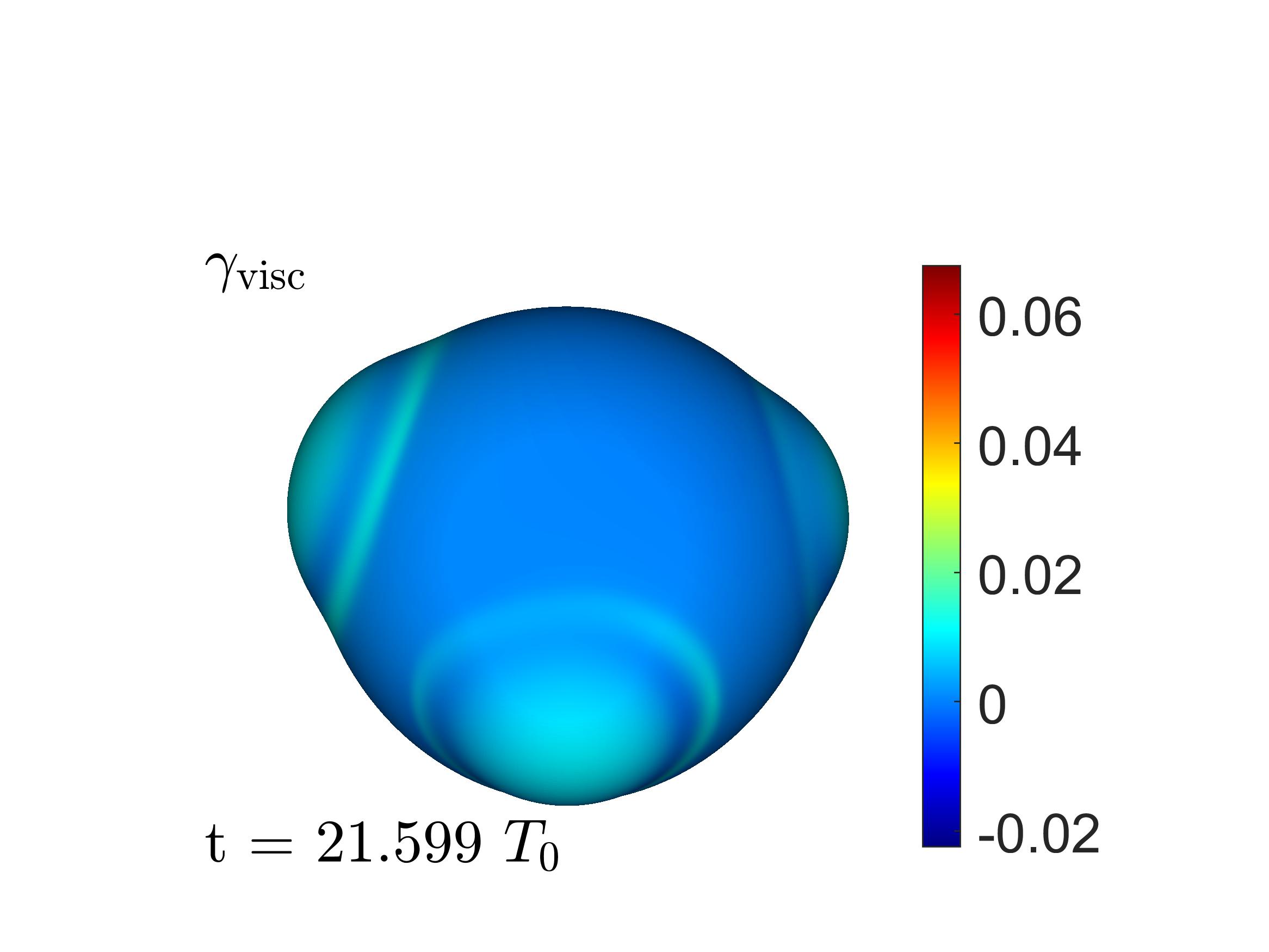}
\put(0,0){\crule[white]{2.1cm}{0.275cm}}
\put(5,80){$\gamma_\mathrm{visc}$}
\put(5,0){$t=21.599\,T_0$}
\end{overpic}
\includegraphics[width=0.064\linewidth, trim = 1595 140 300 430,clip]{ifigs/coup/S/ex_r2/gtot_m8_1650_ts}
\caption{Phase separation on a deforming sphere: Surface tensions $\gamma$, $\gamma_\mathrm{el}$, $\gamma_\mathrm{CH}$ and $\gamma_\mathrm{visc}$ at $t=21.599\,T_0$. The colors show the surface tensions in the units $[N\omega]$. See also supplemental movie file at \href{https://doi.org/10.5446/40803}{\texttt{https://doi.org/10.5446/40803}}.}
\label{fig:s_gam2}
\end{figure}

\begin{figure}[H]
\centering
\begin{overpic}[width=0.21\linewidth, trim = 312 145 715 538,clip]{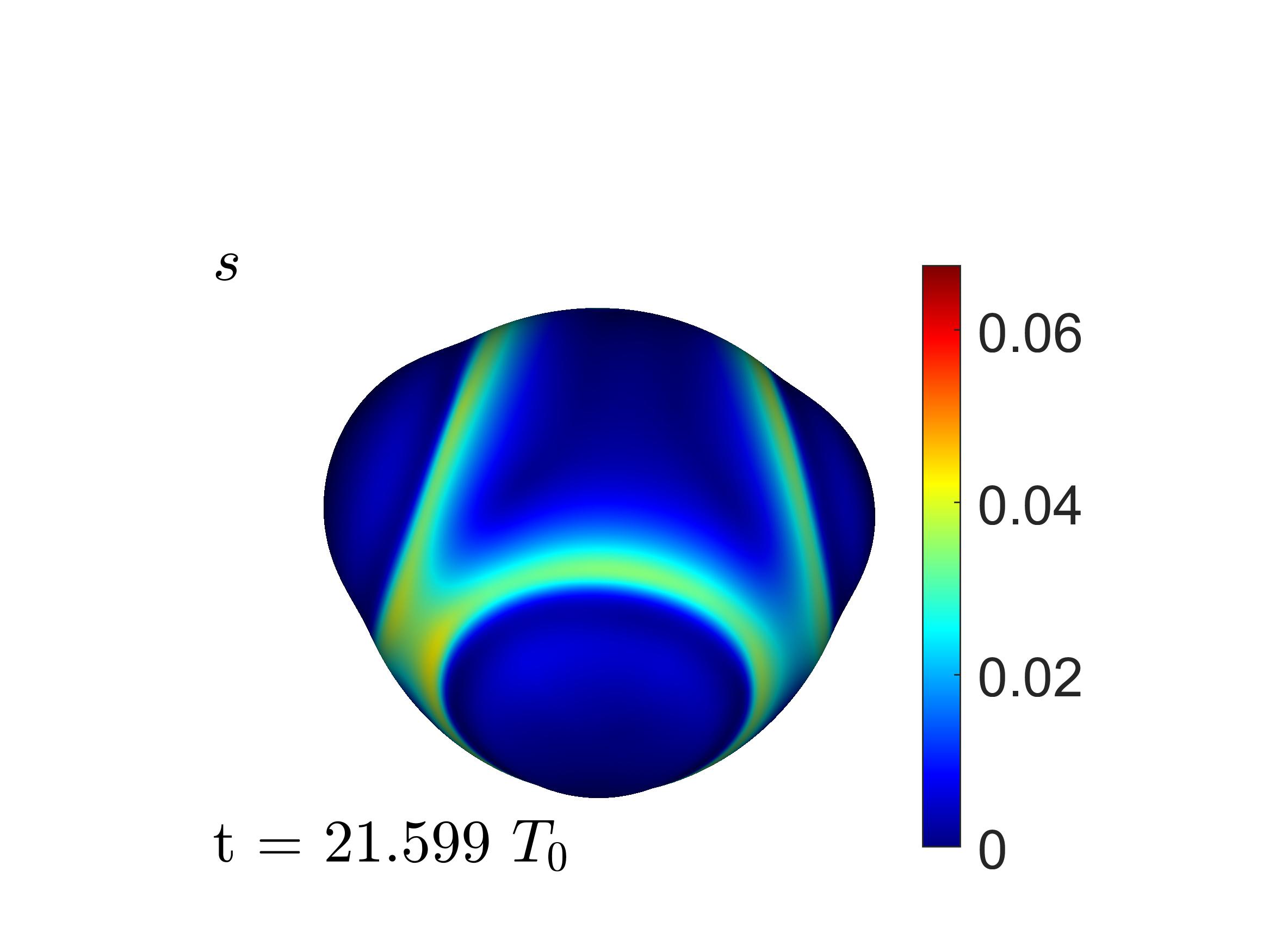}
\put(0,0){\crule[white]{2.1cm}{0.275cm}}
\put(5,80){$s$}
\put(5,0){$t=21.599\,T_0$}
\end{overpic}
\begin{overpic}[width=0.21\linewidth, trim = 312 145 715 538,clip]{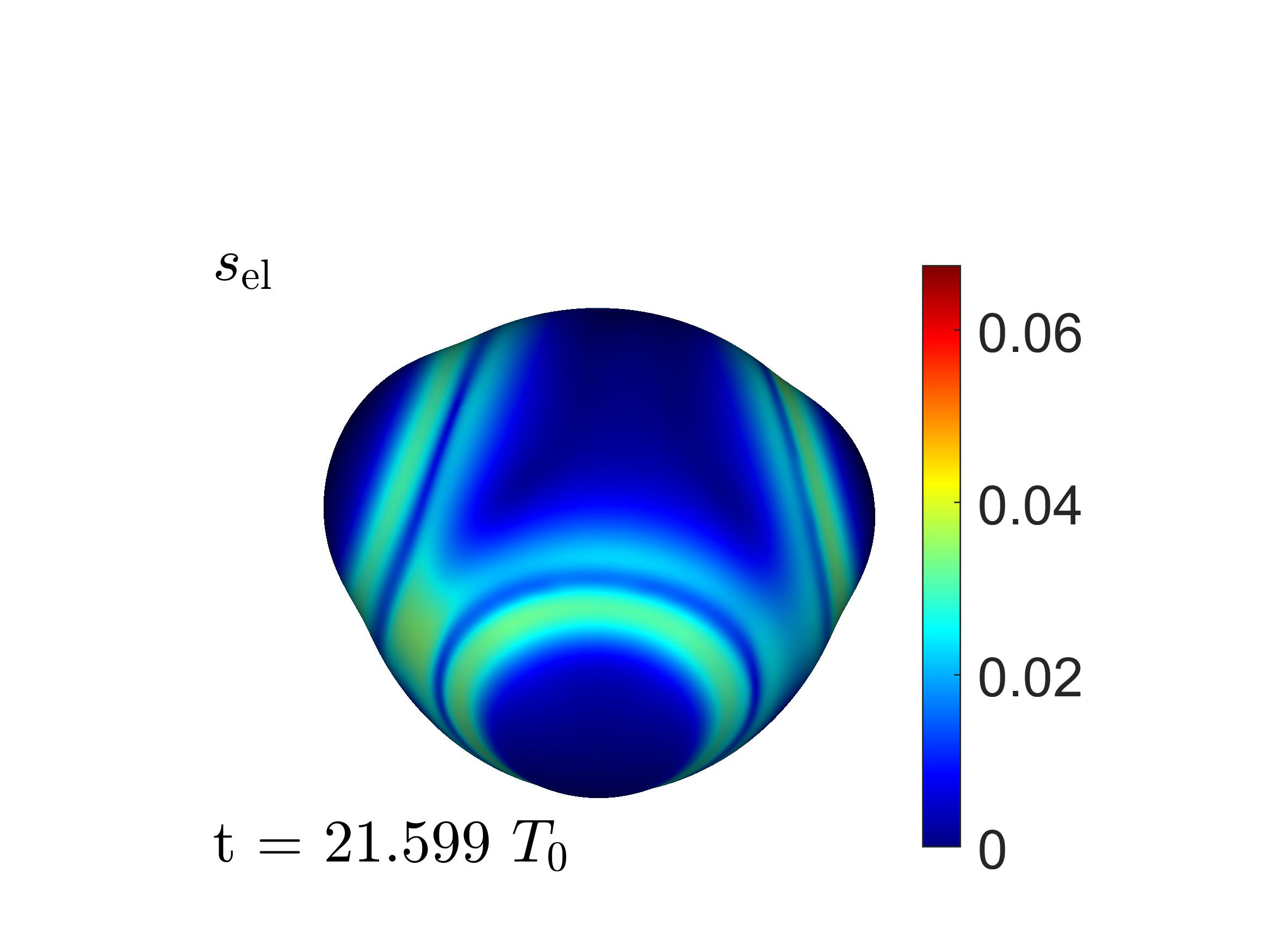}
\put(0,0){\crule[white]{2.1cm}{0.275cm}}
\put(5,80){$s_\mathrm{el}$}
\put(5,0){$t=21.599\,T_0$}
\end{overpic}
\begin{overpic}[width=0.21\linewidth, trim = 312 145 715 538,clip]{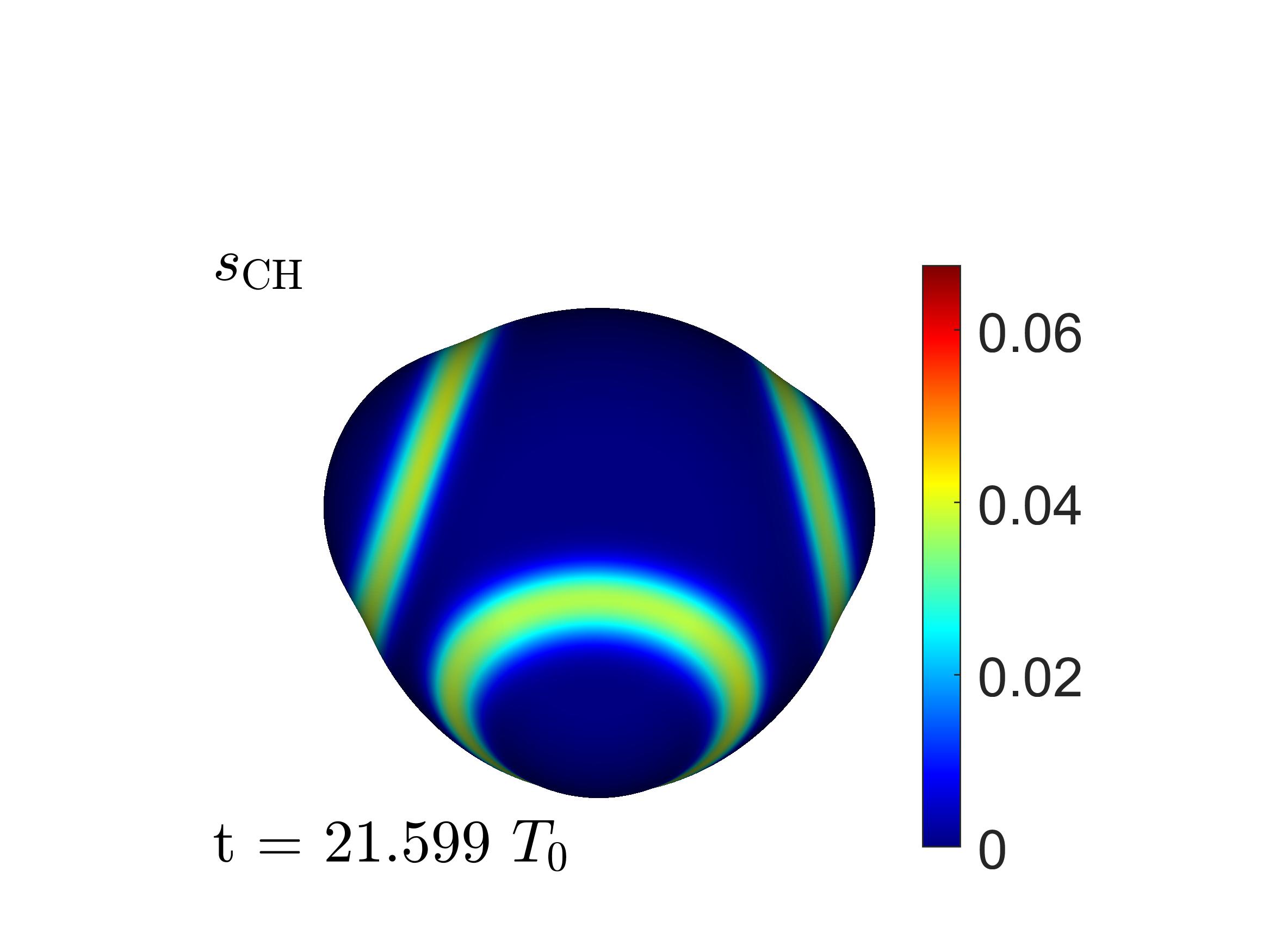}
\put(0,0){\crule[white]{2.1cm}{0.275cm}}
\put(5,80){$s_\mathrm{CH}$}
\put(5,0){$t=21.599\,T_0$}
\end{overpic}
\begin{overpic}[width=0.21\linewidth, trim = 312 145 715 538,clip]{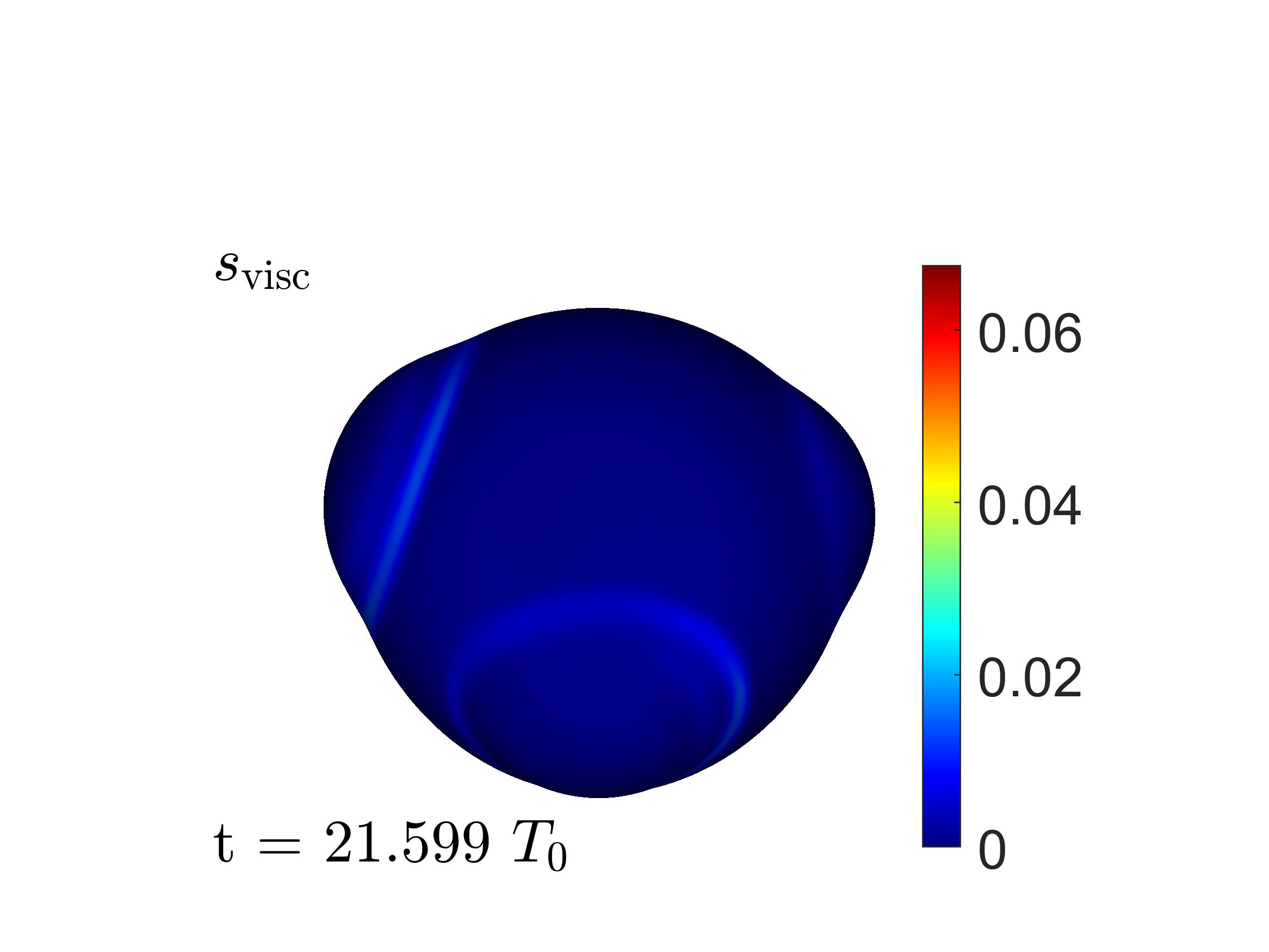}
\put(0,0){\crule[white]{2.1cm}{0.275cm}}
\put(5,80){$s_\mathrm{visc}$}
\put(5,0){$t=21.599\,T_0$}
\end{overpic}
\includegraphics[width=0.060\linewidth, trim = 1620 140 300 430,clip]{ifigs/coup/S/ex_r2/stot_m8_1650_ts}
\caption{Phase separation on a deforming sphere: Stresses $s$, $s_\mathrm{el}$, $s_\mathrm{CH}$ and $s_\mathrm{visc}$ at $t=21.599\,T_0$. The colors show the stresses in the units $[N\omega]$. See also supplemental movie file at \href{https://doi.org/10.5446/40804}{\texttt{https://doi.org/10.5446/40804}}.}
\label{fig:s_str2}
\end{figure}

The left side of Fig.~\ref{fig:s_ts2} shows the evolution of the time step size that results from the adaptive time-stepping procedure of Sec.~\ref{Sec:ATS}. 
The right side of Fig.~\ref{fig:s_ts2} shows the local time truncation errors $\text{err}^\mrp$ and $\text{err}^\mrd$. 
The time step size show similar characteristics as observed in the previous example. 
The maximum time step size is limited to $\Delta t= 0.25\,T_0$ in this example. Fig.~\ref{fig:s_ene2} shows the evolution of the characteristic energies of the system. 
Like before, in the beginning, $\bar\Psi_\mathrm{CH}$ is largest.
Later, $\bar\Psi_\mathrm{CH}$ decreases, while $\bar\Psi_\mathrm{el}$ increases. 
Next, we illustrate and compare two stress measures: 
The surface tension
\eqb{l}
\gamma_\bullet := \ds\frac{1}{2}N_\bullet^{\alpha\beta}\,a_{\alpha\beta}\,,
\label{e:ngc}
\eqe
and the deviatoric stress norm
\eqb{l}
s_\bullet := \ds\sqrt{N_\mathrm{dev}^{\alpha\beta}\,N^\mathrm{dev}_{\alpha\beta}}\,,\quad
N_\mathrm{dev}^{\alpha\beta} := N^{\alpha\beta}_\bullet - \gamma_\bullet\,a^{\alpha\beta}\,,
\label{e:Mises}
\eqe
that follow from the elastic, viscous and Korteweg stresses introduced in \eqref{e:sigtot} and \eqref{e:Ntot}. Note that $N_{\mathrm{visc}}^{\alpha\beta} = \sigma_{\mathrm{visc}}^{\alpha\beta}$ and $N_{\mathrm{CH}}^{\alpha\beta} = \sigma_{\mathrm{CH}}^{\alpha\beta}$.
In theory $\gamma_\mathrm{CH}=0$ (according to Eq.~\eqref{e:sigCH}), while $\gamma_\mathrm{visc} \neq 0$ (unless area-incompressibility is assumed).
The two stress measures are shown in Fig.~\ref{fig:s_gam2} and \ref{fig:s_str2}.\footnote{To avoid numerical round-off errors in the evaluation of Eq.~\eqref{e:Mises}, the various terms should be multiplied out analytically before implementation.} 
The Korteweg stress is largest around bulges at the phase interface. 
The viscous stress is small in comparison to the Korteweg and elastic stresses. In order to resolve the stress at the phase interface, at least 2 elements should be used per $\sqrt{\lambda}$.

\subsection{Phase separation on a deforming double torus}

\begin{figure}[h]
\centering
\includegraphics[width=0.32\linewidth, trim = 250 380 250 470,clip]{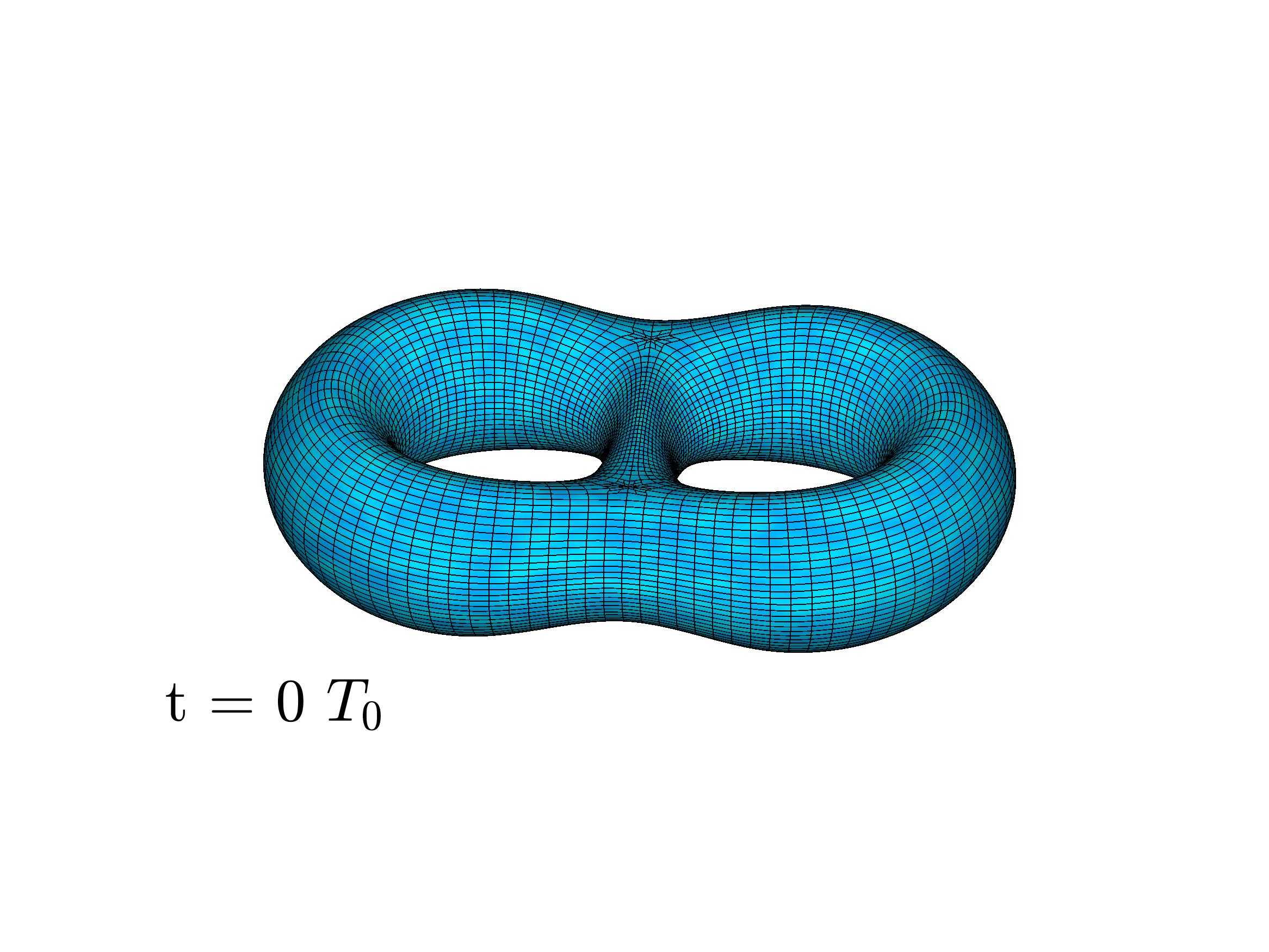}
\includegraphics[width=0.32\linewidth, trim = 250 380 250 470,clip]{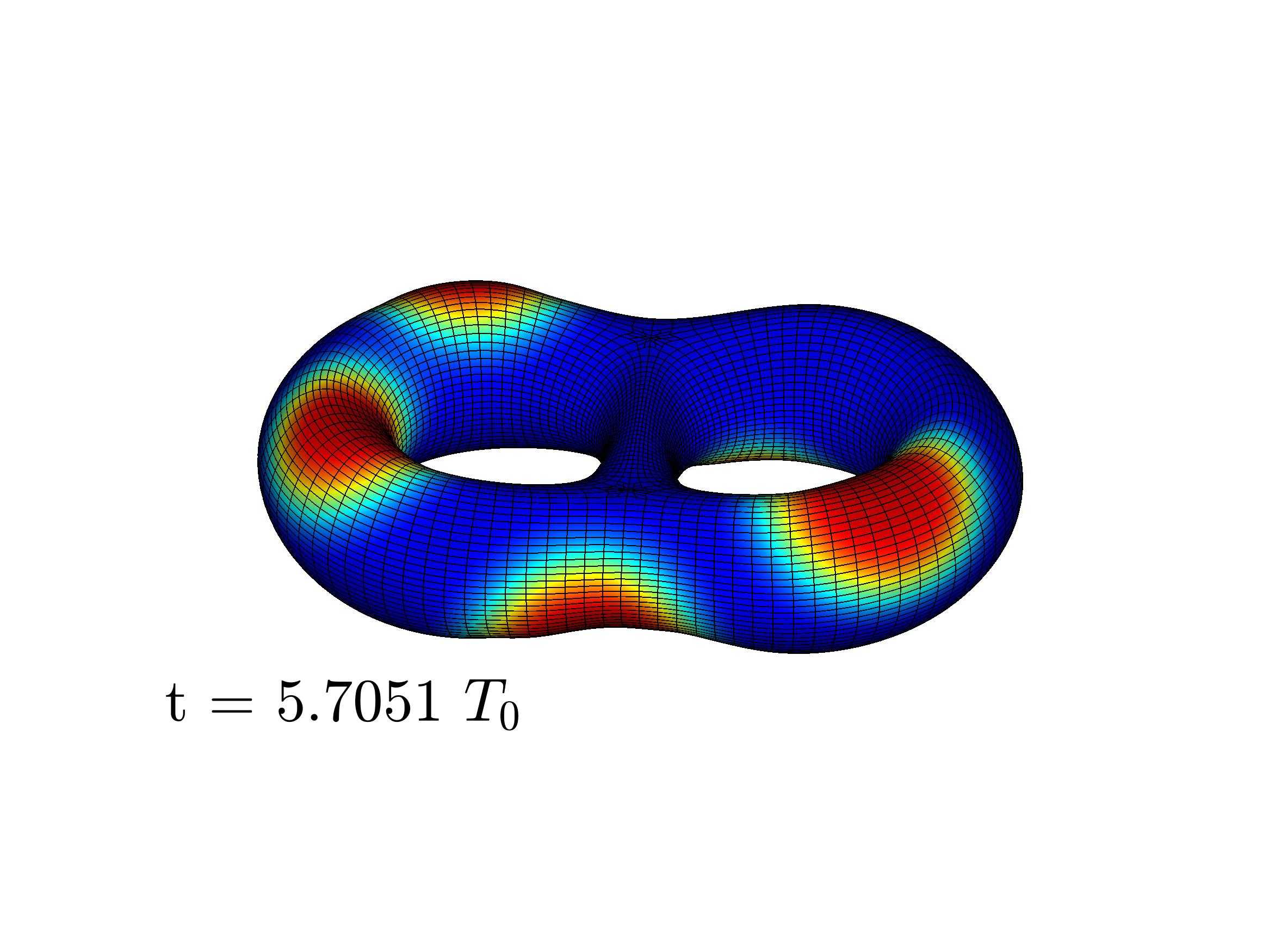}
\includegraphics[width=0.32\linewidth, trim = 250 380 250 470,clip]{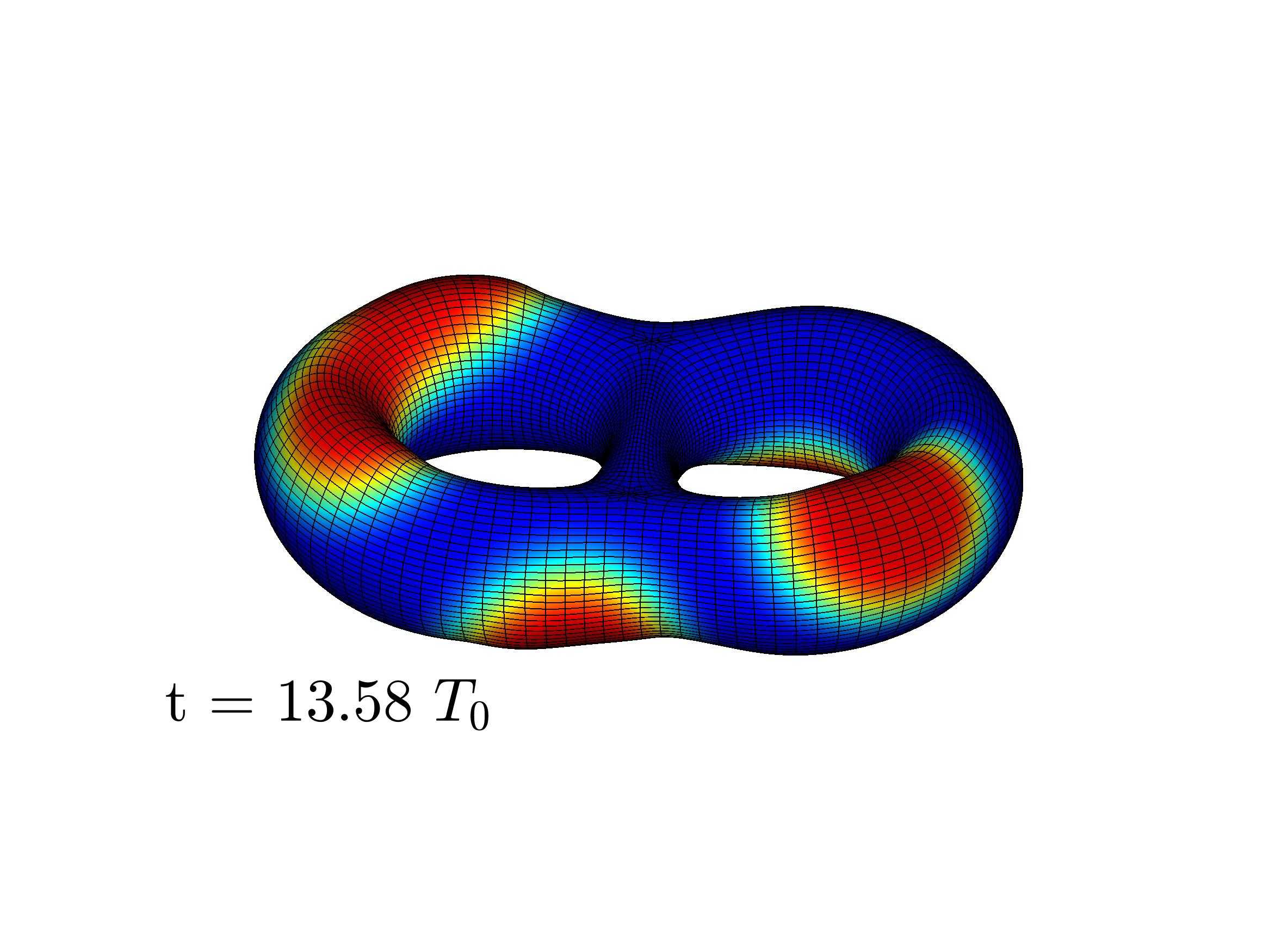}\\
\includegraphics[width=0.32\linewidth, trim = 250 380 250 460,clip]{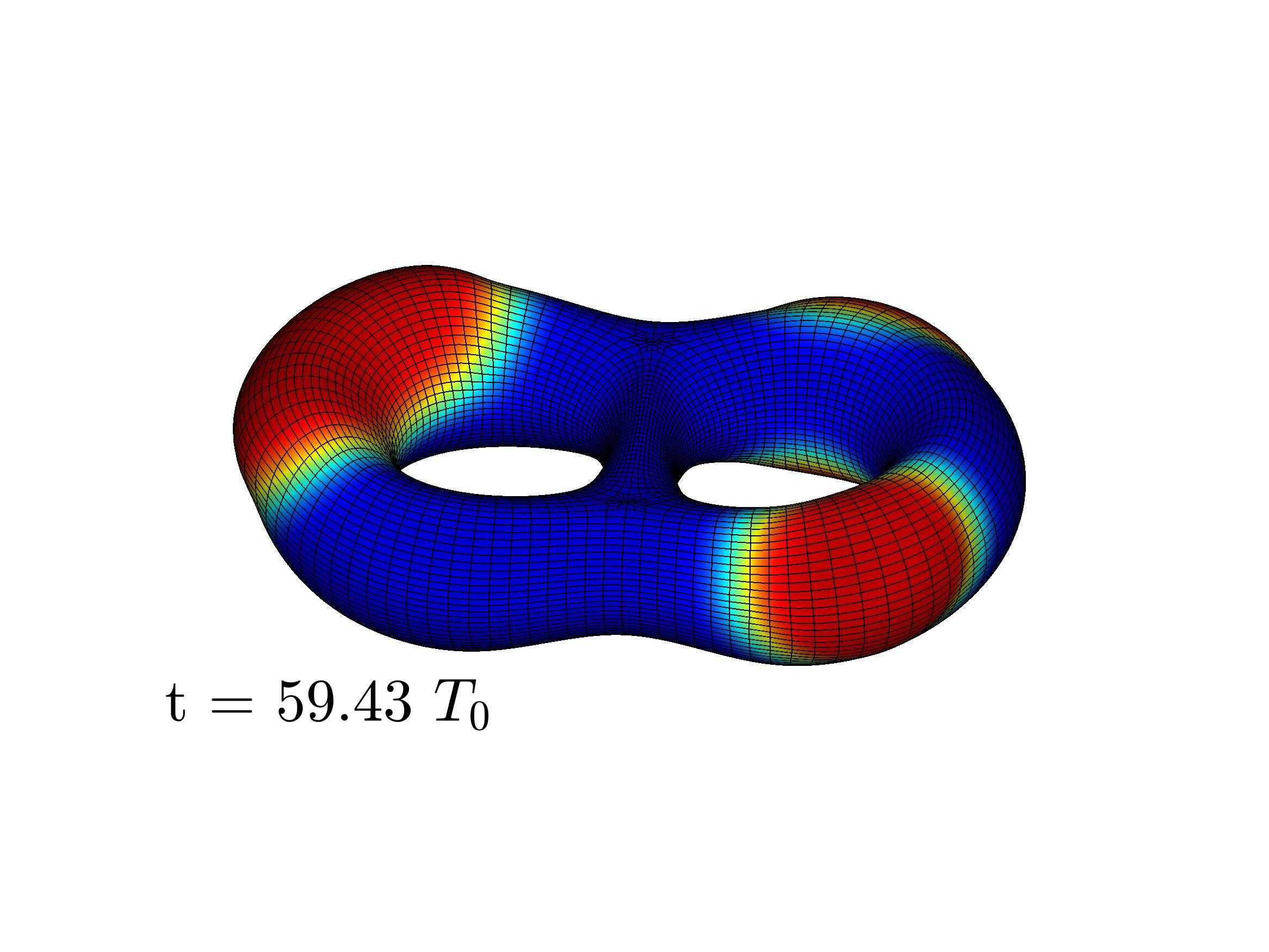}
\includegraphics[width=0.32\linewidth, trim = 250 380 250 460,clip]{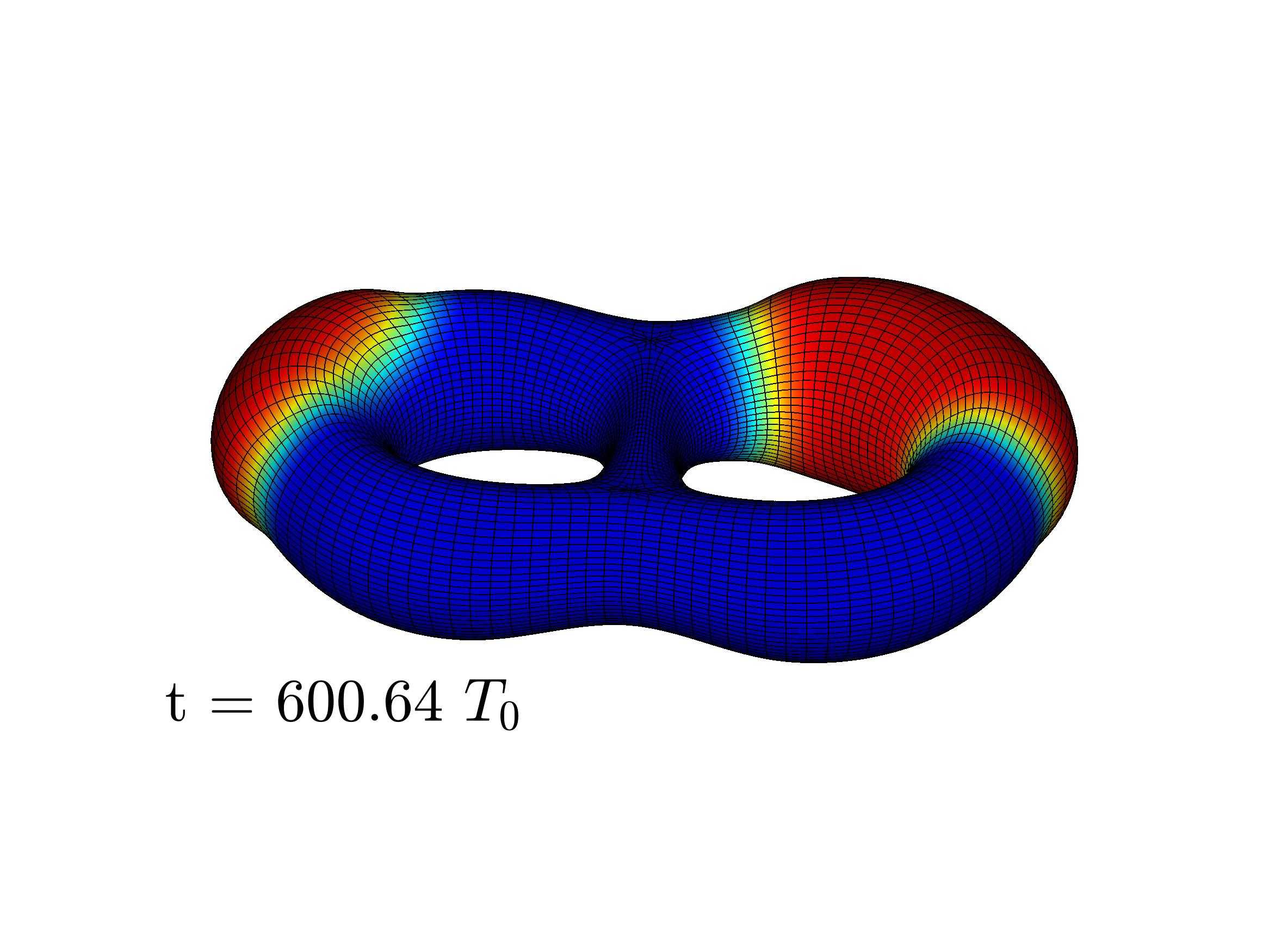}
\includegraphics[width=0.32\linewidth, trim = 250 380 250 460,clip]{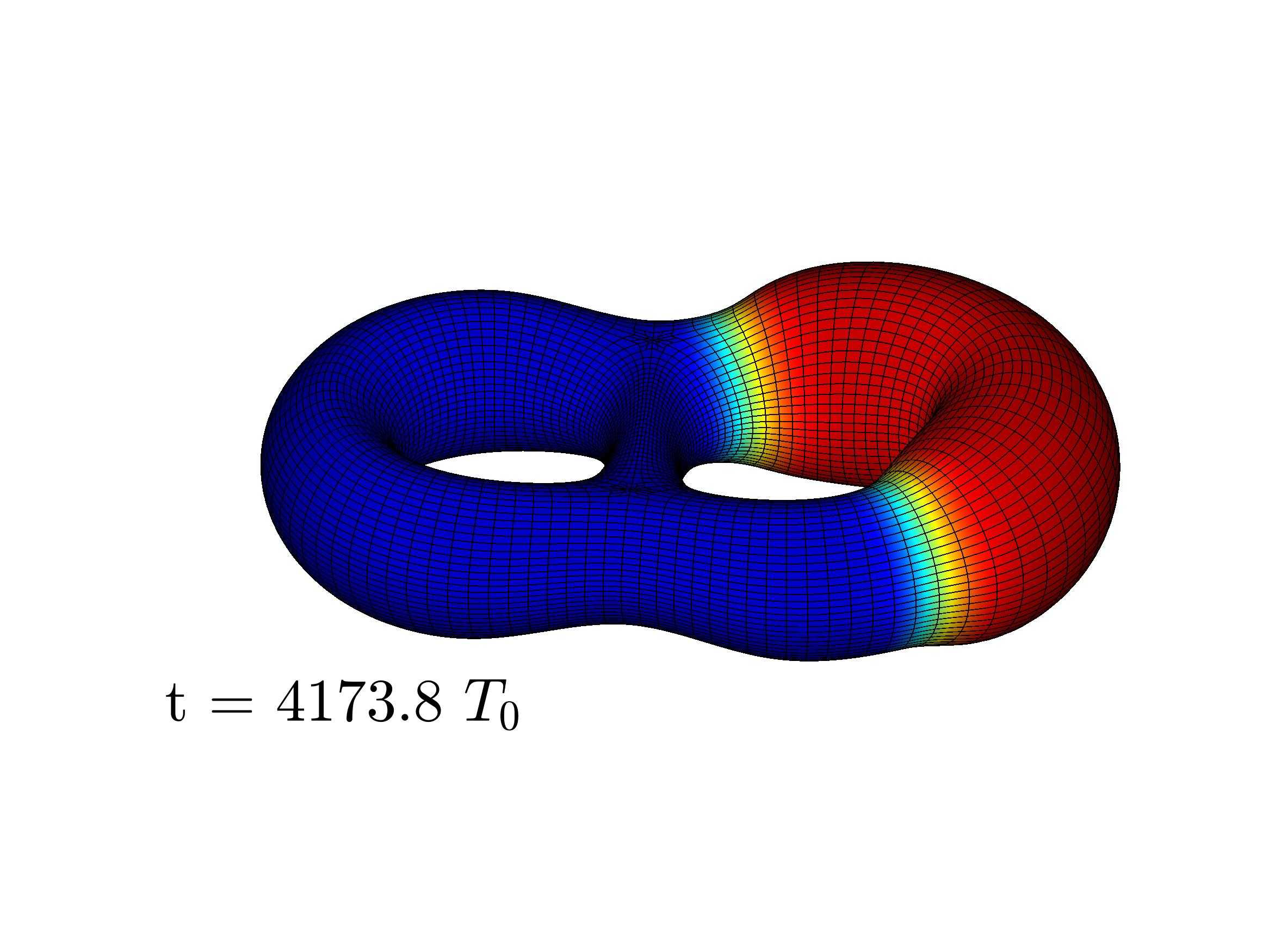}\\
\caption{Phase separation on a deforming double torus: Evolution of the solution with $\lambda=0.025\,L_0^2$ on an unstructured mesh containing $8264$ cubic elements, with $\sqrt{\lambda}\approx 1.6\,h$. The colors follow Table~\ref{tab:mat}.}
\label{fig:d_evo1}
\end{figure}

The last example studies phase separation on a deforming double torus, which is discretized by the unstructured splines from Sec.~\ref{sec:spaceDisc_splines}. 
The parameters are $D=4\,T_0$ and $\lambda=0.025\,L_0^2$ along with the parameters in Table~\ref{tab:mat}. 
The constant internal pressure $p_{\mathrm{int}}=0.03\,EL_0^{-1}$ is prescribed for all $t$. The unstructured mesh consists of $8264$ cubic elements and has $4$ extraordinary points. 
As in the previous example, the discretization is $C^2$-continuous except for the extraordinary points. 
Rigid body deformations are prevented by analogous boundary conditions to those shown in Fig.~\ref{fig:t_IC1}. 
Fig.~\ref{fig:d_evo1} shows the evolution of the phase separation at various times. 
The mechanical deformation and the phase field evolve simultaneously and affect each other. 

\begin{figure}[H]
\centering
\includegraphics[width=0.485\linewidth, trim = 0 0 0 0,clip]{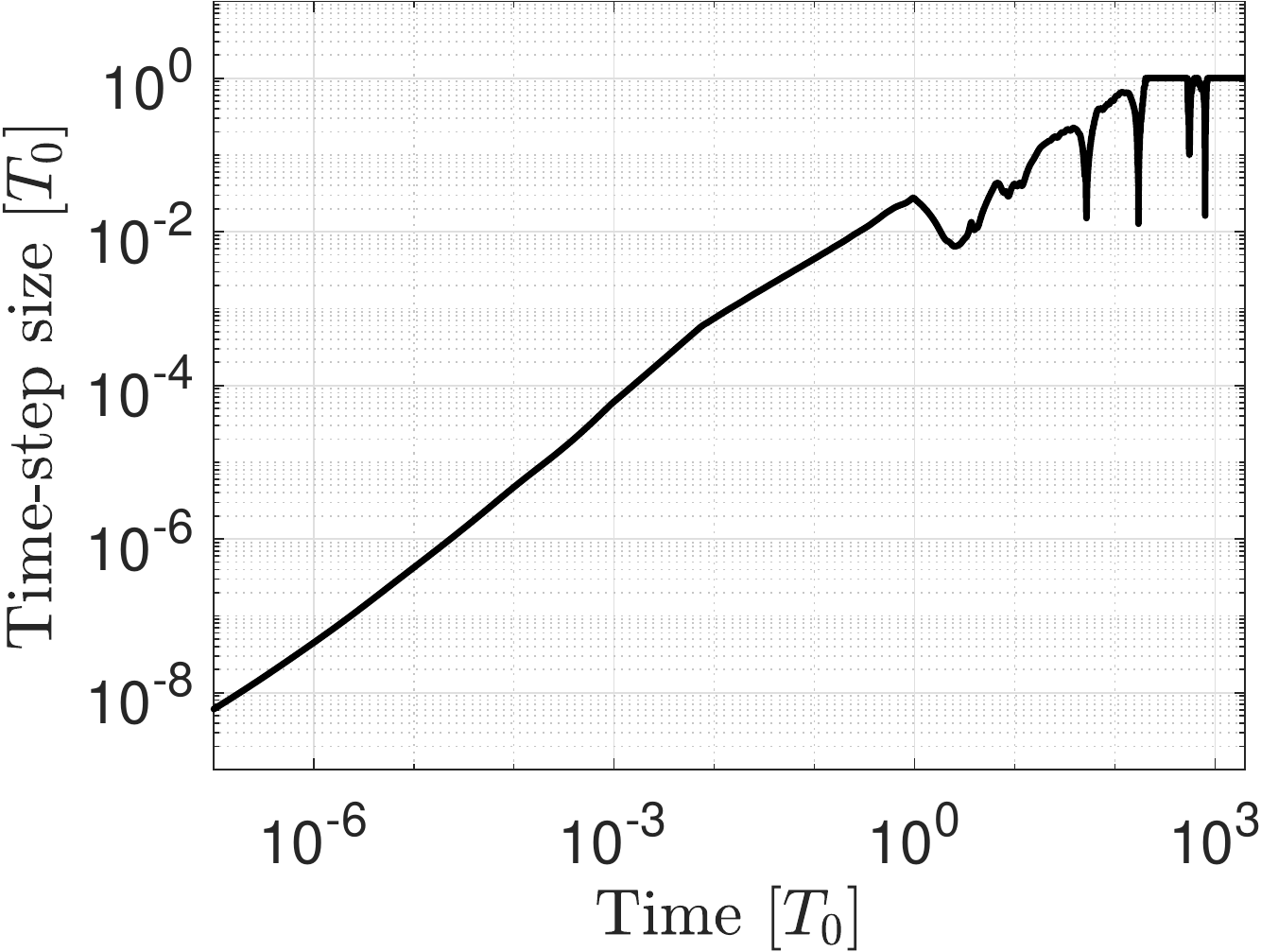}
\includegraphics[width=0.485\linewidth, trim = 0 0 0 0,clip]{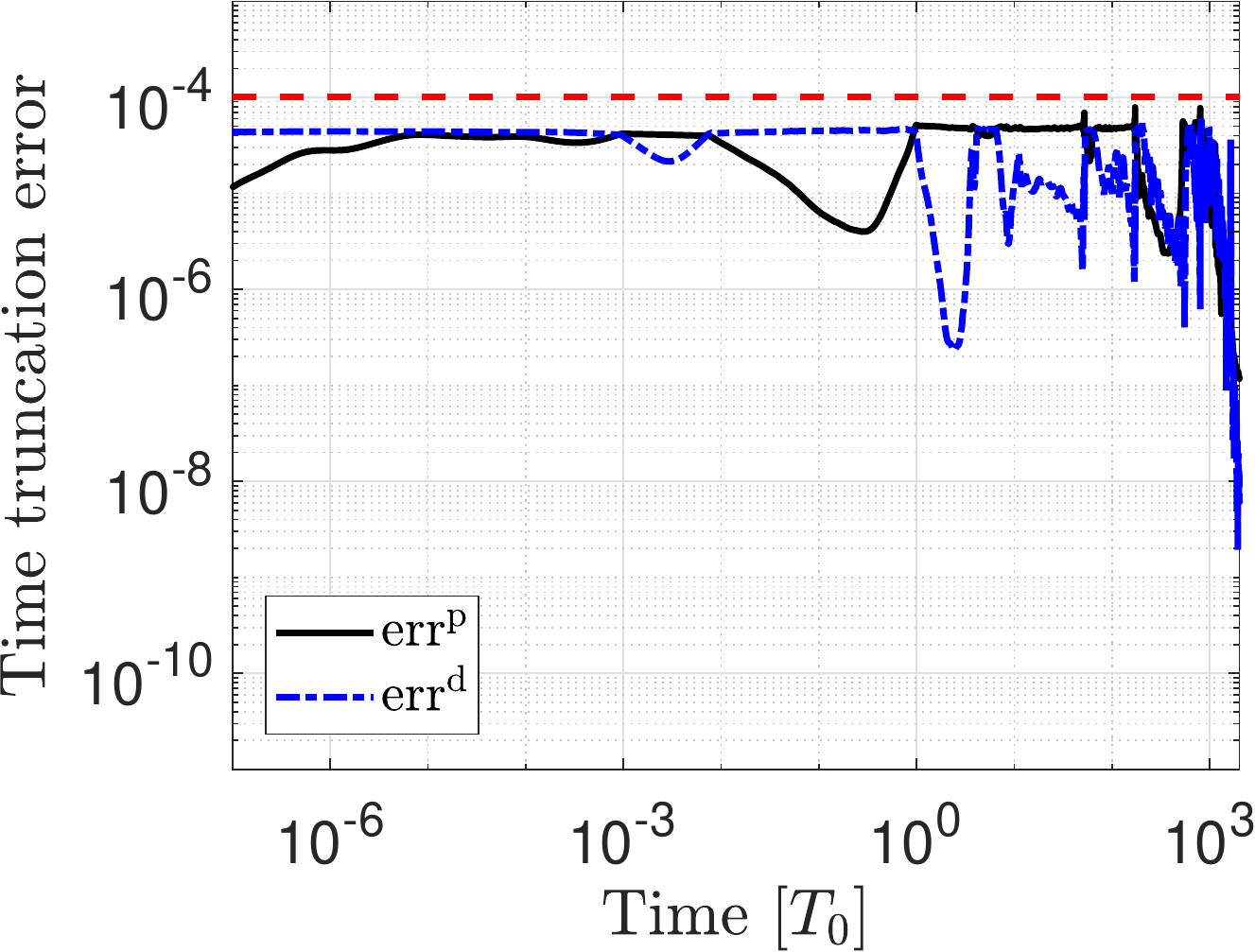}
\caption{Phase separation on a deforming double torus: Left: Adaptive time step size. Right: Evolution of the local time truncation errors of the phase field, $\mathrm{err}^\mathrm{p}$, and mechanical field, $\mathrm{err}^\mathrm{d}$. 
The chosen temporal error bound is shown by a dashed red line.}
\label{fig:d_ts1}
\end{figure}

\begin{figure}[H]
\centering
\includegraphics[width=0.485\linewidth, trim = 0 0 0 0,clip]{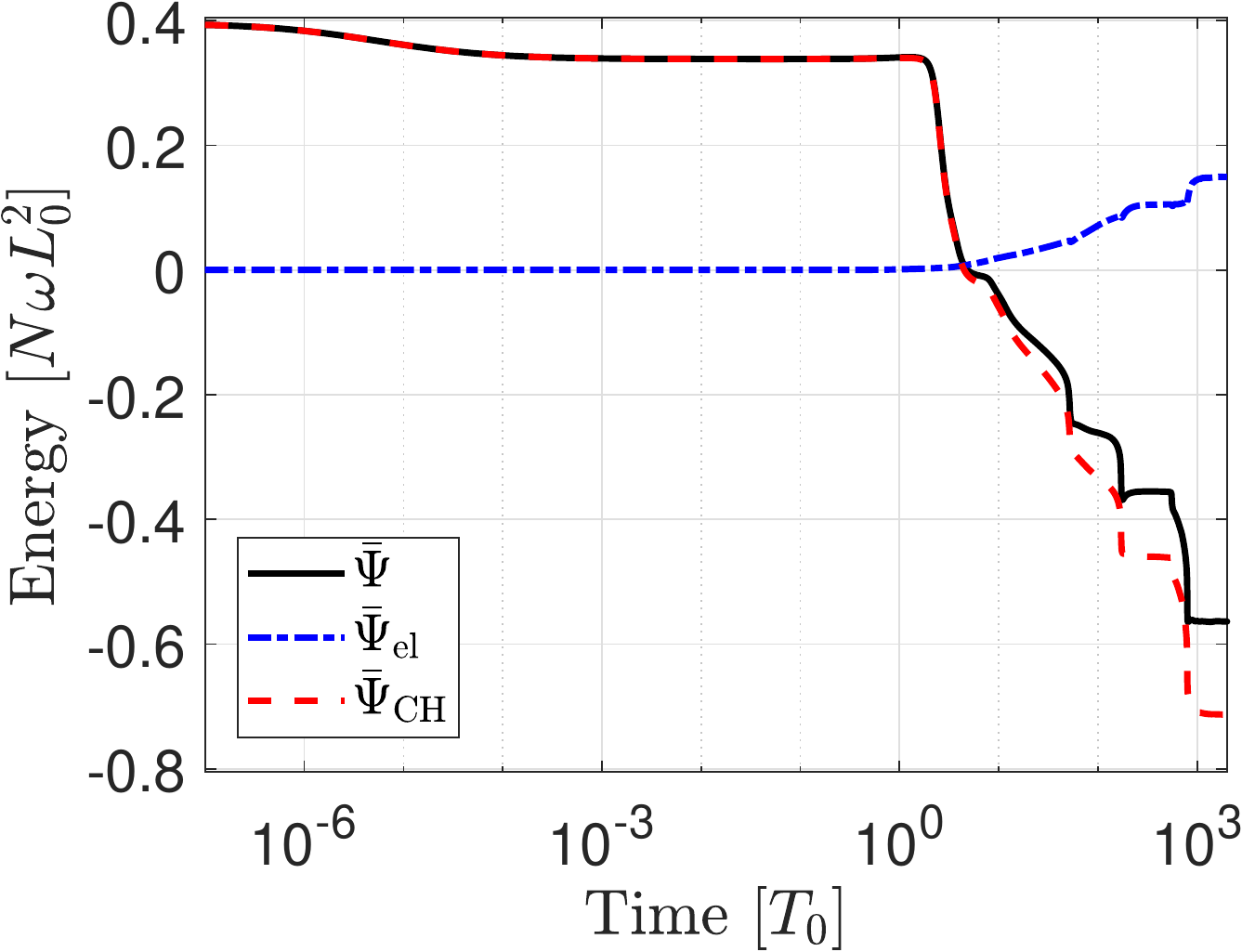}
\includegraphics[width=0.485\linewidth, trim = 0 0 0 0,clip]{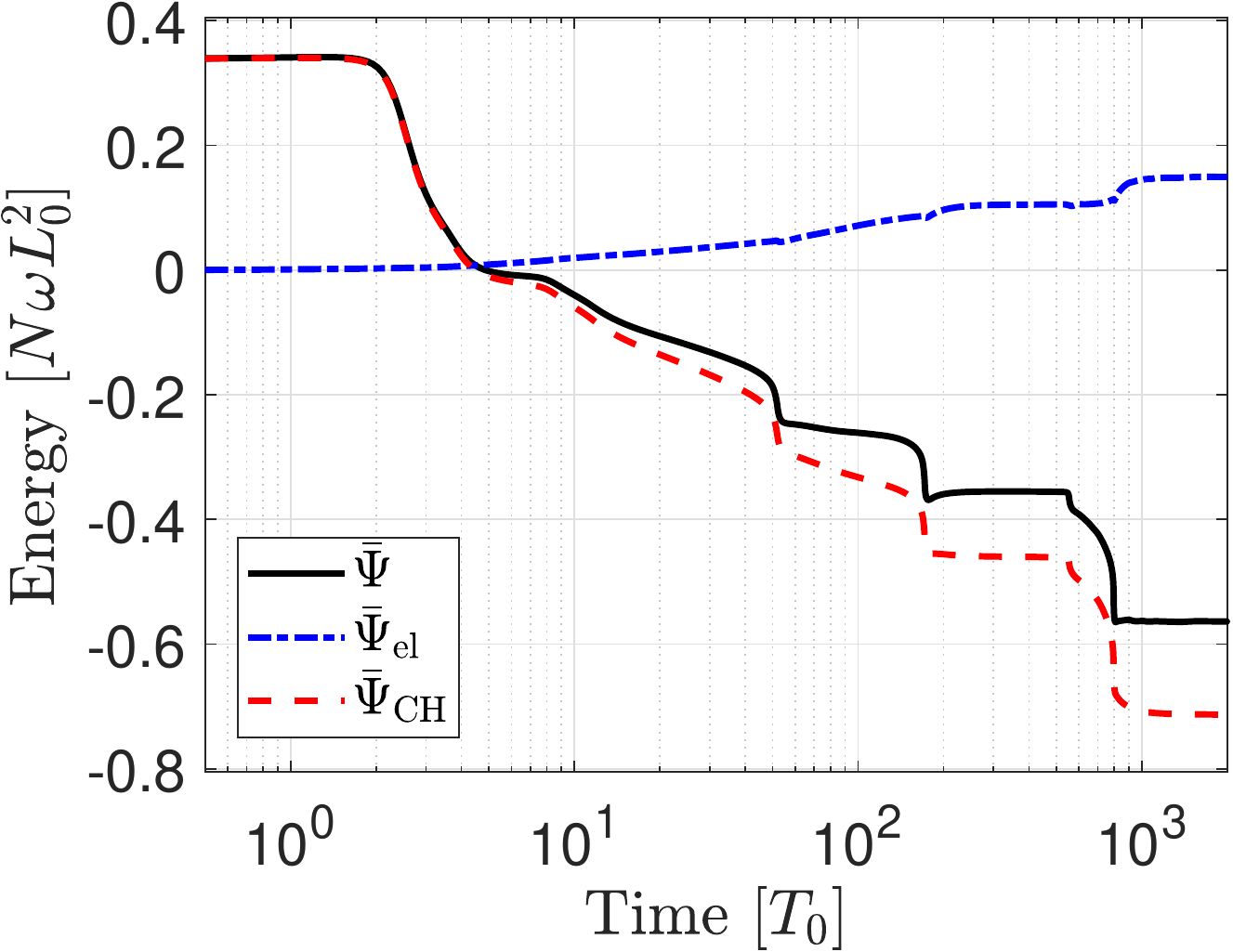}
\caption{Phase separation on a deforming double torus:  Evolution of the characteristic energies of the system. Close-up on the right side.}
\label{fig:d_ene1}
\end{figure}

Fig.~\ref{fig:d_ts1} shows the evolution of the time step size, which is limited to $\Delta t= T_0$, and the local time truncation error. 
The error restriction of the time step size is alternating similar to the torus case in Sec.~\ref{Sec:ex2}.
Fig.~\ref{fig:d_ene1} shows the evolution of the characteristic energies of the system. 

\section{Conclusion}
\label{s:concl}

This work presents a novel coupled formulation for the modeling of phase fields on deforming shell surfaces within the framework of isogeometric finite elements. 
The phase changes are described by the Cahn-Hilliard phase field theory, which is coupled to nonlinear thin shell theory. 
A phase-dependent material model is presented to describe mixtures. 
A monolithic and fully implicit time integration scheme is used to solve the coupled system simultaneously. 
An adaptive time-stepping approach is formulated to adjust the time step size. 
For the numerical examples, bi-quadratic NURBS discretizations and bi-cubic unstructured quadrilateral spline discretizations are used. 
Both provide global $C^1$-continuity.

The examples presented in Sec.~\ref{s:num_ex} demonstrate the direct coupling of phase transitions and mechanical deformations. 
The simultaneous evolution of both fields can be observed for the chosen parameters. 
Other parameters have been observed to produce little or no coupling and they are not reported here for this reason.
The adaptive time-stepping approach allows an automatic control of the time step size. 
The evolution of the phase separation process appears at both small and large time scales. 
In the absence of fast phase separation and large deformation, the time integration error estimation leads to an almost steady increase of the time step size. 
Suitable material behavior is required to allow for large deformations and an appropriate interaction of both fields. 
Due to the direct interaction of mechanical and phase field, the coupled system needs to be damped by viscosity to avoid the build-up of surface oscillations from phase-separation induced deformations. The condition numbers of the tangent matrices of the examples indicate similar observations as in \citet{bartezzaghi16}: They increase with mesh refinement and with spline order. 
They decrease with the time step size $\Delta t$. 
Examining the Newton-Raphson accuracy shows that the present simulation results are not affected by any ill-conditioning (as the Newton-Raphson accuracy reaches machine precision).

Possible extensions of this work include studying applications such as battery systems, liquid droplets and lipid bilayers. 
The presented shell formulation also applies to liquid menisci \citep{droplet} and lipid bilayers \citep{Katira2016,liquidshell}, but additional numerical tools are needed to handle the coupling of surface flows and phase fields. 
Another possible extension is the modeling of contact, since large deformations can lead to self-contact. 
The development of adaptive spatial refinement strategies in order to resolve very thin phase interfaces would also be beneficial. In the future, experiments are also called for in order to calibrate and validate the proposed formulation.

\section*{Acknowledgments}

Thomas J.R. Hughes and Deepesh Toshniwal were partially supported by the Office of Naval Research (Grant Nos. N00014-17-1-2119 and N00014-13-1-0500). Kranthi K. Mandadapu acknowledges support from the University of California Berkeley, from the National Institutes of Health Grant R01-GM110066 and from the Department of Energy (contract DE-AC02-05CH11231, FWP no. CHPHYS02). Roger A. Sauer acknowledges the support from a J. Tinsley Oden fellowship and funding from the German Research Foundation (DFG) through project GSC 111.

\bigskip


\appendix
\section*{Appendix}

\section{On the constitutive relations}\label{s:thermo}

This section briefly summarizes the derivation of the constitutive equations in Sec.~\ref{s:flux} and \ref{s:sigM} following \citet{sahu17}. 
The local form of the energy balance on a curved surface can be written as 
\begin{equation} \label{e:energy}
\rho\,\dot{u} = \rho\,r - q^\alpha_{; \alpha} + \dfrac{1}{2} \sigma^{\alpha \beta} \dot{a}_{\alpha \beta} + M^{\alpha \beta} \dot{b}_{\alpha \beta}\,,
\end{equation}
where $u$ denotes the internal energy density, $r$ is a heat source, $\bq=q^\alpha\ba_\alpha$ is the heat flux on the surface, and $\sigma^{\alpha\beta}$ and $M^{\alpha\beta}$ are stress and bending moment components, respectively. 
The local form of the entropy balance is
\begin{equation} \label{e:entropy}
\rho\,\dot{s} = -{J_s}^{\alpha}_{;\alpha} + \rho\,s_e + \rho\,s_i \,.
\end{equation}
Here, $s$ is the entropy density per unit mass, $J_{s}^{\alpha}$ is the total entropy flux, $s_e$ is the total external entropy rate, and $s_i$ is the total internal entropy production rate. 
The second law of thermodynamics dictates that 
\begin{equation}\label{e:2law}
\rho\,s_i \geq 0\,.
\end{equation}
Let the Helmholtz free energy (per unit mass) be defined as 
\begin{equation} \label{e:psi}
\psi = u - T s
\end{equation}
and assume that it depends on the kinematic variables as
\begin{equation}
\psi = \psi\big(a_{\alpha \beta}, b_{\alpha \beta}, T, \phi, \phi_{;\alpha}\big)\,.
\end{equation}
Taking a time derivative then gives 
\eqb{l}
\dot\psi = \ds\pa{\psi}{a_{\alpha\beta}}\dot a_{\alpha\beta} + \pa{\psi}{b_{\alpha\beta}}\dot b_{\alpha\beta} +  \pa{\psi}{T}\dot T +  \pa{\psi}{\phi}\dot\phi + \pa{\psi}{\phi_{;\alpha}}\dot\phi_{;\alpha}\,,
\eqe
which can be rewritten into
\eqb{l}
\rho\,\dot\psi = \rho \ds\pa{\psi}{a_{\alpha\beta}}\dot a_{\alpha\beta} + \rho \pa{\psi}{b_{\alpha\beta}}\dot b_{\alpha\beta} + \rho \pa{\psi}{T}\dot T + 
\rho\,\tilde\mu_\mrc\,\dot\phi
+ \bigg(\rho \pa{\psi}{\phi_{;\alpha}}\dot\phi\bigg)_{\!\!;\alpha},
\label{e:psidot}\eqe
where
\begin{equation}\label{e:chempot}
\tilde\mu_\mrc  := \pa{\psi}{\phi}-\frac{1}{\rho}\bigg(\rho\pa{\psi}{\phi_{;\alpha}}\bigg)_{\!\!;\alpha}
\end{equation}
introduces the chemical potential (per unit mass).
From \eqref{e:psi} we get 
\begin{equation} \label{e:psi2}
\dot{\psi} = \dot{u} - T \dot{s} - \dot{T}s\,,
\end{equation}
which can be combined with \eqref{e:energy} and \eqref{e:psidot} to yield 
\eqb{lll}\label{e:entropy2}
\rho\,\dot{s}  \is\ds  \dfrac{1}{T} \big(\rho r - q^\alpha_{; \alpha}\big) -\frac{\rho}{T}\bigg(s+\frac{\partial\psi}{\partial T} \bigg)\dot{T}- \rho\frac{\tilde\mu_\mrc}{T} \dot\phi - \frac{1}{T}\bigg(\rho \pa{\psi}{\phi_{;\alpha}}\dot\phi\bigg)_{\!\!;\alpha}\\[4mm]
\plus \dfrac{\dot{a}_{\alpha \beta}}{T} \bigg( \dfrac{1}{2} \sigma^{\alpha \beta} - \rho \dfrac{\partial \psi}{\partial a_{\alpha \beta}}\bigg)
+ \dfrac{\dot{b}_{\alpha \beta}}{T} \bigg( M^{\alpha \beta} - \rho \dfrac{\partial \psi}{\partial b_{\alpha \beta}} \bigg).
\eqe
Let us now assume isothermal conditions such that there are no in-plane temperature gradients ($T_{;\alpha}=0$), and define the entropy as
\eqb{l}
s = \ds - \frac{\partial\psi}{\partial T}\,.
\eqe
Substituting the mass balance equation for $\phi$ \eqref{e:sf} into \eqref{e:entropy2}, we then get 
\eqb{lll}\label{eq:local-entropy-balance-2}
\rho\,\dot{s}
\is - \ds\bigg(\frac{q^\alpha}{T} - \frac{j^\beta_{;\beta}}{T} \pa{\psi}{\phi_{;\alpha}}  -\frac{j^\alpha\tilde\mu_\mrc}{T} \bigg)_{\!\!; \alpha} + \dfrac{\rho\,r}{T} \\[4mm]
\plus \ds\frac{\dot{a}_{\alpha \beta}}{T} \bigg(\frac{1}{2} \sigma^{\alpha \beta} - \rho \dfrac{\partial \psi}{\partial a_{\alpha \beta}}\bigg)
	  + \frac{\dot{b}_{\alpha \beta}}{T} \bigg( M^{\alpha \beta} - \rho \frac{\partial \psi}{\partial b_{\alpha \beta}} \bigg) 
	  - \ds\frac{j^\alpha\tilde\mu_{\mrc;\alpha}}{T}\,.
\eqe
Comparing the above equation with the entropy balance \eqref{e:entropy} lets us identify the entropy flux 
\begin{equation}
J_s^{\alpha} = \frac{q^\alpha}{T}  - \frac{j^\beta_{;\beta}}{T} \pa{\psi}{\phi_{;\alpha}} -\frac{j^\alpha\tilde\mu_\mrc}{T}\,,
\end{equation}
the external entropy rate
\begin{equation}
\rho\,s_e = \frac{\rho\,r}{T}
\end{equation}
and the total entropy production
\begin{equation}\label{eq:ent-prod}
\rho\,s_i = \ds\frac{\dot{a}_{\alpha \beta}}{T} \bigg(\frac{1}{2} \sigma^{\alpha \beta} - \rho \dfrac{\partial \psi}{\partial a_{\alpha \beta}}\bigg)
	  + \frac{\dot{b}_{\alpha \beta}}{T} \bigg( M^{\alpha \beta} - \rho \frac{\partial \psi}{\partial b_{\alpha \beta}} \bigg) 
	  - \frac{j^\alpha\tilde\mu_{\mrc;\alpha}}{T}\,,
\end{equation}
which has to be positive according to \eqref{e:2law}.
Since $T \geq 0$ and since the quantities $\dot{a}_{\alpha \beta}$, $\dot{b}_{\alpha \beta}$ and 
$\tilde\mu_{\mrc; \alpha}$ 
can be varied independently, this implies that
\begin{equation} \label{e:eprod1}
\dot{a}_{\alpha \beta}\bigg(\dfrac{1}{2} \sigma^{\alpha \beta} - \rho \dfrac{\partial \psi}{\partial a_{\alpha \beta}}\bigg) \geq 0 
\quad \forall~\dot{a}_{\alpha \beta}
\end{equation}
and 
\begin{equation}\label{e:eprod2}
\dot{b}_{\alpha \beta} \bigg( M^{\alpha \beta} - \rho \dfrac{\partial \psi}{\partial b_{\alpha \beta}} \bigg) \geq 0 
\quad \forall~\dot{b}_{\alpha \beta}
\end{equation}
and
\begin{equation}\label{e:eprod3}
- j^\alpha \tilde\mu_{\mrc;\alpha} \geq 0\,.
\end{equation}
If we assume that the bending behavior is purely elastic (such that the bending moments do not depend on the curvature rate $\dot b_{\alpha \beta}$), Eq.~\eqref{e:eprod2} implies that the bending moments are given by
\begin{equation} \label{e:Mphi}
M^{\alpha \beta} =   \rho \dfrac{\partial \psi}{\partial b_{\alpha \beta} }\,.
\end{equation}
If we assume that the membrane behavior contains an elastic part and a linear viscosity part, Eq.~\eqref{e:eprod1} implies that the membrane stresses are given by
\begin{equation} \label{e:sigphi}
\sigma^{\alpha \beta} = 2 \rho \dfrac{\partial \psi}{\partial a_{\alpha \beta}} - \eta\,\dot{a}^{\alpha \beta}\,,
\end{equation}
for $\eta\geq0$.
The simplest model satisfying inequality \eqref{e:eprod3} is the linear flux relationship 
\begin{equation} \label{e:jphi}
j_\alpha = - \tilde M \, \tilde\mu_{\mrc;\alpha}\,, 
\end{equation}
as long as $\tilde M\geq0$. These linear models for viscosity and flux are the simplest cases. Alternatively, one can also use nonlinear relationships satisfying \eqref{e:eprod2} and \eqref{e:eprod3}.
Introducing the initial density $\hat\rho = J\rho$ of the undeformed initial mixture, which is considered to be uniform (such that $\hat\rho_{;\alpha}=0$ and $\partial\hat\rho/\partial\phi=0$), we can define the chemical potential per reference area,
\eqb{l}
\mu_\mrc := \hat\rho\,\tilde\mu_\mrc = \ds\pa{\Psi}{\phi} - J\ds\bigg(\frac{1}{J}\pa{\Psi}{\phi_{;\alpha}}\bigg)_{\!\!;\alpha},
\eqe
where $\Psi:=\hat\rho\,\psi$ is the Helmholtz free energy per reference area, and use this to rewrite
\begin{equation} \label{e:jphi2}
j_\alpha = - \ds\frac{M}{J} \, \mu_{\mrc;\alpha}\,,
\end{equation}
where $M := \tilde M/\rho$. 
Similarly, $\sig^{\alpha\beta}$ and $M^{\alpha\beta}$ can be rewritten into the expressions of Eq.~\eqref{e:tauM0}.

\section{Extraction operator initialization} \label{App:ExtractionInitialization}
Initializing the extraction operator for a bi-cubic spline function $N$ on an element $\Omega$ is equivalent to defining the polynomial that $N|_\Omega$ equals. Choosing the element-local polynomial basis as tensor product Bernstein polynomials $B_{ij}$, we only need to specify coefficients $a_{ij}$ such that,
\begin{equation}
	N|_\Omega = \sum_{i,j=0}^{3} a_{ij} B_{ij}\;.
\end{equation}
\begin{figure}[H]
	\centering
	\includegraphics[width=0.25\linewidth]{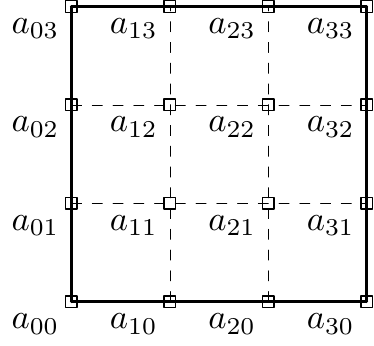}
	\caption{A graphical representation of the extraction operator coefficients on a particular B\'ezier element in the mesh. The coefficient $a_{ij}$ corresponds to the tensor product Bernstein polynomial $B_{ij}$.}
	\label{fig:bezierElement}
\end{figure}
In the following, the extraction coefficients $a_{ij}$ will be denoted graphically as shown in Figure \ref{fig:bezierElement}. Then, the initialization of extraction operators for functions $N_i^D$ and $N_i^A$ spanning spline spaces $\splSpace{D}$ and $\splSpace{A}$, respectively, is done as follows. (Note that the following assumes all elements to be of uniform size in the parametric domain; this is true for all the numerical results presented in this paper. Please see \citep{toshniwal2017smooth} for a more general case.)
\begin{itemize}
	\item For $\splSpace{D}$, a basis function is assigned to each vertex of the mesh (black and red disks in Figure \ref{fig:splineDOFConfig}). The top-left figure in Figure \ref{fig:extractionInitialization} shows the extraction coefficients for the basis function corresponding to the gray disk; $\mu_i$ are the number of edges incident on the corners of the element.
	\item For $\splSpace{A}$, a basis function is assigned to each regular vertex of the mesh (black disks in Figure \ref{fig:splineDOFConfig}), and $4$ additional basis functions are assigned to each element containing an extraordinary vertex (blue squares in Figure \ref{fig:splineDOFConfig}). We call the former vertex-based basis, and the latter face-based basis. The extraction coefficients for them are initialized in a two-step process:
	\begin{enumerate}
		\item For each vertex-based basis, the extraction coefficients are initialized as per the top-left figure in Figure \ref{fig:extractionInitialization}. For each face-based basis, depending on its location w.r.t. the extraordinary point, the extraction coefficients are initialized as per the top-right, mid-left and mid-right figures in Figure \ref{fig:extractionInitialization}. (The particular face-based basis being initialized corresponds to the blue square in these figures.)
		\item In order to retain partition of unity, the extraction coefficients for some of the vertex-based basis functions are truncated. Depending on their location w.r.t. the extraordinary point, the extraction coefficients that are set equal to $0$ have been crossed out in the bottom-left and bottom-right figures in Figure \ref{fig:extractionInitialization}. (The particular vertex-based basis being truncated corresponds to the gray disk in these figures.)
	\end{enumerate}
\end{itemize}
\begin{figure}
	\centering
	\includegraphics[width=0.35\linewidth]{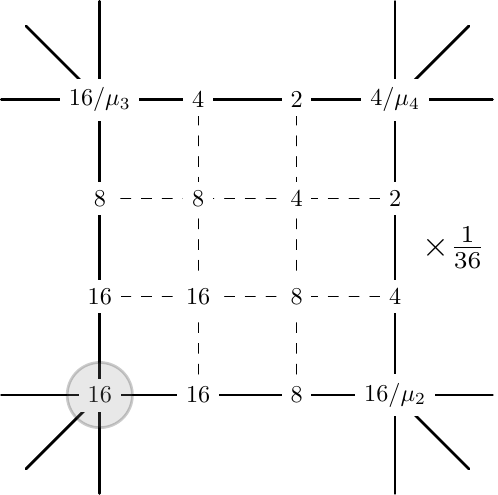} \hspace{1cm}
	\includegraphics[width=0.35\linewidth]{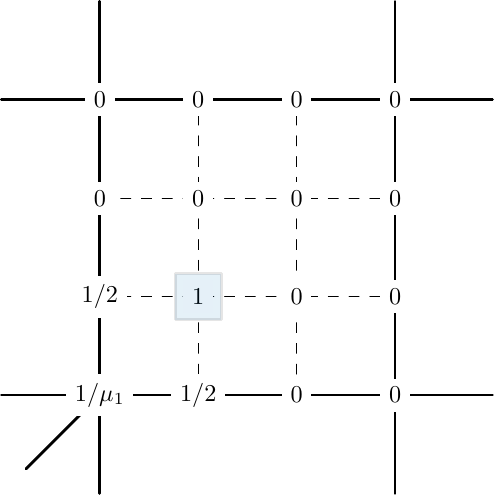} \\ \vspace{0.5cm}
	\includegraphics[width=0.35\linewidth]{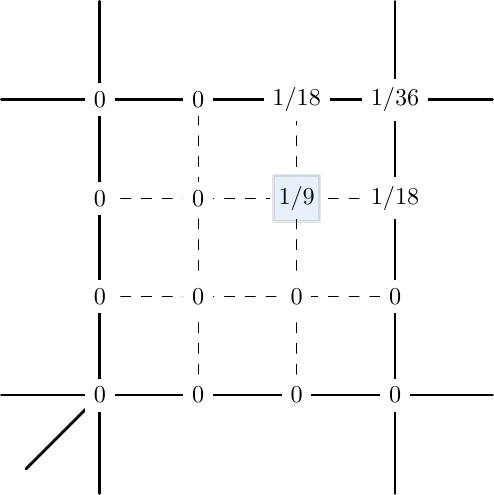} \hspace{1cm}
	\includegraphics[width=0.35\linewidth]{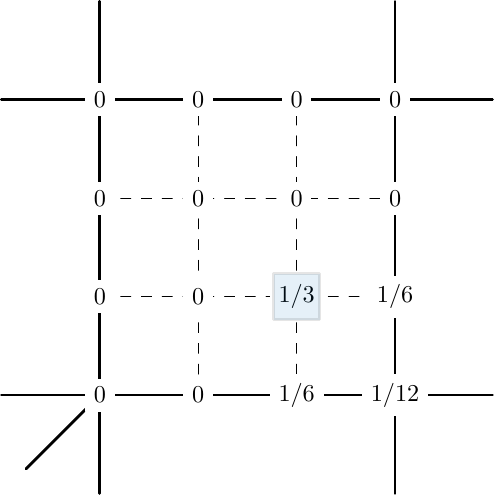} \\ \vspace{0.5cm}
	\includegraphics[width=0.35\linewidth]{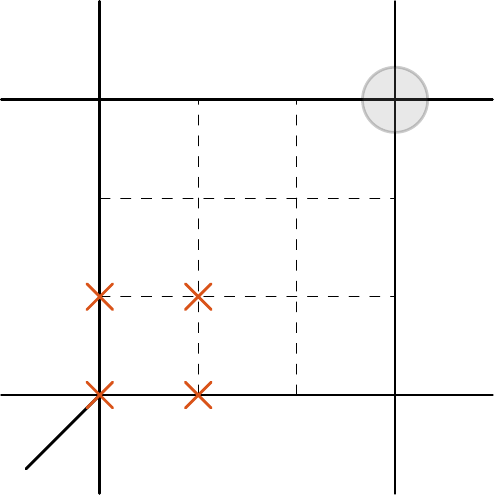} \hspace{1cm}
	\includegraphics[width=0.35\linewidth]{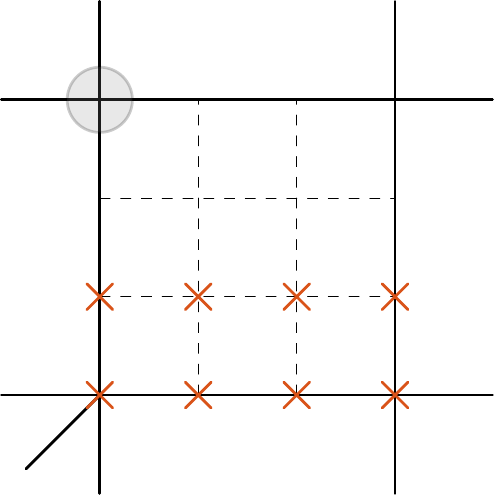}
	\caption{A graphical depiction of the extraction operator initialization. The particular degrees-of-freedom for which the initialization is being performed have been denoted as gray disks or blue squares. The number of edges incident on an element corner has been denoted with $\mu_i$; $\mu_i \neq 4$ implies an extraordinary point. In all but the top-left figure, the extraordinary point is assumed to coincide with the bottom-left corner of B\'ezier element. Additionally, it is assumed that an element contains at most one extraordinary point for its corners.}
	\label{fig:extractionInitialization}
\end{figure}


\section{Temporal discretization and Newton's method}
Details on the temporal discretization and the iterative solution procedure of Newton's method are given in this section. 
\subsection{Generalized-$\alpha$ method}\label{sec:genA}
The intermediate quantities, and the quantities at time step $n+1$, are evaluated for the generalized-$\alpha$ method as,
\eqb{lll}
\mx_{n+1} \is \mx_n+\Delta t_{n+1}\, \dot{\mx}_n+\big(\big(0.5 -\beta \big){\Delta t^2_{n+1}} \big)\ddot{\mx}_n+\beta\Delta t^2_{n+1} \ddot{\mx}_{n+1}\,, \\[2mm]
\dot\mx_{n+1} \is \dot\mx_n+\big(\big(1 -\gamma \big)\Delta t_{n+1} \big)\ddot{\mx}_n+\gamma\Delta t_{n+1}\ddot{\mx}_{n+1}\,, \\[2mm]
\mx_{n+\alpha_\mrf} \is \big(1-\alpha_\mrf \big)\mx_n+\alpha_\mrf\mx_{n+1}\,, \\[2mm]
\dot\mx_{n+\alpha_\mrf} \is \big(1-\alpha_\mrf \big)\dot\mx_n+\alpha_\mrf\dot\mx_{n+1}\,, \\[2mm]
\ddot\mx_{n+\alpha_\mrm} \is \big(1-\alpha_\mrm \big)\ddot\mx_n+\alpha_\mrm\ddot\mx_{n+1}\,, \\[2mm]
\bphi_{n+1} \is \bphi_n+\Delta t_{n+1}\,\dot{\bphi}_n+\gamma \Delta t_{n+1} \left(\dot{\bphi}_{n+1}-\dot{\bphi}_{n} \right)\,, \\[2mm]
\bphi_{n+\alpha_\mrf} \is \big(1-\alpha_\mrf \big)\bphi_n+ \alpha_\mrf {\bphi}_{n+1}\,, \\[2mm]
\dot{\bphi}_{n+\alpha_\mrm} \is \big(1-\alpha_\mrm \big)\dot{\bphi}_n+\alpha_\mrm\dot{\bphi}_{n+1}\,,
\label{e:gen_a2}
\eqe
where $\Delta t_{n+1} = t_{n+1}-t_n$ is the time step. 
The algorithmic parameters $\gamma$, $\beta$, $\alpha_\mrf$ and $\alpha_\mrm$ in Eqs.~\eqref{e:gen_a1} and \eqref{e:gen_a2} control numerical dissipation.
They can be expressed in terms of $\rho_{\infty}\in [0,1]$, which is an algorithmic parameter corresponding to the spectral radius of the amplification matrix as $\Delta t_{n+1} \rightarrow \infty$, i.e.
\eqb{lll}
\alpha_\mrf = \ds\frac{1}{1+\rho_{\infty}}\,,\quad
\alpha_\mrm = \ds\frac{2-\rho_{\infty}}{1+\rho_{\infty}}\,,\\[5mm]
\gamma = \ds\frac{1}{2}+\alpha_\mrm-\alpha_\mrf\,,\quad
\beta = \ds\frac{1}{4}\,(1+\alpha_\mrm-\alpha_\mrf)^2\,,
\eqe
(see \citet{chung93} for further details). The choice $\rho_\infty = 0.5$ shows good performance in the subsequent numerical examples.
\subsection{Newton-Raphson iteration}
\label{sec:NR}
Following \citet{bazilevs13}, the initial guess for the Newton-Raphson iteration is set to
\eqb{lll}
\mx_{n+1}^0 = \mx_n+\Delta t_{n+1}\, \dot{\mx}_n+\big(\big(0.5 -\beta \big)\Delta t^2_{n+1} \big)\ddot{\mx}_n+\big(\beta\Delta t^2_{n+1} \big)\ddot{\mx}^0_{n+1}\,, \\[2mm]
\dot\mx_{n+1}^0 = \dot\mx_n\,,\\[2mm]
\ddot\mx_{n+1}^0 = \ddot\mx_n\ds\frac{\gamma -1}{\gamma}\,, \\[2mm]
\bphi_{n+1}^0 = \bphi_n~,\\[2mm]
\dot\bphi_{n+1}^0 = \dot\bphi_n \ds\frac{\gamma -1}{\gamma}\,,
\label{e:init_g}
\eqe
and then updated from interation step $i\rightarrow i+1$ by
\eqb{lll}
\mx_{n+1}^{i+1} = \mx_{n+1}^{i} + \Delta\mx_{n+1}^{i+1}~,\\[3mm]
\dot\mx_{n+1}^{i+1} = \dot\mx_{n+1}^{i} + \Delta\mx_{n+1}^{i+1} \ds\frac{1}{\gamma\,\Delta t_{n+1}}\,, \\[3mm]
\ddot\mx_{n+1}^{i+1} = \ddot\mx_{n+1}^{i} + \Delta\mx_{n+1}^{i+1} \ds\frac{1}{\beta\,\Delta t^2_{n+1}}\,, \\[3mm]
\bphi_{n+1}^{i+1} = \bphi_{n+1}^{i}+\Delta \bphi_{n+1}^{i+1}\,, \\[3mm]
\dot\bphi_{n+1}^{i+1} = \dot\bphi_{n+1}^{i} + \bphi_{n+1}^{i+1} \ds\frac{1}{\gamma\,\Delta t_{n+1}}\,,
\label{e:upd_g}
\eqe
until convergence is achieved. The stopping criterion for the Newton-Raphson iteration is chosen as
\eqb{l}
\text{max}\ds \left\lbrace \frac{\norm{\mf^i_{n+1}}}{\norm{\mf^0_{n+1}}},\frac{\norm{\bar\mf^i_{n+1}}}{\norm{\bar\mf^0_{n+1}}} \right\rbrace \leq tol^\mathrm{NR}~,
\eqe
where $\norm{...}$ denotes the Euclidean norm. 
The value $tol^\mathrm{NR}=10^{-4}$ is observed to be sufficient for all examples to ensure convergence. This algorithm is also known as a \textit{predictor-multicorrector algorithm}, with \eqref{e:init_g} as the \textit{prediction} and \eqref{e:upd_g} as the \textit{multicorrection}.\\
\\
\textbf{Remark:} The choices \eqref{e:init_g}, \eqref{e:upd_g} show good convergence of Newton's method for fluid structure interaction applications \citep{bazilevs13}, and good convergence is also achieved for the examples in this work.

\section{Linearization}\label{sec:lin}
The linearization of the mechanical force vector $\mf^e:=\mf^e_\mathrm{in}+\mf^e_\mathrm{int}-\mf^e_\mathrm{ext}$ of finite element $\Omega^e$ \eqref{e:finext} with respect to the nodal positions of $\Omega^e$, $\mx_e$, can be found in \citet{solidshell2}. 
The linearization of $\mf^e$ with respect to the nodal phase variables of $\Omega^e$, $\bphi_e$, according to \eqref{e:finext}, is
\eqb{lll}
\Delta_\phi \mf^e = \big[{\mk}^e_{\sig\phi}+{\mk}^e_{M\phi}\big]\,\Delta\bphi_e\,,
\eqe
with
\eqb{lll}
{\mk}^e_{\sig\phi} \dis \ds\int_{\Omega_0^e}\frac{\partial\tau^{\alpha\beta}}{\partial \phi}\,\mN^\mrT_{\!,\alpha}\,\ba^h_\beta\,\bar\mN\,\dif A
+ \ds\int_{\Omega_0^e}\frac{\partial\tau^{\alpha\beta}}{\partial \phi_{;\gamma}}\,\mN^\mrT_{\!,\alpha}\,\ba^h_\beta\,\bar\mN_{\!,\gamma}\,\dif A\,,  \\[4mm]
{\mk}^e_{M\phi} \dis \ds \int_{\Omega_0^e}\frac{\partial M_0^{\alpha\beta}}{\partial \phi}\,\mN^\mrT_{\!;\alpha\beta}\,\bn^h\,\bar\mN\,\dif A\,,
\eqe
where $\tau^{\alpha\beta}:=J\sig^{\alpha\beta}$ and $M^{\alpha\beta}_0:=JM^{\alpha\beta}$.
According to Sec.~\ref{s:sigM}, we find
\eqb{lll}
\ds\pa{\tau^{\alpha\beta}}{\phi} \is \ds \frac{K'}{2}\left(J^2-1 \right)a^{\alpha\beta}+\frac{G'}{2\,J}\left(2\,A^{\alpha\beta}-I_1\,a^{\alpha\beta}\right)-J\,\eta'\,\dot{a}^{\alpha\beta} \,, \\[4mm]
\ds\pa{\tau^{\alpha\beta}}{\phi_{;\gamma}} \is J\lambda \left(a^{\alpha\beta}\,a^{\gamma\delta}-a^{\alpha\gamma}\,a^{\beta\delta}-a^{\alpha\delta}\,a^{\beta\gamma} \right)\phi_{;\delta} \,, \\[4mm]
\ds\pa{M_0^{\alpha\beta}}{\phi} \is \ds c'\left(A^{\alpha\gamma}\,b_{\gamma\delta}\,A^{\beta\delta}-B^{\alpha\beta} \right)\,.
\eqe
According to Eq.~\eqref{e:ODEphi}, the linearization of $\bar\mf^e$  with respect to the nodal positions of $\Omega^e$, $\mx_e$, is
\eqb{lll}
\Delta_\mrx \bar\mf^e = \Delta_\mrx \bar\mf^e_\mathrm{int}\,.
\eqe
Since we can write
\eqb{lll}
\bar\mf^e_\mathrm{int} 
= \ds\int_{\Omega^e_0} \Big[\bar\mN_{\!,\alpha}^\mrT\,a^{\alpha\beta}\Big(\!M \mu'_\phi - M'\big(\mu_\mri+\mu_\mathrm{el}\big)\!\Big)\,\phi_{;\beta} - \Delta_\mrs\bar\mN^\mrT M\big(\mu_\mri + \mu_\mathrm{el} \Big)\Big]\,\dif A
\eqe
where $\mu_\bullet$ are the contributions given in \eqref{e:mubi2}, we obtain
\eqb{lll}
\Delta_\mrx \bar\mf^e = \Delta_\mrx\bar{\mf}^e_{\mrx 1} + \Delta_\mrx\bar{\mf}^e_{\mrx 2} + \Delta_\mrx\bar{\mf}^e_{\mrx 3}\,,
\eqe
with
\eqb{lll}
\Delta_\mrx\bar\mf^e_\mathrm{x1} \is \ds\int_{\Omega^e_0} \bar\mN_{\!,\alpha}^\mrT\Big(\!M \mu'_\phi - M'\big(\mu_\mri+\mu_\mathrm{el}\big)\!\Big)\,\phi_{;\beta}\,\Delta a^{\alpha\beta}\,\dif A\,, \\[4mm]
\Delta_\mrx\bar\mf^e_\mathrm{x2} \is -\ds\int_{\Omega^e_0} \Big(\bar\mN_{\!,\alpha}^\mrT\,a^{\alpha\beta}\,\phi_{;\beta}\,M'+\Delta_\mrs\bar\mN^\mrT M\Big)\big(\Delta_\mrx\mu_\mri + \Delta_\mrx\mu_\mathrm{el}\big)\,\dif A\,, \\[4mm]
\Delta_\mrx\bar\mf^e_\mathrm{x3} \is -\ds\int_{\Omega^e_0}M\big(\mu_\mri + \mu_\mathrm{el}\big)\, \Delta\Delta_\mrs\bar\mN^\mrT\,\dif A\,. \\[4mm]
\eqe
According to \citet{membrane} and \citet{shelltheo2} we have
\eqb{lll}
\Delta \ba_\alpha \is \mN_{\!,\alpha}\,\Delta\mx_e\,, \\[2mm]
\Delta \ba_{\alpha,\beta} \is \mN_{\!,\alpha\beta}\,\Delta\mx_e\,, \\[2mm]
\Delta J \is J \ba^\alpha\cdot\mN_{\!,\alpha}\,\Delta\mx_e\,, \\[2mm]
\Delta a^{\alpha\beta} \is - \big(\ba^\alpha a^{\beta\gamma} + \ba^\beta a^{\alpha\gamma} \big) \cdot\mN_{\!,\gamma}\,\Delta\mx_e\,,\\[2mm]
\Delta \ba^\alpha \is \big(a^{\alpha\beta}\bn\otimes\bn - \ba^\beta \otimes\ba^\alpha) \cdot\mN_{\!,\beta}\,\Delta\mx_e\,.
\eqe
With this we find
\eqb{l}
\Delta\Gamma^\gamma_{\alpha\beta} = \big(\ba^\gamma\cdot\mN_{;\alpha\beta} + a^{\gamma\delta}\ba_{\alpha;\beta}\cdot\mN_{\!,\delta}\big)\,\Delta\mx_e\,,
\eqe
and thus
\eqb{l}
a^{\alpha\beta}\,\Delta\Gamma^\gamma_{\alpha\beta} = \mN^{\gamma}\,\Delta\mx_e\,,
\eqe
where
\eqb{l}
\mN^\gamma := a^{\alpha\beta}\,\ba^\gamma\cdot\mN_{;\alpha\beta} + 2H\,a^{\gamma\delta}\,\bn\cdot\mN_{\!,\delta}\,.
\eqe
From this follows
\eqb{l}
\Delta\Delta_\mrs\bar\mN^\mrT = \mM_{\Delta_\mrs}\,\Delta\mx_e\,,
\eqe
where 
\eqb{l}
\mM_{\Delta_\mrs} := \ds -2\,\bar\mN^\mrT_{;\alpha\beta}\,\ba^\alpha\,a^{\beta\gamma}\cdot\mN_{\!,\gamma} - \bar\mN^\mrT_{\!,\gamma}\,\mN^\gamma \,.
\eqe
Hence $\Delta_\mrx\Delta_\mrs\phi = \bphi_e^\mrT\,\mM_{\Delta_\mrs}\,\Delta\mx_e$ and
\eqb{l}
\Delta_\mrx\mu_\mri = \mN_\mri\,\Delta\mx_e\,,
\eqe
with
\eqb{l}
\mN_\mri := \mu_\mri\,\ba^\alpha\cdot\mN_{\!,\alpha} - J\lambda\,\bphi_e^\mrT\,\mM_{\Delta_\mrs}\,.
\label{e:mNi}\eqe
Similarly, we have
\eqb{l}
\Delta_\mrx\mu_\mathrm{el} = \mN_\mathrm{el}\,\Delta\mx_e\,,
\eqe
with
\eqb{l}
\mN_\mathrm{el} := {\tau^{\alpha\beta}_\mathrm{el}}'\ba_\alpha\cdot\mN_{\!,\beta} + {M^{\alpha\beta}_0}'\bn\cdot\mN_{\!;\alpha\beta} \,,
\eqe
since
\eqb{lll}
\ds\pa{\mu_\mathrm{el}}{a_{\alpha\beta}} \is \ds\pa{\tau^{\alpha\beta}_\mathrm{el}}{\phi} := {\tau^{\alpha\beta}_\mathrm{el}}'\,,\\[4mm]
\ds\pa{\mu_\mathrm{el}}{b_{\alpha\beta}} \is \ds\pa{M^{\alpha\beta}_0}{\phi} := {M^{\alpha\beta}_0}'\,.
\eqe
We can thus write
\eqb{lll}
\Delta_\mrx \bar\mf^e = \big[\bar{\mk}^e_{\mrx 1} + \bar{\mk}^e_{\mrx 2} + \bar{\mk}^e_{\mrx 3}\big]\,\Delta\mx_e\,,
\eqe
where
\eqb{lll}
\bar{\mk}^e_{\mrx 1} \dis -\ds \int_{\Omega^e_0} \left(M\,\mu'_\phi- M'\,\left(\mu_i\,+\mu_\mathrm{el}\right)\right)\bar\mN^\mrT_{\!,\alpha}\,\phi_{;\beta}\big(\ba^\alpha a^{\beta\gamma} + \ba^\beta a^{\alpha\gamma} \big) \cdot\mN_{\!,\gamma}\,\dif A\,,\\[4mm]
\bar{\mk}^e_{\mrx 2} \dis -\ds \int_{\Omega^e_0} \Big(\bar\mN_{\!,\alpha}^\mrT\,a^{\alpha\beta}\,\phi_{;\beta}\,M'+\Delta_\mrs\bar\mN^\mrT M\Big)\big(\mN_\mri + \mN_\mathrm{el}\big)\,\dif A\,,\\[4mm]
\bar{\mk}^e_{\mrx 3} \dis \ds -\int_{\Omega^e_0}M\big(\mu_\mri + \mu_\mathrm{el}\big)\,\mM_{\Delta_\mrs}\, \dif A\,.
\eqe

The linearization of $\bar{\mf}^e$ with respect to the phase variables of $\Omega^e$, $\bphi_e$, according to \eqref{e:ODEphi}, is
\eqb{lll}
\Delta_\phi\bar{\mf}_\mathrm{int}^e = \ds\big[\bar\mk^e + \bar\mk_{\phi 1}^e + \bar\mk_{\phi 2}^e + \bar\mk_{\phi 3}^e\big]\,\Delta\bphi_e~,
\eqe
with $\bar\mk^e$ given in \eqref{e:ODEphi}, and 
\eqb{lll}
\bar\mk_{\phi 1}^e = \ds \int_{\Omega^e_0} \bar\mN_{,\alpha}^\mrT\,a^{\alpha\beta}\,\phi_{,\beta}\left(M'\left(\mu'_\phi - \mu'_\mathrm{el} \right)+M\,\mu''_\phi-M''\left(\mu_\mri+\mu_\mathrm{el}\right) \right)\,\bar\mN\, \dif A\,,\\[4mm]
\bar\mk_{\phi 2}^e = \ds \int_{\Omega^e_0} \bar\mN_{,\alpha}^\mrT\,a^{\alpha\beta}\,\phi_{,\beta}\,J\lambda\,M'\,\Delta_\mrs\bar\mN\, \dif A\,,\\[4mm]
\bar\mk_{\phi 3}^e = -\ds \int_{\Omega^e_0} \Delta_\mrs\bar\mN^\mrT\,\left( M'\,(\mu_\mri+\mu_\mathrm{el})+M\,\mu'_\mathrm{el}\right) \bar\mN\,\dif A\,.
\eqe

\section{Error estimation}
\subsection{{Error estimates for the mechanical field}}
\label{App:B}

As proposed in \cite{hulbert95}, the local error of the deformation and the velocity is given by
\eqb{lll}
\be^d_{n+1} = \mx_{n+1}-\hat\mx_{n+1}\,,\qquad \be^v_{n+1} = \dot\mx_{n+1}-\dot{\hat{\mx}}_{n+1}\,,
\label{eq:1}
\eqe
where $\hat\mx_{n+1}:=\hat\mx_{n+1}(t_{n+1})$ and $\dot{\hat{\mx}}_{n+1}:=\dot{\hat{\mx}}_{n+1}(t_{n+1})$ are the solutions of the local problem. Expressions for $\hat\mx_{n+1}$ and $\dot{\hat{\mx}}_{n+1}$ are obtained by a Taylor series with finite remainder about $t_n$ (Appendix \ref{App:A}, Eq.~\eqref{eq:tay}). At time $t_n$, $\mx_{n}=\hat\mx_{n}$ and $\dot{\mx}_{n}=\dot{\hat{\mx}}_{n}$ holds. Using the Newmark formulae (Eq.~\eqref{eq:NM}) for $\mx_{n+1}$ and $\dot\mx_{n+1}$, the local errors can be expressed as
\eqb{lll}
\be^d_{n+1} = \ds \Delta t_{n+1}^2 \left(\beta \Delta \ddot{\mx}_n -\frac{1}{6}\Delta t_{n+1}\dddot{\hat\mx} (\xi_u)+\frac{1}{2}\left(\ddot{\mx}_n - \ddot{\hat\mx}_n \right) \right)\,,\\[4mm]
\be^v_{n+1} = \ds \Delta t_{n+1} \left(\gamma \Delta \ddot{\mx}_n -\frac{1}{2}\Delta t_{n+1}\dddot{\hat\mx} (\xi_v)+\left(\ddot{\mx}_n - \ddot{\hat\mx}_n \right) \right)\,,
\label{eq:2}
\eqe
where $\Delta \ddot{\mx}_n := \ddot{\mx}_{n+1}-\ddot{\mx}_n$.
Values for $\dddot{\hat\mx} (\xi_u)$ and $\dddot{\hat\mx} (\xi_v)$ with $\xi_u,\xi_v \in [t_n,t_{n+1}]$ are obtained by the following approximation
\eqb{lll}
\dddot{\hat\mx} (\xi_u) = \dddot{\hat\mx} (\xi_v) \approx \Delta t_{n+1}^{-1} \left(\ddot{\hat\mx}_{n+1}-\ddot{\hat\mx}_n \right)\,.
\eqe
Substituting Eq.~\eqref{eq:2} into~\eqref{eq:1} and employing the results in the basic form of the generalized-$\alpha$ method
\eqb{lll}
\mM \ddot\mx_{n+1-\alpha_\mrm} + \mf_{\mathrm{int}}(\mx_{n+1-\alpha_\mrf},\dot{\mx}_{n+1-\alpha_\mrf}) - \mf_{\mathrm{ext}}(\mx_{n+1-\alpha_\mrf},\dot{\mx}_{n+1-\alpha_\mrf},t_{n+1-\alpha_\mrf})=\boldsymbol{0}\,,
\label{eq:gov}
\eqe
results in
\eqb{lll}
\mM (\ddot{\mx}_{n+1} - \ddot{\hat\mx}_{n+1}) = \alpha_\mrm\mM\Delta\ddot{\mx}_n-\alpha_\mrf\mM\Delta\ddot{\hat\mx}_n-(1-\alpha_\mrf)\mK_\mrx\be^d_{n+1}-(1-\alpha_\mrf)\mC\be^v_{n+1}\,.
\eqe
Note, that Eq.~\eqref{eq:gov} is solved at intermediate time steps and we use
\eqb{lll}
\mx_{n+1-\alpha_\mrf}=(1-\alpha_\mrf)\mx_{n+1}+\alpha_\mrf\mx_n\,,\\[4mm]
\dot\mx_{n+1-\alpha_\mrf}=(1-\alpha_\mrf)\dot\mx_{n+1}+\alpha_\mrf\dot\mx_n\,,\\[4mm]
\ddot\mx_{n+1-\alpha_\mrm}=(1-\alpha_\mrm)\ddot\mx_{n+1}+\alpha_\mrm\ddot\mx_n\,. 
\eqe
Replacing $n$ by $n-1$, and after some algebraic manipulations, we obtain
\eqb{lll}
(1-\alpha_\mrf)\mM\Delta \ddot{\hat\mx}_n = \mM\left((1-\alpha_\mrm)\Delta \ddot{\mx}_n+\alpha_\mrm\Delta \ddot{\mx}_{n-1}-\alpha_\mrf\Delta \ddot{\hat\mx}_{n-1}\right)+(1-\alpha_\mrf)(\mK_\mrx\Delta\be^d_n+\mC\Delta\be^v_{n})\,,\\[4mm]
\mM (\ddot{\mx}_n - \ddot{\hat\mx}_n )=(1-\alpha_\mrf)\mM\Delta \ddot{\hat\mx}_n -(1-\alpha_\mrm)\mM\Delta \ddot{\mx}_n-(1-\alpha_\mrf)(\mK_\mrx\be^d_{n+1}-\mC\be^v_{n+1})\,.
\label{eq:4}
\eqe
Multiplying Eq.~\eqref{eq:2} by $\mM$ and using Eq.~\eqref{eq:4} results in
\eqb{lll}
\begin{bmatrix}
\ds\mM+\frac{1}{6}\Delta t^2_{n+1}\mK_\mrx & \ds\frac{1}{6}\Delta t^2_{n+1}\bC\\[4mm]
\ds\frac{1}{2}\Delta t^2_{n+1}\mK_\mrx & \ds\mM+\frac{1}{2}\Delta t^2_{n+1}\mC
\end{bmatrix}
\begin{bmatrix}
\ds\be^d_{n+1}\\[4mm]
\ds\Delta t_{n+1}\be^v_{n+1}
\end{bmatrix}
=\\[12mm]
\Delta t_{n+1}^2
\begin{bmatrix}
\ds\left(\left(\beta-\frac{1-\alpha_\mrm}{6\left(1-\alpha_\mrf \right)} \right)\Delta\ddot\mx_n + \left(\frac{1}{6\left(1-\alpha_\mrf \right)}-\frac{1}{2} \right)\bw_n\right)\\[4mm]
\ds\left(\left(\gamma-\frac{1-\alpha_\mrm}{2\left(1-\alpha_\mrf \right)} \right)\Delta\ddot\mx_n + \left(\frac{1}{2\left(1-\alpha_\mrf \right)}-1 \right)\bw_n\right)
\end{bmatrix}
+\\[12mm]\Delta t_{n+1}^2
\begin{bmatrix}
\ds \left(\frac{1}{2}\alpha_\mrf-\frac{1}{3} \right)\left(\mK_\mrx\be^d_{n}+\mC\be^v_n \right)\\[4mm]
\ds \left(\alpha_\mrf-\frac{1}{2} \right)\left(\mK_\mrx\be^d_{n}+\mC\be^v_n \right)
\end{bmatrix}\,,
\label{eq:5}
\eqe
with
\eqb{lll}
\bw_n = \alpha_\mrm\Delta\ddot\mx_{n-1}-\alpha_\mrf\Delta\ddot{\hat{\mx}}_{n-1}\,.
\eqe
Dropping the higher order terms, $\Delta t_{n+1}^2\mK_\mrx$ and $\Delta t_{n+1}\mC$, the local errors are expressed as
\eqb{lll}
\be_{n+1}^d = \ds \Delta t_{n+1}^2\left(\left(\beta-\frac{1-\alpha_\mrm}{6\left(1-\alpha_\mrf \right)} \right)\Delta\ddot\mx_n + \left(\frac{1}{6\left(1-\alpha_\mrf \right)}-\frac{1}{2} \right)\bw_n\right)\,,\\[4mm]
\be_{n+1}^v = \ds \Delta t_{n+1}\left(\left(\gamma-\frac{1-\alpha_\mrm}{2\left(1-\alpha_\mrf \right)} \right)\Delta\ddot\mx_n + \left(\frac{1}{2\left(1-\alpha_\mrf \right)}-1 \right)\bw_n\right)\,.
\label{eq:6}
\eqe
Rewriting $\bw_n$ into
\eqb{lll}
\bw_n = \ds \frac{\alpha_\mrm-\alpha_\mrf}{(1-\alpha_\mrf)^2}\Delta\ddot\mx_{n}-\frac{\alpha_\mrf}{1-\alpha_\mrf}\bw_{n-1},\quad \text{with}\quad \bw_0=\boldsymbol{0}\,,
\eqe
following the idea of \cite{hulbert95} and interpreting $\bw_n$ as a history vector, Eq.~\eqref{eq:6} can be expressed as Eq.~\eqref{eq:ede} with constants
\eqb{lll}
c_1^\mrd = \ds \beta - \frac{1-\alpha_\mrm}{6\left(1-\alpha_\mrf \right)}\,, \qquad c_2^\mrd = \ds \left(1+\rho_\infty\right)\left(1-\rho_\infty\right)\left(\frac{1}{6}-\frac{1}{2}\left(1-\alpha_\mrf \right) \right)\,.
\label{eq:cd}
\eqe
\subsection{{Error estimate for the phase field}}
\label{App:B2}
The local error of the phase field is given by
\eqb{lll}
\be^p_{n+1} = \bphi_{n+1}-\hat\bphi_{n+1}\,.
\label{eq:ep}
\eqe
By following the same approach as for the mechanical field, Eq.~\eqref{eq:ep} can be expressed as
\eqb{lll}
\be^p_{n+1} = \ds \Delta t_{n+1} \left(\gamma \Delta \dot{\bphi}_n -\frac{1}{2}\Delta t_{n+1}\ddot{\hat\bphi} (\xi_\phi)+\left(\dot{\bphi}_n - \dot{\hat\bphi}_n \right) \right)\,,
\label{eq:ep2}
\eqe
where $\Delta \dot{\bphi}_n = \dot{\bphi}_{n+1}-\dot{\bphi}_n$. At time $t_n$, $\bphi_{n}=\hat\bphi_{n}$ holds. Values for $\ddot{\hat\bphi} (\xi_\phi)$ with $\xi_\phi \in [t_n,t_{n+1}]$ are obtained by the following approximation
\eqb{lll}
\ddot{\hat\bphi} (\xi_\phi) \approx \Delta t_{n+1}^{-1} \left(\dot{\hat\bphi}_{n+1}-\dot{\hat\bphi}_n \right)\,.
\eqe
Substituting Eq.~\eqref{eq:ep} and Eq.~\eqref{eq:ep2} and employing the results in the basic form of the generalized-$\alpha$ method
\eqb{lll}
\bar\mM\dot{\bphi}_{n+1-\alpha_\mrm} + \bar\mf_{\mathrm{int}}(\bphi_{n+1-\alpha_\mrf})=\boldsymbol{0}\,,
\label{eq:gov2}
\eqe
results in
\eqb{lll}
\bar\mM(\dot{\bphi}_{n+1}-\dot{\hat{\bphi}}_{n+1})=\alpha_\mrm\bar\mM\Delta\dot{\bphi}_{n}-\alpha_\mrf\bar\mM\Delta\dot{\hat{\bphi}}_{n}-(1-\alpha_\mrf)\bar{\mK}_\phi\be^p_{n+1}\,.
\eqe
Note that Eq.~\eqref{eq:gov2} is solved at intermediate time steps and we use
\eqb{lll}
\bphi_{n+1-\alpha_\mrf}=(1-\alpha_\mrf)\bphi_{n+1}+\alpha_\mrf\bphi_n\,,\\[4mm]
\dot\bphi_{n+1-\alpha_\mrm}=(1-\alpha_\mrm)\dot\bphi_{n+1}+\alpha_\mrm\dot\bphi_n\,. 
\eqe
Following a similar derivation as in Appendix \ref{App:B}, the local error can be expressed as
\eqb{lll}
\be_{n+1}^p = \ds \Delta t_{n+1}\left(\left(\gamma-\frac{1-\alpha_\mrm}{2\left(1-\alpha_\mrf \right)} \right)\Delta\dot{\bphi}_n + \left(\frac{1}{2\left(1-\alpha_\mrf \right)}-1 \right)\bw_n\right)\,,
\label{eq:ep3}
\eqe
with
\eqb{lll}
\bw_n = \alpha_\mrm\Delta\dot\bphi_{n-1}-\alpha_\mrf\Delta\dot{\hat{\bphi}}_{n-1}\,.
\eqe
Rewriting $\bw_n$ into
\eqb{lll}
\bw_n = \ds \frac{\alpha_\mrm-\alpha_\mrf}{(1-\alpha_\mrf)^2}\Delta\dot\bphi_{n}-\frac{\alpha_\mrf}{1-\alpha_\mrf}\bw_{n-1}\,,\quad \text{with}\quad \bw_0=\boldsymbol{0}\,,
\eqe
Eq.~\eqref{eq:ep3} can be expressed as Eq.~\eqref{eq:ep4} with constants
\eqb{lll}
c_1^\mrp = \ds \gamma - \frac{1-\alpha_\mrm}{2\left(1-\alpha_\mrf \right)}\,, \qquad c_2^\mrp = \ds \left(1+\rho_\infty\right)\left(1-\rho_\infty\right)\left(\frac{1}{2}-1+\alpha_\mrf \right)\,.
\label{eq:cp}
\eqe

\subsection{{Taylor series expansion and approximations}}
\label{App:A}
For the derivation of the error estimates in Sec.~\ref{App:B} the
Taylor series expansion of $\hat\mx_{n+1}$ and $\dot{\hat{\mx}}_{n+1}$ with finite remainder about $t_n$
\eqb{lll}
\hat\mx_{n+1} &= \ds\hat{\mx}_{n}+\Delta t_{n+1}\dot{\hat{\mx}}_{n}+\frac{1}{2}\Delta t^2_{n+1}\ddot{\hat{\mx}}_{n}+\frac{1}{6}\Delta t^3_{n+1}\dddot{\hat{\mx}}_{n}+\mathcal{O}(\Delta t^4_{n+1})\,, \\[4mm]
\dot{\hat{\mx}}_{n+1}&=\ds \dot{\hat{\mx}}_{n}+\Delta t_{n+1} \ddot{\hat{\mx}}_{n} + \frac{1}{2}\Delta t_{n+1}^2\dddot{\hat{\mx}}_{n}+\mathcal{O}(\Delta t^3_{n+1})\,,
\label{eq:tay}
\eqe
are used. The derivation in Sec.~\ref{App:B2} uses the expansion of $\hat\bphi_{n+1}$ in a Taylor series with finite remainder about $t_n$, which gives
\eqb{lll}
\hat\bphi_{n+1} = \ds\hat\bphi_{n}+\Delta t_{n+1}\dot{\hat{\bphi}}_{n}+\frac{1}{2}\Delta t^2_{n+1}\ddot{\hat{\bphi}}_{n}+\mathcal{O}(\Delta t^3_{n+1})\,.
\label{eq:tay2}
\eqe
We also make use of Newmark's formulae for second order systems
\eqb{lll}
\mx_{n+1}=\ds\mx_n+\Delta t_{n+1}\dot{\mx}_n+\Delta t^2_{n+1}\left(\left(\frac{1}{2}-\beta\right)\ddot{\mx}_n+\beta\ddot{\mx}_{n+1}\right)\,,\\[4mm]
\dot{\mx}_{n+1}=\ds\mx_n+\Delta t_{n+1}\left(\left(1-\gamma\right)\ddot{\mx}_n+\gamma\ddot{\mx}_{n+1}\right)\,.
\label{eq:NM}
\eqe
and Newmark's formulae for first order systems 
\eqb{lll}
{\bphi}_{n+1}=\ds\bphi_n+\Delta t_{n+1}\left(\left(1-\gamma\right)\dot{\bphi}_n+\gamma\dot{\bphi}_{n+1}\right)\,.
\label{eq:NM2}
\eqe


\bigskip

\bibliographystyle{apalike}
\bibliography{bibliography}

\begin{thebibliography}{}

\bibitem[Akkerman et~al., 2008]{akkerman2008role}
Akkerman, I., Bazilevs, Y., Calo, V.~M., Hughes, T. J.~R., and Hulshoff, S.
  (2008).
\newblock The role of continuity in residual-based variational multiscale
  modeling of turbulence.
\newblock {\em Computational Mechanics}, {\bf 41}(3):371--378.

\bibitem[Barrett et~al., 1999]{barrett99}
Barrett, J.~W., Blowey, J.~F., and Garcke, H. (1999).
\newblock Finite element approximation of the {C}ahn-{H}illiard equation with
  degenerate mobility.
\newblock {\em SIAM Journal on Numerical Analysis}, {\bf 37}(1):286--318.

\bibitem[Bartezzaghi et~al., 2015]{bartezzaghi15}
Bartezzaghi, A., Ded\`e, L., and Quarteroni, A. (2015).
\newblock Isogeometric analysis of high order partial differential equations on
  surfaces.
\newblock {\em Comput. Meth. Appl. Mech. Engrg.}, {\bf 295}:446--469.

\bibitem[Bartezzaghi et~al., 2016]{bartezzaghi16}
Bartezzaghi, A., Ded\`e, L., and Quarteroni, A. (2016).
\newblock Isogeometric analysis of geometric partial differential equations.
\newblock {\em Comput. Meth. Appl. Mech. Engrg.}, {\bf 311}:625--647.

\bibitem[Baumgart et~al., 2003]{baumgart03}
Baumgart, T., Hess, S.~T., and Webb, W.~W. (2003).
\newblock Imaging coexisting fluid domains in biomembrane models coupling
  curvature and line tension.
\newblock {\em Nature}, 425(6960):821--824.

\bibitem[Bazilevs et~al., 2013]{bazilevs13}
Bazilevs, Y., Takizawa, K., and Tezduyar, T.~E. (2013).
\newblock {ALE} and space-time methods for {FSI}.
\newblock In {\em Computational Fluid-Structure Interaction}, pages 111--137.
  John Wiley \& Sons, Ltd.

\bibitem[Bertal\'{i}mo et~al., 2001]{bertalimo01}
Bertal\'{i}mo, M., Cheng, L.-T., Osher, S., and Sapiro, G. (2001).
\newblock Variational problems and partial differential equations on implicit
  surfaces.
\newblock {\em Journal of Computational Physics}, {\bf 174}(2):759--780.

\bibitem[Borden et~al., 2016]{borden16-1}
Borden, M.~J., Hughes, T. J.~R., Landis, C.~M., Anvari, A., and Lee, I.~J.
  (2016).
\newblock A phase-field formulation for fracture in ductile materials: Finite
  deformation balance law derivation, plastic degradation, and stress
  triaxiality effects.
\newblock {\em Comput. Meth. Appl. Mech. Engrg.}, {\bf 312}:130--166.
\newblock Special Issue on Phase Field Approaches to Fracture.

\bibitem[Borden et~al., 2014]{borden14-1}
Borden, M.~J., Hughes, T. J.~R., Landis, C.~M., and Verhoosel, C.~V. (2014).
\newblock A higher-order phase-field model for brittle fracture: Formulation
  and analysis within the isogeometric analysis framework.
\newblock {\em Comput. Meth. Appl. Mech. Engrg.}, {\bf 273}:100--118.

\bibitem[Borden et~al., 2011]{borden2011isogeometric}
Borden, M.~J., Scott, M.~A., Evans, J.~A., and Hughes, T. J.~R. (2011).
\newblock Isogeometric finite element data structures based on {B}{\'e}zier
  extraction of {NURBS}.
\newblock {\em International Journal for Numerical Methods in Engineering},
  {\bf 87}(1--5):15--47.

\bibitem[Borden et~al., 2012]{borden201277}
Borden, M.~J., Verhoosel, C.~V., Scott, M.~A., Hughes, T. J.~R., and Landis,
  C.~M. (2012).
\newblock A phase-field description of dynamic brittle fracture.
\newblock {\em Comput. Meth. Appl. Mech. Engrg.}, {\bf 217--220}:77--95.

\bibitem[Cahn, 1961]{cahn61}
Cahn, J.~W. (1961).
\newblock On spinodal decomposition.
\newblock {\em Acta Metallurgica}, {\bf 9}(9):795--801.

\bibitem[Cahn and Hilliard, 1958]{cahn58-1}
Cahn, J.~W. and Hilliard, J.~E. (1958).
\newblock {Free Energy of a Nonuniform System. I. Interfacial Free Energy}.
\newblock {\em The Journal of Chemical Physics}, {\bf 28}(2):258--267.

\bibitem[Chung and Hulbert, 1993]{chung93}
Chung, J. and Hulbert, G.~M. (1993).
\newblock {A Time Integration Algorithm for Structural Dynamics With Improved
  Numerical Dissipation: The Generalized-alpha Method}.
\newblock {\em Journal of Applied Mechanics}, {\bf 60}(2):371--375.

\bibitem[Ciarlet, 1993]{ciarlet}
Ciarlet, P.~G. (1993).
\newblock {\em Mathematical Elasticity: Three Dimensional Elasticity}.
\newblock North-Holland.

\bibitem[Cottrell et~al., 2009]{cottrell}
Cottrell, J.~A., Hughes, T. J.~R., and Bazilevs, Y. (2009).
\newblock {\em Isogeometric Analysis}.
\newblock Wiley.

\bibitem[Cottrell et~al., 2007]{cottrell2007studies}
Cottrell, J.~A., Hughes, T. J.~R., and Reali, A. (2007).
\newblock Studies of refinement and continuity in isogeometric structural
  analysis.
\newblock {\em Comput. Meth. Appl. Mech. Engrg.}, {\bf 196}(41):4160--4183.

\bibitem[Cottrell et~al., 2006]{cottrell2006isogeometric}
Cottrell, J.~A., Reali, A., Bazilevs, Y., and Hughes, T. J.~R. (2006).
\newblock Isogeometric analysis of structural vibrations.
\newblock {\em Comput. Meth. Appl. Mech. Engrg.}, {\bf 195}(41):5257--5296.

\bibitem[Ded\`e et~al., 2012]{dede12}
Ded\`e, L., Borden, M.~J., and Hughes, T. J.~R. (2012).
\newblock Isogeometric analysis for topology optimization with a phase field
  model.
\newblock {\em Archives of Computational Methods in Engineering}, {\bf
  19}(3):427--465.

\bibitem[Di~Leo et~al., 2014]{dileo14-1}
Di~Leo, C.~V., Rejovitzky, E., and Anand, L. (2014).
\newblock A {C}ahn-{H}illiard-type phase-field theory for species diffusion
  coupled with large elastic deformations: Application to phase-separating
  li-ion electrode materials.
\newblock {\em Journal of the Mechanics and Physics of Solids}, {\bf 70}:1--29.

\bibitem[Dokken et~al., 2013]{dokken13}
Dokken, T., Lyche, T., and Pettersen, K.~F. (2013).
\newblock Polynomial splines over locally refined box-partitions.
\newblock {\em Computer Aided Geometric Design}, {\bf 30}(3):331--356.

\bibitem[Duong et~al., 2017]{solidshell2}
Duong, T.~X., Roohbakhshan, F., and Sauer, R.~A. (2017).
\newblock A new rotation-free isogeometric thin shell formulation and a
  corresponding continuity constraint for patch boundaries.
\newblock {\em Comput. Meth. Appl. Mech. Engrg.}, {\bf 316}:43--83.
\newblock Special Issue on Isogeometric Analysis: Progress and Challenges.

\bibitem[Dziuk and Elliott, 2007]{dziuk07}
Dziuk, G. and Elliott, C.~M. (2007).
\newblock {Finite elements on evolving surfaces}.
\newblock {\em IMA Journal of Numerical Analysis}, {\bf 27}(2):262--292.

\bibitem[Dziuk and Elliott, 2012]{dziuk12}
Dziuk, G. and Elliott, C.~M. (2012).
\newblock A fully discrete evolving surface finite element method.
\newblock {\em SIAM Journal on Numerical Analysis}, {\bf 50}(5):2677--2694.

\bibitem[Ebner et~al., 2013]{ebner13}
Ebner, M., Marone, F., Stampanoni, M., and Wood, V. (2013).
\newblock Visualization and quantification of electrochemical and mechanical
  degradation in li ion batteries.
\newblock {\em Science}, {\bf 342}(6159):716--720.

\bibitem[Eilks and Elliott, 2008]{eilks08}
Eilks, C. and Elliott, C.~M. (2008).
\newblock Numerical simulation of dealloying by surface dissolution via the
  evolving surface finite element method.
\newblock {\em Journal of Computational Physics}, {\bf 227}(23):9727--9741.

\bibitem[Elliott et~al., 1989]{elliott89}
Elliott, C.~M., French, D.~A., and Milner, F.~A. (1989).
\newblock A second order splitting method for the {C}ahn-{H}illiard equation.
\newblock {\em Numerische Mathematik}, {\bf 54}(5):575--590.

\bibitem[Elliott and Ranner, 2015]{elliott15}
Elliott, C.~M. and Ranner, T. (2015).
\newblock Evolving surface finite element method for the {C}ahn-{H}illiard
  equation.
\newblock {\em Numerische Mathematik}, {\bf 129}(3):483--534.

\bibitem[Elliott and Stinner, 2009]{elliott09}
Elliott, C.~M. and Stinner, B. (2009).
\newblock Analysis of a diffuse interface approach to an advection diffusion
  equation on a moving surface.
\newblock {\em Mathematical Models and Methods in Applied Sciences}, {\bf
  19}(05):787--802.

\bibitem[Elliott and Stinner, 2010]{elliot10}
Elliott, C.~M. and Stinner, B. (2010).
\newblock Modeling and computation of two phase geometric biomembranes using
  surface finite elements.
\newblock {\em Journal of Computational Physics}, {\bf 229}(18):6585--6612.

\bibitem[Embar et~al., 2013]{embar13}
Embar, A., Dolbow, J., and Fried, E. (2013).
\newblock Microdomain evolution on giant unilamellar vesicles.
\newblock {\em Biomechanics and Modeling in Mechanobiology}, {\bf
  12}(3):597--615.

\bibitem[Giannelli et~al., 2012]{giannelli2012thb}
Giannelli, C., J{\"u}ttler, B., and Speleers, H. (2012).
\newblock {THB}-splines: The truncated basis for hierarchical splines.
\newblock {\em Computer Aided Geometric Design}, {\bf 29}(7):485--498.

\bibitem[Gomez et~al., 2008]{gomez08-1}
Gomez, H., Calo, V.~M., Bazilevs, Y., and Hughes, T. J.~R. (2008).
\newblock Isogeometric analysis of the {C}ahn-{H}illiard phase-field model.
\newblock {\em Comput. Meth. Appl. Mech. Engrg.}, {\bf 197}(49--50):4333--4352.

\bibitem[H{\"o}llig, 2003]{hollig2003finite}
H{\"o}llig, K. (2003).
\newblock {\em Finite Element Methods with B-Splines}.
\newblock Frontiers in Applied Mathematics. Society for Industrial and Applied
  Mathematics.

\bibitem[Hughes et~al., 2005]{hughes05}
Hughes, T. J.~R., Cottrell, J.~A., and Bazilevs, Y. (2005).
\newblock Isogeometric analysis: {CAD}, finite elements, {NURBS}, exact
  geometry and mesh refinement.
\newblock {\em Comput. Meth. Appl. Mech. Engrg.}, {\bf 194}:4135--4195.

\bibitem[Hulbert and Jang, 1995]{hulbert95}
Hulbert, G.~M. and Jang, I. (1995).
\newblock Automatic time step control algorithms for structural dynamics.
\newblock {\em Comput. Meth. Appl. Mech. Engrg.}, {\bf 126}(1):155--178.

\bibitem[Johannessen et~al., 2014]{johannessen14}
Johannessen, K.~A., Kvamsdal, T., and Dokken, T. (2014).
\newblock Isogeometric analysis using {LR} {B}-splines.
\newblock {\em Comput. Meth. Appl. Mech. Engrg.}, {\bf 269}:471--514.

\bibitem[K\"{a}stner et~al., 2016]{kastner16-1}
K\"{a}stner, M., Metsch, P., and de~Borst, R. (2016).
\newblock Isogeometric analysis of the {C}ahn-{H}illiard equation - a
  convergence study.
\newblock {\em Journal of Computational Physics}, {\bf 305}(C):360--371.

\bibitem[Katira et~al., 2016]{Katira2016}
Katira, S., Mandadapu, K.~K., Vaikuntanathan, S., Smit, B., and Chandler, D.
  (2016).
\newblock Pre-transition effects mediate forces of assembly between
  transmembrane proteins: The orderphobic effect.
\newblock {\em Biophysical Journal}, {\bf 110}(3, Supplement 1):567a.

\bibitem[Li, 2015]{li2015some}
Li, X. (2015).
\newblock Some properties for analysis-suitable {T}-splines.
\newblock {\em Journal of Computational Mathematics}, {\bf 33}:428--442.

\bibitem[Lipton et~al., 2010]{lipton2010robustness}
Lipton, S., Evans, J.~A., Bazilevs, Y., Elguedj, T., and Hughes, T. J.~R.
  (2010).
\newblock Robustness of isogeometric structural discretizations under severe
  mesh distortion.
\newblock {\em Comput. Meth. Appl. Mech. Engrg.}, {\bf 199}(5):357--373.

\bibitem[Liu et~al., 2013]{liu13-1}
Liu, J., Ded\`e, L., Evans, J.~A., Borden, M.~J., and Hughes, T. J.~R. (2013).
\newblock Isogeometric analysis of the advective {C}ahn-{H}illiard equation:
  Spinodal decomposition under shear flow.
\newblock {\em Journal of Computational Physics}, {\bf 242}:321--350.

\bibitem[Lowengrub et~al., 2009]{lowengrub09}
Lowengrub, J.~S., R\"atz, A., and Voigt, A. (2009).
\newblock Phase-field modeling of the dynamics of multicomponent vesicles:
  Spinodal decomposition, coarsening, budding, and fission.
\newblock {\em Physical Review E}, {\bf 79}:031926.

\bibitem[Lubich et~al., 2013]{lubich13}
Lubich, C., Mansour, D., and Venkataraman, C. (2013).
\newblock {Backward difference time discretization of parabolic differential
  equations on evolving surfaces}.
\newblock {\em IMA Journal of Numerical Analysis}, {\bf 33}(4):1365--1385.

\bibitem[{McWhirter} et~al., 2004]{mcwhirter04}
{McWhirter}, J., {Ayton}, G., and {Voth}, G. (2004).
\newblock {Coupling Field Theory with Mesoscopic Dynamical Simulations of
  Multicomponent Lipid Bilayers}.
\newblock {\em Biophysical Journal}, {\bf 87}:3242--3263.

\bibitem[Mercker et~al., 2012]{mercker12}
Mercker, M., Ptashnyk, M., K\"{u}hnle, J., Hartmann, D., Weiss, M., and
  J\"{a}ger, W. (2012).
\newblock A multiscale approach to curvature modulated sorting in biological
  membranes.
\newblock {\em Journal of Theoretical Biology}, {\bf 301}(Supplement C):67--82.

\bibitem[Morganti et~al., 2015]{morganti2015patient}
Morganti, S., Auricchio, F., Benson, D.~J., Gambarin, F.~I., Hartmann, S.,
  Hughes, T. J.~R., and Reali, A. (2015).
\newblock Patient-specific isogeometric structural analysis of aortic valve
  closure.
\newblock {\em Comput. Meth. Appl. Mech. Engrg.}, {\bf 284}:508--520.

\bibitem[Myles and Peters, 2011]{myles2011c2polar}
Myles, A. and Peters, J. (2011).
\newblock {$C^2$} splines covering polar configurations.
\newblock {\em Computer-Aided Design}, {\bf 43}:1322--1329.

\bibitem[Naghdi, 1973]{naghdi1971theory}
Naghdi, P.~M. (1973).
\newblock The theory of shells and plates.
\newblock In Truesdell, C., editor, {\em Linear Theories of Elasticity and
  Thermoelasticity: Linear and Nonlinear Theories of Rods, Plates, and Shells},
  pages 425--640, Berlin, Heidelberg. Springer.

\bibitem[Nguyen and Peters, 2016]{nguyen2016refinable}
Nguyen, T. and Peters, J. (2016).
\newblock Refinable {$C^1$} spline elements for irregular quad layout.
\newblock {\em Computer Aided Geometric Design}, {\bf 43}:123--130.

\bibitem[Piegl and Tiller, 2012]{piegl2012nurbs}
Piegl, L. and Tiller, W. (2012).
\newblock {\em The {NURBS} Book}.
\newblock Springer-Verlag.

\bibitem[Rangamani et~al., 2013]{rangamani13}
Rangamani, P., Agrawal, A., Mandadapu, K.~K., Oster, G., and Steigmann, D.~J.
  (2013).
\newblock Interaction between surface shape and intra-surface viscous flow on
  lipid membranes.
\newblock {\em Biomechanics and Modeling in Mechanobiology}, {\bf
  12}(4):833--845.

\bibitem[Rangamani et~al., 2014]{rangamani14}
Rangamani, P., Mandadapu, K.~K., and Oster, G. (2014).
\newblock Protein-induced membrane curvature alters local membrane tension.
\newblock {\em Biophysical Journal}, {\bf 107}(3):751--762.

\bibitem[Reif, 1997]{reif1997refineable}
Reif, U. (1997).
\newblock A refineable space of smooth spline surfaces of arbitrary topological
  genus.
\newblock {\em Journal of Approximation Theory}, {\bf 90}:174--199.

\bibitem[Reusken, 2015]{reusken14}
Reusken, A. (2015).
\newblock Analysis of trace finite element methods for surface partial
  differential equations.
\newblock {\em IMA Journal of Numerical Analysis}, {\bf 35}(4):1568--1590.

\bibitem[Sahu et~al., 2017]{sahu17}
Sahu, A., Sauer, R.~A., and Mandadapu, K.~K. (2017).
\newblock Irreversible thermodynamics of curved lipid membranes.
\newblock {\em Physical Review E}, {\bf 96}:042409.

\bibitem[Sauer, 2014]{droplet}
Sauer, R.~A. (2014).
\newblock Stabilized finite element formulations for liquid membranes and their
  application to droplet contact.
\newblock {\em International Journal for Numerical Methods in Fluids}, {\bf
  75}(7):519--545.

\bibitem[Sauer, 2018]{cism}
Sauer, R.~A. (2018).
\newblock On the computational modeling of lipid bilayers using thin-shell
  theory.
\newblock In Steigmann, D., editor, {\em CISM Advanced School `On the role of
  mechanics in the study of lipid bilayers}, pages 221--286. Springer.

\bibitem[Sauer and Duong, 2017]{shelltheo2}
Sauer, R.~A. and Duong, T.~X. (2017).
\newblock On the theoretical foundations of thin solid and liquid shells.
\newblock {\em Mathematics and Mechanics of Solids}, {\bf 22}(3):343--371.

\bibitem[Sauer et~al., 2014]{membrane}
Sauer, R.~A., Duong, T.~X., and Corbett, C.~J. (2014).
\newblock A computational formulation for constrained solid and liquid
  membranes considering isogeometric finite elements.
\newblock {\em Comput. Meth. Appl. Mech. Engrg.}, {\bf 271}:48--68.

\bibitem[Sauer et~al., 2017]{liquidshell}
Sauer, R.~A., Duong, T.~X., Mandadapu, K.~K., and Steigmann, D.~J. (2017).
\newblock A stabilized finite element formulation for liquid shells and its
  application to lipid bilayers.
\newblock {\em Journal of Computational Physics}, {\bf 330}:436--466.

\bibitem[Schillinger et~al., 2012]{schill12}
Schillinger, D., Ded\`e, L., Scott, M.~A., Evans, J.~A., Borden, M.~J., Rank,
  E., and Hughes, T. J.~R. (2012).
\newblock An isogeometric design-through-analysis methodology based on adaptive
  hierarchical refinement of {NURBS}, immersed boundary methods, and {T}-spline
  {CAD} surfaces.
\newblock {\em Comput. Meth. Appl. Mech. Engrg.}, {\bf 249}:116--150.

\bibitem[Scott et~al., 2012]{scott12}
Scott, M., Li, X., Sederberg, T., and Hughes, T. J.~R. (2012).
\newblock Local refinement of analysis-suitable {T}-splines.
\newblock {\em Comput. Meth. Appl. Mech. Engrg.}, {\bf 213}:206--222.

\bibitem[Scott et~al., 2013]{scott2013isogeometric}
Scott, M.~A., Simpson, R.~N., Evans, J.~A., Lipton, S., Bordas, S. P.~A.,
  Hughes, T. J.~R., and Sederberg, T.~W. (2013).
\newblock Isogeometric boundary element analysis using unstructured
  {T}-splines.
\newblock {\em Comput. Meth. Appl. Mech. Engrg.}, {\bf 254}:197--221.

\bibitem[Sethian, 1999]{sethian99}
Sethian, J.~A. (1999).
\newblock Level set methods and fast marching methods. evolving interfaces in
  computational geometry, fluid mechanics, computer vision, and materials
  science.
\newblock {\em Cambridge University Press}, {Vol. \bf 3}.

\bibitem[Steigmann, 1999]{steigmann99b}
Steigmann, D.~J. (1999).
\newblock Fluid films with curvature elasticity.
\newblock {\em Archive for Rational Mechanics and Analysis}, {\bf
  150}:127--152.

\bibitem[Stein and Xu, 2014]{stein2014}
Stein, P. and Xu, B. (2014).
\newblock 3{D} isogeometric analysis of intercalation-induced stresses in
  li-ion battery electrode particles.
\newblock {\em Comput. Meth. Appl. Mech. Engrg.}, {\bf 268}:225--244.

\bibitem[Tang et~al., 2010]{tang10}
Tang, M., Carter, W.~C., and Chiang, Y.-M. (2010).
\newblock Electrochemically driven phase transitions in insertion electrodes
  for lithium-ion batteries: Examples in lithium metal phosphate olivines.
\newblock {\em Annual Review of Materials Research}, {\bf 40}(1):501--529.

\bibitem[Taylor et~al., 1997]{taylor97}
Taylor, M., Tribbia, J., and Iskandarani, M. (1997).
\newblock The spectral element method for the shallow water equations on the
  sphere.
\newblock {\em Journal of Computational Physics}, {\bf 130}(1):92--108.

\bibitem[Toshniwal et~al., 2017a]{toshniwal2017multi}
Toshniwal, D., Speleers, H., Hiemstra, R.~R., and Hughes, T. J.~R. (2017a).
\newblock Multi-degree smooth polar splines: A framework for geometric modeling
  and isogeometric analysis.
\newblock {\em Comput. Meth. Appl. Mech. Engrg.}, {\bf 316}:1005--1061.

\bibitem[Toshniwal et~al., 2017b]{toshniwal2017smooth}
Toshniwal, D., Speleers, H., and Hughes, T. J.~R. (2017b).
\newblock Smooth cubic spline spaces on unstructured quadrilateral meshes with
  particular emphasis on extraordinary points: Geometric design and
  isogeometric analysis considerations.
\newblock {\em Comput. Meth. Appl. Mech. Engrg.}, {\bf 327}:411--458.

\bibitem[Veatch and Keller, 2003]{veatch03}
Veatch, S.~L. and Keller, S.~L. (2003).
\newblock Separation of liquid phases in giant vesicles of ternary mixtures of
  phospholipids and cholesterol.
\newblock {\em Biophysical journal}, 85(5):3074--3083.

\bibitem[Wang and Du, 2008]{wang08}
Wang, X. and Du, Q. (2008).
\newblock Modelling and simulations of multi-component lipid membranes and open
  membranes via diffuse interface approaches.
\newblock {\em Journal of Mathematical Biology}, {\bf 56}:347--371.

\bibitem[Wells et~al., 2006]{wells06}
Wells, G.~N., Kuhl, E., and Garikipati, K. (2006).
\newblock A discontinuous {G}alerkin method for the {C}ahn-{H}illiard equation.
\newblock {\em Journal of Computational Physics}, {\bf 218}(2):860--877.

\bibitem[Xia et~al., 2007]{xia07}
Xia, Y., Xu, Y., and Shu, C.-W. (2007).
\newblock Local discontinuous {G}alerkin methods for the {C}ahn-{H}illiard type
  equations.
\newblock {\em Journal of Computational Physics}, {\bf 227}(1):472--491.

\bibitem[Xu et~al., 2016]{Xu2016}
Xu, B.-X., Zhao, Y., and Stein, P. (2016).
\newblock Phase field modeling of electrochemically induced fracture in li-ion
  battery with large deformation and phase segregation.
\newblock {\em GAMM-Mitteilungen}, {\bf 39}(1):92--109.

\bibitem[Zhao et~al., 2015]{zhao15-1}
Zhao, Y., Stein, P., and Xu, B.-X. (2015).
\newblock Isogeometric analysis of mechanically coupled {C}ahn-{H}illiard phase
  segregation in hyperelastic electrodes of {L}i-ion batteries.
\newblock {\em Comput. Meth. Appl. Mech. Engrg.}, {\bf 297}:325--347.

\bibitem[Zhao et~al., 2016]{zhao2016}
Zhao, Y., Xu, B.-X., Stein, P., and Gross, D. (2016).
\newblock Phase-field study of electrochemical reactions at exterior and
  interior interfaces in li-ion battery electrode particles.
\newblock {\em Comput. Meth. Appl. Mech. Engrg.}, {\bf 312}:428--446.

\bibitem[Zimmermann and Sauer, 2017]{zimmermann17}
Zimmermann, C. and Sauer, R.~A. (2017).
\newblock Adaptive local surface refinement based on {LR NURBS} and its
  application to contact.
\newblock {\em Computational Mechanics}, {\bf 60}(6):1011--1031.

\end{thebibliography}

\end{document}